\documentclass[12pt]{iopart}

\usepackage{subfig}

\newcommand{\dd}{ \,\textrm{d}}

\usepackage{graphicx,color}
\usepackage{iopams}

\usepackage{latexsym}
\usepackage{amsfonts}

\newtheorem{algorithm}{Algorithm}

\newtheorem{lemma}{Lemma}[section]

\newtheorem{remark}{Remark}[section]

\newtheorem{proposition}{Proposition}[section]

\def\n {\mathbf{n}}

\def\S{\mathbf{S}}

\def\0{\mathbf{0}}

\def\n {\mathbf{n}}

\begin{document}
\title[An iterative regularizing ensemble Kalman method]{A regularizing iterative ensemble Kalman method for PDE-constrained inverse problems}
\author{Marco A. Iglesias}
\address{School of Mathematical Sciences, The University of Nottingham,
University Park, Nottingham, NG7 2RD, United Kingdom.}
\begin{abstract}
We introduce a derivative-free computational framework for approximating solutions to nonlinear PDE-constrained inverse problems. The general aim is to merge ideas from iterative regularization with ensemble Kalman methods from Bayesian inference to develop a derivative-free stable method easy to implement in applications where the PDE (forward) model is only accessible as a black box (e.g. with commercial software). The proposed regularizing ensemble Kalman method can be derived as an approximation of the regularizing Levenberg-Marquardt (LM) scheme \cite{Hanke} in which the derivative of the forward operator and its adjoint are replaced with empirical covariances from an ensemble of elements from the admissible space of solutions. The resulting ensemble method consists of an update formula that is applied to each ensemble member and that has a regularization parameter selected in a similar fashion to the one in the LM scheme. Moreover, an early termination of the scheme is proposed according to a discrepancy principle-type of criterion.  The proposed method can be also viewed as a regularizing version of standard Kalman approaches which are often unstable unless ad-hoc fixes, such as covariance localization, are implemented.

 The aim of this paper is to provide a detailed numerical investigation of the regularizing and convergence properties of the proposed regularizing ensemble Kalman scheme; the proof of these properties is an open problem. By means of numerical experiments, we investigate the conditions under which the proposed method inherits the regularizing properties of the LM scheme of \cite{Hanke} and is thus stable and suitable for its application in problems where the computation of the Fr\'echet derivative is not computationally feasible. More concretely, we study the effect of ensemble size, number of measurements, selection of initial ensemble and tunable parameters on the performance of the method. The numerical investigation is carried out with synthetic experiments on two model inverse problems: (i) identification of conductivity on a Darcy flow model  and (ii) electrical impedance tomography with the complete electrode model. We further demonstrate the potential application of the method in solving shape identification problems that arises from the aforementioned forward models by means of a level-set approach for the parameterization of unknown geometries.  

\end{abstract}
\submitto{\IP}

\section{Introduction}\label{Intro}


We propose a computational derivative-free regularization method for approximating solutions to nonlinear PDE-constrained inverse problems. More precisely, the aim of the method is to identify parameters in PDE models given data/observations of the solution of the PDE. Inverse problems of this kind are ill-posed in the sense of stability; small perturbations of the data may have uncontrolled effects on the approximation of the unknown parameters in the PDE. Therefore, the computation of stable approximations of solutions to these inverse problems requires \textit{regularization}. Classical regularization methods (e.g. Tikhonov regularization) \cite{Engl} reformulate the inverse problem so that the regularized version can be solved with, for example, standard optimization methods.  In contrast, \textit{iterative regularization} approaches regularize while computing a stable approximation to the inverse problem \cite{Iterative}. The aim of these methods is to compute an estimate/approximation controlled by the noise level. As the noise in the data goes to zero, this approximation converges to a solution of the identification problem. 

\subsection{Contribution of this work}
Most existing iterative regularization methods \cite{Iterative} require the implementation of the Fr\'echet derivative of the forward map as well as the corresponding adjoint operator. In various applications, however, standard software for forward simulation does not provide numerical approximations of the Fr\'echet derivative and/or its associated adjoint. In some cases, even if the linearization of the forward map is computed within the forward simulation (e.g. in a Newton-type solver for nonlinear PDEs) this may not be accessible in a modular fashion suitable for an iterative regularization framework. In this paper we present a derivative-free ensemble Kalman-based iterative regularization technique for the approximation of PDE-constrained identification problems. While standard ensemble Kalman methods are aimed at approximating an inverse problem posed in a Bayesian inference framework, the objective of the present work is to merge ideas from iterative regularization with ensemble Kalman methods to develop a computational framework with the regularization properties needed to solve classical (deterministic) identification problems. The proposed framework offers the flexibility of typical implementations of ensemble Kalman methods often used for large-scale data assimilation. In particular, it uses the forward map in a black-box fashion which makes it easy to implement and thus ideal for applications where modeling and simulation are performed with complex computer codes. 

The proposed method can be derived as an approximation of the regularizing Levenberg-Marquard (LM) scheme developed by Hanke in \cite{Hanke}. More specifically, we construct an iterative ensemble method from the regularizing LM scheme by replacing operators involving the Fr\'echet derivative of the forward map and its adjoint by empirical covariances from an ensemble of elements from the parameter space. Members from this ensemble are iteratively updated according to an expression that resembles the Kalman update from standard ensemble Kalman filter/smoother methods \cite{evensen2009data}. However, in contrast to those standard methods, we propose an update formula with a regularization parameter whose selection is made similar to the one of the regularization parameter in the LM scheme \cite{Hanke} but with the aforementioned ensemble approximations of the forward map and its adjoint. Similarly, we propose an early termination of the scheme motivated by the discrepancy principle used in iterative regularization methods. The proposed ensemble Kalman regularizing scheme can then be regarded as (i) a derivative-free approximation of the regularizing LM scheme and/or (ii) a regularizing version of standard iterative Kalman methods \cite{EnKF_US}. 

An iterative regularizing ensemble Kalman method of the type presented here has been recently introduced in \cite{Yo} for the solution of Bayesian inverse problems in reservoir modeling applications. The aim of \cite{Yo} was to show that importing ideas from iterative regularization may improve the performance of Kalman methods for approximating the Bayesian posterior. The classical approach to the inverse problem is now pursued in the present work. More concretely, we wish to assess the convergence and regularizing properties of ensemble Kalman methods derived from iterative regularization techniques. In contrast to  \cite{Yo} where the framework was Bayesian and focused on reservoir applications, our objective here is to study the performance of regularizing ensemble methods for solving classical identification problems in generic PDE-constrained applications.

The contribution of the present work is threefold.  First, we show that the proposed method addresses the ill-conditioning typically exhibited by existing implementations of ensemble methods. We show that importing ideas from iterative regularization such as the discrepancy principle can offer stability that, in existing implementations are often treated with ad-hoc methodologies such as covariance localization and covariance inflation. The second contribution of this paper is to showcase the potential application of the proposed methods for addressing a wide range of parameter identification problems. We show that, with reasonable computational cost, the proposed method can be not only computationally advantageous but also robust and accurate. We consider two model inverse problems: (i) identification of hydraulic head in a Darcy flow model and (ii) electrical impedance tomography (EIT) with a complete electrode model (CEM).  In addition, we display the capabilities of the proposed method to solve shape identification problems where the computation of the shape derivative of the forward map may be cumbersome. In concrete, we combine the proposed ensemble method with the level-set approach of \cite{Level_set_US} to estimate shapes whose boundaries determine regions of sharp discontinuities of parameters in the PDE models under consideration. Since our methodology does not require derivatives, our level-set based  formulation is considerably more simple than standard approaches where shape derivatives are needed \cite{BurgerSurvey, BurgerGB,Santosa}. Moreover, the proposed ensemble level-set approach does not use the level-set equation as in standard formulations; our iterative scheme induces a stable evolution of the unknown interface. In addition, by initializing the ensemble according to the ideas recently proposed in \cite{Level_set_US} we avoid computational issues such as the flattening of the level-set function which is often observed with standard methods and typically addressed with ad-hoc techniques. The third contribution of the present work is to provide an extensive numerical investigation in order to assess the convergence and regularization properties of the proposed ensemble scheme with respect to ensemble size, number of measurements, initial ensemble and tunable parameters. Although the convergence theory of the proposed regularizing ensemble Kalman method is an open problem, our numerical study offers deep insight into the potential application of iterative regularization for the development of derivative-free ensemble methods thereby opening up a whole new field for application and theory.

The rest of the manuscript is organized as follows. In \Sref{Sec:PDE} we introduce the general framework for the identification problem. The proposed ensemble Kalman method is introduced in subsection \ref{subsec:EnKF}. The test models used for the validation of the proposed scheme are introduced in subsection \ref{test}. Preliminary numerical examples that show the regularizing properties of the method are presented in subsection \ref{unreg}. In Section \ref{sec:EnKF} we discuss general aspects and properties of the proposed scheme including the derivation of the method as an ensemble approximation of the regularizing LM scheme of \cite{Hanke}. In Section \ref{Numerics} we then conduct an extensive numerical investigation of the proposed scheme. In particular, we study the relation between the ensemble size and the number of measurements in terms of their effect on the regularization and convergence properties of the scheme. We investigate the effect of the tunable parameters of the scheme and the regularization properties with respect to the ensemble size in the small noise limit. In Section \ref{Sec:Applications} we discuss some potential applications of the scheme for the solution of geometrical inverse problems. Conclusions and future research are provided in \Sref{Conclusions}.

\section{PDE-constrained inverse problems}\label{Sec:PDE}

\subsection{General framework}\label{prel}

Let us denote by $G: X \to Y$ the (nonlinear) {\em forward operator} that arises from the PDE-constrained model under consideration. In other words, $G$ maps the space $X$ of PDE parameters (e.g. coefficients, source terms and/or boundary conditions) to the observation space $Y$ defined in terms of observable quantities related to the solution of the PDE (e.g. pointwise measurements of the solution). We assume that $X$ and $Y$ separable Hilbert spaces with norms denoted by $\vert\vert \cdot \vert\vert_{X}$ and $\vert\vert \cdot \vert\vert_{Y}$ respectively. For the applications under consideration, the dimension of the observation space is often small (i.e. order of $10^3$ - $10^4$). Therefore, we consider the case where $Y$ is finite dimensional. However, upon discretization, the dimensions of $X$ could be large (i.e. greater than $10^6$). Therefore, in order to derive algorithms robust under grid refinement we consider the case where $\dim(X)=\infty$; our aim is thus to keep our formulation within a functional analytical framework. 

The PDE-constrained model has an unknown parameter $u^{\dagger}\in X$ that we wish to identify from noisy measurements $y^{\eta}\in Y$ defined by
\begin{equation}
\label{eq:data}
y^{\eta}\equiv G(u^{\dagger})+ \xi
\end{equation}
where $\xi \in Y$ is noise which for the purpose of this exposition will be considered deterministic. In addition to $y^{\eta}$, we assume we have knowledge of the (weighted) noise level $\eta$ defined by
\begin{equation}\label{eq:nl}
 \eta \equiv    \vert\vert \Gamma^{-1/2}(y^{\eta} - G(u^{\dagger})) \vert\vert_{Y} 
\end{equation}
where $\Gamma:Y\to Y$ is a self-adjoint positive-definite operator which for the present analysis can be understood as a weighting/scaling operator that enable us, for example, to include information concerning the precision of our measurement device.

The identification problem is the following: \textit{Given $y^{\eta}$, find  $u\in X$ such that $G(u)=y^{\eta}$}. Since the data may not be in the range of the forward operator we may alternatively formulate the identification as the minimizer of 
\begin{eqnarray}\label{data_misfit}
\Phi(u)\equiv  \vert\vert \Gamma^{-1/2}(y^{\eta}-G(u))\vert\vert_{Y}
\end{eqnarray}
However, forward operators $G$ that arise from PDE models are typically compact and weakly sequentially closed \cite{Ill-posed}. From this property it follows that the inverse problem is ill-posed in the sense of stability. In other words, we may find $u_{n}$ such that $G(u_{n})\to y^{\eta}$ but $u_{n}\nrightarrow u^{\dagger}$. Therefore the computation of a minimizer of (\ref{data_misfit}) with standard (unregularized) iterative minimization approaches may be unstable. This lack of stability can be alleviated by means of iterative regularization methods. As discussed in \Sref{Intro}, iterative regularization provides a computational framework to compute stable approximation to the inverse problem.  In concrete, their aim is to compute an estimate/approximation $u^{\eta}$ controlled by the noise level $\eta$, that converges, in the small noise limit, to a solution of the identification problem, i.e. $u^{\eta}\to u$ as $\eta\to 0$, where the limit $u$ satisfies $G(u)=G(u^{\dagger})$. Most iterative regularization approaches require the computation of the Fr\'echet derivative of the forward operator $G$.  In some applications where commercial software is used for forward simulations Fr\'echet derivatives may not available and so the application of those iterative methods may be limited. The objective of the present work is to introduce a derivative-free regularization ensemble Kalman method for the stable computations of the identification problem.

\subsection{The regularizing ensemble Kalman method}\label{subsec:EnKF}

Let us assume that we are given an ensemble of $N_{e}$ elements $u_{0}^{(j)}$ ($j\in\{1,\dots,N_{e}\}$) of the parameter space $X$. One may think of $\{u_{0}^{(j)}\}_{j=1}^{N_{e}}$ as an ensemble of potential initial guesses for any iterative regularization method applied for the stable identification of the unknown parameter $u^{\dagger}$. We will discuss the selection of such initial ensemble $\{u_{0}^{(j)}\}_{j=1}^{N_{e}}$ in subsection \ref{initial}.

We now propose an iterative scheme where, at every iteration level $n$, each ensemble member $u_{n}^{(j)}$ is updated in such a way that the corresponding ensemble mean
\begin{eqnarray}\label{eq:m7A}
\overline{u}_{n}\equiv \frac{1}{N_{e}}\sum_{j=1}^{N_{e}}u_{n}^{(j)}
\end{eqnarray}
approximates the inverse problem in the small noise limit as described in subsection \ref{prel}. In other words, we need that algorithm stops at finite iteration level $n^{\star}$ producing an estimate $u^{\eta}\equiv \overline{u}_{n}^{\star}$ such that $u^{\eta}\to u$ as $\eta\to 0$, where $G(u)=G(u^{\dagger})$ . The proposed regularizing ensemble Kalman method is presented below.

\begin{algorithm}{Iterative regularizing ensemble Kalman method}\label{Al1}\\
Let $\{u_{0}^{(j)}\}_{j=1}^{N_{e}}\subset X$ be the initial ensemble of $N_{e}$ elements. Let $\rho\in (0,1)$ and $\tau>1/\rho$.\\
For $n=0,1,\dots$
\begin{itemize}
\item[(1)] \textbf{Prediction step.} Evaluate  
\begin{eqnarray}\label{eq:m14}
w_{n}^{(j)}= G(u_{n}^{(j)}),\qquad j\in\{1,\dots, N_{e}\}
\end{eqnarray}
and define $\overline{w}_{n}=\frac{1}{N_{e}}\sum_{j=1}^{N_{e}}w_{n}^{(j)}$
\item[(2)] \textbf{Discrepancy principle}. If
\begin{eqnarray}\label{eq:m15}
 \vert\vert \Gamma^{-1/2}(y^{\eta}-\overline{w}_{n})\vert\vert_{Y} \leq \tau \eta 
\end{eqnarray}
stop. Output $\overline{u}_{n} \equiv \frac{1}{N_{e}}\sum_{j=1}^{N_{e}} u_{n}^{(j)}$.\\
\item[(3)]\textbf{Analysis step.}  Define $C_{n}^{uw}$,  $C_{n}^{ww}$ by 
\begin{eqnarray}\label{eq:m8b}
C_{n}^{ww}(\cdot)=\frac{1}{N_{e}-1}\sum_{j=1}^{N_{e}}(G(u_{n}^{(j)})-\overline{w}_{n})\langle G(u_{n}^{(j)})-\overline{w}_{n}),\cdot\rangle_{Y}\\
C_{n}^{uw}(\cdot) = \frac{1}{N_{e}-1}\sum_{j=1}^{N_{e}} (u_{n}^{(j)}-\overline{u}_{n})\langle G(u_{n}^{(j)})-\overline{w}_{n}),\cdot \rangle_{Y}.\label{eq:m8c}
\end{eqnarray}
Update each ensemble member:
\begin{eqnarray}\label{eq:m16}
\fl u_{n+1}^{(j)} =u_{n}^{(j)}+C_{n}^{uw}(C_{n}^{ww} +\alpha_{n}\Gamma   )^{-1}(y^{\eta}-w_{n}^{(j)}),\qquad j\in\{1,\dots, N_{e}\}
\end{eqnarray}
where $\alpha_{n}$ is chosen by the following sequence
\begin{equation}\label{eq:m17}
\alpha_{n}^{i+1}=2^{i}\alpha_{n}^{0}.
\end{equation}
where $\alpha_{n}^{0}$ is an initial guess. We then define $\alpha_{n}\equiv \alpha_{n}^{N}$ where $N$ is the first integer such that
$$\alpha_{n}^{N}\vert\vert \Gamma^{1/2}(C_{n}^{ww} +\alpha_{n}^{N}\Gamma   )^{-1}(y^{\eta}-\overline{w}_{n})\vert\vert\ge \rho\vert\vert \Gamma^{-1/2}(y^{\eta}-\overline{w}_{n})\vert\vert$$
\end{itemize}
\end{algorithm}

Note that, in contrast to other regularization techniques, the discrepancy principle in (\ref{eq:m15}) is applied to $\overline{w}_{n}$ which (see expression (\ref{eq:m7})) is the average of the model output $G(u^{(j)})$ of each ensemble member. While this quantity approximates $G(\overline{u}_{n})$ to first order (see subsection \ref{deriva}), it is clearly not the data misfit of the proposed estimate.

The proof of convergence of Algorithm \ref{Al1} to a stable solution of the inverse problem is beyond the scope of the present manuscript; our aim is to offer numerical evidence of the convergence and regularization properties of the scheme.

\subsection{Test Models}\label{test}

The proposed ensemble Kalman method will be tested on two PDE-constrained inverse problems. We consider small test inverse problems where we have the computational flexibility to conduct a large amount of numerical experiments in order to understand the effect that the tunable parameters, ensemble size, selection of initial ensemble and number of measurements have on the regularization properties and accuracy of the proposed method. We introduce the test models under consideration below.

\subsubsection{Test model I. Darcy flow.}\label{Darcy}

We consider single-phase steady-state Darcy flow in a two-dimensional confined aquifer whose physical domain is $D=[0,6]\times [0,6]$. $\kappa$ denotes the hydraulic conductivity. The flow is described in term of the piezometric head $h(x)$ ($x\in D$) given by the solution to \cite{Bear}
\begin{eqnarray}\label{eq:a1}
-\nabla\cdot \kappa \nabla h&=f &\qquad\textrm{in}~~D
\end{eqnarray}
where $f$ is the source which for the present work is defined by 
\begin{eqnarray}\label{eq:a2}
f(x_{1},x_{2})=\left\{\begin{array}{ccc}
0 &\textrm{if}& 0< x_{2}\leq 4,\\
137&\textrm{if}& 4< x_{2}< 5,\\
274&\textrm{if}& 5\leq x_{2} < 6.\end{array}\right.
\end{eqnarray}
The following boundary conditions are considered
\begin{eqnarray}\label{eq:a3}
\fl h(x,0)=100, \qquad \frac{\partial h}{\partial x}(6,y)=0,\qquad
-\kappa\frac{\partial h}{\partial x}(0,y)=500,\qquad   \frac{\partial h}{\partial y}(x,6)=0,
\end{eqnarray}
We are interested in recovering the logarithm of the hydraulic conductivity $u\equiv \log{\kappa}$, from noisy pointwise measurements of the piezometric head $h$. In other words, we consider 
\begin{eqnarray}\label{eq:a4}
G(u) = (h(x_{1}),\dots, h(x_{M}))
\end{eqnarray}
where $h$ is the solution to (\ref{eq:a1})-(\ref{eq:a3}) and $\{x_{j}\}_{j=1}^{M}$ are the measurement locations. This groundwater model was used first used as benchmark for inverse modeling in \cite{Carrera}. It has been also used as a test model for the identification of parameters with iterative regularization methods in \cite{Hanke,Repre} and with an ensemble Kalman approach in \cite{EnKF_US}. 

\subsubsection{Test model II. Complete Electrode Model.}\label{EIT}

Our second test model is based on the Complete Electrode Model (CEM) for Electrical Impedance Tomography (EIT). The objective of EIT is to identify the conductivity $\kappa$ of a body $D$ given measurements of voltages from a configuration of $m_{e}$ electrodes on $\{e_{k}\}_{k=1}^{m_{e}}$ places on the boundary $\partial D$. The measured voltage arise from current patterns applied on those electrodes.  The forward model associated to EIT is the CEM which consist of computing the voltage $v$ in $D$ and the voltages $\{V_{k}\}_{k=1}^{m_{e}}$ on  $\{e_{k}\}_{k=1}^{m_{e}}$ that satisfy
\begin{eqnarray}
&~~~~~~~\nabla \cdot \kappa \nabla v&=0\qquad\textrm{in}~~D,\\
&v+z_{k} \kappa \nabla v\cdot \n&=V_{k} \qquad\textrm{on}~~e_k, ~~k=1,\dots,m_{e},\\
&~~~~~~~~~\nabla v\cdot \n &=0 \qquad\textrm{on}~~\partial D\setminus \cup_{k=1}^{m_{e}}e_{k},\\
&\int_{e_{k}}\kappa \nabla v\cdot \n ~ds &= I_{k}\qquad ~~k=1,\dots,m_{e},
\end{eqnarray}
where $\{I_{k}\}_{k=1}^{m_{e}}$ are the currents injected through the electrodes, $\{z_{k}\}_{k=1}^{m_{e}}$ are the contact impedances of the electrodes and $u=\log{\kappa}$ where $\kappa$ is the conductivity. Well posedness of the CEM model requires conservation of charge
\[
\sum_{k=1}^{m_{e}}I_{k}=0
\]
Given $\kappa$ and  $\{z_{k}\}_{k=1}^{m_{e}}$, for each current patter $I=\{I_{k}\}_{k=1}^{m_{e}}$ there exists a unique solution $[v,\{V_{k}\}_{k=1}^{m_{e}}]$ to the CEM \cite{cheney}. The EIT problem consists now of finding $\kappa$ and $\{z_{k}\}_{k=1}^{m_{e}}$ from a set of $n_{p}$ measurements of voltages $V_1=\{V_{1,k}\}_{k=1}^{m_{e}}, \dots , V_{n_p}=\{V_{n_p,k}\}_{k=1}^{m_{e}}$ obtained from $n_{p}$ current paters $I_{1}=\{I_{1,k}\}_{k=1}^{m_{e}},\dots I_{n_{p}}=\{I_{n_p,k}\}_{k=1}^{m_{e}}$. For simplicity we assume that the contact impedances of the electrodes are known. Therefore, the identification problem is to find $\kappa$ given 
\[
G(u)=\big[\{V_{1,k}\}_{k=1}^{m_{e}}, \dots , \{V_{n_p,k}\}_{k=1}^{m_{e}}\big]
\]
For a review of the EIT problem we refer the reader to \cite{EIT_revew}.

\subsection{Regularizing properties of the ensemble Kalman method.}\label{unreg}

Although ensemble Kalman methods have been typically used in the statistical Bayesian framework, several publications \cite{EnKF_US,IterativeEnKF} have explored the use of these approaches for solving (deterministic) inverse problems such as the parameter identification PDE-constrained problems described earlier.  However, most of these approaches are based on update formulas of the form of (\ref{eq:m16}) but with a fixed parameter $\alpha_{n}=1$, i.e. 
\begin{eqnarray}\label{eq:11d}
u_{n+1}^{(j)} =u_{n}^{(j)}+C_{n}^{uw}(C_{n}^{ww} +\Gamma   )^{-1}(y-w_{n}^{(j)})
\end{eqnarray}
 As we discuss in subsection \ref{bayesian}, the standard choice $\alpha_{n}=1$ is motivated by the application of Kalman methods for solving Bayesian inference problems when the model $G$ is linear, and the underlying prior distribution is Gaussian \cite{Tarantola}. For nonlinear forward models, however, the same choice of $\alpha_{n}$ leads to instabilities that are often fixed with ad-hoc methods such as covariance localization. The main objective of the present work is to show numerically that such instabilities can be addressed by the proposed ensemble Kalman method when the ensemble size $N_e$ is sufficiently large.

\begin{figure}[htbp]
\begin{center}
\includegraphics[scale=0.245]{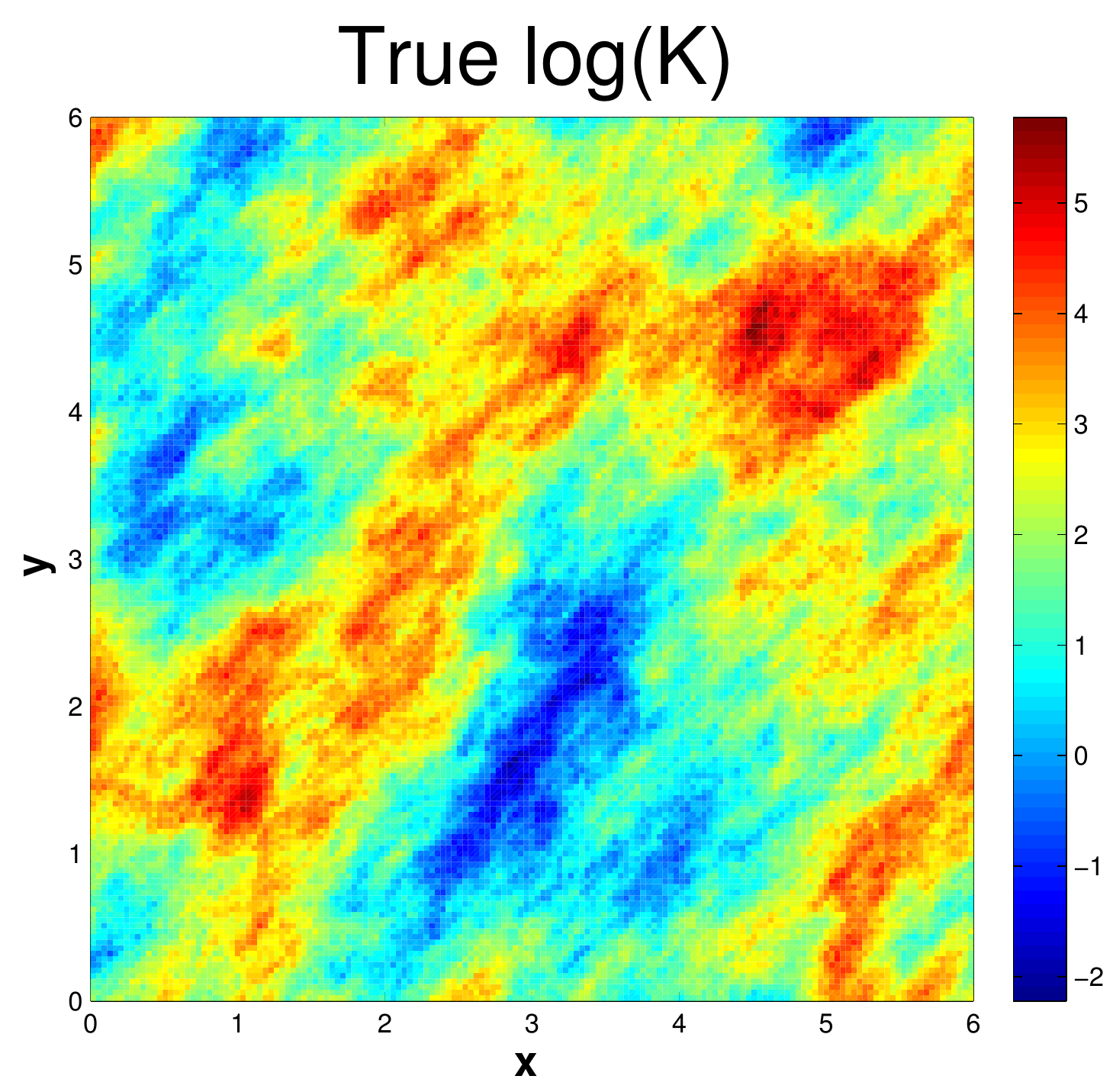}\\
\includegraphics[scale=0.245]{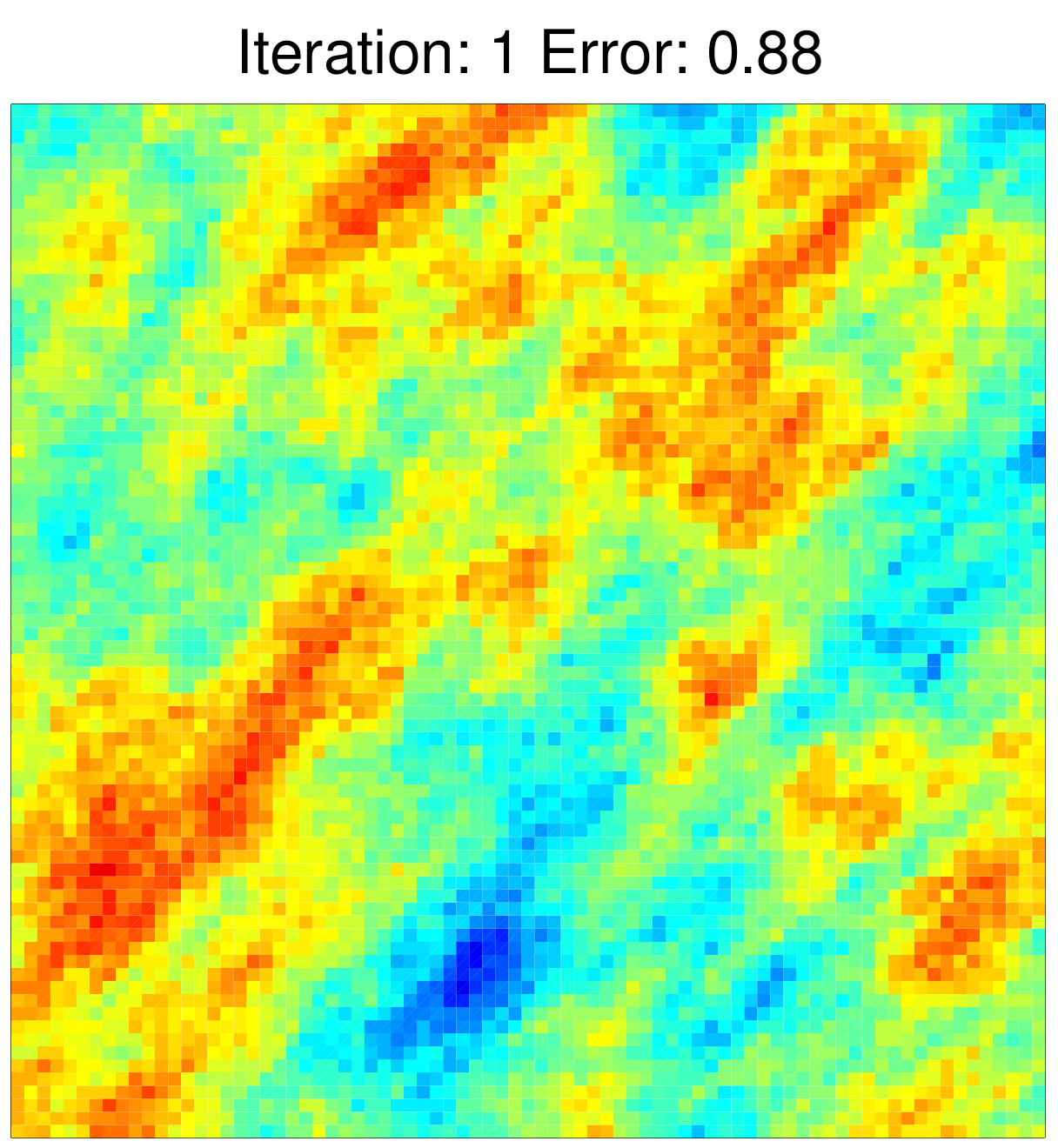}
\includegraphics[scale=0.245]{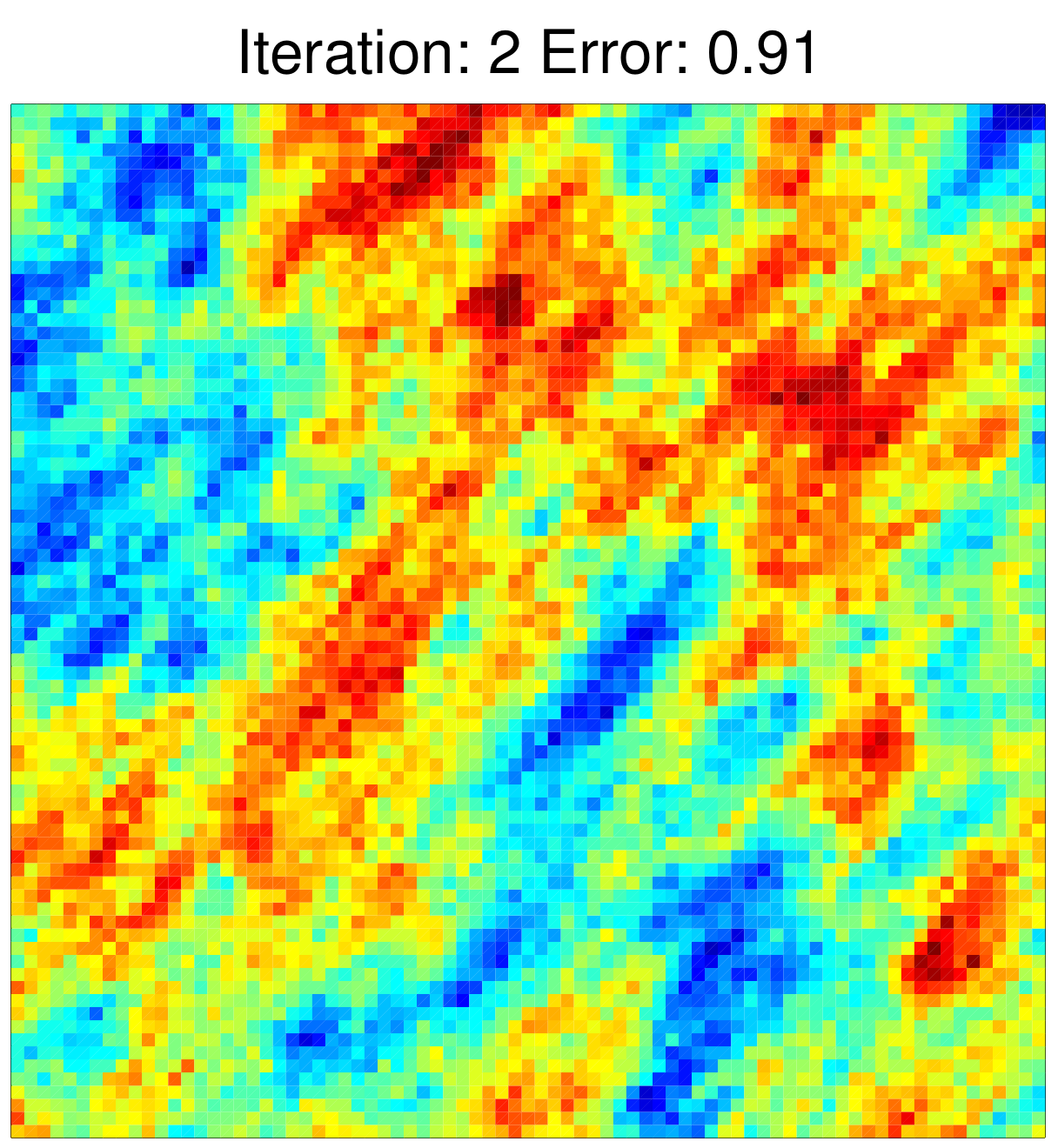}
\includegraphics[scale=0.245]{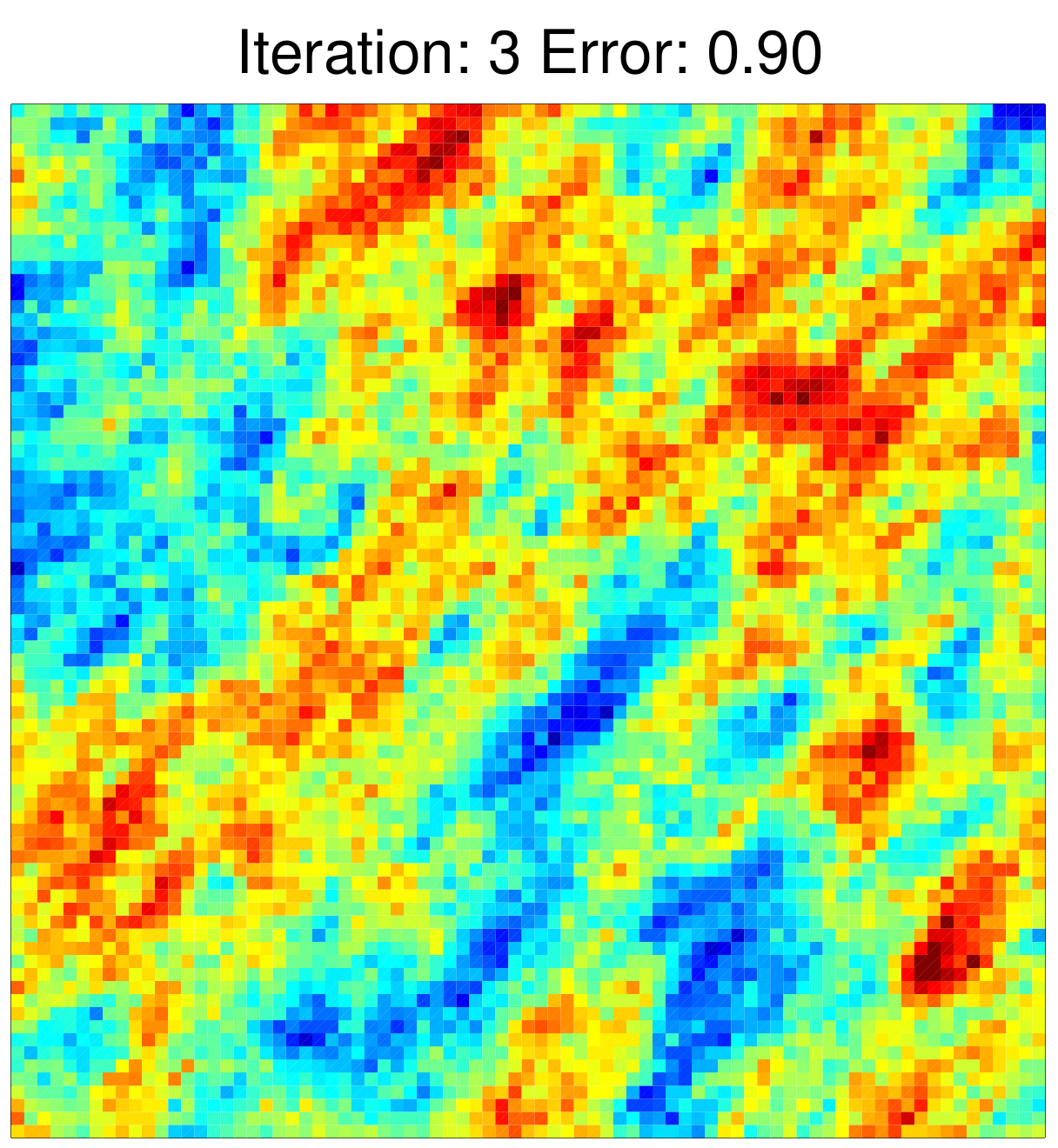}
\includegraphics[scale=0.245]{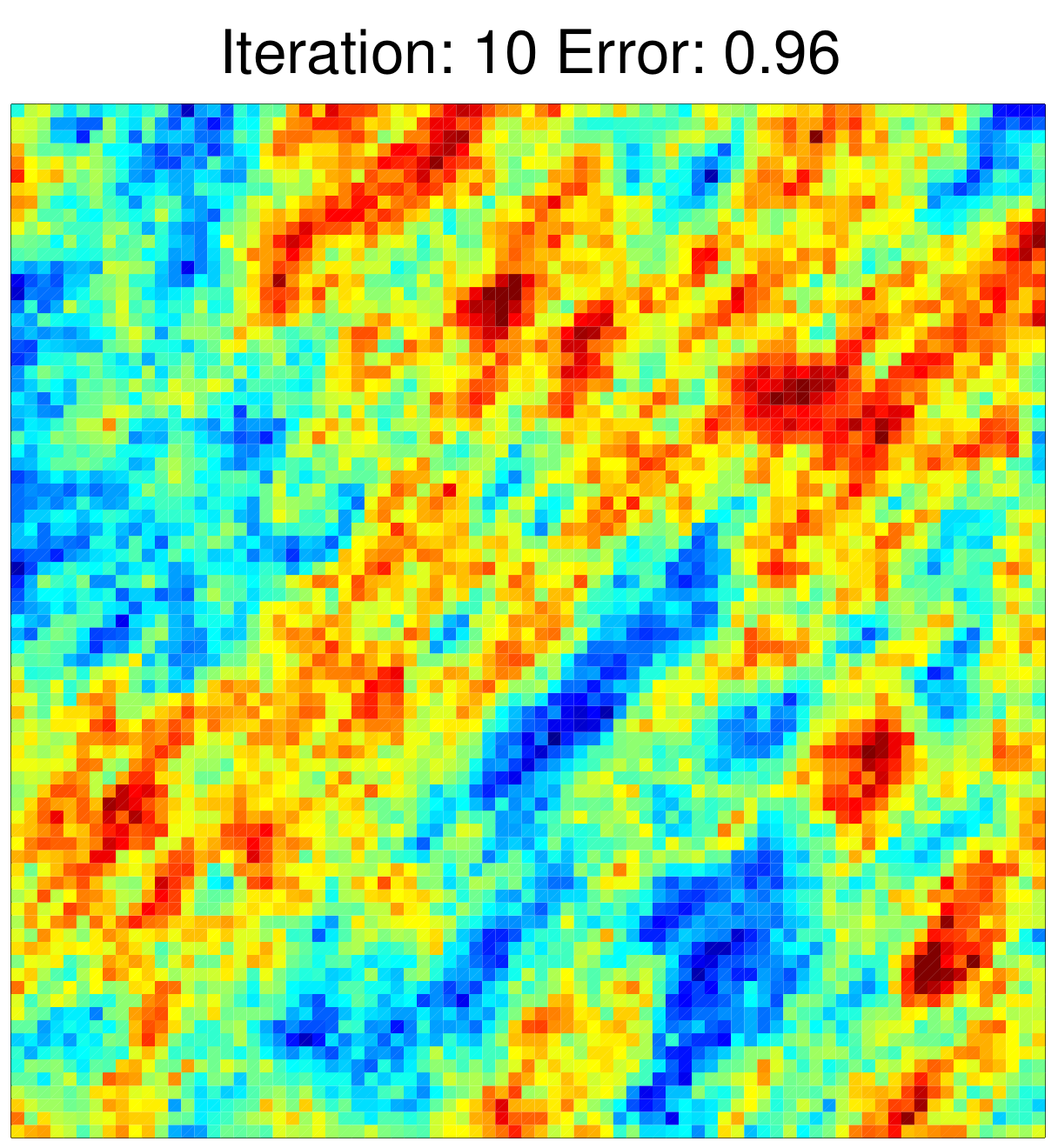}
\includegraphics[scale=0.245]{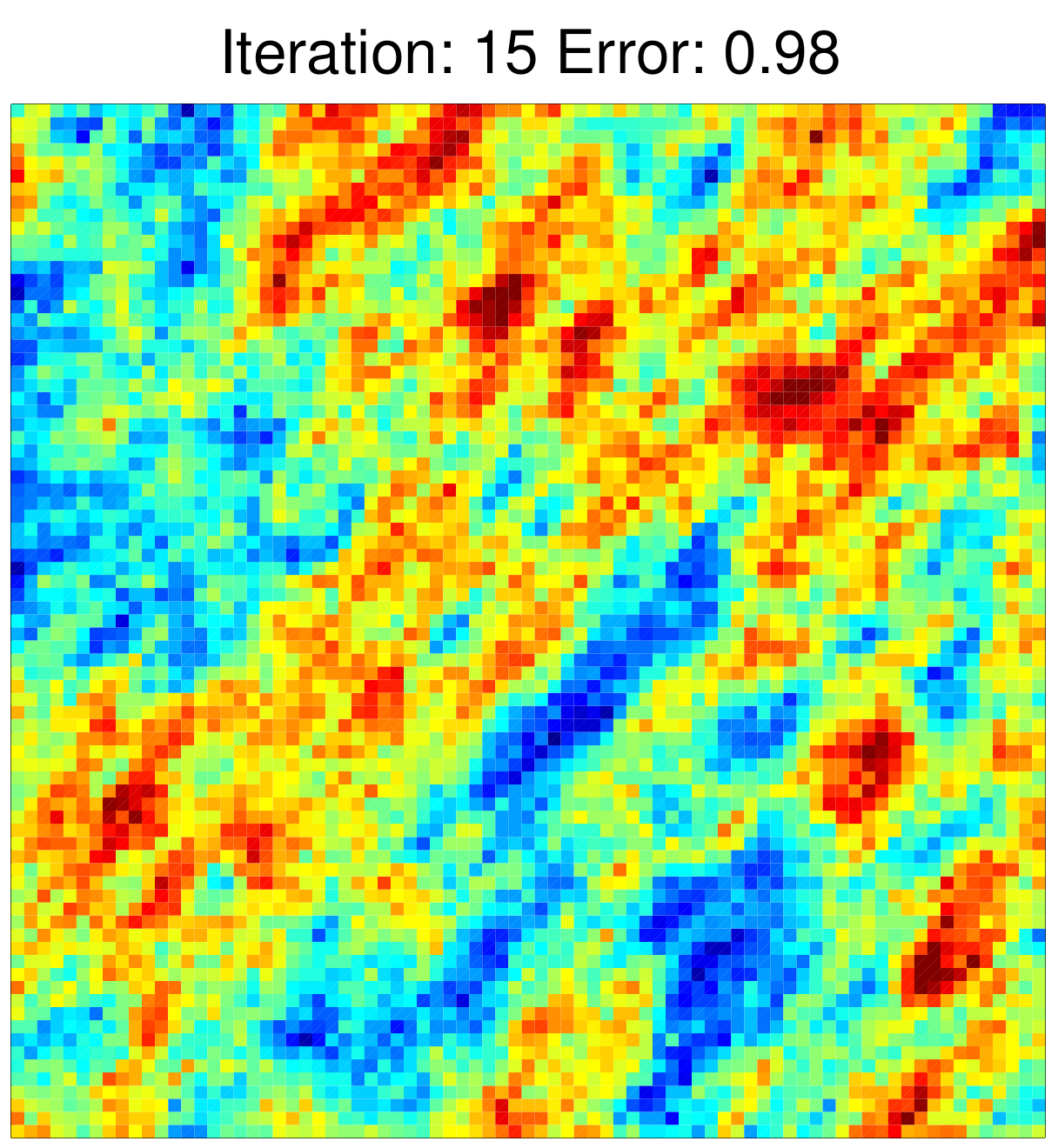}\\
\includegraphics[scale=0.245]{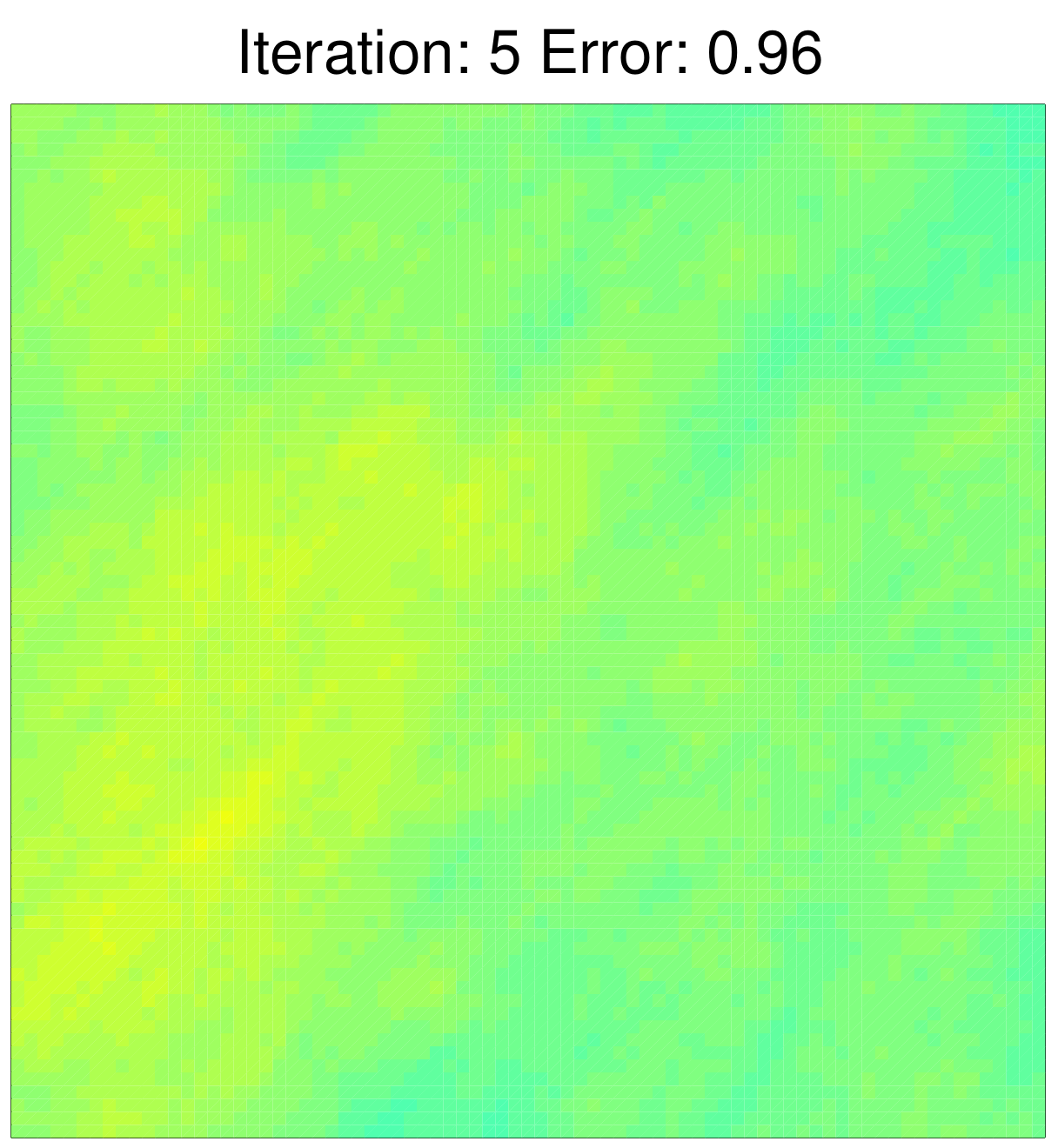}
\includegraphics[scale=0.245]{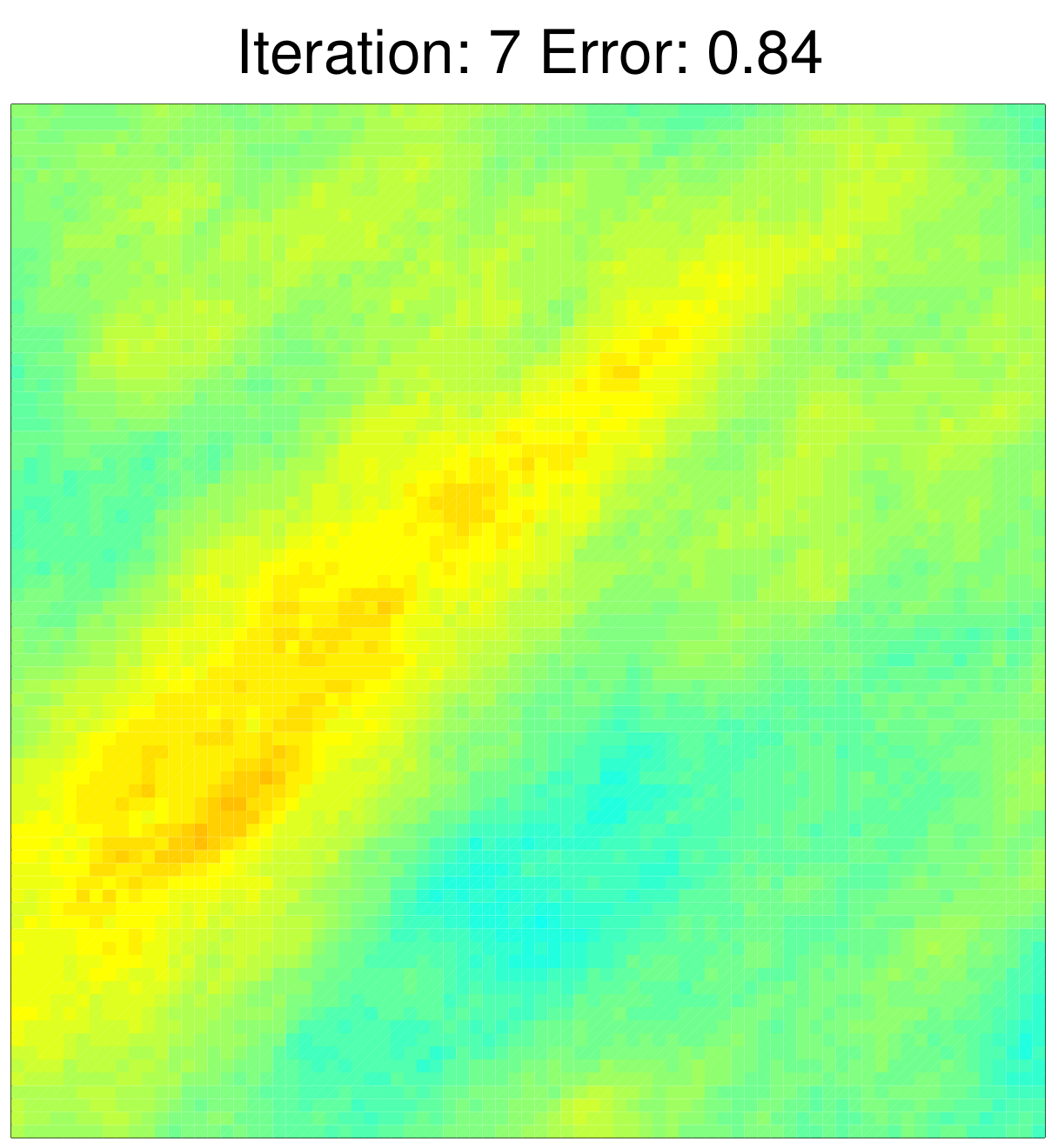}
\includegraphics[scale=0.245]{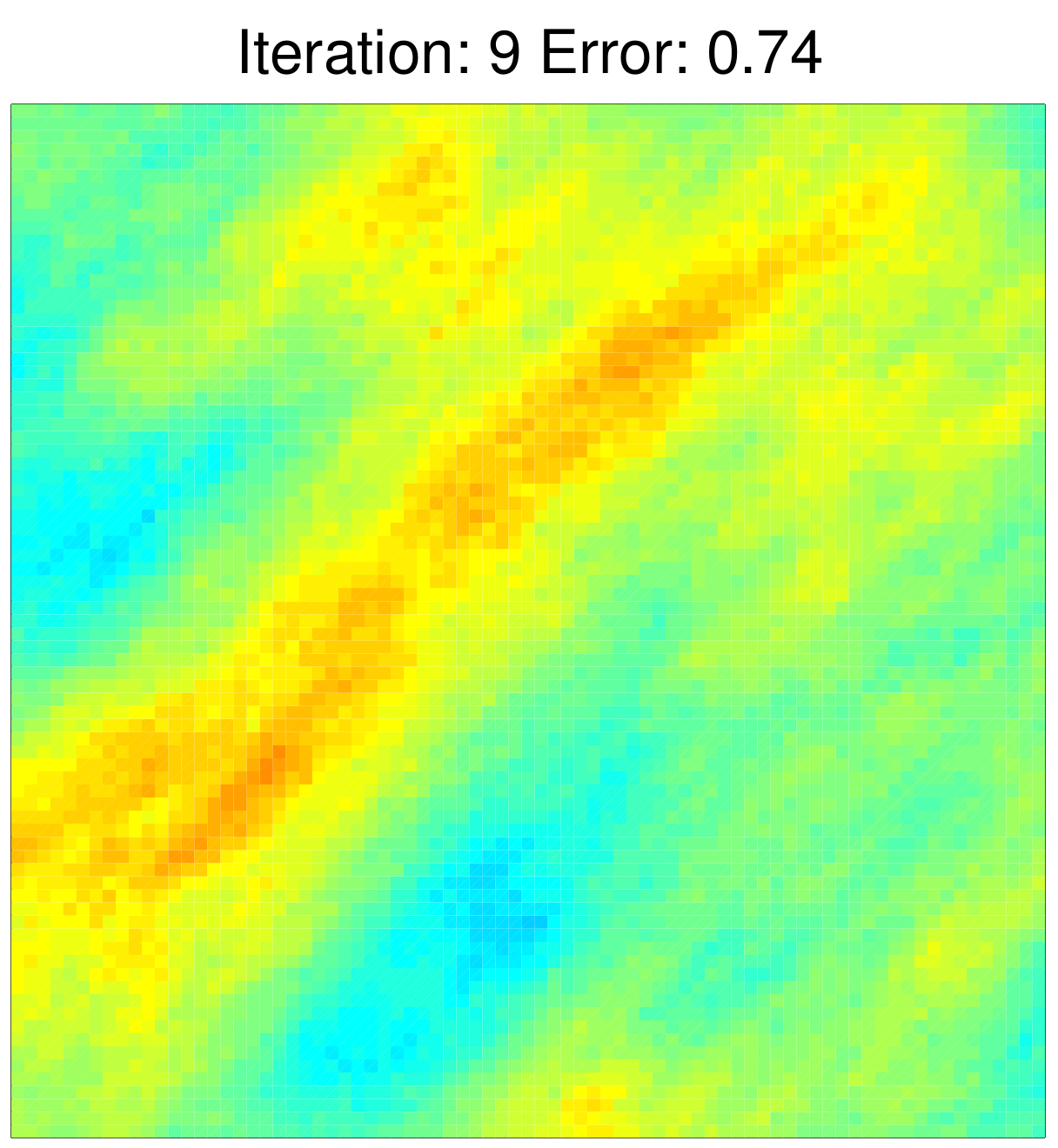}
\includegraphics[scale=0.245]{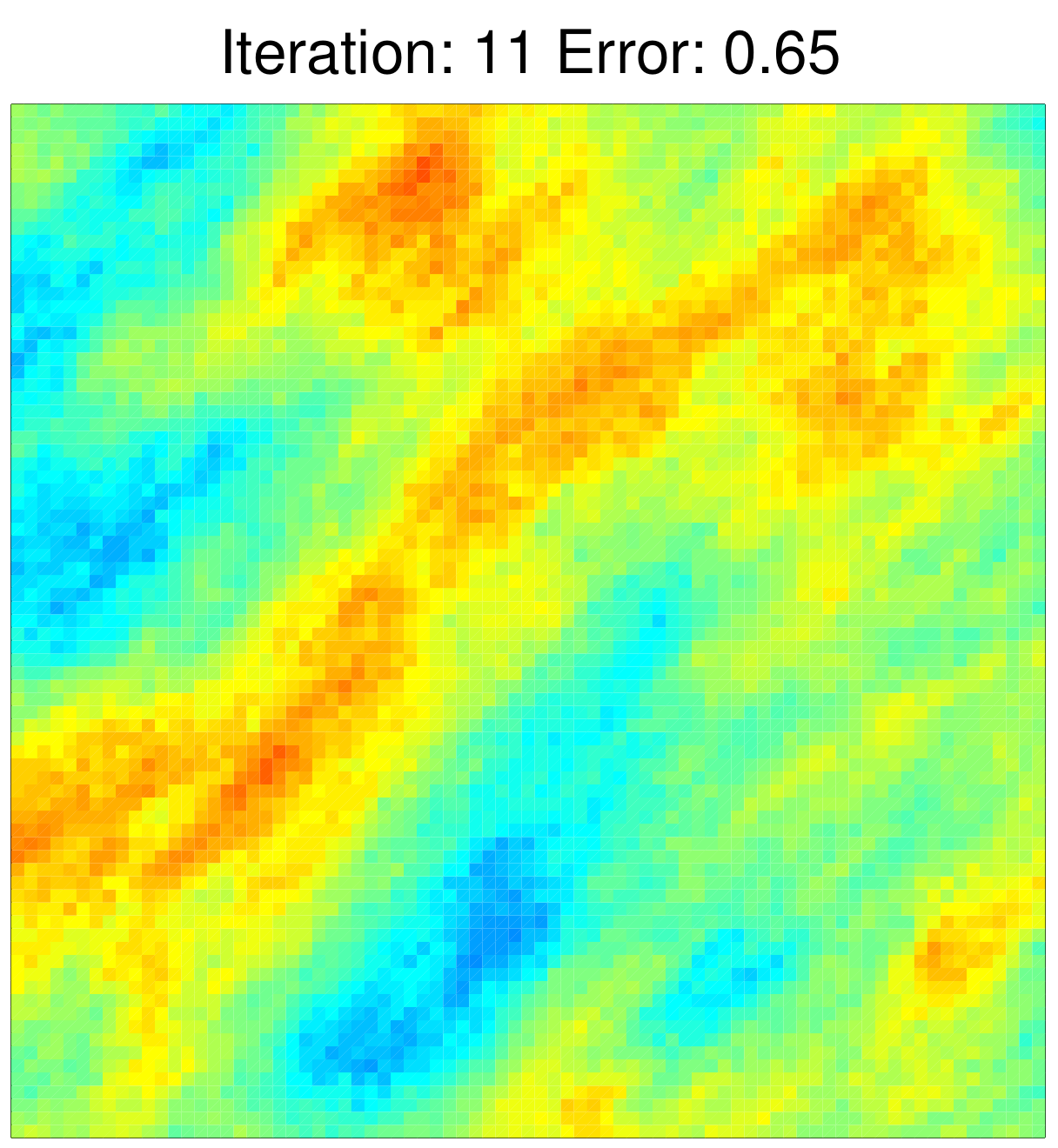}
\includegraphics[scale=0.245]{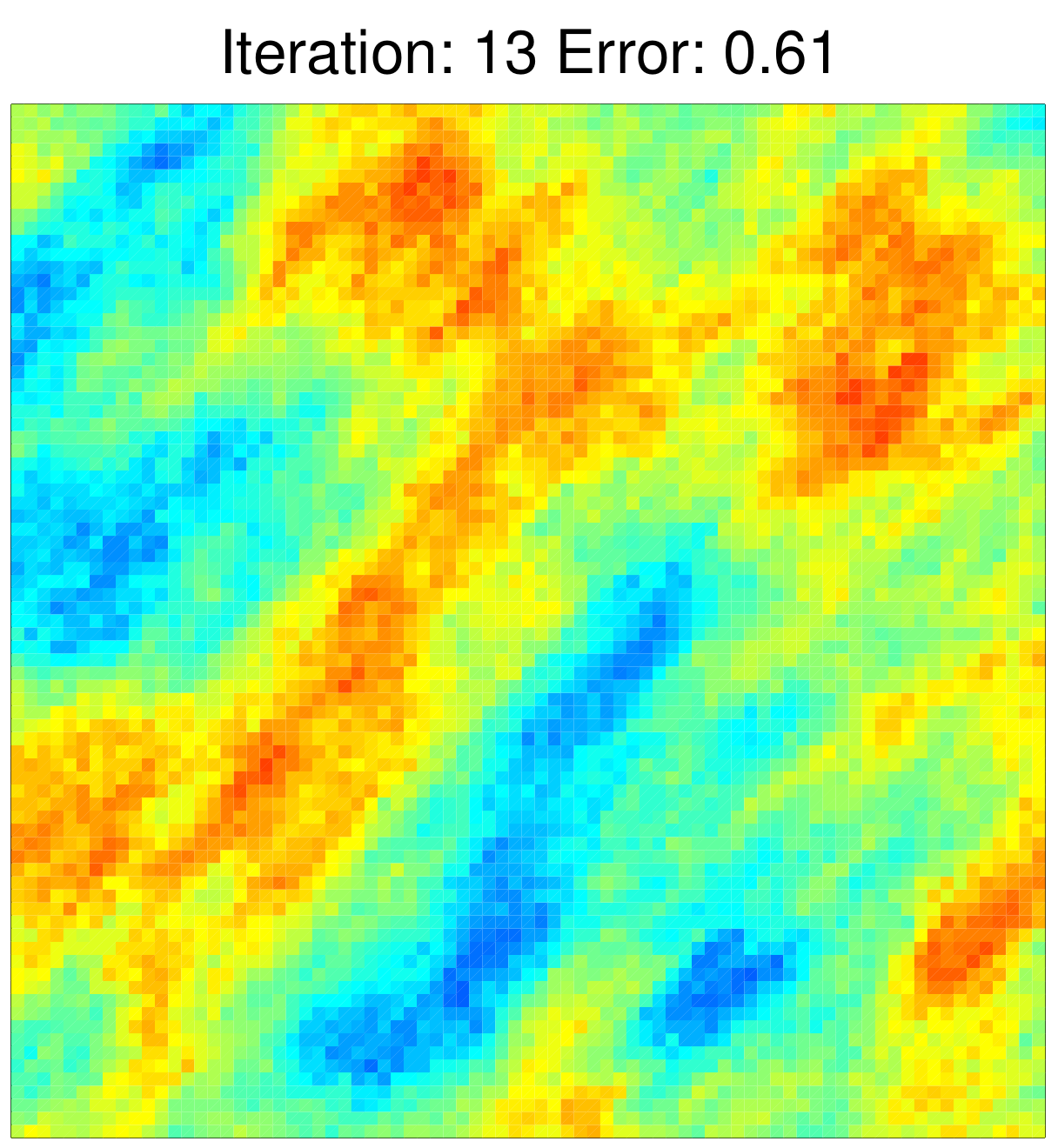}

 \caption{Top: true log-conductivity. Middle row: ensemble mean at some iterations of the standard (unregularized) ensemble Kalman method. Bottom row: ensemble mean at some iterations of the regularizing ensemble Kalman method. }
    \label{Fig10}
\end{center}
\end{figure}

\begin{figure}[htbp]
\begin{center}
\includegraphics[scale=0.325]{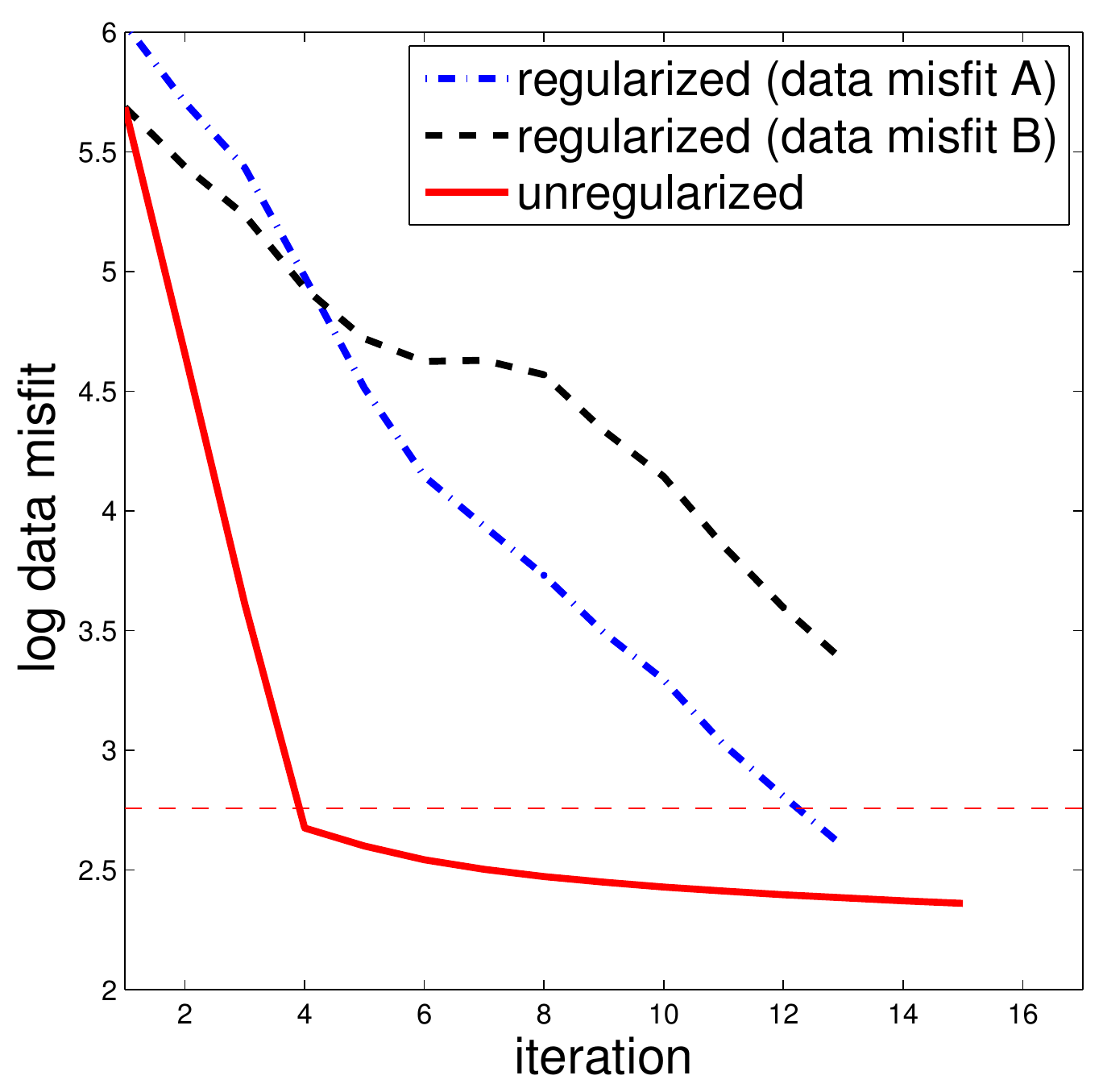}
\includegraphics[scale=0.325]{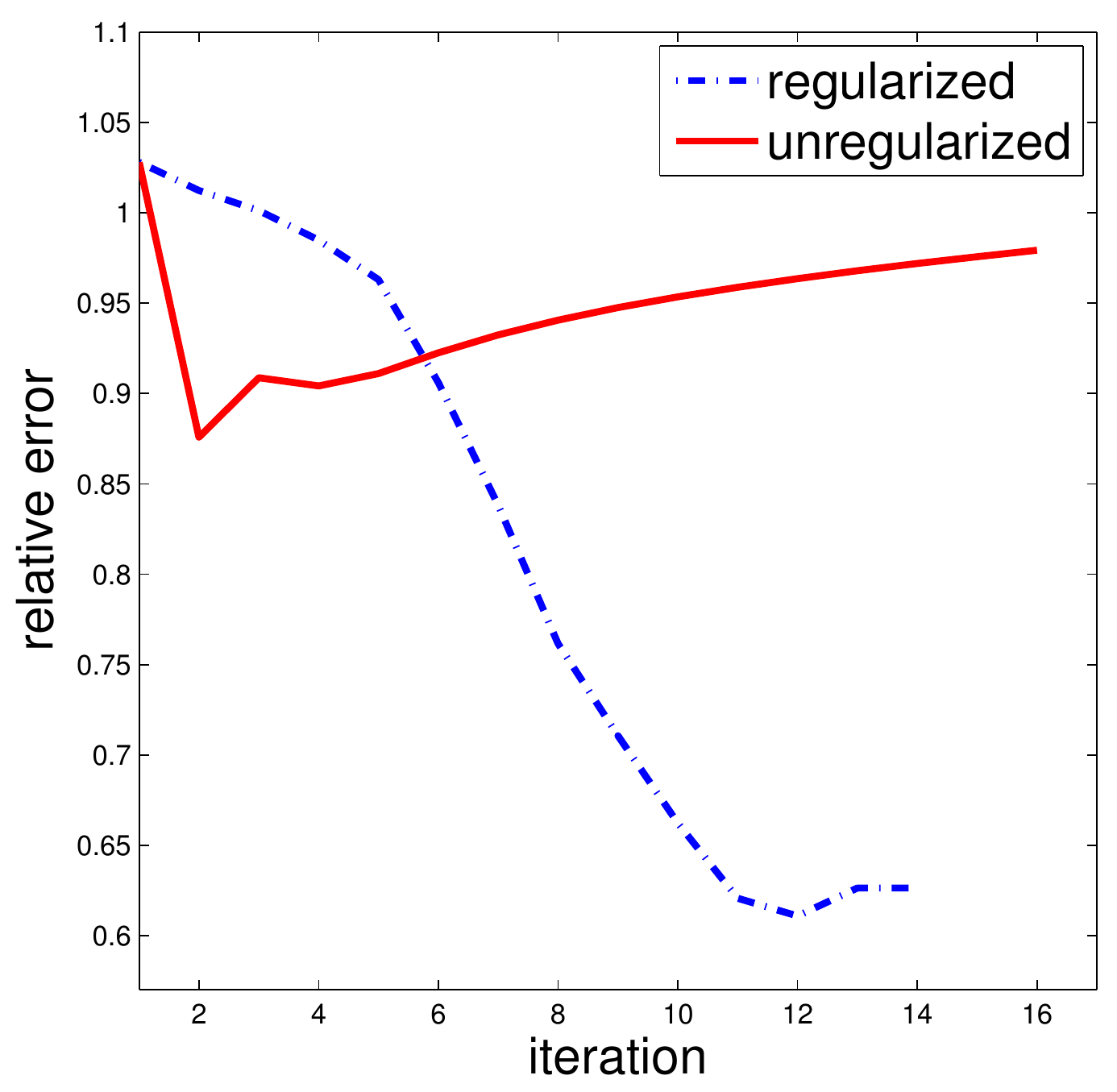}
\includegraphics[scale=0.325]{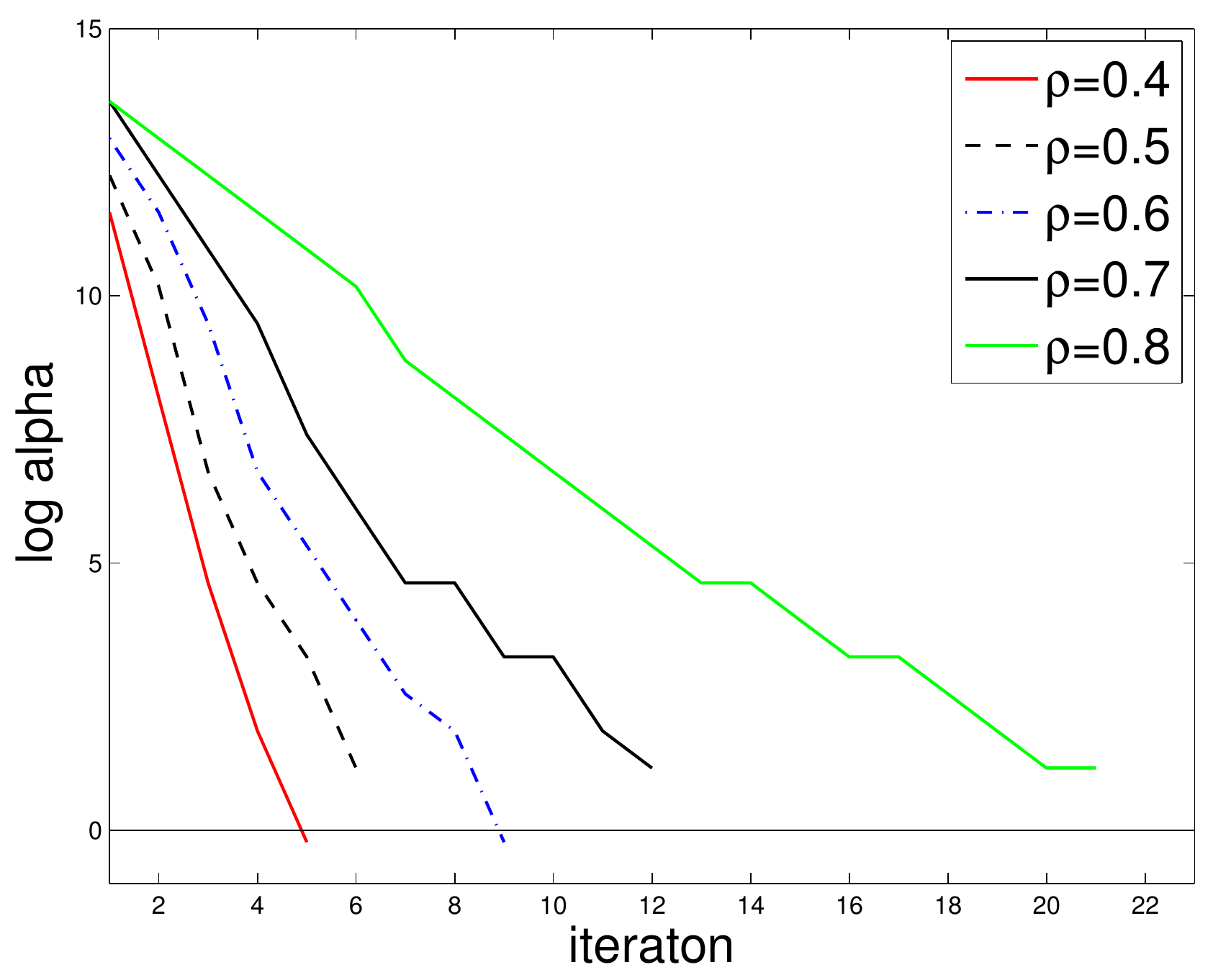}
 \caption{Log-data misfit (left) and relative error w.r.t truth (middle) of the ensemble mean at some iterations obtained with the unregularized ensemble Kalman method (solid red line) and with the proposed regularizing scheme (dotted blue line). Right: $\log (\alpha_{n})$ as a function of the number of iterations of Algorithm \ref{Al1} from 40 experiments with different initial ensembles of size $N_{e}=150$.}
   \label{Fig11}
\end{center}
\end{figure}

In this subsection we briefly present numerical evidence of the regularization properties of the proposed scheme; a detailed numerical investigation will be presented in \Sref{Numerics}. We apply Algorithm \ref{Al1} to the identification of the ``true'' log hydraulic conductivity $u^{\dagger}$ displayed in Figure \ref{Fig10}  (top). We use use synthetic measurements of hydraulic head  from the Darcy model of subsection \ref{Darcy} with the measurement locations displayed in  \Fref{Fig1} (left). Both the truth and the elements of the initial ensemble are generated from a Gaussian distribution. For details on the generation of synthetic data (avoiding inverse crimes) and the generation of initial ensemble,  we refer the reader to subsection \ref{Darcy_Num}. In the middle row (resp. bottom row) of Figure \ref{Fig10}  we display the estimate obtained from the ensemble mean at different iterations computed with the standard unregularized approach (resp. the proposed regularized method). For the regularized method we select $\rho=0.7$ and $\tau=1/\rho$. For the unregularized method we simply apply (\ref{eq:11d}) with no stopping criterion. Both methods are applied with the same (fixed) initial ensemble of $N_{e}=150$ members generated as described  in subsection \ref{Darcy_Num}. \Fref{Fig11} (left) shows the data misfit (\ref{data_misfit}) for the ensemble mean $\overline{u}_{n}$ computed with the standard unregularized approach (solid red line) and our regularized method (dotted-black line). Additionally, in \Fref{Fig11} (left) we also display (the dotted-blue line ) the data misfit with respect to the averaged ensemble data predictions (expression (\ref{eq:m15})) which we monitor for the termination of Algorithm \ref{Al1}. 


In \Fref{Fig11} (middle) we show the relative $L^2$-error with respect to the truth for the ensemble mean $\overline{u}_{n}$, (i.e. $\vert\vert \overline{u}_{n}-u^{\dagger}\vert\vert_{L^2}/\vert\vert u^{\dagger}\vert\vert_{L^2}$), computed with the standard unregularized approach (solid red line) and our regularized method (dotted-blue line). For the standard method we see that a fast decrease of the data misfit and relative error are observed at the early iterations. While the data misfit keeps decreasing, the error starts to increase when apparently the data misfit drops below the noise level whose value is indicated with the horizontal line in \Fref{Fig11} (left). This increase in the error comes as no surprise since the corresponding estimates overfit the data. Iterative methods applied for the solution of ill-posed inverse problems often display such behavior \cite{Engl,Kirsch}. In fact, a similar semiconvergent  behavior was reported in the ensemble Kalman method of \cite{EnKF_US} and inspired the present work where not only a regularization of the estimates need to be introduced (here by means of the parameter $\alpha_{n}$) but also the proper early termination of the scheme that avoids fitting the noise. Indeed, we note that the regularizing algorithm (see \Fref{Fig11}), both relative error and data misfit decreases very slowly. Once the (averaged) data misfit (of expression (\ref{eq:m15})) is close to the noise level, the algorithm is stopped producing a considerably more accurate estimate of the truth than the one obtained with the unregularized algorithm which blows up the update at the early iterations. However, if the proposed regularized algorithm is not stopped, the error with respect to the truth could potentially increase as we observe in subsequent experiments (see \Fref{Fig7}). From the log-conductivity estimates ( \Fref{Fig10} middle and bottom rows), we see that the unregularized method produces uncontrolled estimates of the unknown while the proposed method results in small incremental transitions which eventually captures the truth more accurately. As we will discuss in Section \ref{Numerics}, the ability of the proposed method for regularizing the computations of the inverse problem depends on the size of the ensemble $N_e$.

The selection of $\alpha_{n}$ according to (\ref{eq:m12}) is crucial for the regularizing of the proposed method.  We monitor numerically these values of $\alpha_{n}$ obtained from the application of Algorithm \ref{Al1} to the identification problem described in the preceding paragraph. We apply the method for several choices of the parameter $\rho$ ($\rho=0.4, 0.5, 0.6, 0.7, 0.8$). In \Fref{Fig11} (right) we display the values of $\alpha_{n}$ that we obtained from the proposed selection (\ref{eq:m12}) as the number of iterations increases. At the early iterations $\alpha_{n}$ is large hence controlling the updates of the approximation. As the iteration progresses and the data misfit decreases, $\alpha_{n}$ also decreases. For larger $\rho$ the decay of $\alpha_{n}$ is slower; this provides more regularization to the expense of a more costly algorithm (see Remark \ref{rema1}). From \Fref{Fig11} (right) we can observe that when the data misfit has dropped to a value close to the noise level by which the algorithm will be terminated according to (\ref{eq:dis2}), $\alpha_{n}$ has dropped down to a value close to $\alpha_{n}=1$ which is, in turn, the intrinsic value for the standard Kalman method. Interestingly, there is indeed an apparent intrinsic parameter $\alpha_{n}=1$ once the scheme arrives at the optimal approximation of the inverse problem. However, choosing this parameter at early iterations can be substantially detrimental to the performance of the scheme. In the following section, the selection of $\alpha_{n}$ is derived by using ensemble approximations of the Fr\'echet derivative of the forward map in the selection of the parameter $\alpha_{n}$ that has been proven to provide stability in the regularizing LM scheme of Hanke \cite{Hanke}. As stated earlier, for the present methodology, such a theory is beyond the scope of the present work.

\section{Properties and computational aspects of the regularizing ensemble Kalman method}\label{sec:EnKF}

In this subsection we discuss some properties and general computational aspects of the proposed regularizing ensemble Kalman method presented in subsection \ref{subsec:EnKF}.

\subsection{The initial ensemble}\label{initial}

One of the main properties of the proposed ensemble Kalman method is that the estimate $\overline{u}_{n}$ obtained from Algorithm  \ref{Al1} lives in the subspace generated by the initial ensemble $\{u_{0}^{(j)}\}_{j=1}^{N_{e}}$. 
\begin{proposition}[Invariance Subspace Property]\label{Prp2}
At every iteration of the scheme, the ensemble mean $\overline{u}_{n}$ defined by (\ref{eq:m7A}) satisfies  $\overline{u}_{n}\in  \textrm{span} \{u_{0}^{(j)} \}_{j=1}^{N_{e}}$
\end{proposition}
\textbf{Proof:}  Once we rewrite Algorithm \ref{Al1} in the augmented version of subsection \ref{3_2}, the proof follows directly with the same argument as \cite[Theorem 2.1]{EnKF_US} $\Box$. 

From the invariance subspace property it follows that the selection of the initial ensemble is a design parameter crucial to the performance of the proposed scheme; we study this numerically in subsection \ref{num_EIT}. Prior knowledge of the space of admissible solutions $X$ can be used for such selection. For example, $\{u_{0}^{(j)}\}_{j=1}^{N_{e}}$ can be a truncated basis for $X$. Another example of prior knowledge that we may use for the construction of the initial ensemble is the regularity of the elements in $X$. More concretely, we may construct an ensemble by drawing its members from some probability distribution with the desired regularity. While the problem under consideration is deterministic, the use of a probability distribution is made for the sake of the generation of the initial ensemble with the regularity of the parameter space. The aforementioned invariance subspace property then ensures that the estimate produced by the proposed method inherits the regularity of the space of admissible solutions. Clearly, Gaussian distributions are desirable since sampling from them is relatively easy. Nonetheless, other priors such as Besov \cite{Besov} could also be considered. In subsections \ref{num_EIT} we provide examples where probability distributions are considered in order to generate an initial ensemble that we use with our computational approach for the EIT problem described in subsection \ref{EIT}. In general, the construction of the prior ensemble is based on prior knowledge of the problem under consideration.

\subsection{The Regularizing LM scheme}\label{LM}

In the following subsection we derive the proposed method as an approximation of the regularizing LM scheme \cite{Hanke}. For the subsequent derivation we consider $X$ completed with the norm $\vert\vert C^{-1/2} \cdot \vert\vert_{X}$  where $C^{-1}:D(C^{-1})\subset X\to X$ is a densely-defined unbounded self-adjoint operator with compact resolvent and $C^{-1/2}$ is defined in terms of the spectral decomposition of $C^{-1}$ (see \cite[Section 2.1]{EnKF_US}). For the purpose of this work, $C^{-1}$ is an operator selected a priori  that enforces regularity on the functional space $X$. Introducing an operator $C^{-1}$ in our formulation enable us to (i) derive the proposed method from the regularizing LM scheme of Hanke \cite{Hanke} and (ii) establish a connection between the classical (deterministic) and the Bayesian approach for inverse problems. In the Bayesian framework, $C^{-1}$ is the inverse of the covariance operator $C$ from a prior distribution \cite{Andrew}. It is important to remark that the definition of $C^{-1}$ does not appear in the proposed scheme Algorithm \ref{Al1}.

The regularizing LM scheme of \cite{Hanke} is an iterative method that possess the regularizing properties need for the stable computation of solution to the identification problem described in the preceding section. In concrete Hanke in \cite{Hanke} proposes an iterative scheme where the $u_{n+1}$ iteration level is given by 
\begin{eqnarray}\label{eq:m1}
\fl u_{n+1}=u_{n}+\arg\min_{v\in X}\vert\vert \Gamma^{-1/2}(y^{\eta}-G(u_{n})-DG(u_{n})v\vert\vert_{Y}^2+\alpha_{n}\vert\vert C^{-1/2} v\vert\vert_{X}^{2}
\end{eqnarray}
where $\alpha>0$ is a regularization parameter chosen below and $DG$ denotes the Fr\'echet derivative of $G$. From the first order optimality conditions associated to the right hand side of (\ref{eq:m1}) it follows that the previous expression is equivalent to
\begin{eqnarray}\label{eq:m2}
\fl u_{n+1} =u_{n}+(DG^{\ast}(u_{n})\,\Gamma^{-1}\,DG(u_{n}) +\alpha_{n} C^{-1}  )^{-1}C\,DG^{\ast}(u_{n})(y^{\eta}-G(u_{n})).
\end{eqnarray}
The convergence and regularizing properties of the method of \cite{Hanke} requires that $\alpha_{n}$ in (\ref{eq:m1}) satisfies
\begin{eqnarray}\label{eq:m3}
\fl \rho \vert\vert \Gamma^{-1/2}(y^{\eta}-G(u_{n}))\vert\vert_{Y} \leq \alpha_{n} \vert\vert \Gamma^{-1/2}( y^{\eta}-G(u_{n})-DG(u_{n})(u_{n+1} -u_{n}) )\vert\vert_{Y}
\end{eqnarray}
for $\rho\in (0,1)$ selected a priori and that the scheme is terminated whenever the $n$th iteration level satisfies
\begin{eqnarray}\label{eq:m4}
 \vert\vert \Gamma^{-1/2}(y^{\eta}-G(u_{n}))\vert\vert_{Y} \leq \tau \eta  <  \vert\vert \Gamma^{-1/2}(y^{\eta}-G(u_{n-1}))\vert\vert_{Y}
\end{eqnarray}
for some $\tau>0$ that satisfies $\tau>1/\rho$.  

The theory of Hanke provides assumptions on the forward operator under which the regularizing LM scheme terminates after a finite number of iterations and the corresponding approximation is a solution of the inverse problem in the small noise limit as described in subsection \ref{prel}. For full details of the theoretical framework for the regularizing LM scheme the reader is referred to the work of Hanke in \cite{Hanke}. An application of this method to inverse problems in the geosciences can be found in \cite{LM}. 

In order to derive the proposed ensemble method as an approximation of the LM scheme, we need the following lemma where we assume that $Y=\mathbb{R}^{M}$ (recall that $\textrm{dim}(Y)<\infty$) with the standard Euclidean inner product.  
\begin{lemma}\label{lema:Rep1}
Assume that for any $u\in X$, the linear functionals $D_{m}G(u):X\to \mathbb{R}$ ($m=1,\dots, M$) are linearly independent. Then, expression (\ref{eq:m2}) is equivalent to 
\begin{eqnarray}\label{eq:m6}
\fl u_{n+1} =u_{n}+C\,DG^{\ast}(u_{n})(DG(u_{n})\,C\,DG^{\ast}(u_{n}) +\alpha_{n}\Gamma   )^{-1}(y-G(u_{n}))
\end{eqnarray}
\end{lemma}
\textbf{Proof:} See appendix.\\

Note that (\ref{eq:m6}) can be used to rewrite (\ref{eq:m3}) as follows
\begin{eqnarray}\label{eq:m5}
\fl \rho \vert\vert \Gamma^{-1/2}(y^{\eta}-G(u_{n}))\vert\vert_{Y} \leq \alpha_{n} \vert\vert \Gamma^{1/2}(DG(u_{n})\,CDG^{\ast}(u_n)+\alpha_{n}\Gamma)^{-1}(y^{\eta}-G(u_{n}))\vert\vert_{Y}
\end{eqnarray}
Expressions (\ref{eq:m6})-(\ref{eq:m5}) are computationally more convenient since the operator inversion of $(DG(u_{n})\,C\,DG^{\ast}(u_{n}) +\alpha_{n}\Gamma   )$ is conducted on a finite dimensional space. Upon discretization the aforementioned inversion has often negligible cost since the number of observations is typically small. Note that when $X$ is finite dimensional the relation between (\ref{eq:m2}) and (\ref{eq:m6}) follows simply from matrix Lemmas which cannot be applied in the present case.

\subsection{The proposed method as an approximation of the regularizing LM scheme}\label{deriva}

Let us consider a first order approximation of $G(u_{n}^{(j)})$ around the ensemble mean, i.e.
$$G(u_{n}^{(j)})\approx G(\overline{u}_{n})+DG(\overline{u}_{n})(u_{n}^{(j)}-\overline{u}_{n}) $$
Then,
\begin{eqnarray}\label{eq:m7}
\fl  \overline{w}_{n}\equiv \frac{1}{N_{e}}\sum_{j=1}^{N_{e}}G(u_{n}^{(j)})\approx G(\overline{u}_{n})\quad \textrm{and}\quad
\langle DG(\overline{u}_{n})(u_{n}^{(j)}-\overline{u}_{n}),w \rangle_{Y} \approx \langle G(u_{n}^{(j)})- G(\overline{u}_{n}),w\rangle_{Y}.\nonumber\\
\end{eqnarray}
Define the following covariance operator
\begin{eqnarray}\label{eq:m8a}
 C_{n}^{uu}(\cdot) = \frac{1}{N_{e}-1}\sum_{j=1}^{N_{e}} (u_{n}^{(j)}-\overline{u}_{n})\langle u_{n}^{(j)}-\overline{u}_{n},\cdot\rangle_{X}
\end{eqnarray}
From (\ref{eq:m7})-(\ref{eq:m8a}) and (\ref{eq:m8b}) -(\ref{eq:m8c}) we find that
\begin{eqnarray}\label{eq:m9}
\fl C_{n}^{uu}DG^{*}(\overline{u}_{n})v=  \frac{1}{N_{e}-1}\sum_{j=1}^{N_{e}}(u_{n}^{(j)}-\overline{u}_{n})\langle(u_{n}^{(j)}-\overline{u}_{n}),DG^{*}(\overline{u}_{n})v\rangle_{X}\nonumber\\
\approx \frac{1}{N_{e}-1}\sum_{j=1}^{N_{e}}(u_{n}^{(j)}-\overline{u}_{n})\langle G(u_{n}^{(j)})- G(\overline{u}_{n}),v\rangle =C_{n}^{uw}v
\end{eqnarray}
and from similar arguments we obtain
\begin{eqnarray}\label{eq:m10}
 DG(\overline{u}_{n})C_{n}^{uu}DG(\overline{u}_{n})^{*}v \approx C_{n}^{ww}v
\end{eqnarray}
In expression (\ref{eq:m6}) we now replace the following terms and the corresponding approximations 
\begin{eqnarray}\label{eq:m11}
 ~~~~~~~~~~~ ~~~~~~~~~~~u_{n} &\Longrightarrow \overline{u}_{n},\\
 ~~~~~~~~~~~C\,DG^{*}(u_{n})&\Longrightarrow C_{n}^{uu}DG^{*}(\overline{u}_{n})\approx C_{n}^{uw}\label{eq:m11bb}\\
 DG(u_{n})\,C\,DG^{*}(u_{n})&\Longrightarrow DG(\overline{u}_{n})C_{n}^{uu}DG^{*}(\overline{u}_{n}) \approx C_{n}^{ww}\label{eq:m11cc}
\end{eqnarray}
thereby obtaining 
\begin{eqnarray}\label{eq:11b}
\overline{u}_{n+1} =\overline{u}_{n}+C_{n}^{uw}(C_{n}^{ww} +\alpha_{n}\Gamma   )^{-1}(y^{\eta}-\overline{w}_{n})
\end{eqnarray}
that we use as the update formula for the ensemble mean of our iterative scheme. Similarly, from (\ref{eq:m7}) - (\ref{eq:m11cc}), the selection of $\alpha$ in (\ref{eq:m5}) and the stopping criteria (\ref{eq:m4}) become 
\begin{eqnarray}\label{eq:m12}
\fl  \rho \vert\vert \Gamma^{-1/2}( y^{\eta}-\overline{w}_{n}))\vert\vert_{Y} \leq \alpha_{n} \vert\vert \Gamma^{1/2}(C_{n}^{ww}+\alpha_{n}\Gamma)^{-1}(y^{\eta}-\overline{w}_{n})\vert\vert_{Y}.
\end{eqnarray}
and (\ref{eq:m15}), respectively. The sequence defined in (\ref{eq:m17}) satisfies (\ref{eq:m12}). In addition, note that the sample mean obtained from the ensemble updated according to formula (\ref{eq:m16}) satisfies indeed (\ref{eq:11b}). We then obtain the proposed ensemble Kalman scheme presented in  Algorithm \ref{Al1}.

\subsection{The proposed method as a sequence of linear inverse problems}\label{3_2}

As state earlier, the proof of convergence of Algorithm \ref{Al1} is beyond the scope of the present manuscript. Nonetheless, in this subsection we show some relevant properties which shed light on its regularizing effect. In order to study the aforementioned properties it is advantageous to rewrite the proposed algorithm in terms of the following augmented variables \cite{EnKF_US},
\begin{eqnarray}\label{eq:m18}
z=\left(\begin{array}{c}
u\\
w\end{array}\right),\qquad \Xi(z)=\left(\begin{array}{c}
u\\
G(u)\end{array}\right),
\end{eqnarray}
the space $Z\equiv X\times Y$ and the projection operators $H$ and $H^{\perp}$ defined by $Hz=w$ and $H^{\perp}z=u$, respectively. It is not difficult to see that, in terms of these new variables, Algorithm \ref{Al1} becomes

\begin{algorithm}{~~}\label{Al2}\\
Let $\{u_{0}^{(j)}\}_{j=1}^{N_{e}}\subset X$ be the initial ensemble of $N_{e}$ elements. Let $\rho\in (0,1)$ and $\tau>1/\rho$. Define 
\begin{eqnarray}\label{eq:fwd}
z_{0}^{(j,a)}\equiv \left(\begin{array}{c}
u_{0}^{(j)}\\
G(u_{0}^{(j)})\end{array}\right)
\end{eqnarray}
For $n=1,\dots$
\begin{itemize}
\item[(1)] \textbf{Prediction step.} Evaluate the forward map (\ref{eq:m18}),
\begin{equation}
\label{eq:m19}
z_{n}^{(j,f)}=\Xi(z_{n-1}^{(j,a)}).
\end{equation}
and define sample mean:
\begin{eqnarray}
 \overline{z}_{n}^{f}=\frac{1}{N_{e}}\sum_{j=1}^{N_{e}} z_{n}^{(j,f)}
\end{eqnarray}

\item[(2)] \textbf{Discrepancy principle}. If
\begin{eqnarray}\label{eq:dis2}
 \vert\vert \Gamma^{-1/2}(y-H\overline{z}_{n}^{f})\vert\vert_{Y}\leq \tau \eta 
\end{eqnarray}
stop. Output
\begin{eqnarray}\label{eq:update3}
\overline{u}_{n} \equiv \frac{1}{N_{e}}\sum_{j=1}^{N_{e}} H^{\perp} z_{n}^{(j,f)} =\frac{1}{N_{e}}\sum_{j=1}^{N_{e}} u_{n}^{(j)}
\end{eqnarray}
Define sample covariance:
\begin{eqnarray}
C_{n} z =\frac{1}{N_{e}-1}\sum_{j=1}^{N_{e}} (z_{n}^{(j,f)}-\overline{z}_{n}^{f})\langle z_{n}^{(j,f)}-\overline{z}_{n}^{f},z\rangle_{Z}.\label{eq:cov}
\end{eqnarray}
\item[(3)]\textbf{Analysis step.}  
Update each ensemble member as follows
\begin{eqnarray}
\label{eq:updateB}
z_{n}^{(j,a)}&=&z_{n}^{(j,f)}+C_{n}H^{\ast} (HC_{n}H^{\ast}+\alpha\Gamma)^{-1}(y^{\eta}-Hz_{n}^{(j,f)})
\end{eqnarray}
where $\alpha$ satisfies
\begin{eqnarray}\label{eq:m13}
 \rho \vert\vert \Gamma^{-1/2}(y^{\eta}-H\overline{z}_{n}^{f} )\vert\vert_{Y}\leq  \vert\vert \Gamma^{-1/2}( y^{\eta}-H\overline{z}_{n}^{a}(\alpha))\vert\vert_{Y}
\end{eqnarray}
with 
\begin{eqnarray}\label{eq:new1}
 \overline{z}_{n}^{a}=\frac{1}{N_{e}}\sum_{j=1}^{N_{e}} z_{n}^{(j,a)}
\end{eqnarray}


\end{itemize}
\end{algorithm}

It is not difficult to see \cite{EnKF_US} that $C_{n}$ in (\ref{eq:cov}) can be written as 
\begin{eqnarray*}
C_{n} =\left( \begin{array}{cc}
C_{n}^{uu} &C_{n}^{uw}\\
C_{n}^{wu}  & C_{n}^{ww} \end{array}\right)
\end{eqnarray*}

While the regularization properties of the proposed method will be studied numerically in \Sref{Numerics}, it is clear that the essence of the proposed regularizing Kalman method is the selection of $\alpha$ by means of (\ref{eq:m13}). Since $0<\rho<1$, the existence of such $\alpha$ is ensured by the following proposition
\begin{proposition}\label{Prp2B}
At every iteration of the scheme, the ensemble mean $\overline{z}_{n}^{a}(\alpha)$ defined by (\ref{eq:new1}) satisfies 
\begin{itemize}
\item[(i)] $\vert\vert \Gamma^{-1/2}(y-H\overline{z}_{n}^{a}(\alpha))\vert\vert_{Y}\to\vert\vert \Gamma^{-1/2}(y-H\overline{z}_{n}^{f})\vert\vert_{Y}$ as $\alpha\to \infty$.
\item[(ii)]$\vert\vert \Gamma^{-1/2}(y-H\overline{z}_{n}^{a}(\alpha))\vert\vert_{Y}\to 0$ as $\alpha\to 0$.
\item[(iii)] The map $\alpha\to \vert\vert\Gamma^{-1/2}( y-H\overline{z}_{n}^{a}(\alpha))\vert\vert_{Y}$ is monotonously nondecreasing
\end{itemize}
\end{proposition}
\textbf{Proof:}  The proof can be carried out as in \cite[Theorem 2.16]{Kirsch} $\Box$.

Similar to the regularizing LM scheme where each iterate solves a linear inverse problem (see expression (\ref{eq:m1})), our ensemble Kalman scheme can be posed as Tikhonov-type regularization. This equivalence motivates the need to further regularize standard Kalman methods. 

Note that the ensemble mean of the analysis step (\ref{eq:updateB}) is given by
\begin{eqnarray}\label{eq:mean_update}
\overline{z}_{n}^{a}=\overline{z}_{n}^{f}+C_{n}H^{\ast} (HC_{n}H^{\ast}+\alpha\Gamma)^{-1}(y^{\eta}-H\overline{z}_{n}^{f})
\end{eqnarray}
We now show that (\ref{eq:mean_update}) can be posed as the solution of a regularized linear inverse problem
\begin{proposition}\label{Prp3}
Assume that at each iteration, the ensemble $\{z_{n}^{(j,f)}\}_{j=1}^{N_{e}}$ is linearly independent. Then, the ensemble mean of the analysis step (\ref{eq:mean_update}) of the ensemble Kalman algorithm satisfies
\begin{eqnarray}\label{eq:m19}
\fl  \overline{z}_{n}^{a } =\textrm{argmin}_{z\in \mathcal{Z}_{n}}\Big(\vert\vert \Gamma^{-\frac{1}{2}}(y^{\eta}-H\overline{z}_{n}^{f }-H(z-\overline{z}_{n}^{f }))\vert\vert_{Y}^2+\alpha\vert\vert z-\overline{z}_{n}^{f}\vert\vert_{ C_{n}} ^2\Big)
\end{eqnarray}
where $\mathcal{Z}_{n}$ is the completion of $\mathcal{R}(C_{n})$ with respect to the norm induced by $\langle\cdot ,C_{n}^{-1}\cdot \rangle_{Z}$ and denoted by $\vert\vert \cdot \vert\vert_{ C_{n}}$.
\end{proposition}
\textbf{Proof:}   From the definition of $C_{n}$ in (\ref{eq:cov}) it is clear that the rank of $C_{n}$ is given by $\mathcal{R}(C_{n})=span\{ z_{n}^{(j,f}\}$. Also note that $C_{n}$ is a sum of rank-one operators, hence compact. If  $N_{e}<\infty$, then $\textrm{dim}(\mathcal{R}(C_{n}))=N_{e}<\infty$; thus  from the linearly independence of the initial ensemble it follows that the restriction $C_{n}:\mathcal{R}(C_{n}) \to \mathcal{R}(C_{n})$ is a positive definite operator. Hence its inverse $C_{n}^{-1}:\mathcal{R}(C_{n}) \to\mathcal{R}(C_{n})$ exists. Moreover, we may consider the restriction of $H$ into $\mathcal{R}(C_{n})$ as well as its adjoint $H^{*}:\mathcal{R}(C_{n})\to Y$. Thus, the operators $C_{n}$, $C_{n}^{-1}$ and $H$ can be represented with matrices and so standard matrix algebra can be directly applied to show that (\ref{eq:m1}) is equivalent to \cite{Tarantola}
\begin{eqnarray}\label{eq:ape9}
 ( H^{*}\Gamma^{-1}H + \alpha C_{n}^{-1}  )(\overline{z}_{n}^{a }-\overline{z}_{n}^{f}))=H^{*}\Gamma^{-1}(y-H\overline{z}_{n}^{f})
\end{eqnarray}
which, from simple standard arguments is the solution to the normal equations associated to the minimization of 
\begin{eqnarray}\label{eq:ape10}
 J(z)\equiv \vert\vert \Gamma^{-\frac{1}{2}}(y-H\overline{z}_{n}^{f }-H(z-\overline{z}_{n}^{f }))\vert\vert^2+\alpha \vert\vert z-\overline{z}_{n}^{f}\vert\vert_{ C_{n}} ^2
\end{eqnarray}
in the space $\mathcal{Z}_{n}$. For the case $N_{e}\to \infty$  we may consider $C_{n}^{-1}:\mathcal{R}(C_{n})  \to \mathcal{R}(C_{n})$ as a densely defined unbounded operator. Then a proof based on representers can then be carried out similarly to the proof of Lemma \ref{lema:Rep1}. $\Box$.

From Proposition \ref{Prp3} we notice that the ensemble mean at each iteration of the ensemble Kalman algorithm is the Tikhonov-regularized solution (in the subspace $\mathcal{Z}_{n}$ with norm $\vert\vert \cdot\vert\vert_{C_{n}}$) of the following ``artificial'' linear inverse problem:
\begin{eqnarray}\label{eq:ip}
\fl  \textit{given\quad $\tilde{y}\equiv y^{\eta}-H\overline{z}_{n}^{f }$\quad find \quad $w\equiv z-\overline{z}_{n}^{f}\in \mathcal{Z}_{n}$ \quad such that\quad }
\tilde{y}= Hw
\end{eqnarray}
Let us define $z^{\dagger}=(u^{\dagger}, G(u^{\dagger}))^{T}$ and $w^{\dagger}\equiv  z^{\dagger}-\overline{z}_{n}^{f}$ where $u^{\dagger}$, we recall, is the truth. 
If at the $n$th iteration we were to solve the original identification problem (i.e find/approximate $u^{\dagger}$) that would imply to find $w^{\dagger}$ defined above. In other words, $w^{\dagger}$ is the solution to the inverse problem (\ref{eq:ip}) and so the corresponding exact data is $ \tilde{y}^{\dagger}\equiv  Hw^{\dagger}$ which is 
\begin{eqnarray}\label{eq:m20}
\fl  \tilde{y}^{\dagger}\equiv  Hw^{\dagger}= H(z^{\dagger}-\overline{z}_{n}^{f})=G(u^{\dagger})-H\overline{z}_{n}^{f}=y^{\eta}-\xi-H\overline{z}_{n}^{f}=\tilde{y}-\xi
\end{eqnarray}
where we have used (\ref{eq:data}) and the definition of $H$. From the previous expression and (\ref{eq:nl}) it follows that 
\begin{eqnarray}\label{eq:m21}
  \vert\vert \Gamma^{-1/2}(\tilde{y}^{\dagger}-\tilde{y} )=\vert \vert \Gamma^{-1/2}\xi\vert\vert_{Y}=\eta
\end{eqnarray}
which implies that the noise level for the linearized inverse problem solved by $\overline{z}_{n}^{a}$, is the same as the noise level of the original nonlinear inverse problem that we aim at solving by means of the proposed method. Moreover, note that before convergence is achieved, i.e. when $\tau \eta <\vert\vert \Gamma^{-1/2}(y^{\eta}-H\overline{z}_{n}^{f} )\vert\vert_{Y}$ from (\ref{eq:m13}) we then have that
\begin{eqnarray}\label{eq:m22}
\fl  \eta< \frac{1}{\tau}\vert\vert \Gamma^{-1/2}(y^{\eta}-H\overline{z}_{n}^{f})\vert\vert_{Y}\leq \rho \vert\vert \Gamma^{-1/2}(y-H\overline{z}_{n}^{f})\vert\vert_{Y}\leq \vert\vert \Gamma^{-1/2}(y^{\eta}-H\overline{z}_{n}^{a}(\alpha))\vert\vert_{Y}.
\end{eqnarray}
Thus,  at each iteration, the selection of $\alpha$ in (\ref{eq:m13}) honors the discrepancy principle for the artificial inverse problem that the ensemble mean solves at each iteration of the scheme.

\begin{remark}\label{rema1}
From the augmented version of the proposed ensemble Kalman method  we may appreciate more clearly the role of the parameter $\rho$. Indeed, note from Proposition \ref{Prp2B} that a value of $\rho\approx 1$ ($\rho<1$) in (\ref{eq:m13}) will results in a very large $\alpha$ and so the estimate $\overline{z}_{n}^{a}(\alpha)$ will move slightly from the previous estimate $\overline{z}_{n}^{f}(\alpha)$. Then, a selection of $\tau\approx 1/\rho$, ($\tau>1/\rho$) for $\rho\approx 1$ will enable us to allow the algorithm to progress until the data misfit is close to the noise level. In contrast, a smaller $\rho$ results in a smaller $\alpha$ and thus a less controlled estimate $\overline{z}_{n}^{f}(\alpha)$. Therefore, choosing $\tau\approx 1/\rho$, ($\tau>1/\rho$) for smaller $\rho$ implies that we need to stop the algorithm via (\ref{eq:dis2}) much earlier to prevent the data overfitting as previously discussed. While it seems clear that a choice of $\rho\approx 1$ will result in a more stabilized scheme, it may be also detrimental to the performance of the scheme. In \Sref{Numerics} we investigate, with numerical experiments, practical choices for the parameter $\rho$.
\end{remark}

\subsection{The connection with the Bayesian framework}\label{bayesian}

While here we propose an ensemble Kalman method as a derivative-free regularizing method for classical inverse problems, most  Kalman-based approaches are posed in a Bayesian inference framework where the objective is to produce an ensemble from which statistical information from the Bayesian posterior can be computed. In order to fully understand the underlying motivation of our method it is instructive to consider the role of the proposed method in the context of Bayesian inversion. Let us then assume that the noise $\xi$ in (\ref{eq:data}) centered Gaussian noise with covariance $\Gamma$. Additionally, assume that there is an underlying Gaussian prior distribution $N(\overline{u},C)$ of the unknown $u$, where the mean and covariance of such distribution are $\overline{u}$ and $C$ respectively. Choose an initial ensemble generated from samples of this distribution, i.e. $u_{0}^{(j)}\sim N(\overline{u}, C)$. Let us now consider the first iteration $n=0$ of Algorithm \ref{Al1} with a fixed parameter $\alpha=1$ and replace $y^{\eta}$ in (\ref{eq:m16}) by $y^{(j)}=y^{\eta}+\eta^{(j)}$ where $\eta^{(j)}\sim N(0,\Gamma)$. Then, Algorithm \ref{Al1} becomes
\begin{eqnarray}\label{eq:m23}
u_{1}^{(j)} =u_{0}^{(j)}+C_{0}^{uw}(C_{0}^{ww} +\Gamma   )^{-1}(y^{(j)}-G(u_{0}^{(j)}))
\end{eqnarray}
which is also the standard formula for the so-called Ensemble Smoother \cite{evensen2009data}. If the forward operator $G$ is linear then it can be shown \cite{EnKF_US} that the updated ensemble $\{u_{1}^{(j)}\}_{j=1}^{N_{e}}$ fully characterizes, as $N_{e}\to \infty$, the conditional probability of $u$ given the data $y^{\eta}$. In other words, (\ref{eq:m23}) provides an ensemble approximation of the Bayesian posterior. Moreover, the mean of the ensemble converges, in the limit $N_{e}\to \infty$, to the maximum a posterior estimate, i.e. the maximizer of the  aforementioned posterior distribution. The linear-Gaussian case is trivial in the sense that the posterior is Gaussian and it can be easily characterized with its mean and covariance. However, when the forward operator $G$ is nonlinear, such as the ones described in subsection \ref{test}, the Bayesian posterior is in general non-Gaussian and the ensemble smoother provides only an approximation whose convergence theory, to our best knowledge, is nonexistent. It has been often reported, however, than when straight forward application of (\ref{eq:m23}) are used in the general nonlinear case, inaccurate approximations of the Bayesian posterior may be obtained \cite{Evaluation}. In a recent publication \cite{Yo}, we have demonstrated that the application of iterative regularization methods like the ones proposed here may improve the accuracy of ensemble methods for approximation of the statistical inverse problem. The present work  follows the classical approach; we endow an intrinsically Bayesian method with the regularizing properties needed to address deterministic inverse problem in a derivative-free easy-to-implement computational framework.


\subsection{Small ensemble size ill-conditioning}

As stated above, issues of stability of ensemble Kalman methods for data assimilation have been often reported. Ad-hoc fixes such as covariance localization and covariance inflation have become standard practice to address such lack of stability which have been typically associated and more noted with implementations of small ensemble size with respect to the number of measurements \cite{EnKF_Review}. Indeed, note that since $\textrm{dim}(Y)=M<\infty$,then $\textrm{dim}(\mathcal{R}(C_{n}^{ww}))=\min\{N_{e},M\}$. Clearly, when $N_{e}<M$, $C_{n}^{ww}$ is singular and, even though $(C_{n}^{ww}+\Gamma)^{-1}$ is strictly positive definite, it could potentially have a large conditioning number for some choices of $\Gamma$, which is in turn selected according to measurement information. This type ill-conditioning when $N_{e}<M$ has the similar effect of the ill-posedness of the Fr\'echet derivative of the forward map in the regularizing LM scheme. Thus, it comes as no surprise that the application of ideas from iterative regularization techniques can be useful for addressing the ill-conditioning of the method when used with a small ensemble size with respect to the number of measurements. From either sources of ill-conditioning, it is clear that ensemble Kalman methods need regularization and, to our best knowledge, the standard fixes indicated above do not possess a mathematical framework which provides general guidelines for developing robust implementations. In contrast, iterative regularization has a rigorous general mathematical ground; our goal here is to use it as the tool to empower ensemble Kalman methods with the stabilization/regularization required to applied them as accurate robust but computationally tractable solvers for large-scale PDE-constrained inverse problems.

\subsection{Computational aspects of the scheme}\label{Sec:cost}

Under our assumption that the number of measurements is small compared to the dimension of the (discretized) parameters space, the main computational cost of the proposed scheme is in the evaluation of the forward map in the prediction step (\ref{eq:m14}). The construction of $C_{n}^{uw}$,  $C_{n}^{ww}$ as well as the inversion of $(C_{n}^{ww} +\alpha_{n}\Gamma   )^{-1}$ are negligible for the number of measurement considered for the present applications.Therefore, the cost of Algorithm \ref{Al1} is mainly $n^{*}N_{e}$ forward model evaluations where $n^{\star}$ is the stopping iteration determined by (\ref{eq:m13}). It is clear that Algorithm \ref{Al1} is computationally feasible provided that, with a few $N_{e}$ ensemble members, a reasonably accurate estimate can be achieved in only a few iterations. In Section \ref{Numerics} we provide an extensive numerical investigation of the performance and computational feasibility of the scheme. In some cases, a relatively large number of ensembles are required to achieve stable computations with reasonable accuracy.

\section{Numerical Experiments}\label{Numerics}

By means of numerical experiments, in this section we investigate the convergence and regularizing properties of the iterative regularizing ensemble Kalman method. We consider the identification problems on the test models from subsection \ref{test} and conduct synthetic experiments where we apply the proposed scheme with synthetic data generated with a true parameter. Algorithm \ref{Al1} then produces an estimate (the ensemble mean (\ref{eq:m7A})) which we compare against the truth in order to assess the accuracy of the scheme. By conducting synthetic experiments, our aim is to study the accuracy, convergence and regularizing properties of the proposed methods with respect to 
\begin{itemize}
\item ensemble size $N_{e}$
\item tunable parameters $\rho$ and $\tau$
\item number of measurements $M$
\item noise level $\eta$
\item selection of initial ensemble $\{u_{0}^{(j)}\}_{j=1}^{N_{e}}$
\end{itemize}

\subsection{Test Model I. Darcy flow}\label{Darcy_Num}

\begin{figure}[htbp]
\begin{center}
\includegraphics[scale=0.30]{True_GW}
\includegraphics[scale=0.30]{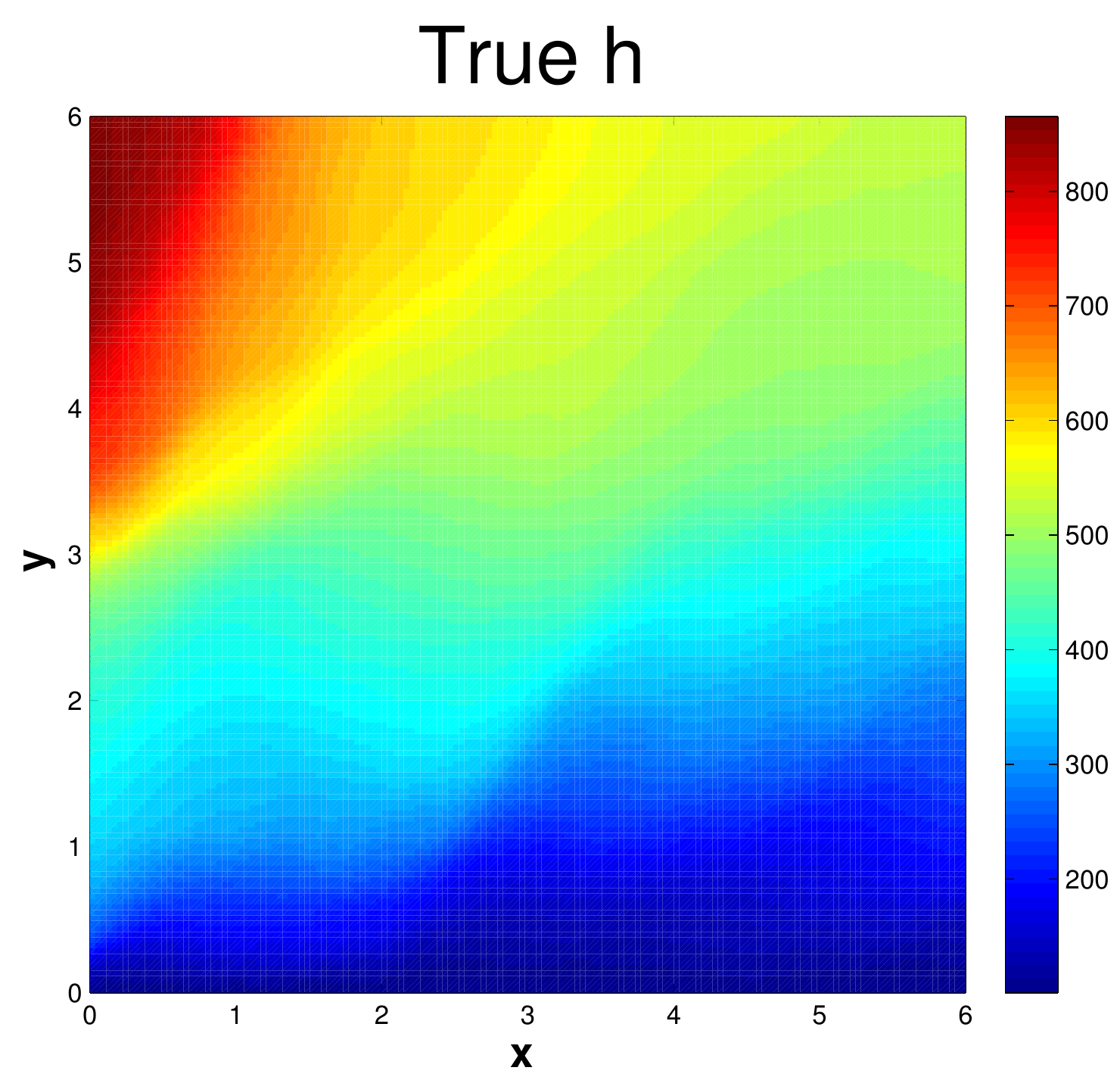}
\includegraphics[scale=0.30]{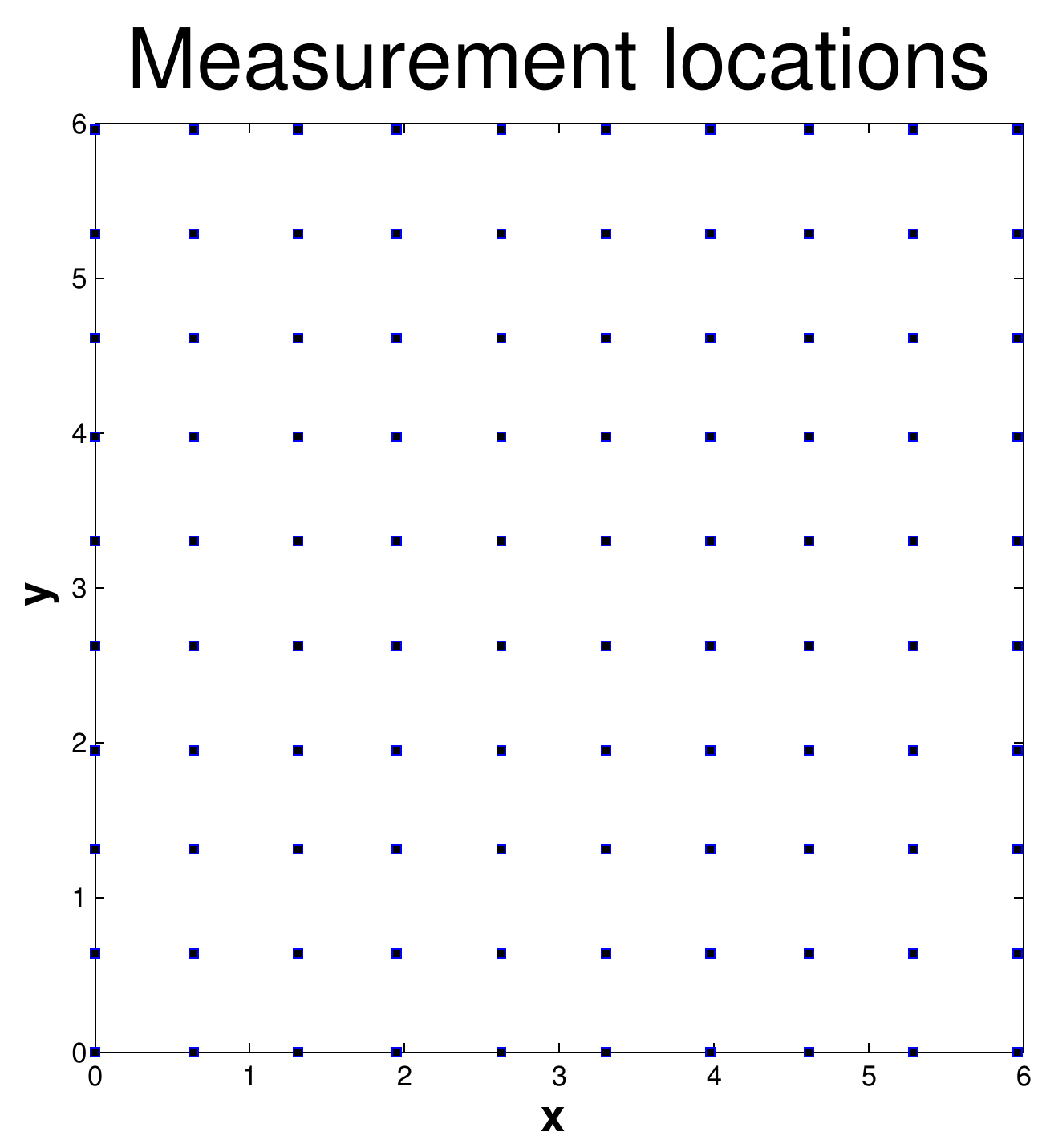}
 \caption{Left: true log-conductivity. Middle: log-data misfit.  Right: relative error with respect to the truth }
    \label{Fig1}
\end{center}
\end{figure}

We now describe the generation of synthetic data from the Darcy flow model of subsection \ref{Darcy}. We consider the identification of the log of the ``true''  hydraulic conductivity $u^{\dagger}$ of Figure \ref{Fig1}  (left).
The true conductivity is a random field drawn from a  Gaussian measure $N(0,C)$ with covariance $C:L^{2}(D)\to L^{2}(D)$ defined by 
  \begin{equation}\label{eq:cova2}
C \phi(x) = \int_{D}c(x, y)\phi (y)\dd y 
\end{equation}
with a spherical covariance function
  \begin{equation}\label{eq:cova2B}
 c(x, y) = \left\{\begin{array}{cc}
c_{0}\Big[1-\frac{3}{2}\frac{h}{a}+\frac{1}{2}\frac{h^3}{a^3}\Big]  &\textrm{if}~~h\equiv \vert x-y\vert<a \\
0 &\textrm{if}~~h>a \end{array}\right.
\end{equation}
where $a>0$ is the correlation length and $c_{0}>0$ is the variance of the field. Gaussian distributions with covariance operators of this type are often used to model the geologic properties of some formations \cite{geos}. The true log conductivity field is used in equations (\ref{eq:a1})-(\ref{eq:a3}) to find the true head field displayed in Figure \ref{Fig1} (middle). The latter is then used in (\ref{eq:a4}) where the measurement locations $\{x_{m}\}_{m=1}^{M}$ will be specified below. Then, Gaussian noise is added to the synthetic data; the norm of the noise is $1\%$ of the norm of the data. The Darcy flow model is solved numerically with cell-centered finite differences \cite{Mixed}. The matrix $\Gamma$ is simply a diagonal matrix with entries proportional to the size of the corresponding entry of $G(u^{\dagger})$. In order to avoid inverse crimes,  Algorithm \ref{Al1} is applied on a coarser grid (of $80\times 80$ cells) than the one (of $160\times 160$ cells) used for the generation of the synthetic data. 

\begin{figure}[htbp]
\begin{center}
\includegraphics[scale=1]{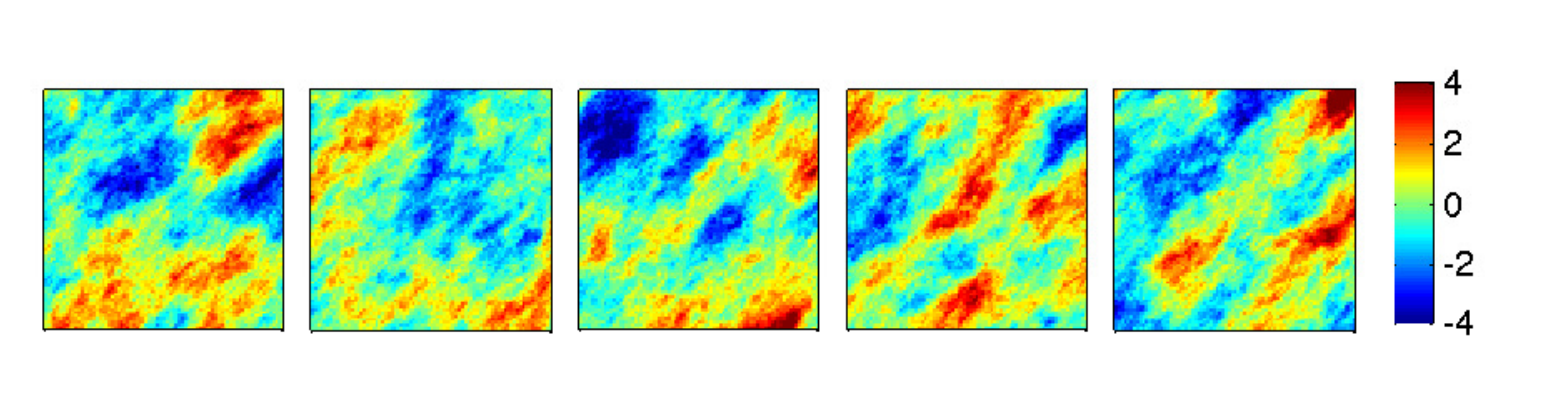}

 \caption{Some elements from the initial ensemble.}
    \label{Fig2}
\end{center}
\end{figure}

For the experiments of this section, we apply Algorithm \ref{Al1} with an initial ensemble that consists of $N_{e}$ (to be defined below) random fields drawn from the  same Gaussian measure $N(0,C)$ that we use to generate the truth. Some members of the initial ensemble are displayed in Figure \ref{Fig2}. The generation of these fields (as well as the truth) is conducted by means of the Karhunen-Loeve (KL) expansion of a random field distributed according to $N(0,C)$. Although the truth and the prior ensemble are samples from the same probability distribution, note that (i) the truth is generated on a much finer grid and (ii) the truth (or its projection on the coarser grid) is not part of the initial ensemble. 

The motivation for selecting the truth as a realization from the same distribution from the initial ensemble is twofold. On the one hand, in subsurface applications there is often prior knowledge of the geologic properties in terms of a prior probability distribution for the truth. It is then natural, as in the Bayesian approach, to use this distribution to generate the initial ensemble.  On the other hand, by choosing the truth and the initial ensemble form the same probability distribution we expect to reduce (although not completely eliminate) the effect that the selection of an initial ensemble has on the performance of the proposed scheme. More precisely, from the invariance subspace property (Proposition \ref{Prp2}) we know that the approximation provided by our method lives in the subspace generated by the initial ensemble. Therefore, each initial ensemble will produce a different estimate. However, by generating the initial ensemble from the same distribution that we use to generate the truth, our initial ensembles, while different from the truth, have the same regularity and spatial structure (e.g. correlation) of the true field. In the following subsection we study an inverse problem where the truth is a  prescribed field with no association whatsoever to the initial ensemble and there we study the effect that the spatial correlation of the initial ensemble has on the accuracy of the proposed ensemble method to identify the truth.

While the initial ensemble generated as described above will provide estimates with the same spatial features of the truth, from the aforementioned subspace property we certainly expect that the results form Algorithm \ref{Al1} will vary with the selection of the ensemble. Therefore, in order to fully assess the numerical properties of this algorithm without the influence of the initial ensemble, for the present work we will conduct multiple experiments corresponding to different choices of the initial ensemble and provide averages over these experiments and thus reliable quantities related to the average performance of the scheme. It is worth reiterating that the initial ensemble is a design parameter that can be chosen in various ways according to available information.

\subsubsection{Effect of the ensemble size $N_{e}$.}

\begin{figure}[htbp]
\begin{center}
\includegraphics[scale=0.215]{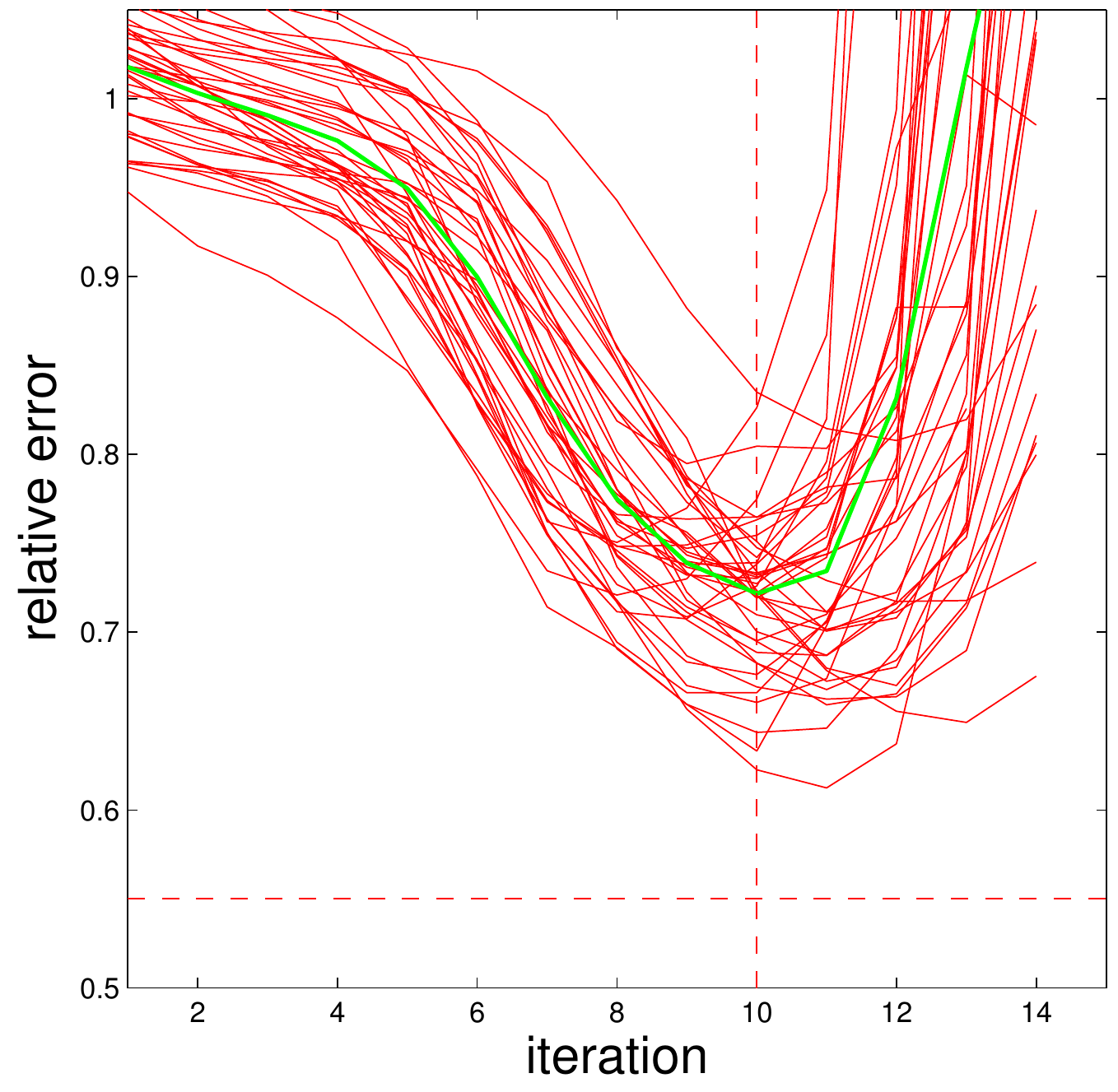}
\includegraphics[scale=0.215]{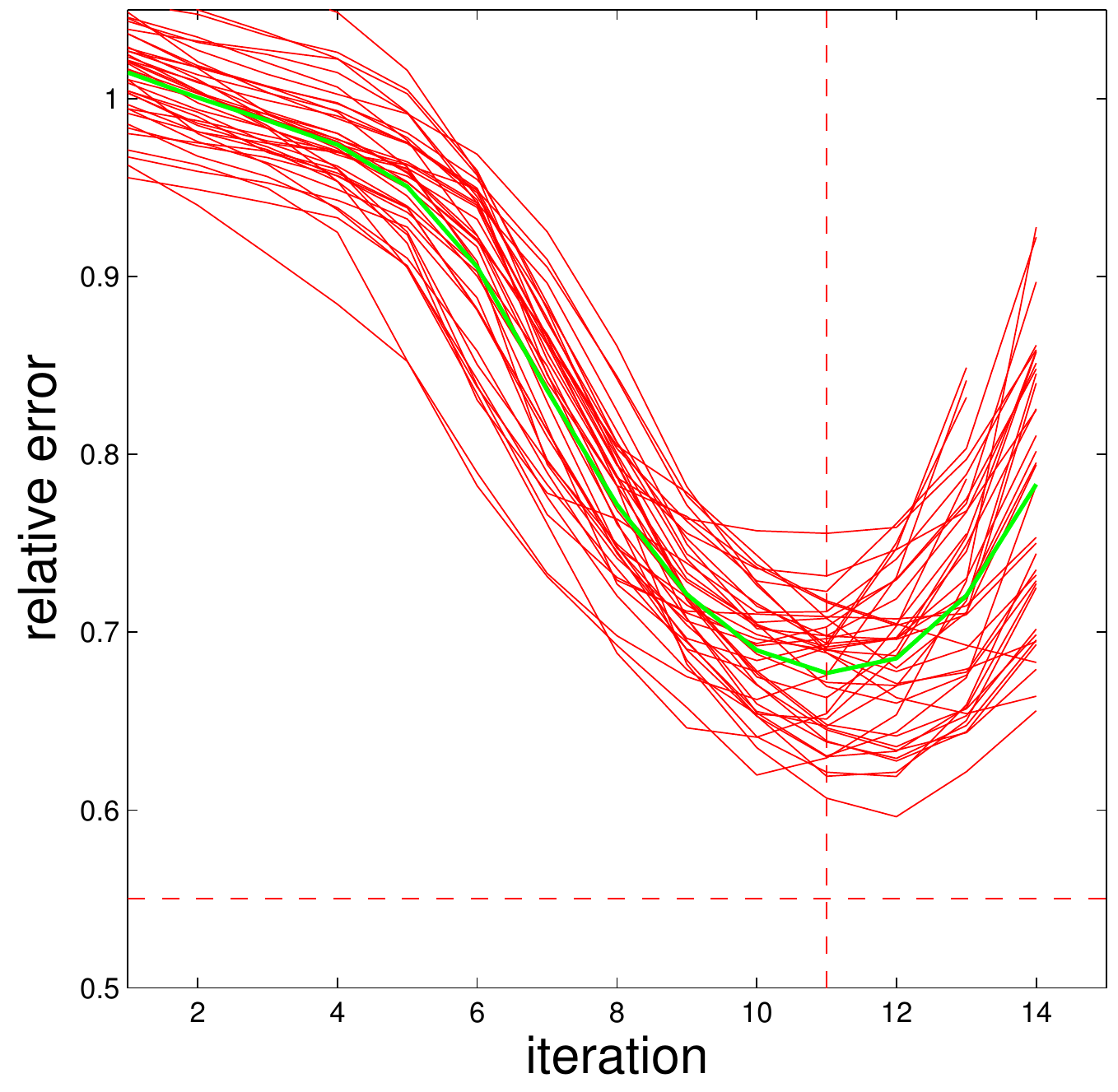}
\includegraphics[scale=0.215]{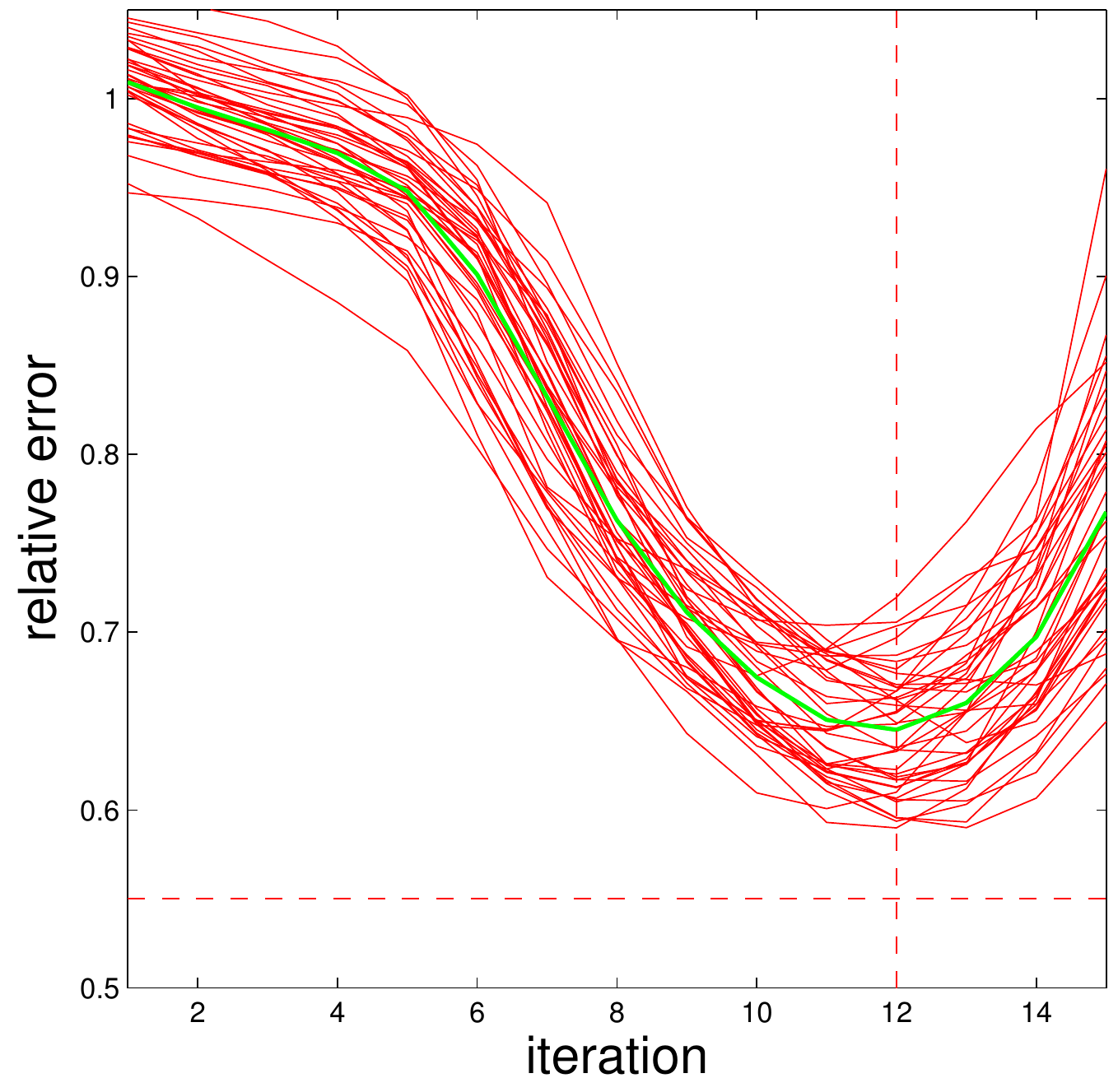}
\includegraphics[scale=0.215]{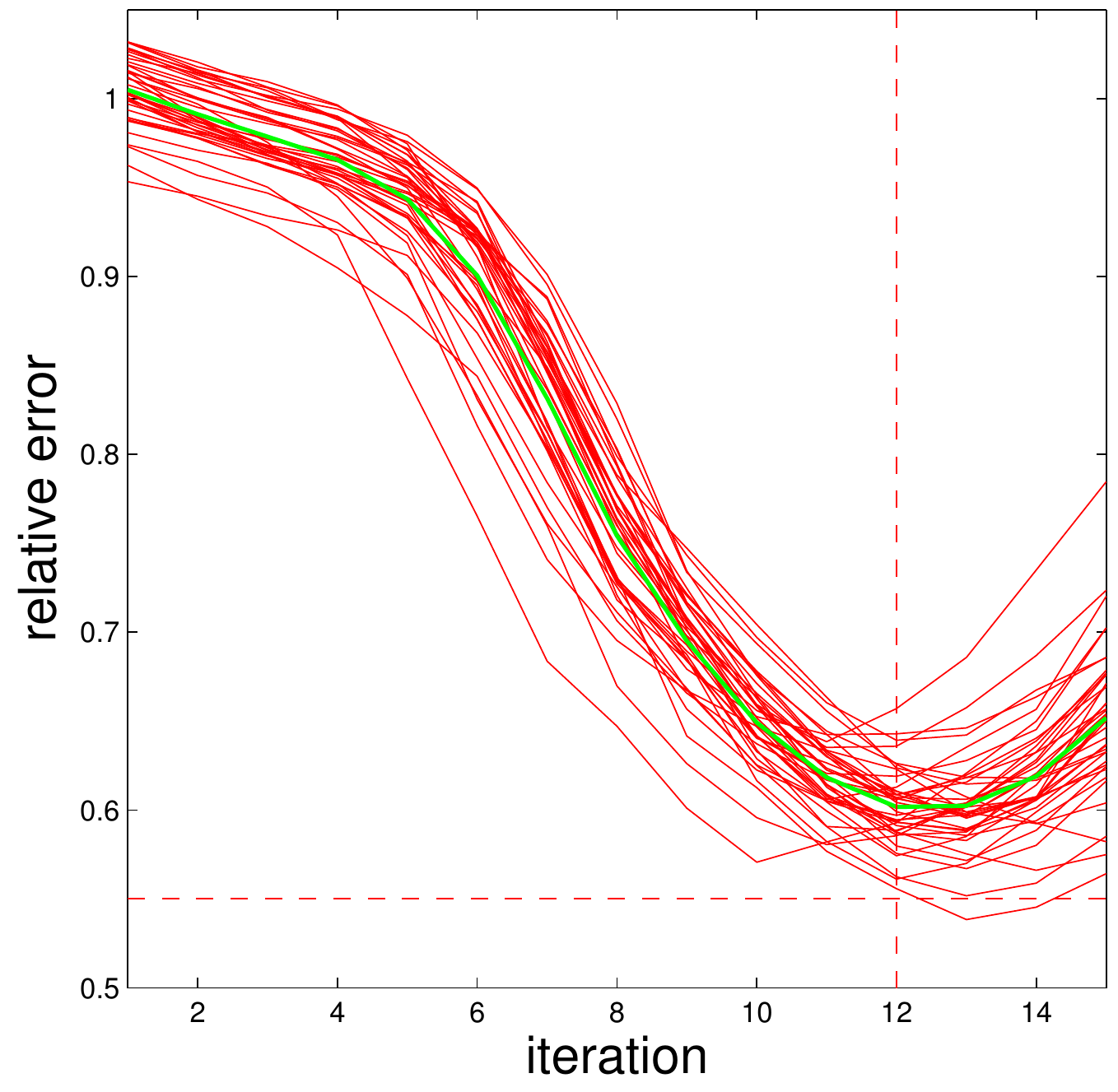}
\includegraphics[scale=0.215]{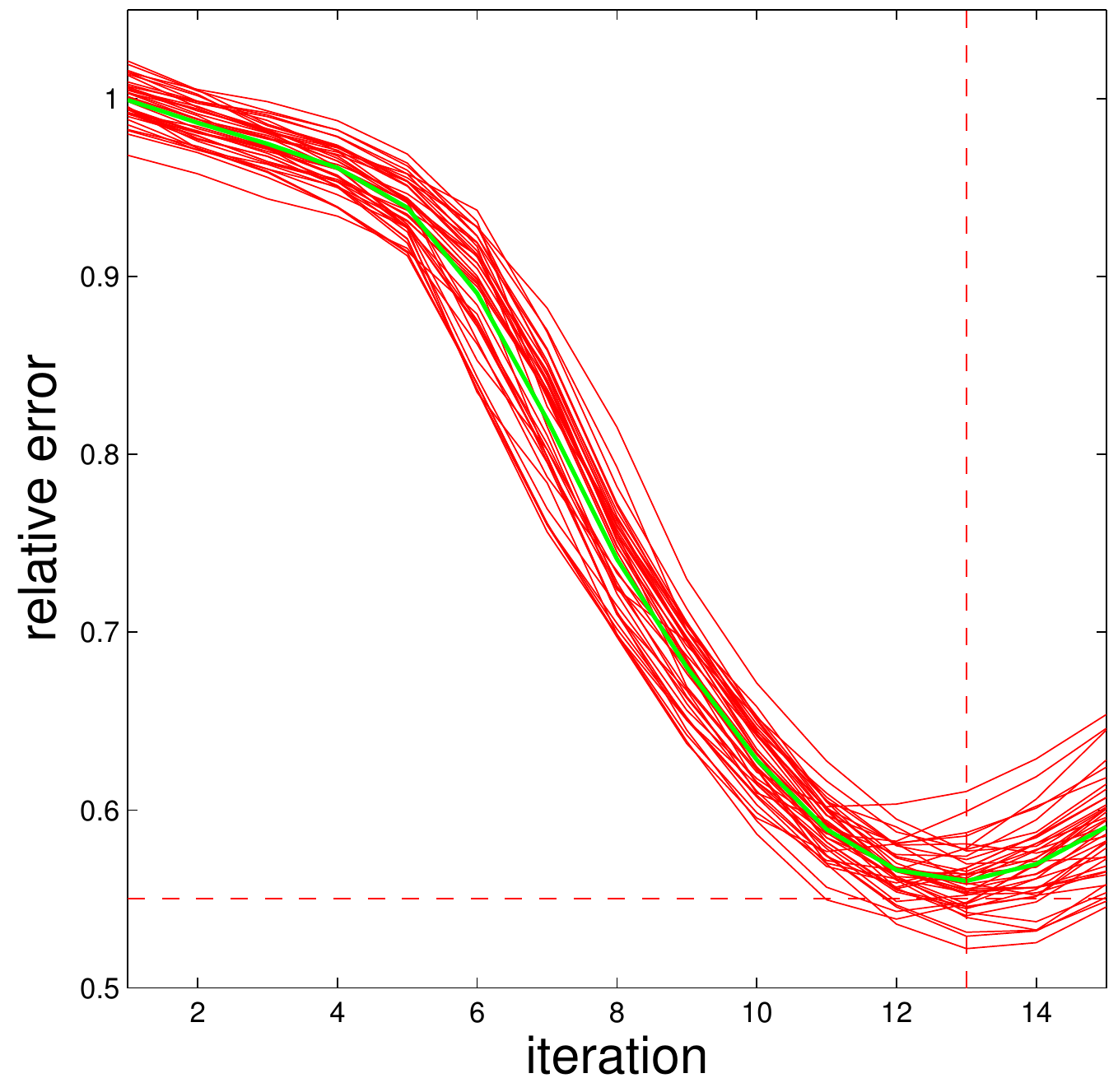}\\
\includegraphics[scale=0.215]{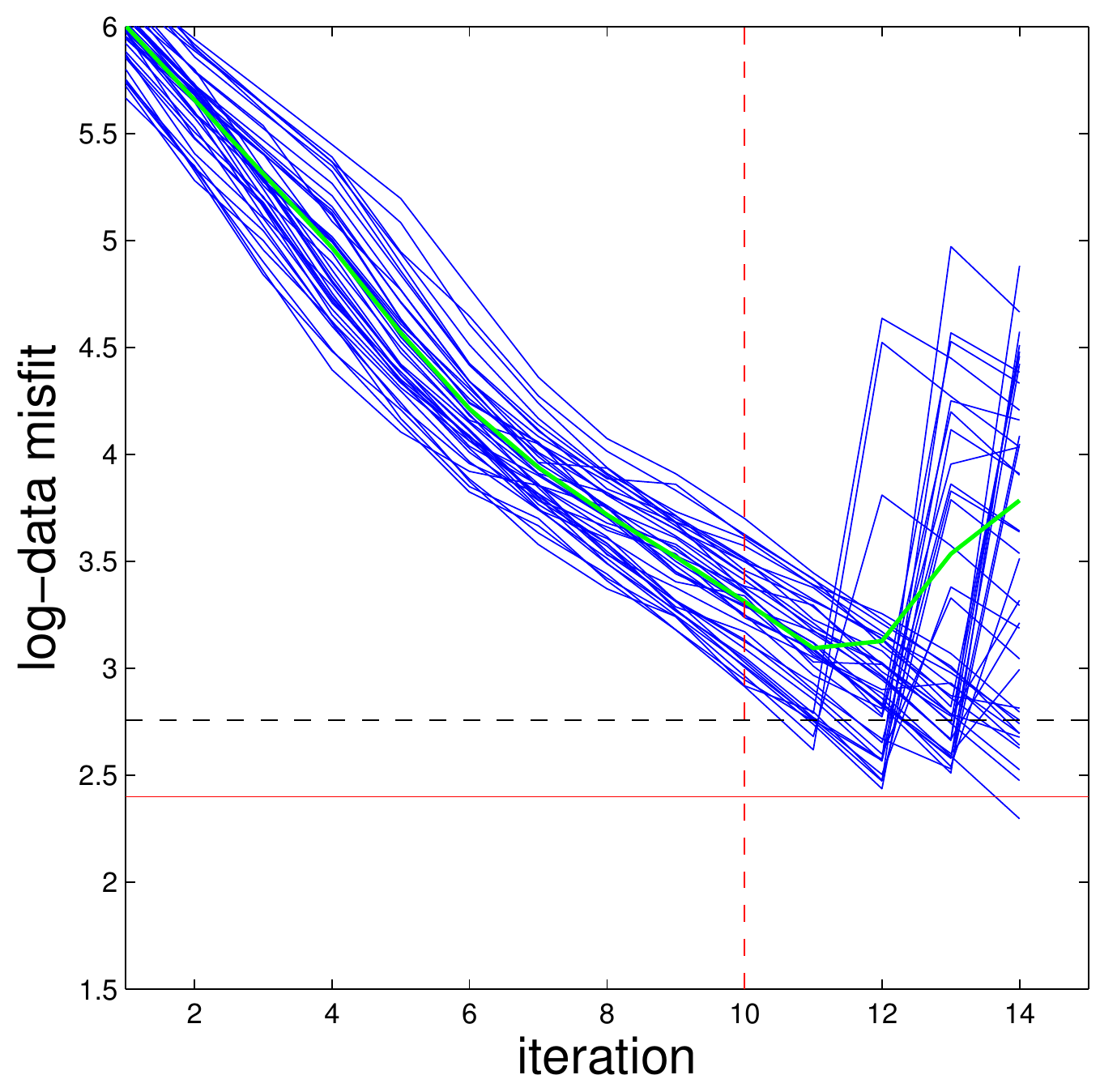}
\includegraphics[scale=0.215]{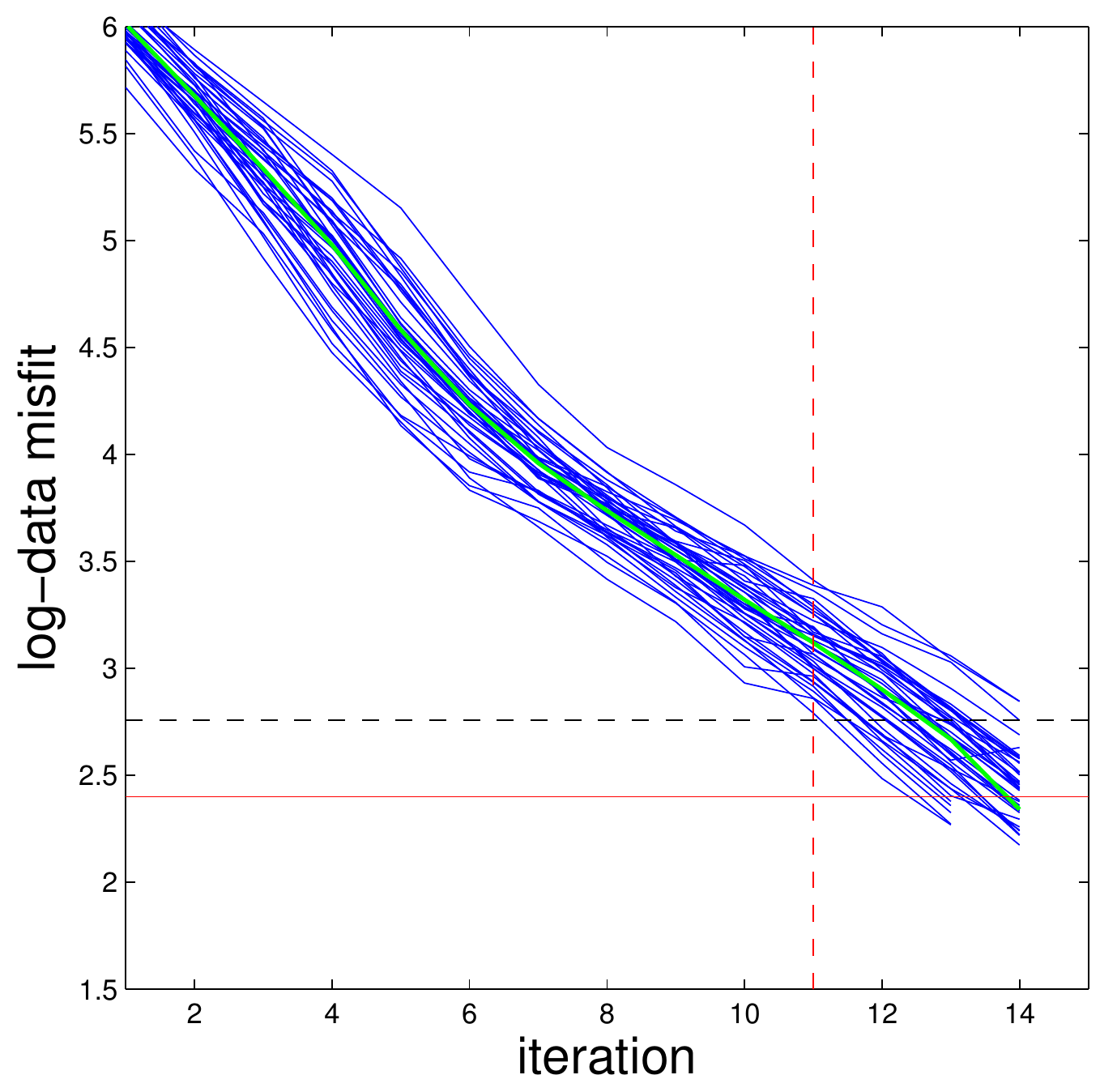}
\includegraphics[scale=0.215]{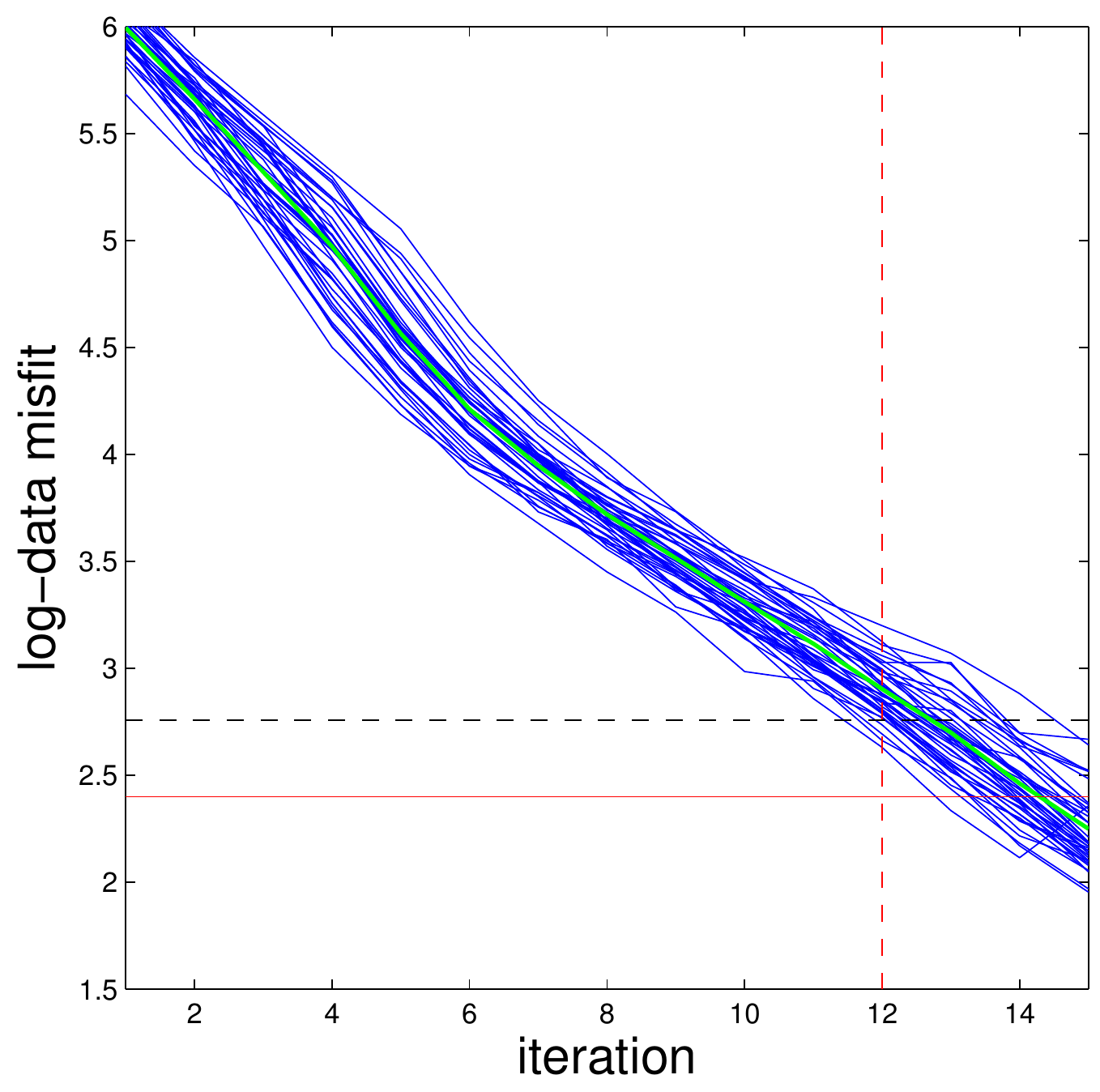}
\includegraphics[scale=0.215]{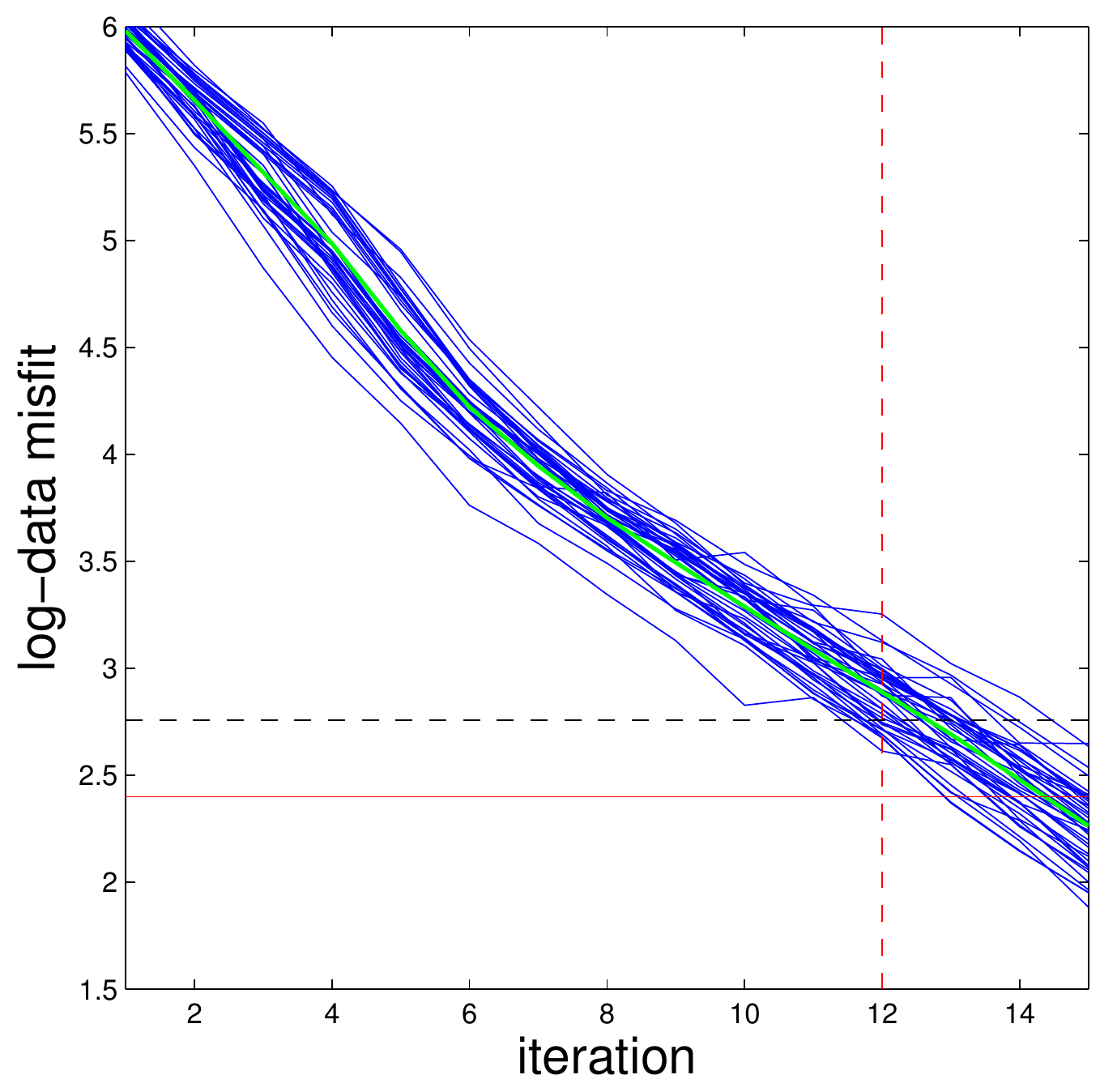}
\includegraphics[scale=0.215]{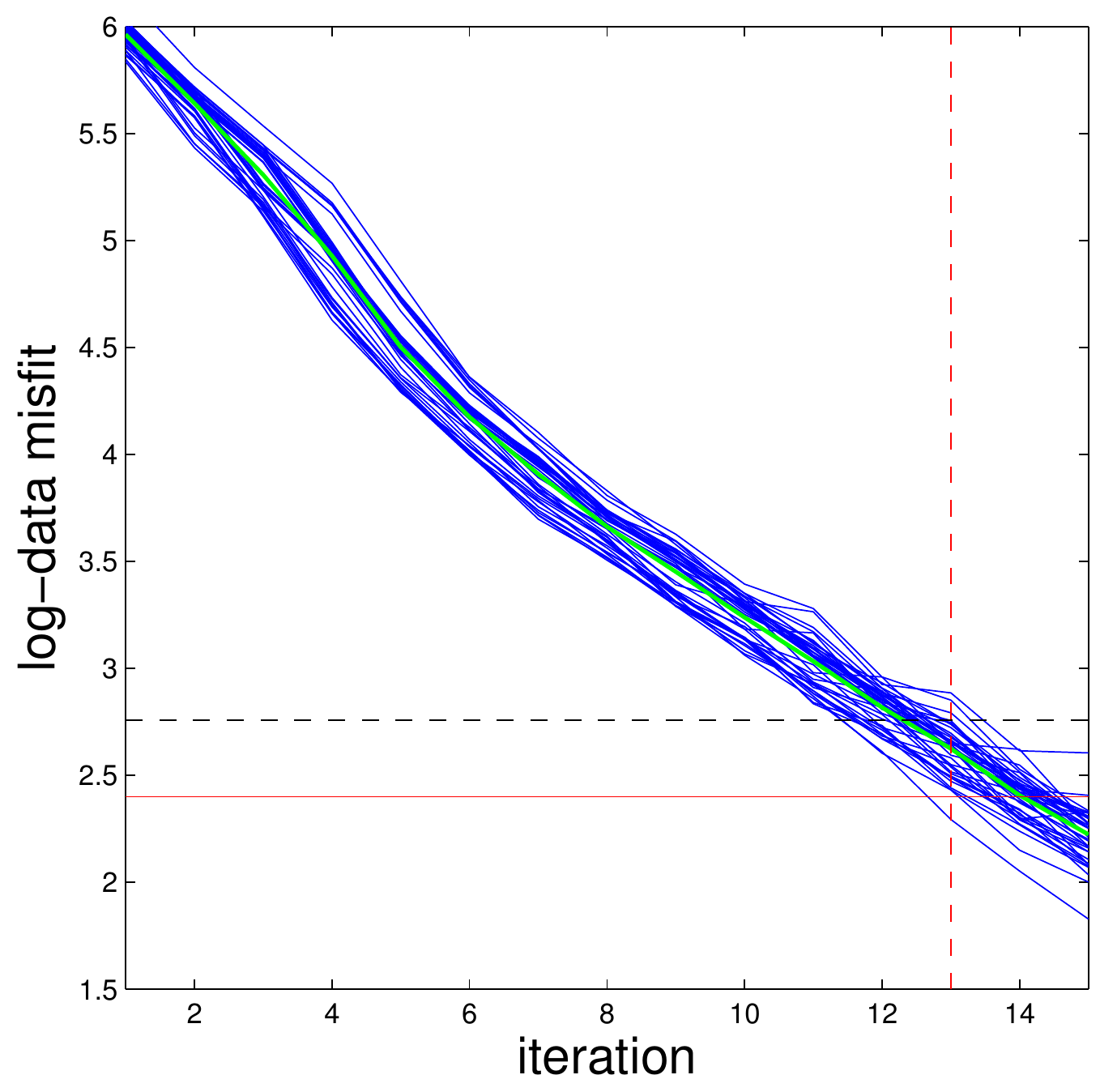}
 \caption{Relative error w.r.t. the truth (top) and log - data misfit (bottom) from 40 different experiments with $\rho=0.7$ associated to different initial ensembles of size (from left to right) $N_{e}=75,100,125,200,400$}
   \label{Fig3}
\end{center}
\end{figure}

We generate synthetic data as described below with the array of 100 measurement locations displayed in Figure \ref{Fig1} (right). We use these data in Algorithm \ref{Al1} with a (fixed) parameter $\rho=0.7$. In order to appreciate the effect of the early termination of the scheme, the algorithm is allowed to progress even when the data misfit goes below $\eta /\rho$ (recall the stopping criteria (\ref{eq:m15}) requires $\tau>1/\rho$). In the top (resp. bottom)  of Figure \ref{Fig3} we plot the relative error with respect to the truth $\vert\vert u_{n} -u^{\dagger}\vert\vert_{L^{2}(D)} /\vert\vert u^{\dagger}\vert\vert_{L^{2}(D)}$ (resp. log - data misfit) from 40 different experiments corresponding to different selection of the initial ensemble (recall each initial ensemble is a set of $N_{e}$ draws from a the Gaussian distribution defined above). Different panels corresponds to different selections of ensemble size $N_{e}$. Each of the blue curves (resp. red curves) represents the log-data misfit (resp. error with respect to the truth) that we obtain from the ensemble mean $\overline{u}_{n}$ (expression (\ref{eq:m7A})) computed with Algorithm \ref{Al1} initialized with each of the 40 initial ensembles mentioned earlier. We reiterate that by ``data misfit'' we mean the misfit between data and the average of the model outputs from the ensemble defined in (\ref{eq:m15}) which is, in turn, the quantity monitored for the convergence of the scheme.

The green curve in the bottom (resp. top) of \Fref{Fig3} represents the log-data misfit (resp. error w.r.t truth) averaged, at each iteration, over the 40 experiments from different initial ensembles. The dotted vertical line defines the iteration number after which the averaged relative error w.r.t. the truth starts increasing. The dotted horizontal line in \Fref{Fig3} (bottom) indicates the log of the value of $\eta/\rho$. The solid red line indicates the value of the noise level $\eta$. As we expected, these results confirm that Algorithm \ref{Al1} reduces the relative error w.r.t. truth. However, after some number of iterations, this error will start increasing unless the algorithm is stopped. The results from \Fref{Fig3} reveals that there is a critical ensemble size for which, on average (over several experiments with different initial ensembles), the discrepancy principle (\ref{eq:m15}), with $\tau\approx 1/\rho$ is a reliable stopping criteria which terminates the algorithm before the error w.r.t. the truth increases due to data overfitting. In other words, the stability that the scheme inherits from the regularizing LM scheme depends on the ensemble size. For this experiment the critical ensemble size is $N_{e}=125$. In average, for $N_{e}<125$ the error with respect to the truth of the ensemble mean will increase before the data misfit has dropped below $\eta /\rho$. For $N_{e}\ge 125$ we observe that the error will increase when the data misfit has reached the value in $\eta/\rho$ as we would expect if the regularization properties were inherited from the regularizing LM scheme. For larger ensemble sizes $N_{e}\ge 200$ the relative error with respect to the truth will not increase even when the data misfit takes values below $\eta/\rho$ but a slow increase starts showing when the data misfit drops below the noise level. 

\begin{figure}[htbp]
\begin{center}

\includegraphics[scale=0.245]{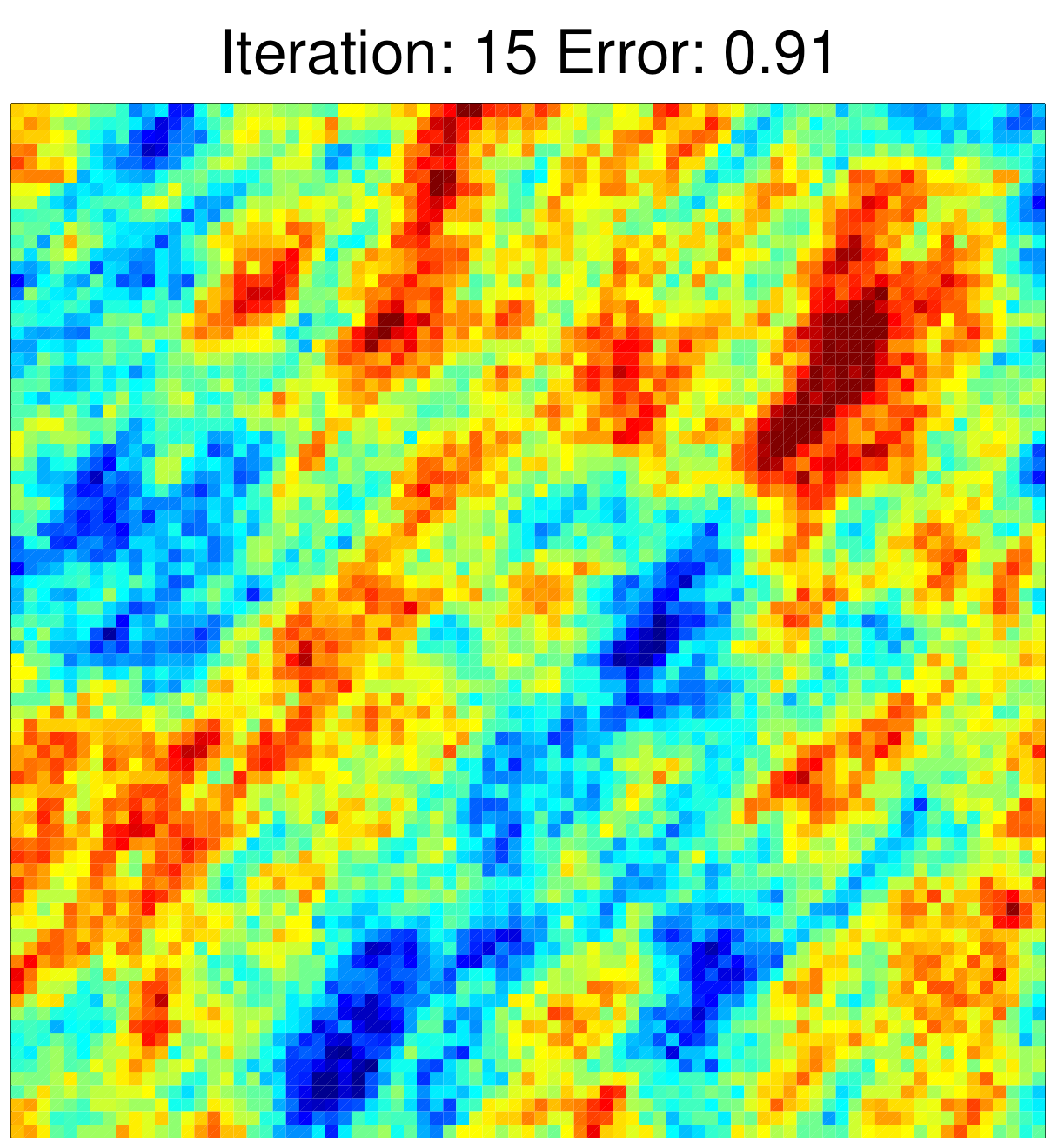}
\includegraphics[scale=0.245]{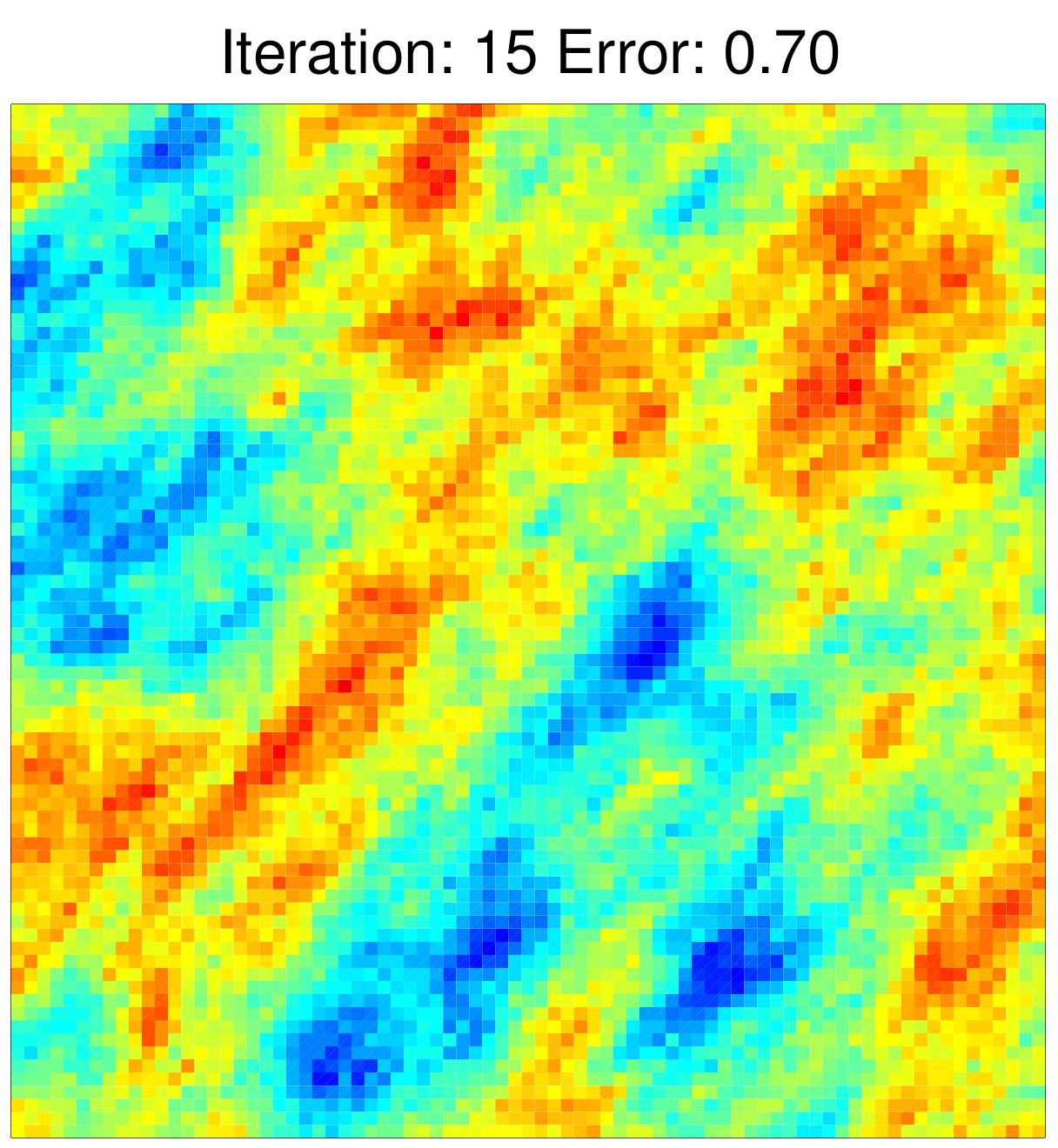}
\includegraphics[scale=0.245]{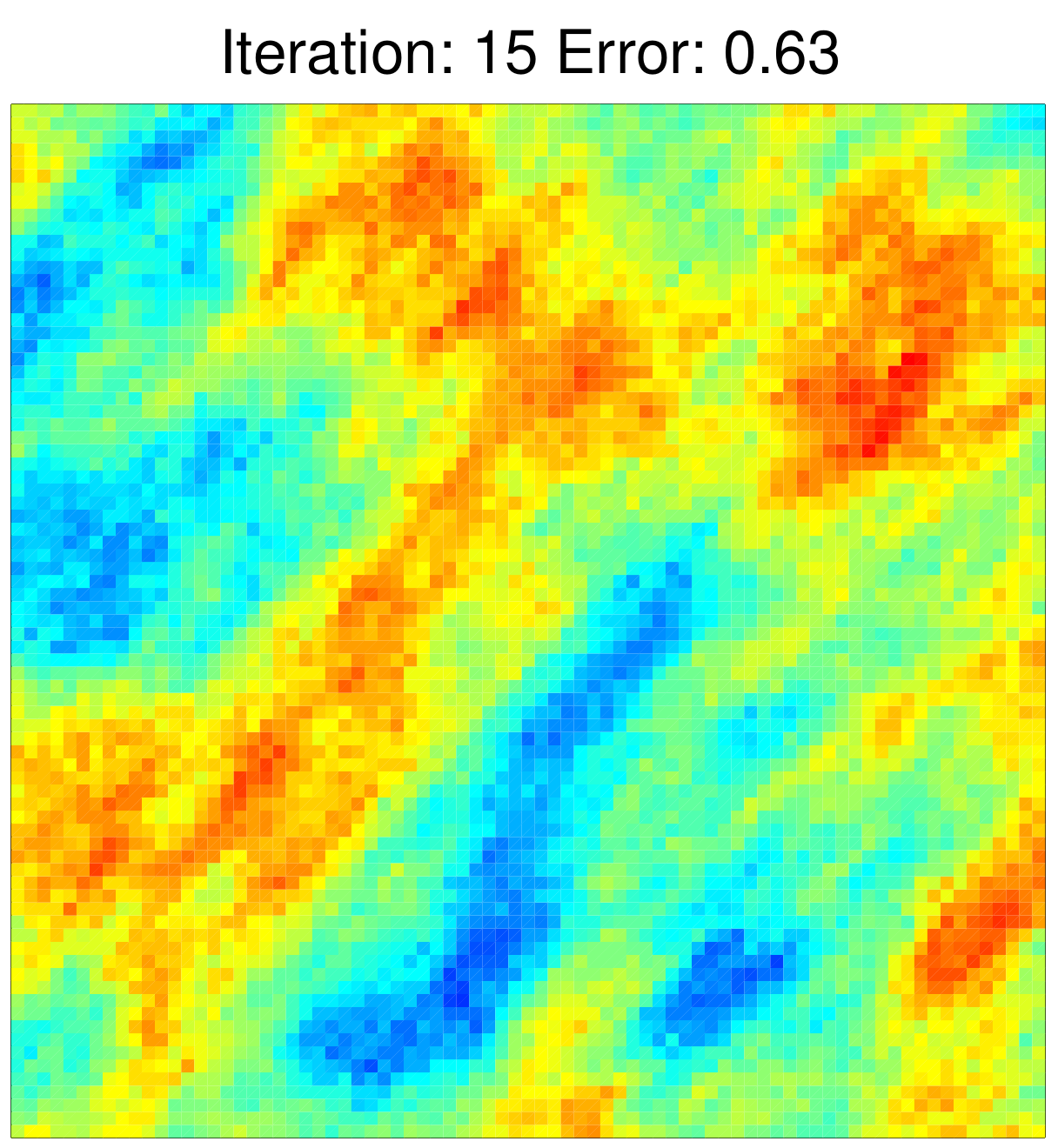}
\includegraphics[scale=0.245]{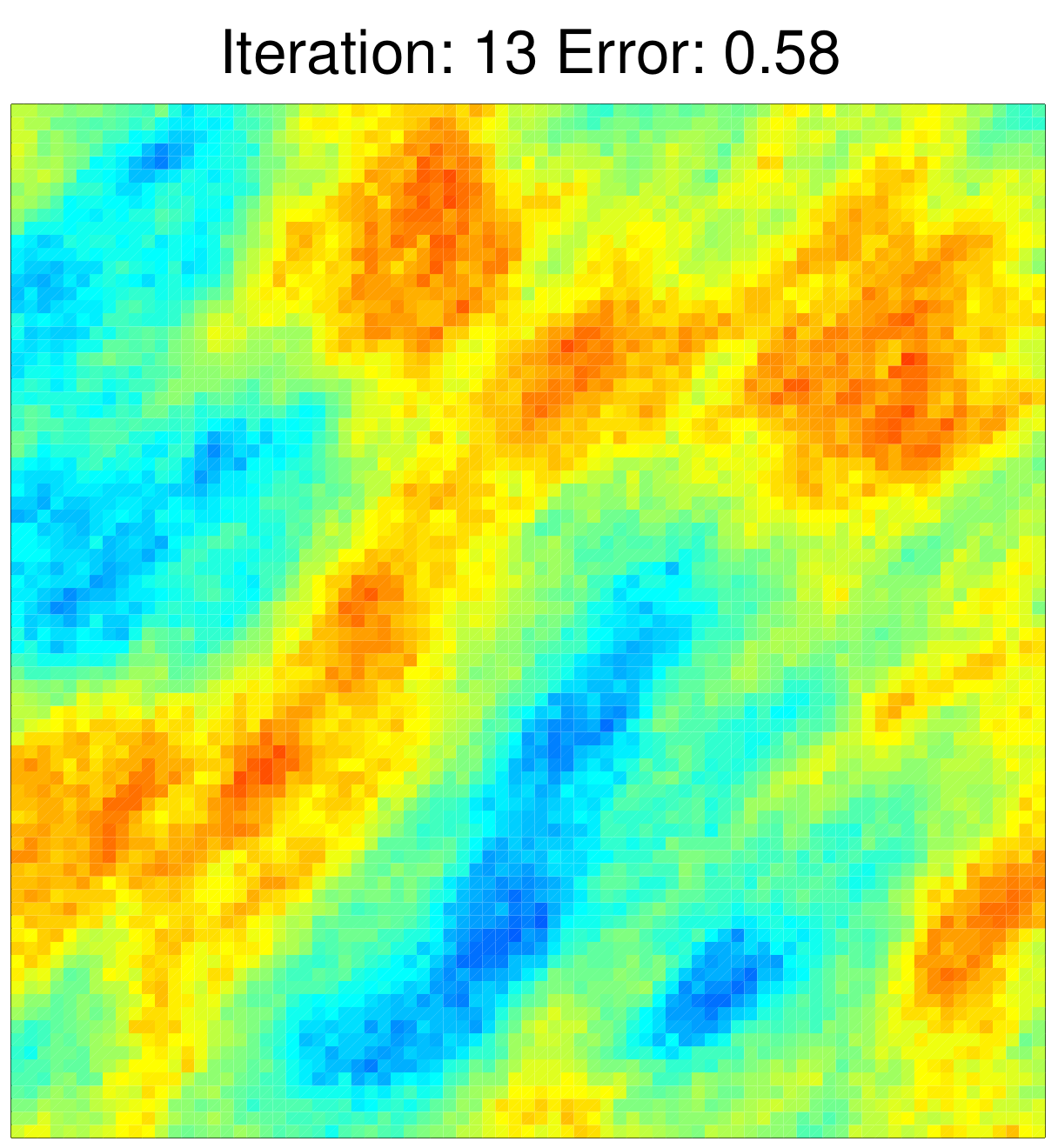}
\includegraphics[scale=0.245]{True_GW}
 \caption{Log - conductivity estimates obtained from an initial ensemble size of (from left to right) 75, 100, 150, 200 (the larger ensembles contains the smaller ones). These estimates are the ensemble mean obtained from Algorithm \ref{Al1} with $\tau=1/\rho$ in (\ref{eq:m15}).}
   \label{Fig4}
\end{center}
\end{figure}

We also note from \Fref{Fig3} that, as we increase the ensemble size, the ensemble mean $\overline{u}_{n}$ provides a more accurate approximation of the truth (provided the scheme is properly stopped). However, a further increase in the ensemble size results in higher computational cost of the scheme. In Figure \ref{Fig4} we show the ensemble mean that we obtain, when Algorithm \ref{Al1} is initialized with different ensemble sizes $N_{e}$ and the scheme stopped according to (\ref{eq:m15}) with $\tau=1/\rho$. For this result, the members of the smaller initial ensembles are contained in the larger ones. We can clearly appreciate that for smaller ensemble sizes, the stabilization of the scheme may be lost and thus its accuracy.


\subsubsection{Effect of the number of measurements.}

We now investigate the relation between the number of measurements and the critical size that enable us to observe the regularization properties inherited by the regularizing LM scheme. As discussed in subsection \ref{unreg} stability issues for $N_{e}<M$ may arise unless the parameter $\alpha_{n}$ controls the inversion in (\ref{eq:m16}). We recall that the selection of $\alpha_{n}$ according to (\ref{eq:m12}) is inspired by the regularizing LM scheme which addresses a different type of ill-posedness (although is reflected in a similar fashion). Thus, this selection of $\alpha_{n}$ may not necessarily address the stability issues that arise in the ensemble method when $N_{e}<M$. In the following experiments we show that under some conditions, the selection of $\alpha_{n}$ and the stopping criteria does indeed resolve the stability issues even when $N_{e}<M$.

\begin{figure}[htbp]
\begin{center}
\includegraphics[scale=0.32]{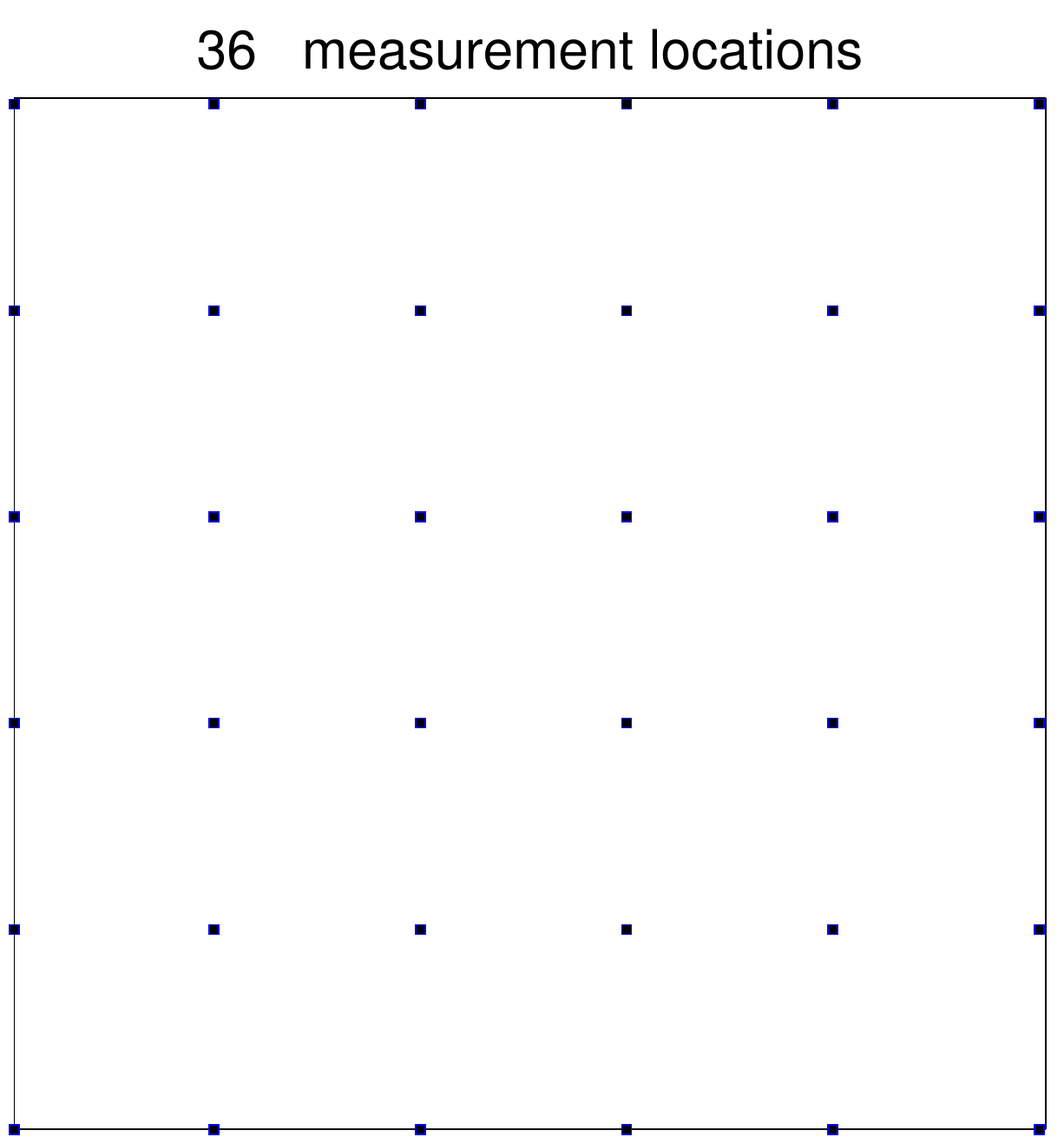}
\includegraphics[scale=0.3]{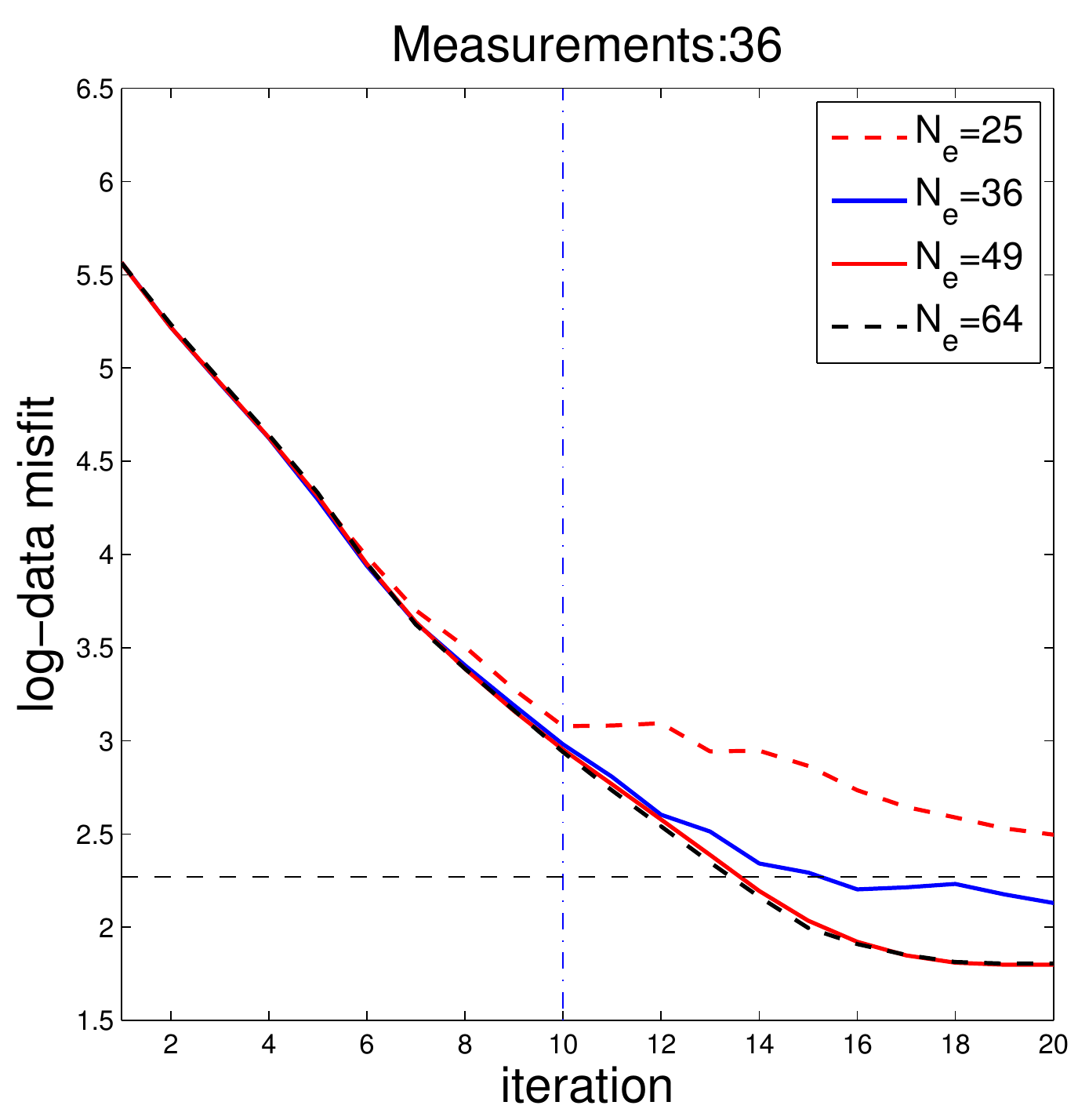}
\includegraphics[scale=0.3]{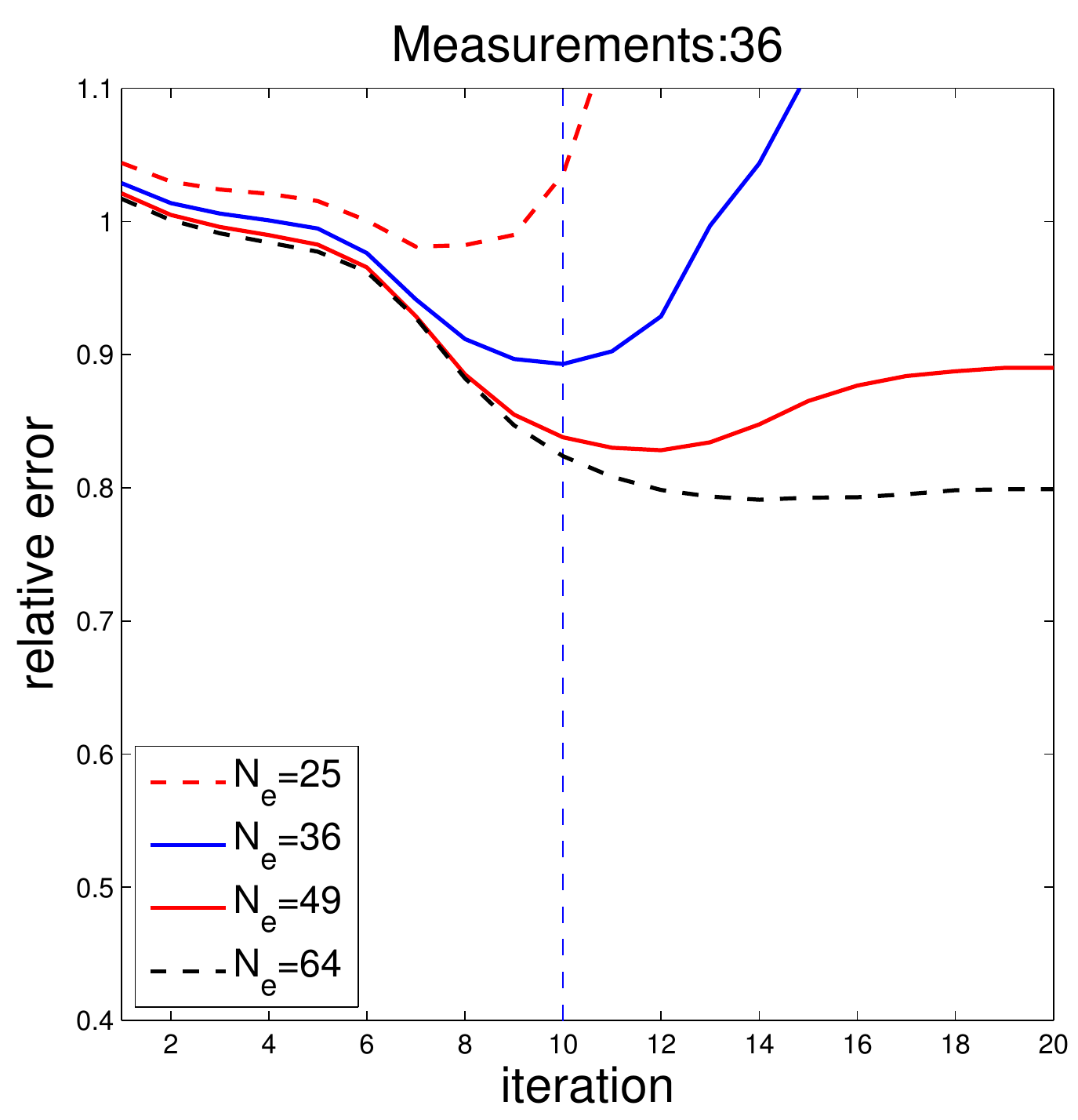}\\
\includegraphics[scale=0.32]{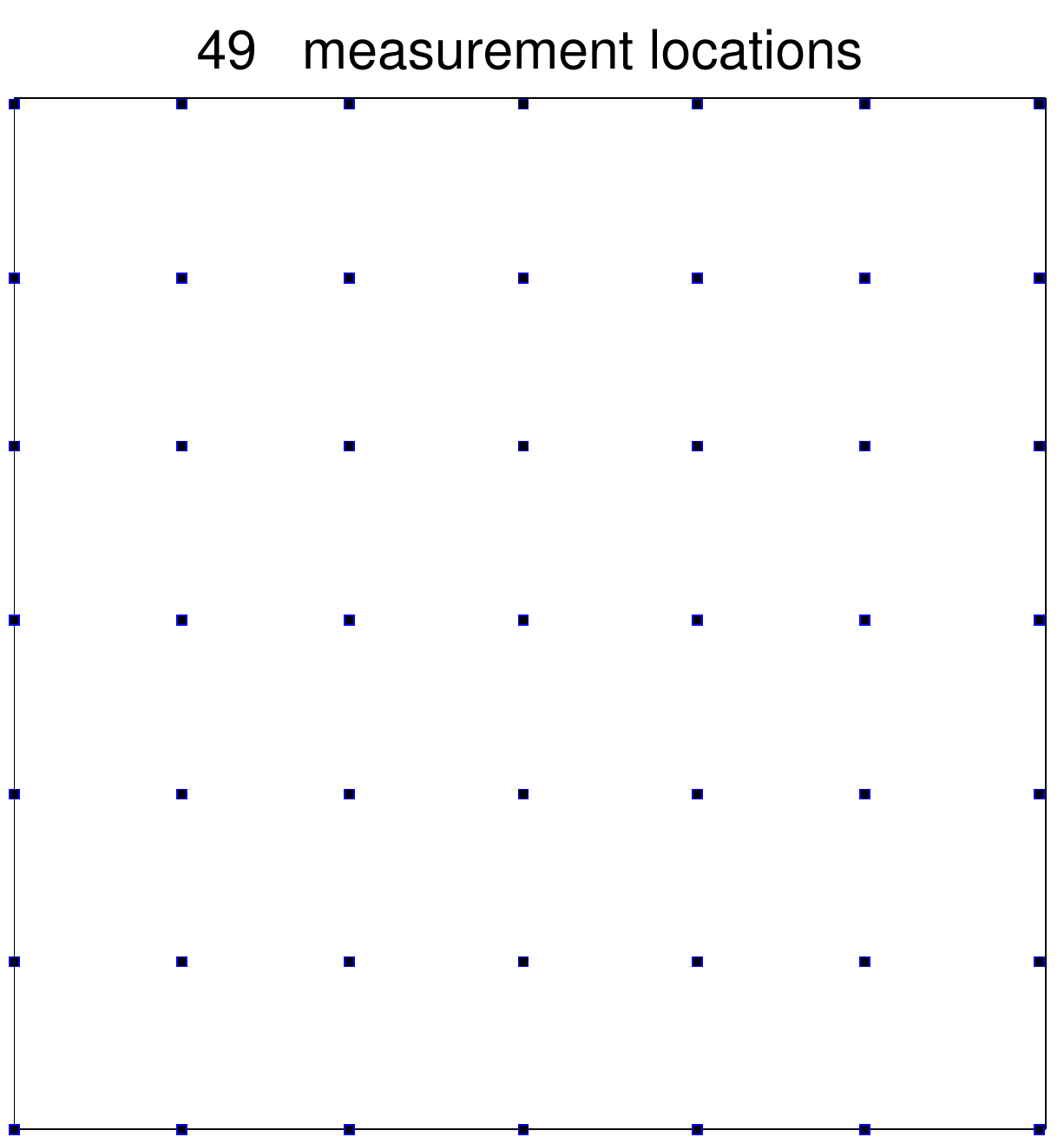}
\includegraphics[scale=0.3]{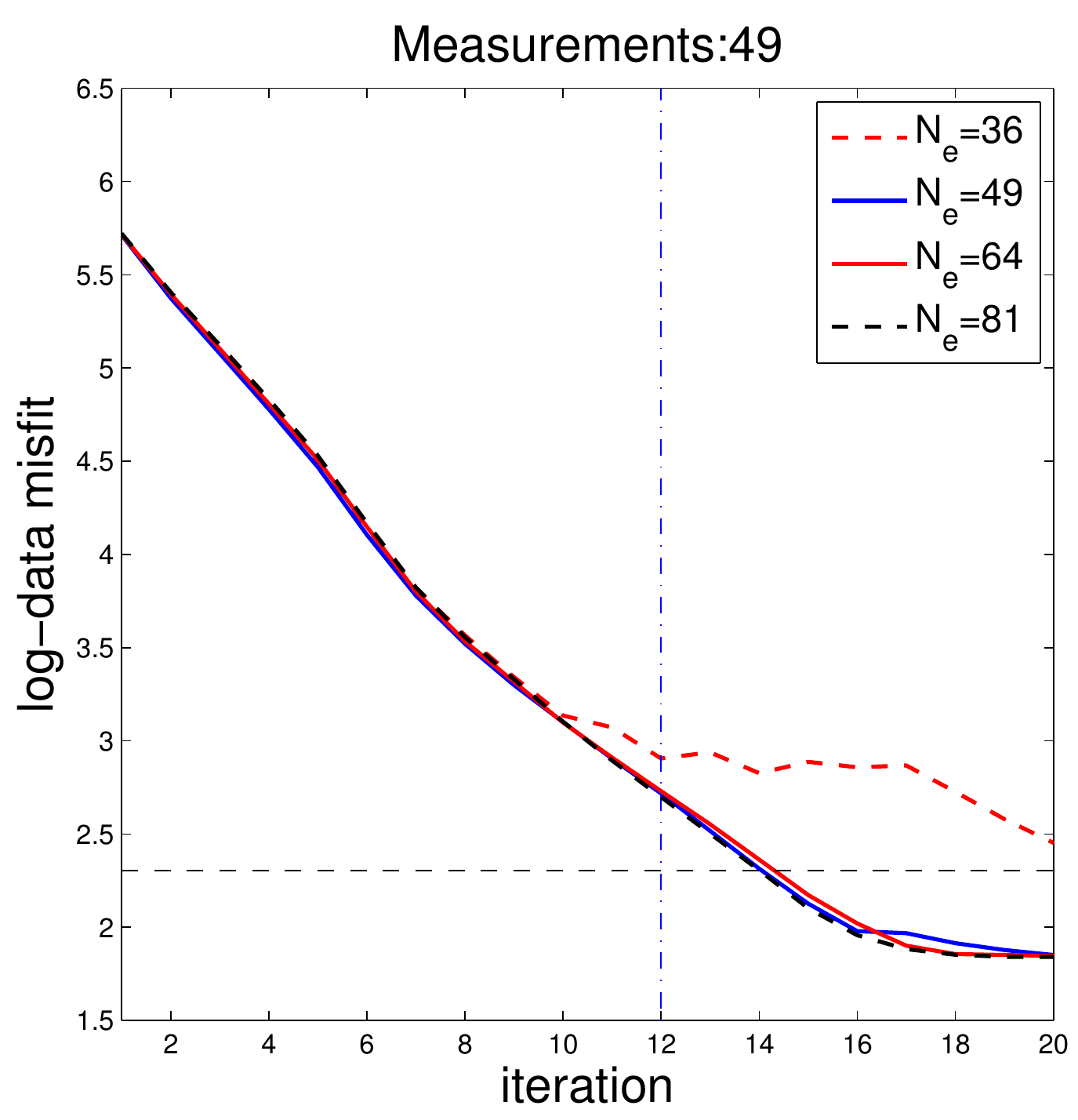}
\includegraphics[scale=0.3]{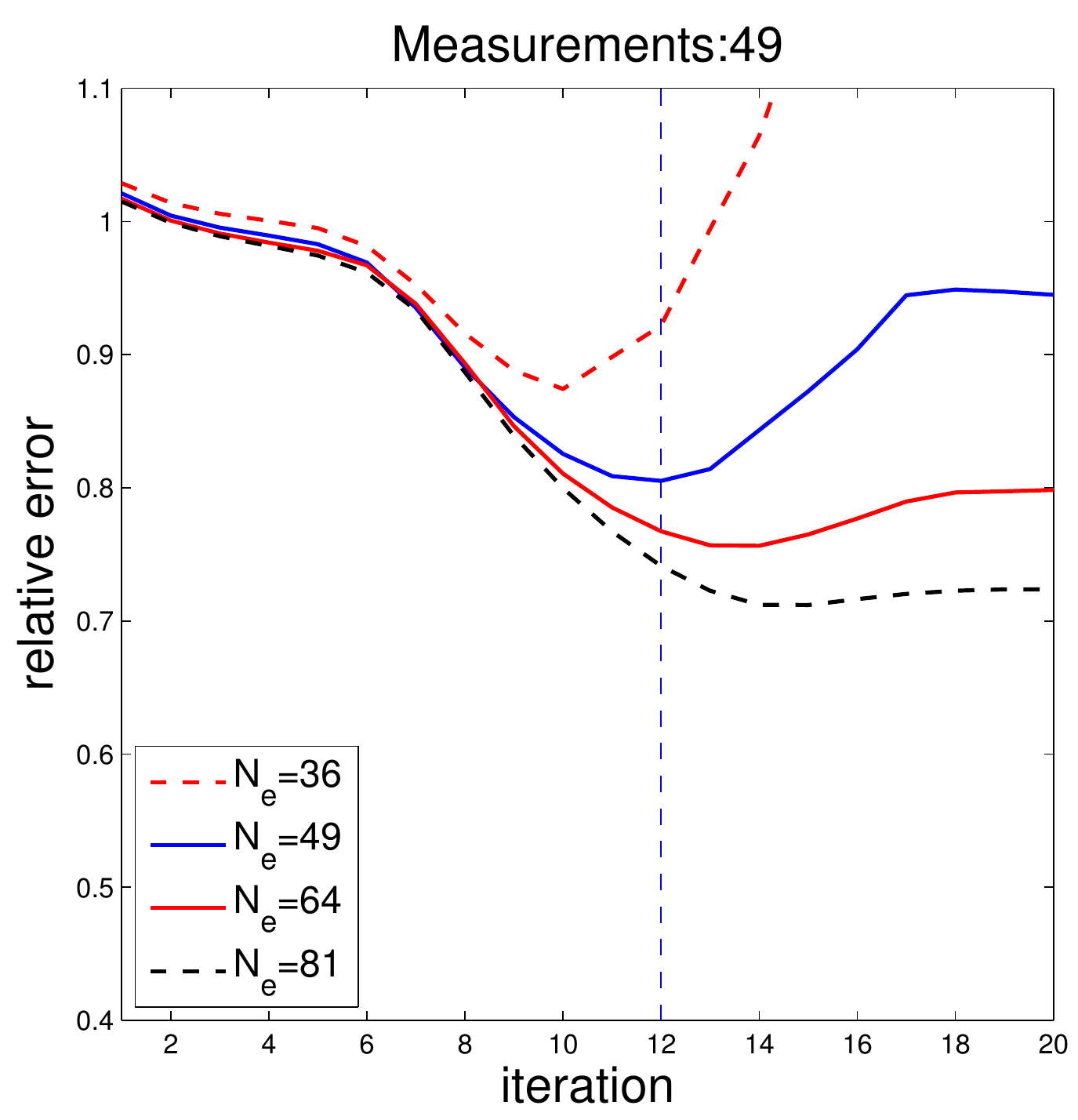}\\
\includegraphics[scale=0.32]{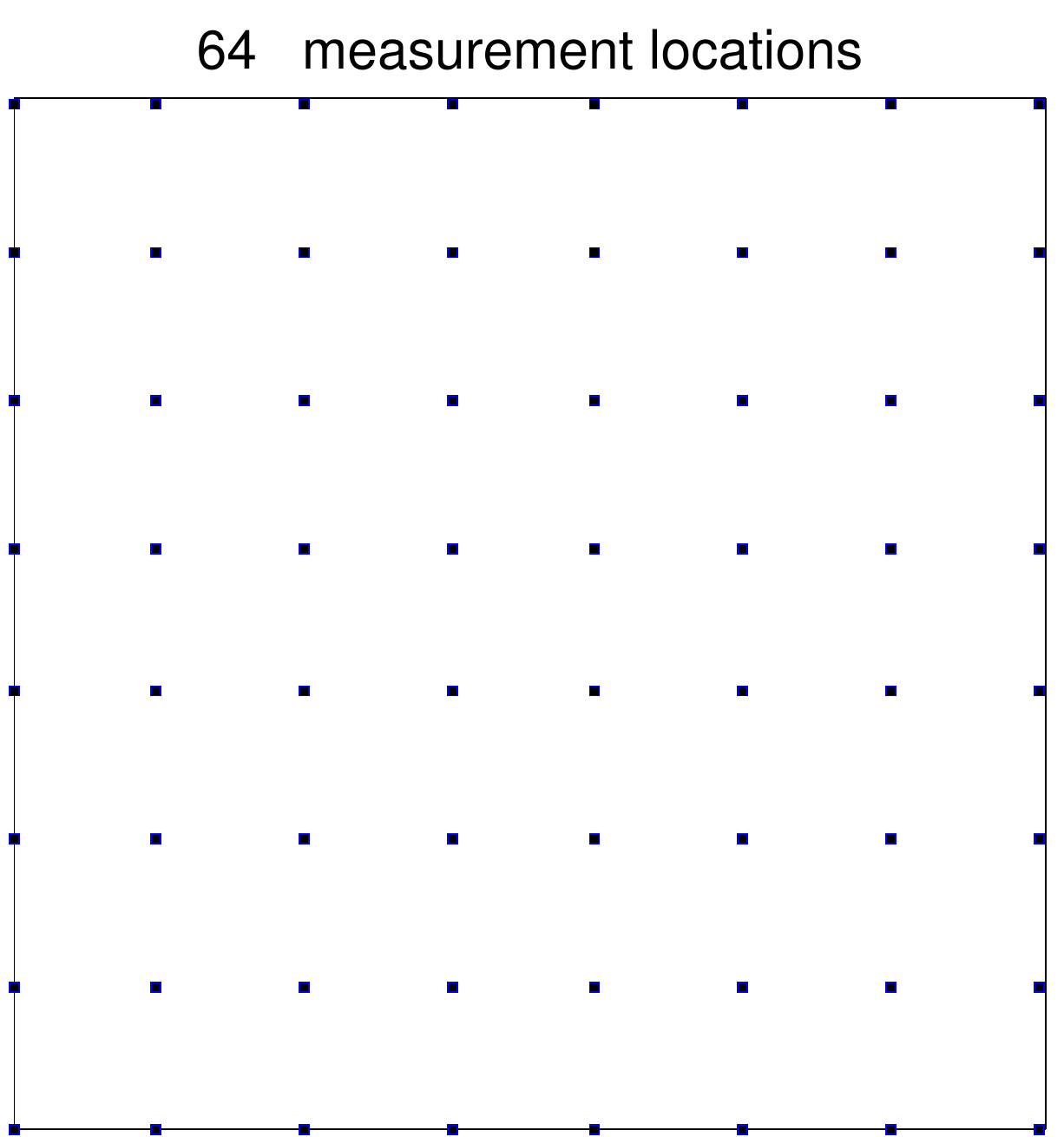}
\includegraphics[scale=0.3]{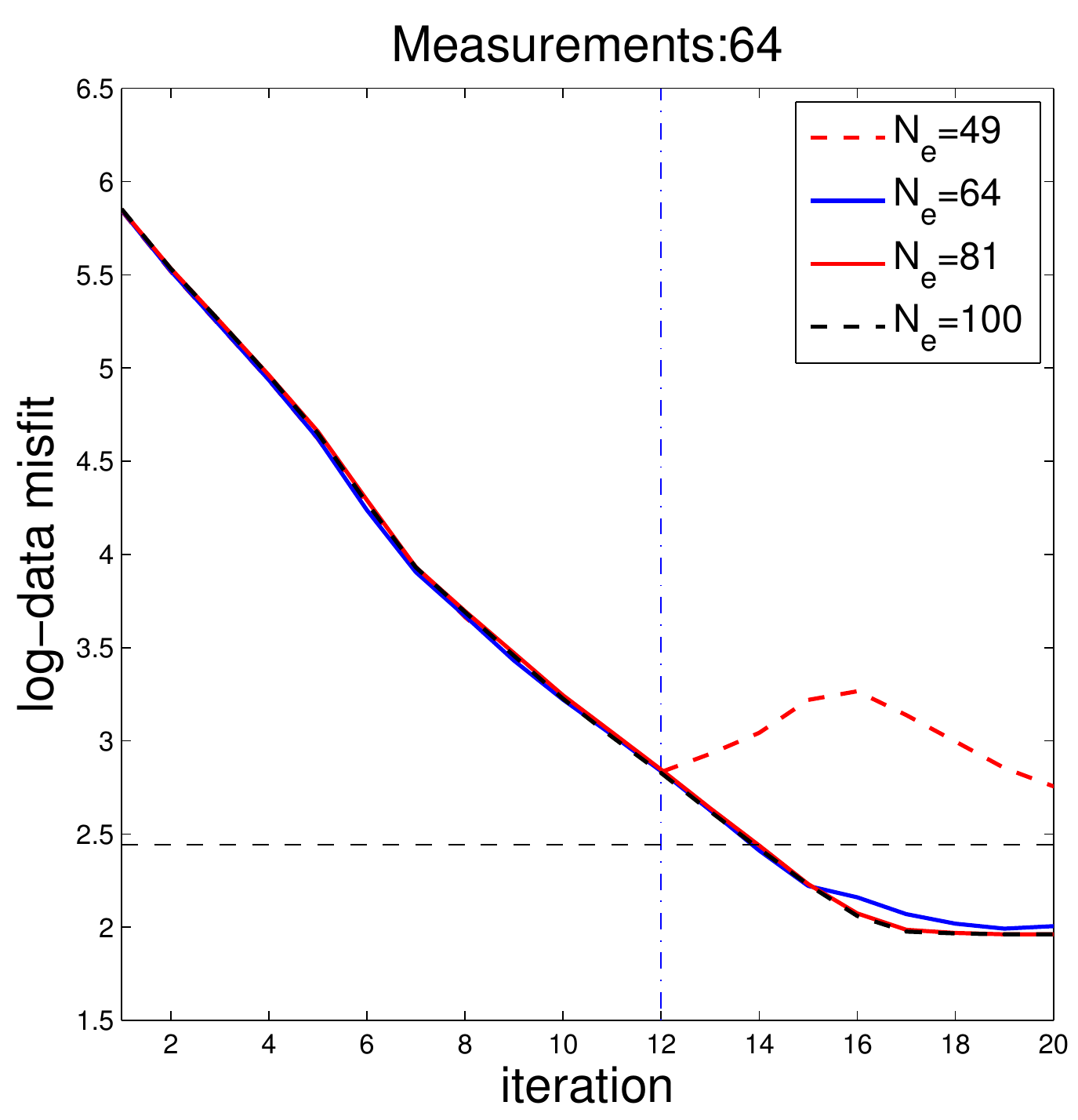}
\includegraphics[scale=0.3]{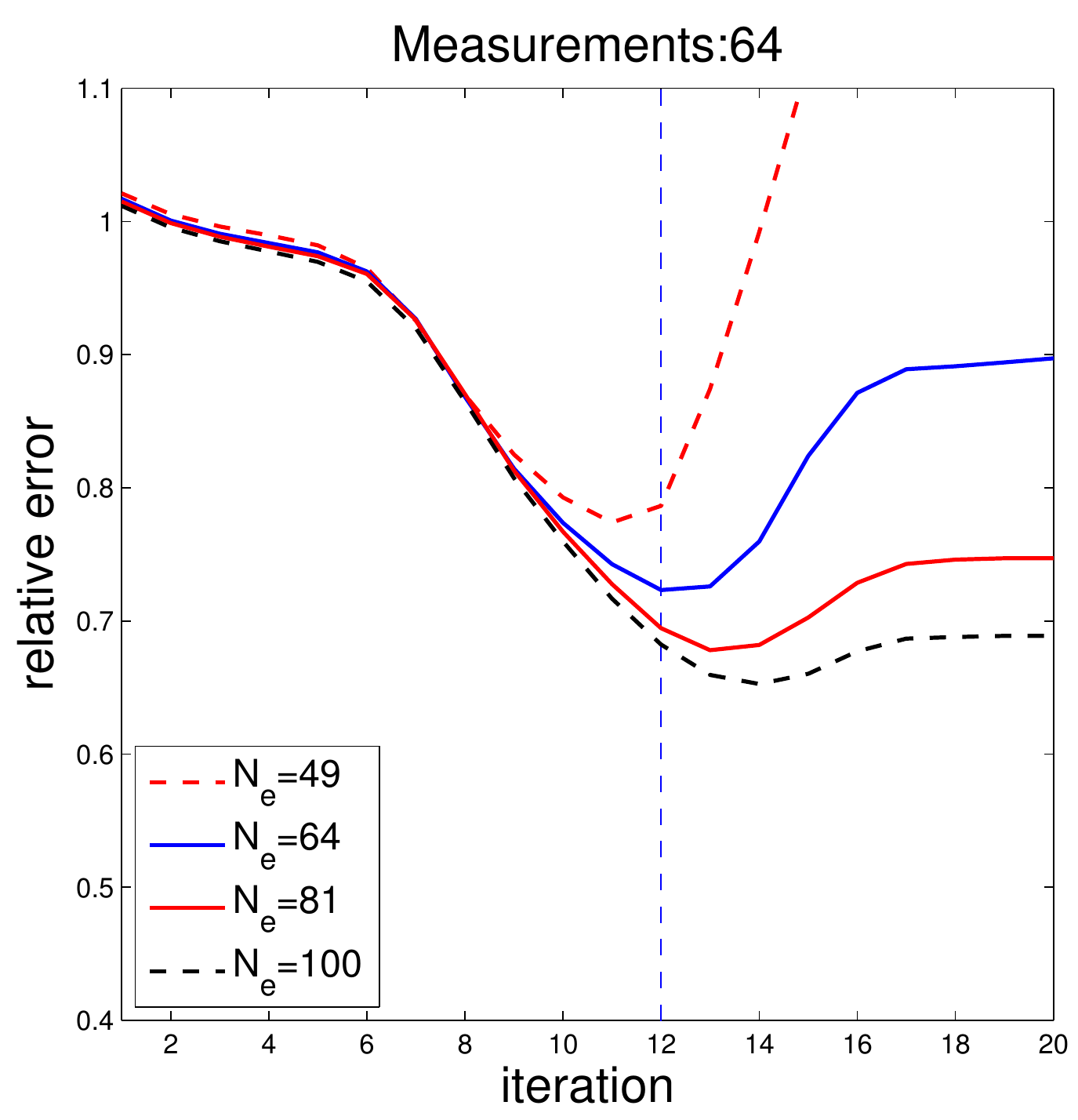}\\
\includegraphics[scale=0.32]{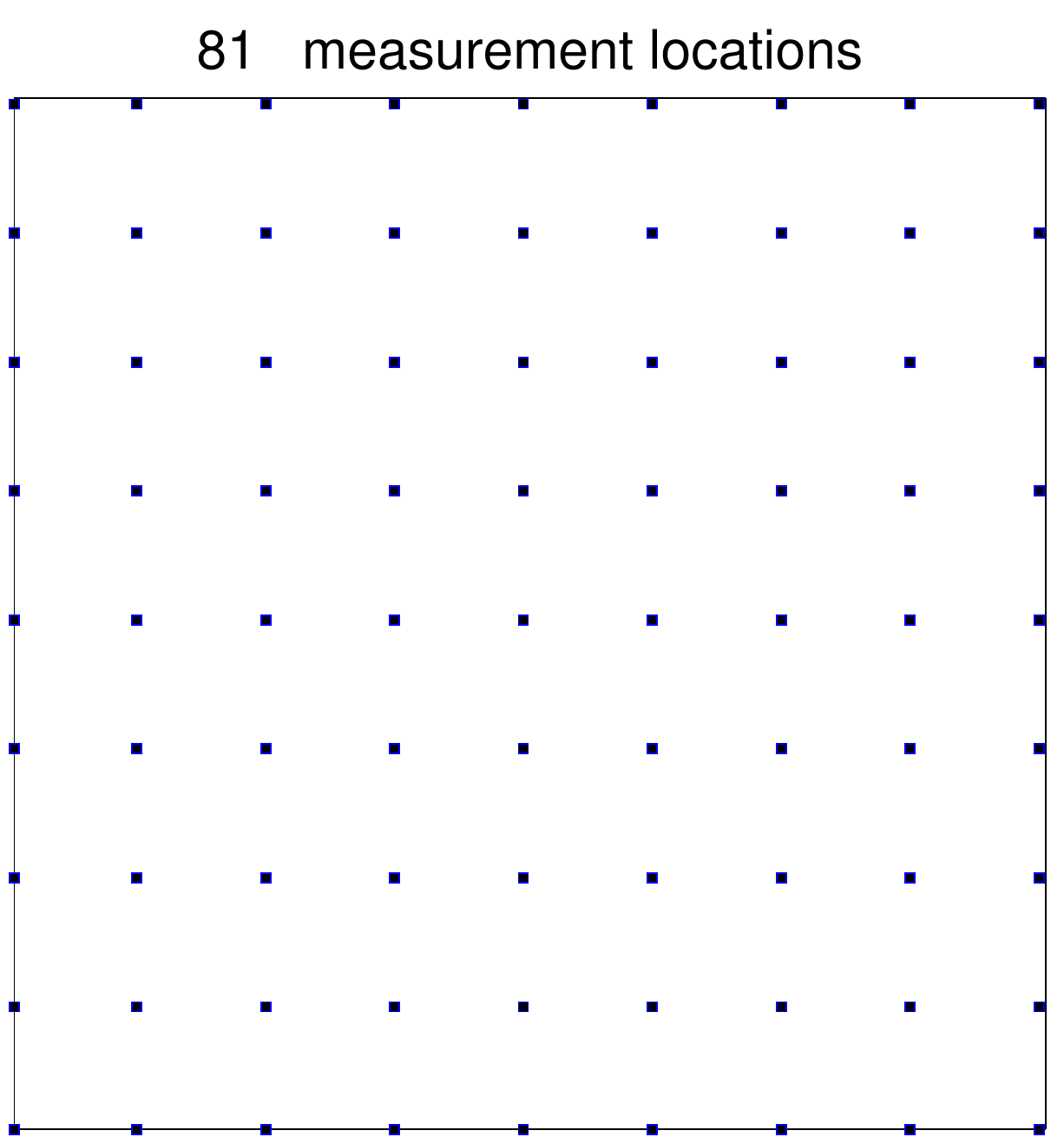}
\includegraphics[scale=0.3]{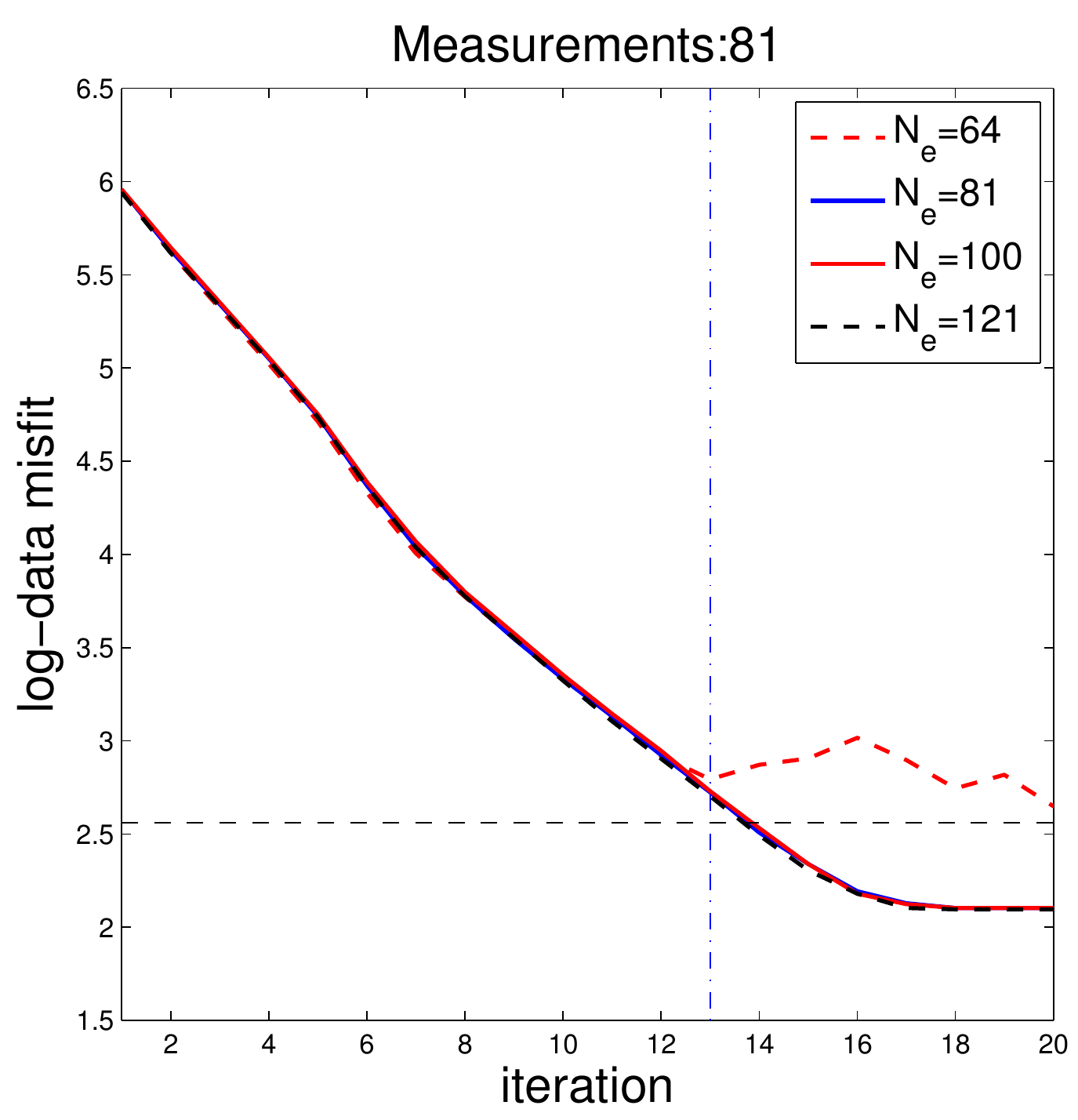}
\includegraphics[scale=0.3]{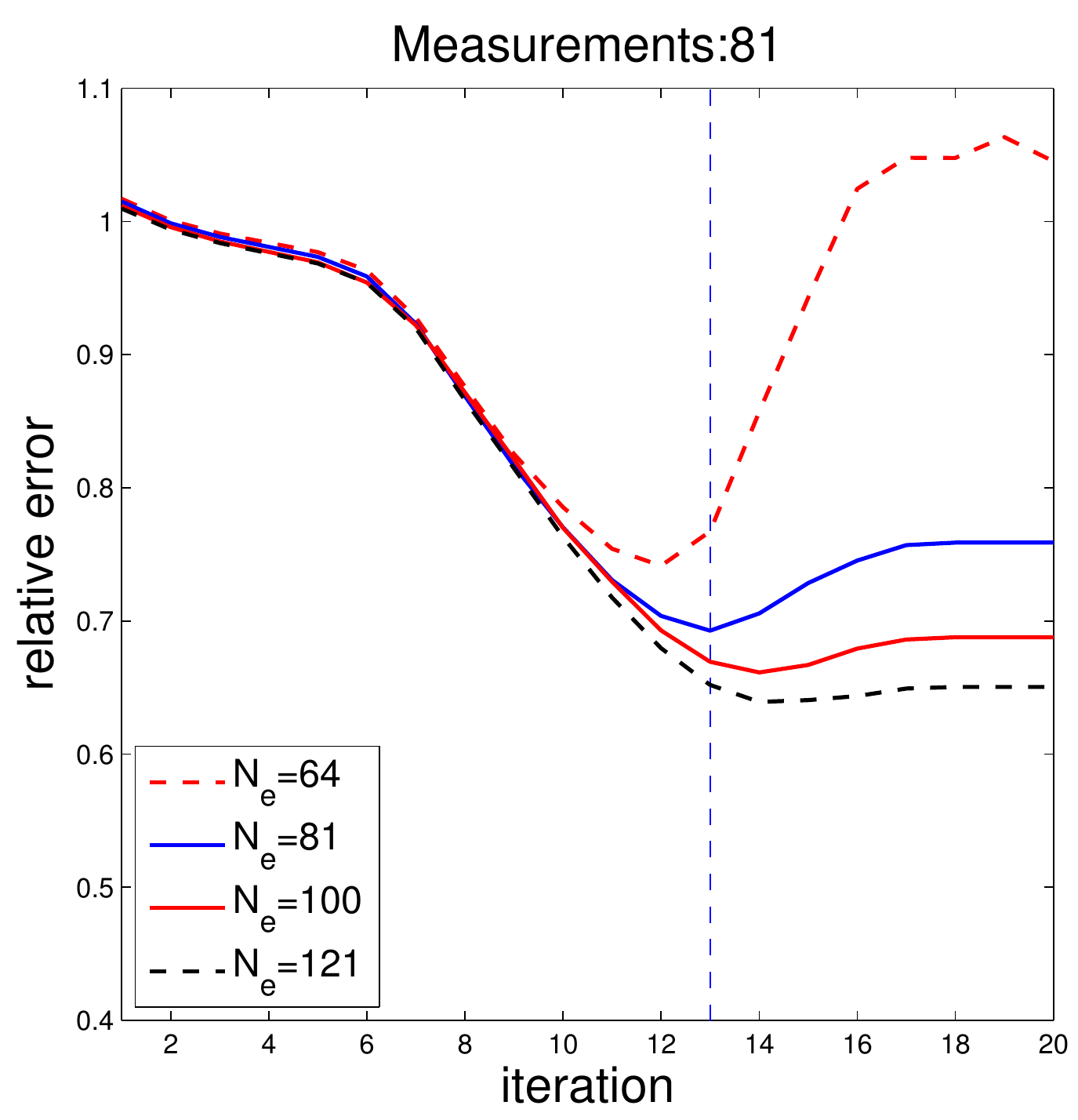}
\includegraphics[scale=0.32]{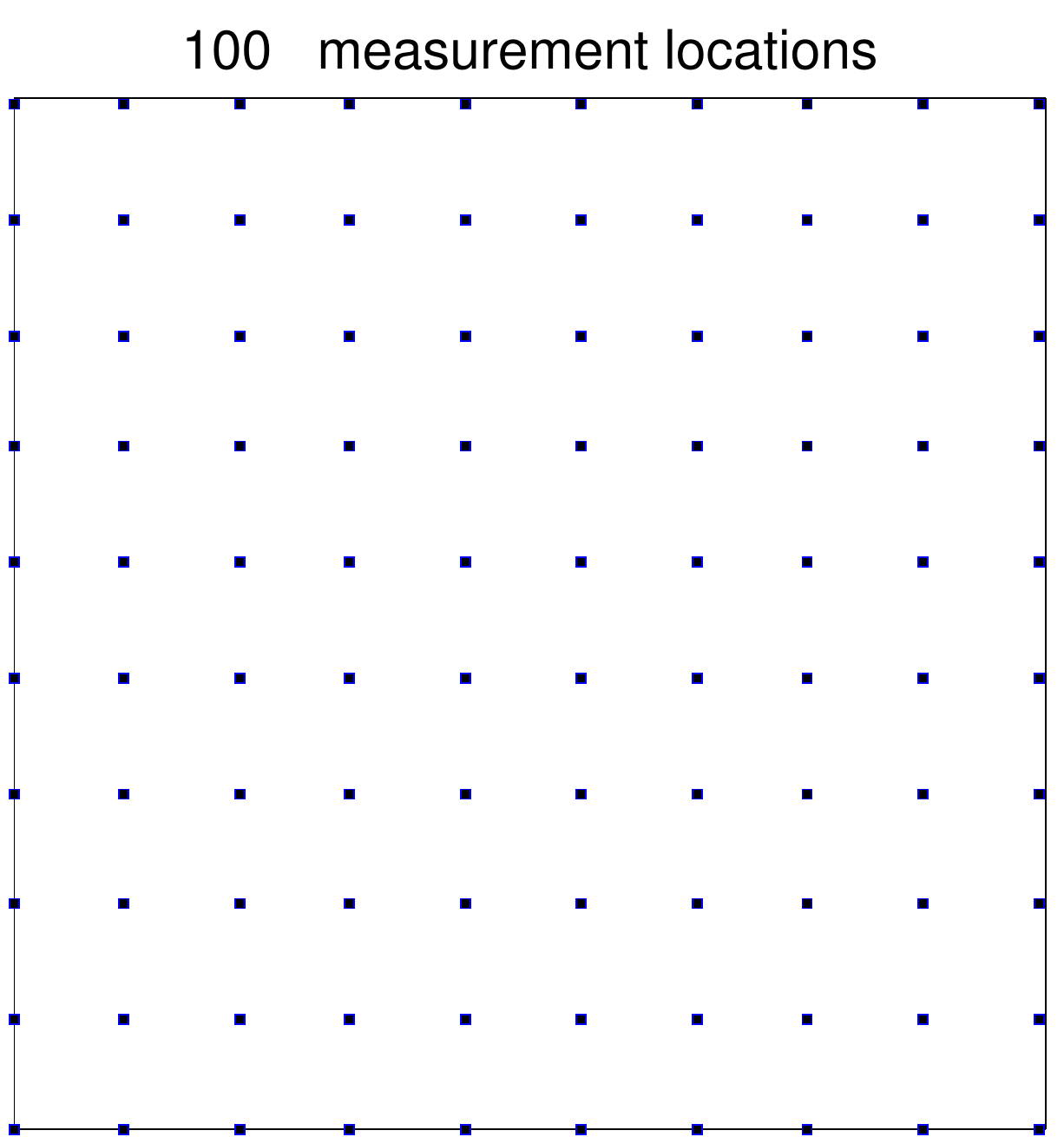}
\includegraphics[scale=0.3]{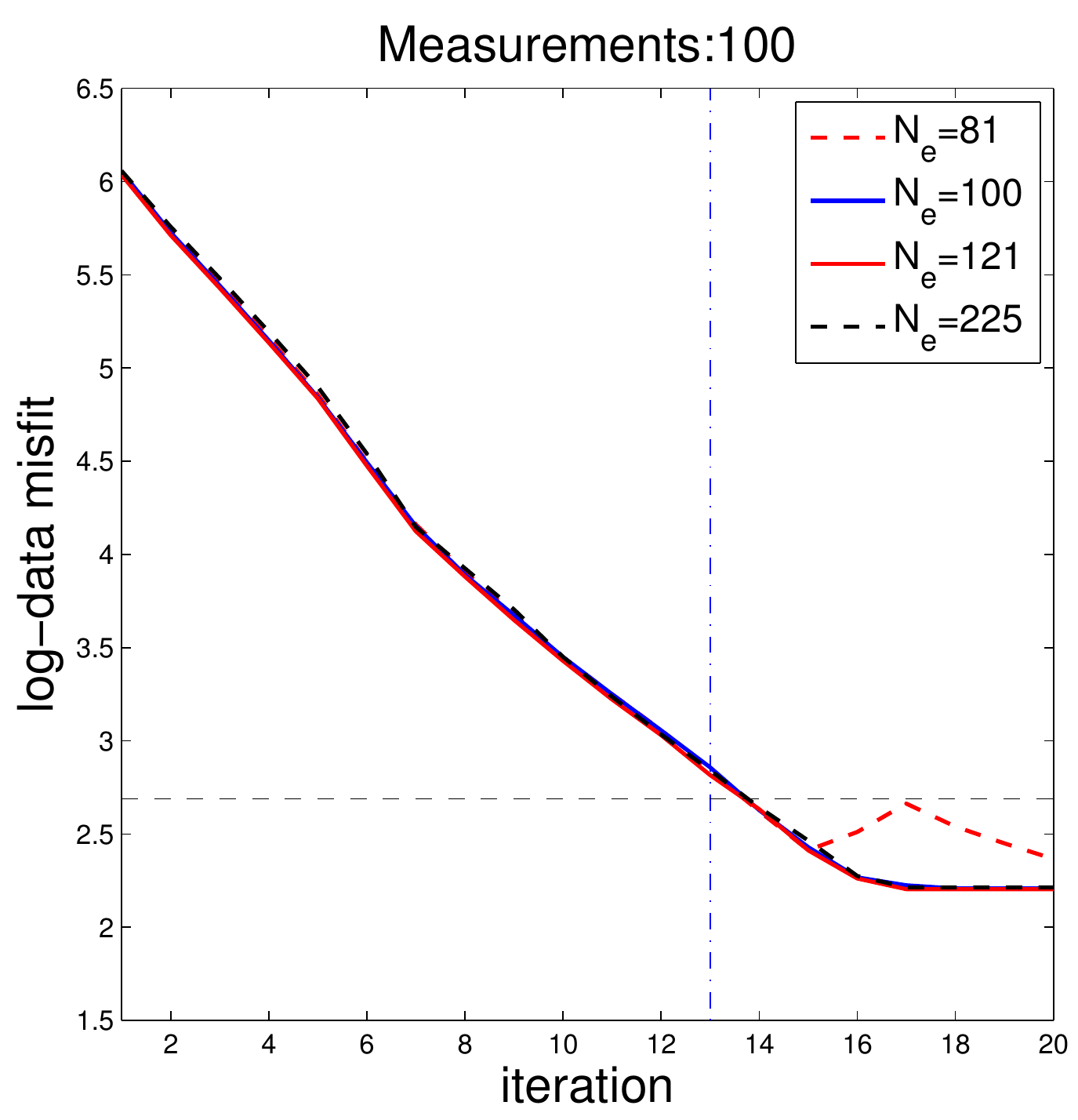}
\includegraphics[scale=0.3]{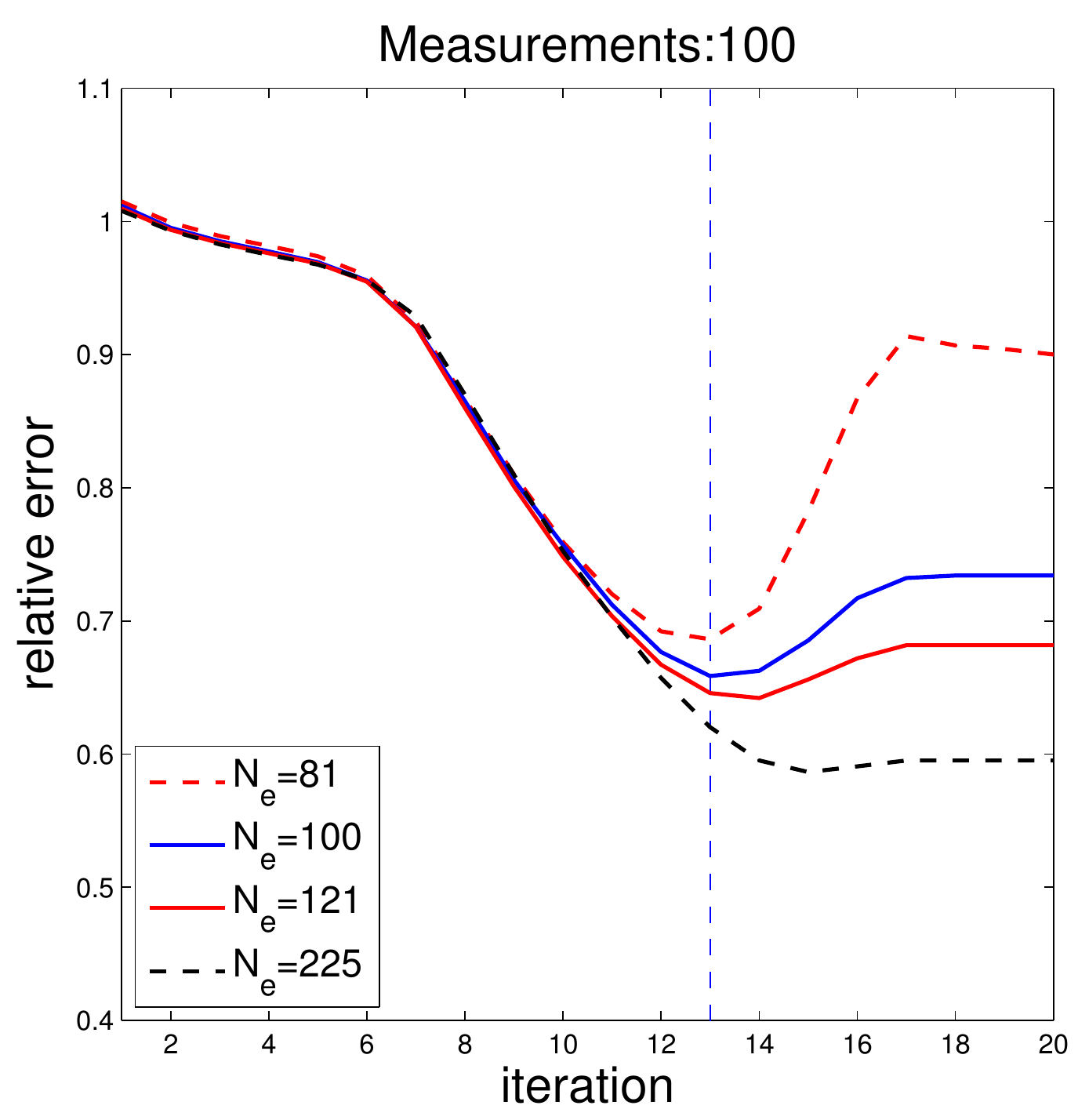}\\  \vskip2pt

 \ \caption{Numerical results with the proposed method for $\rho=0.7$ and different measurement configurations and ensemble sizes. Left column: measurement configuration. Middle column: relative error w.r.t. the truth. Right column: Log - data misfit.   Quantities are averaged at each iteration, over 40 experiments corresponding to different experiment.}
   \label{Fig5}
\end{center}
\end{figure}

\begin{figure}[htbp]
\begin{center}
\includegraphics[scale=0.32]{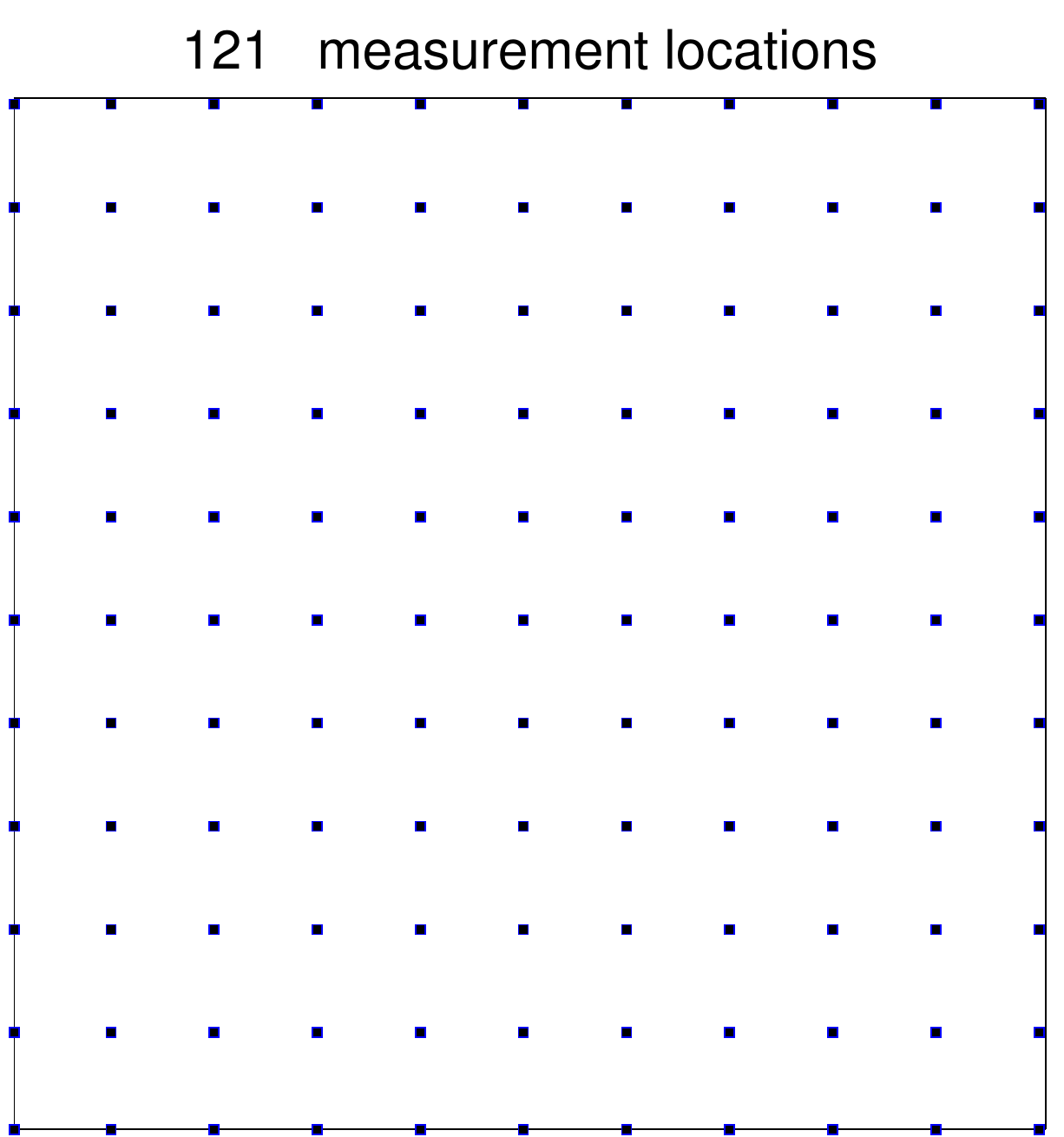}
\includegraphics[scale=0.3]{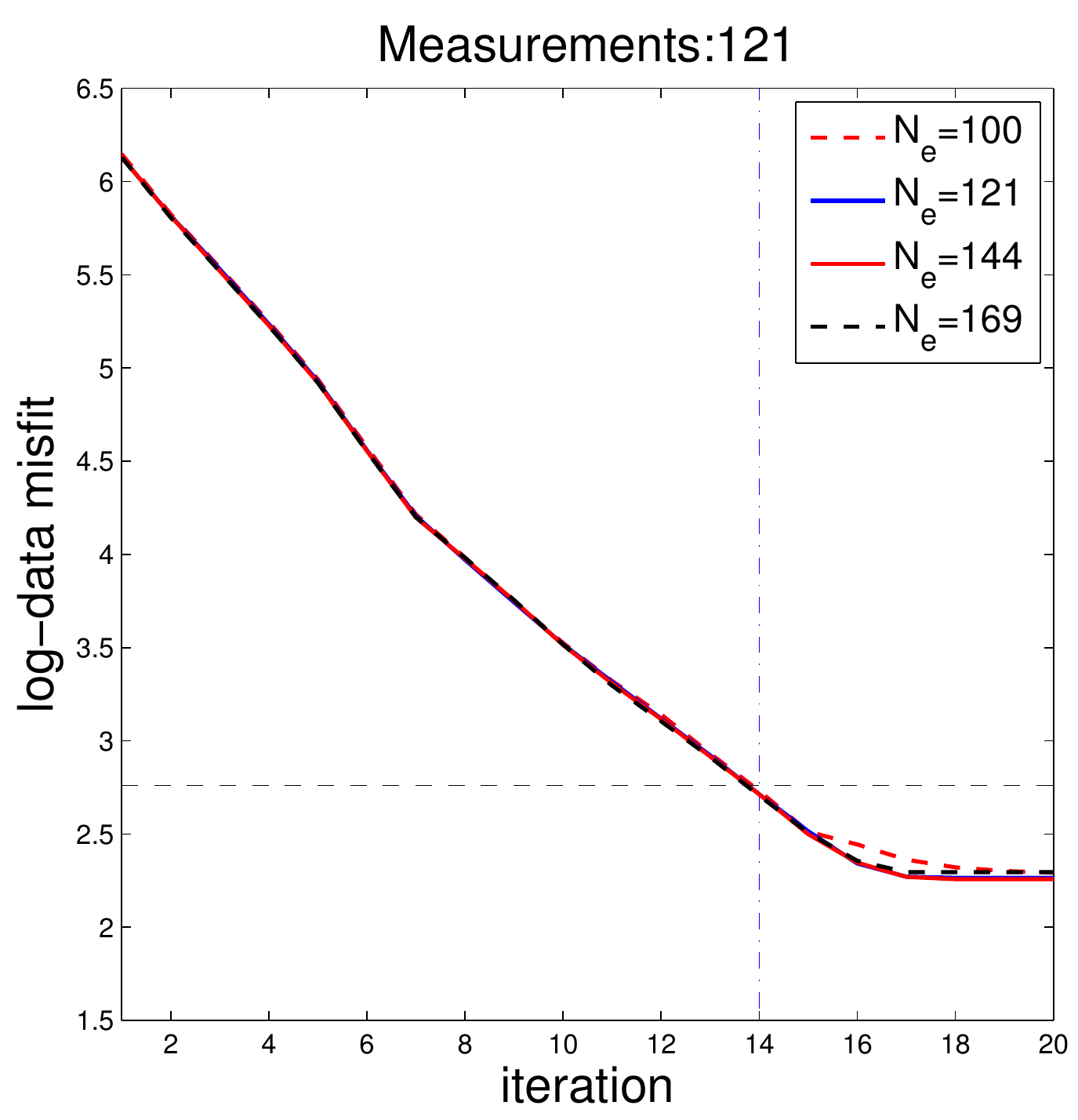}
\includegraphics[scale=0.3]{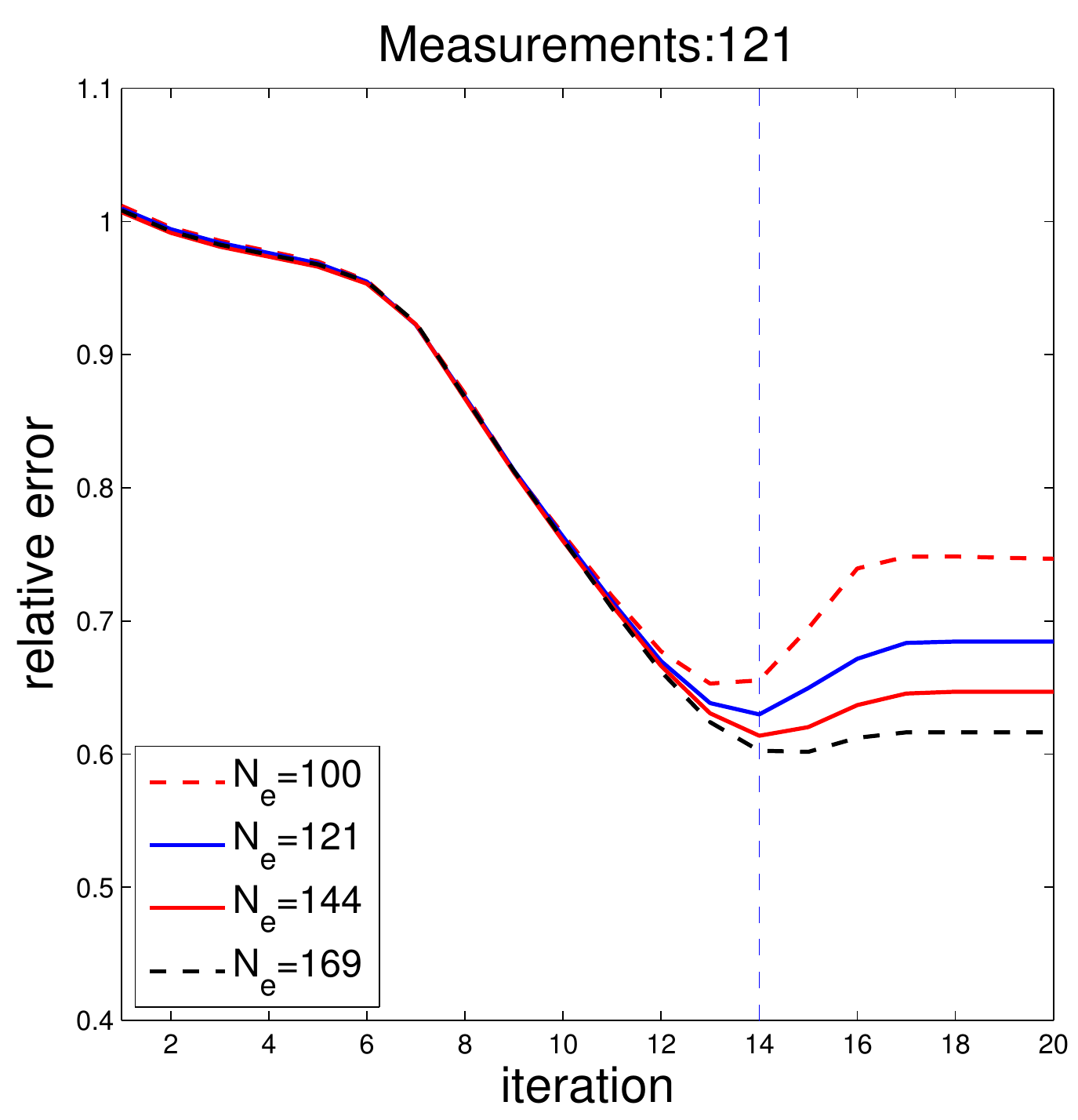}\\  \vskip2pt
\includegraphics[scale=0.32]{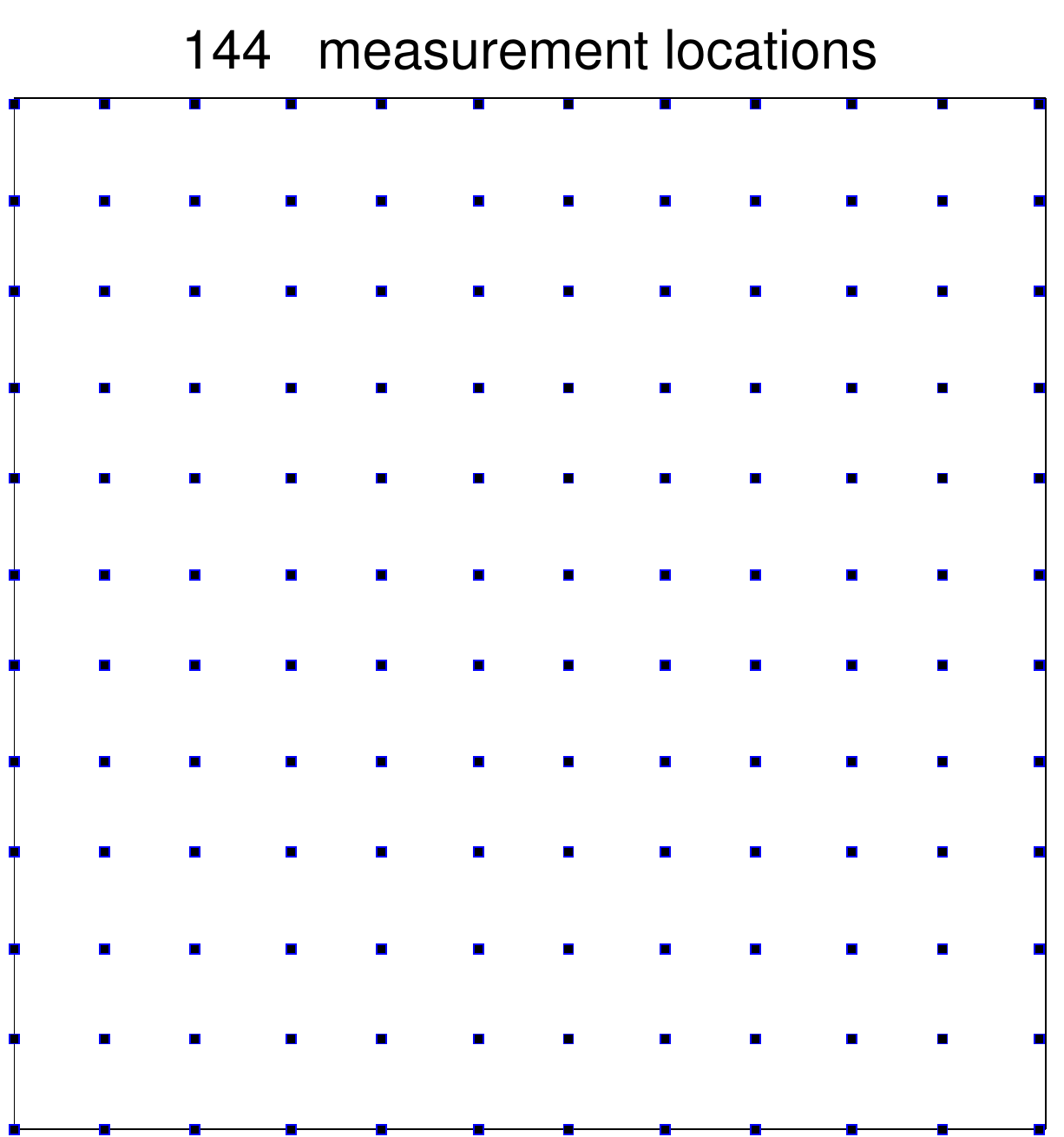}
\includegraphics[scale=0.3]{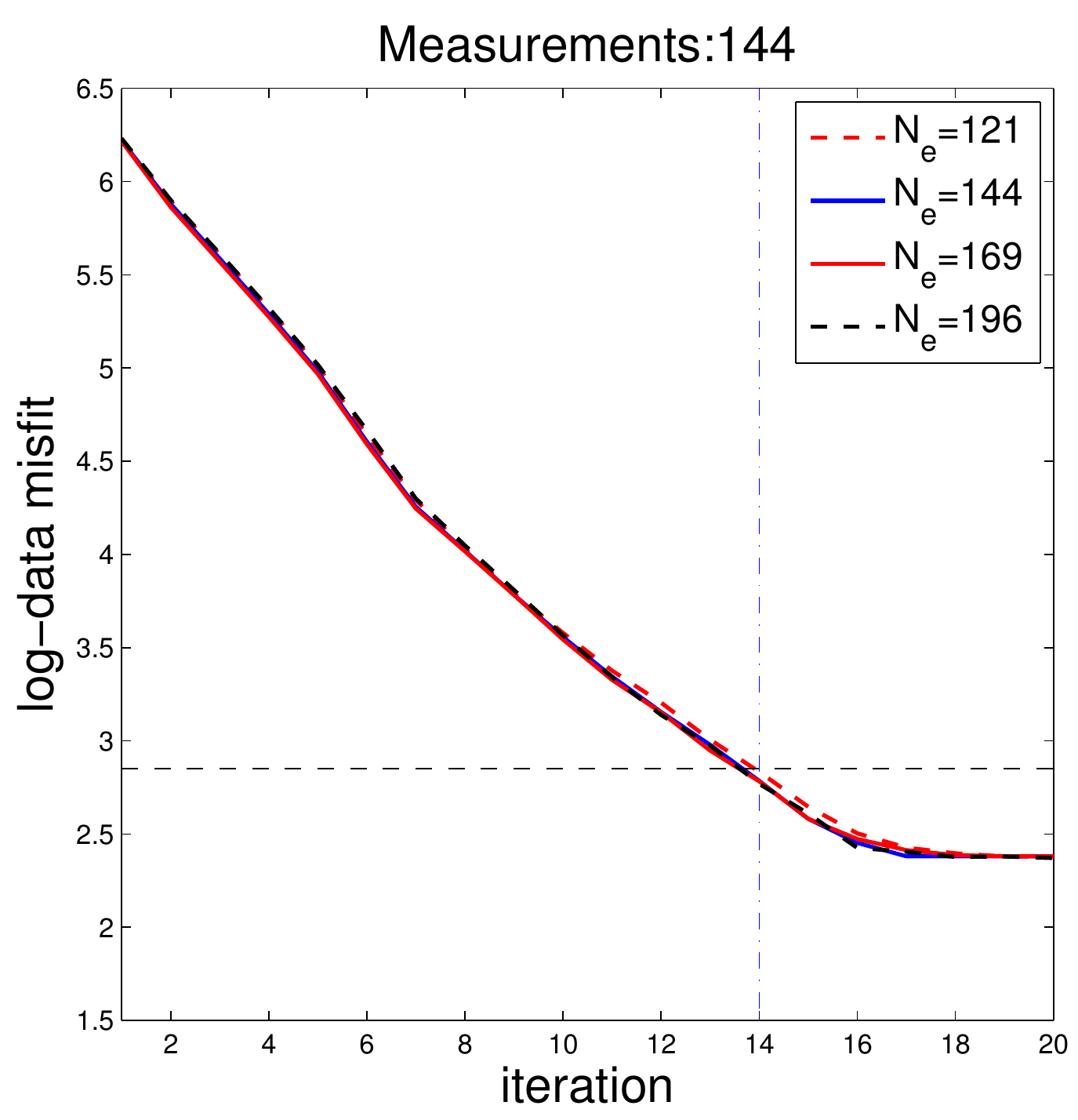}
\includegraphics[scale=0.3]{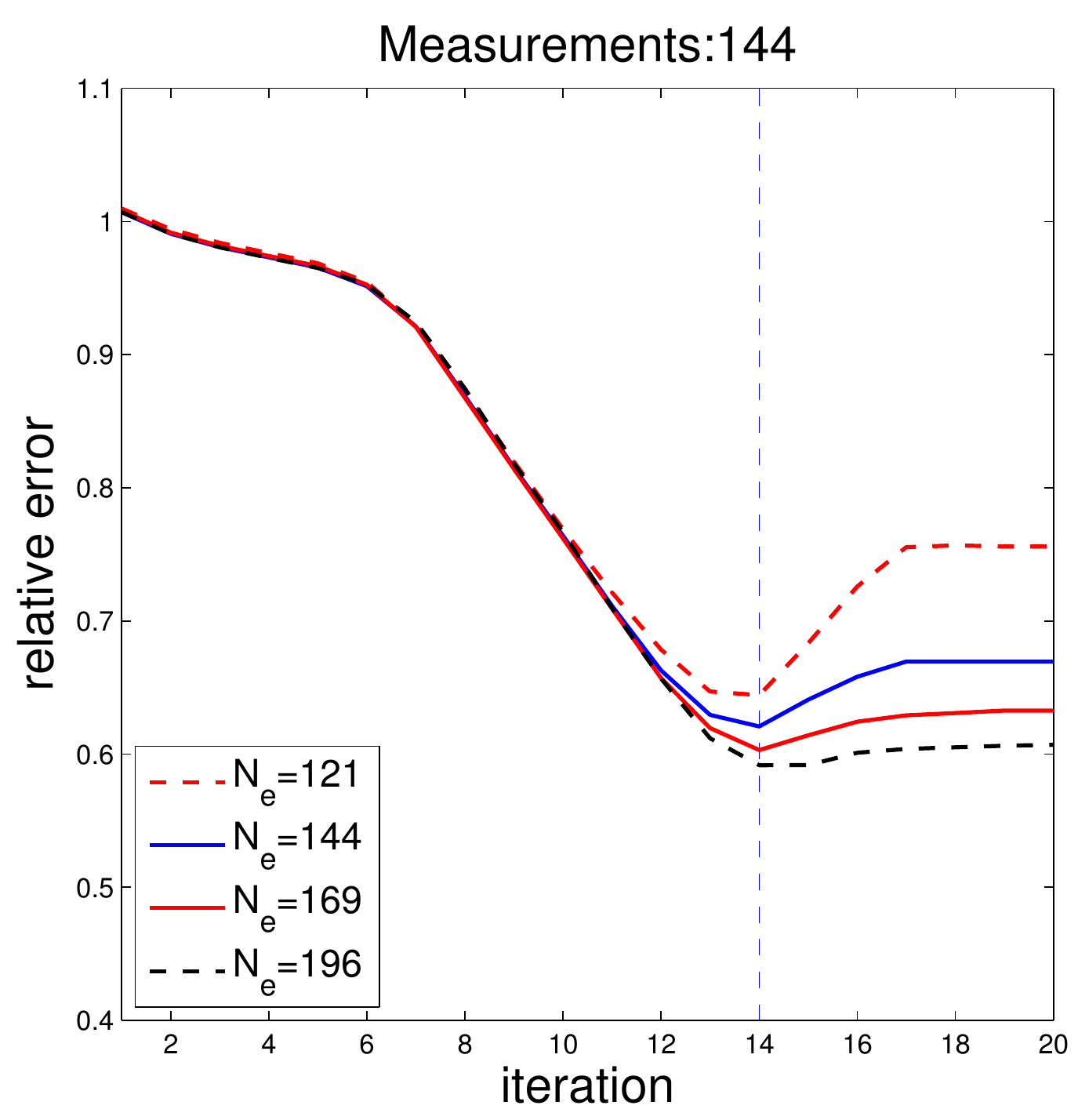}\\
\includegraphics[scale=0.32]{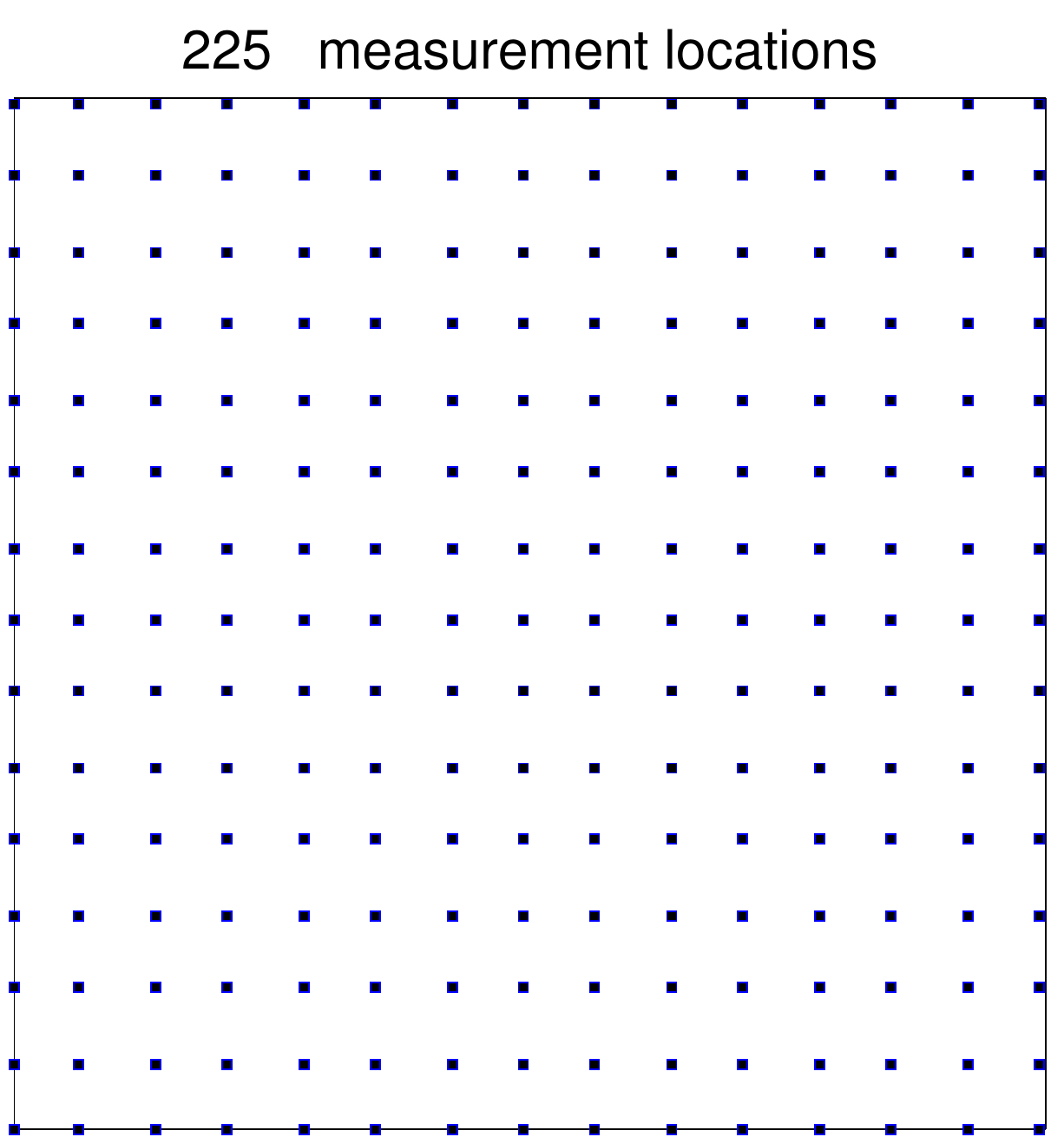}
\includegraphics[scale=0.3]{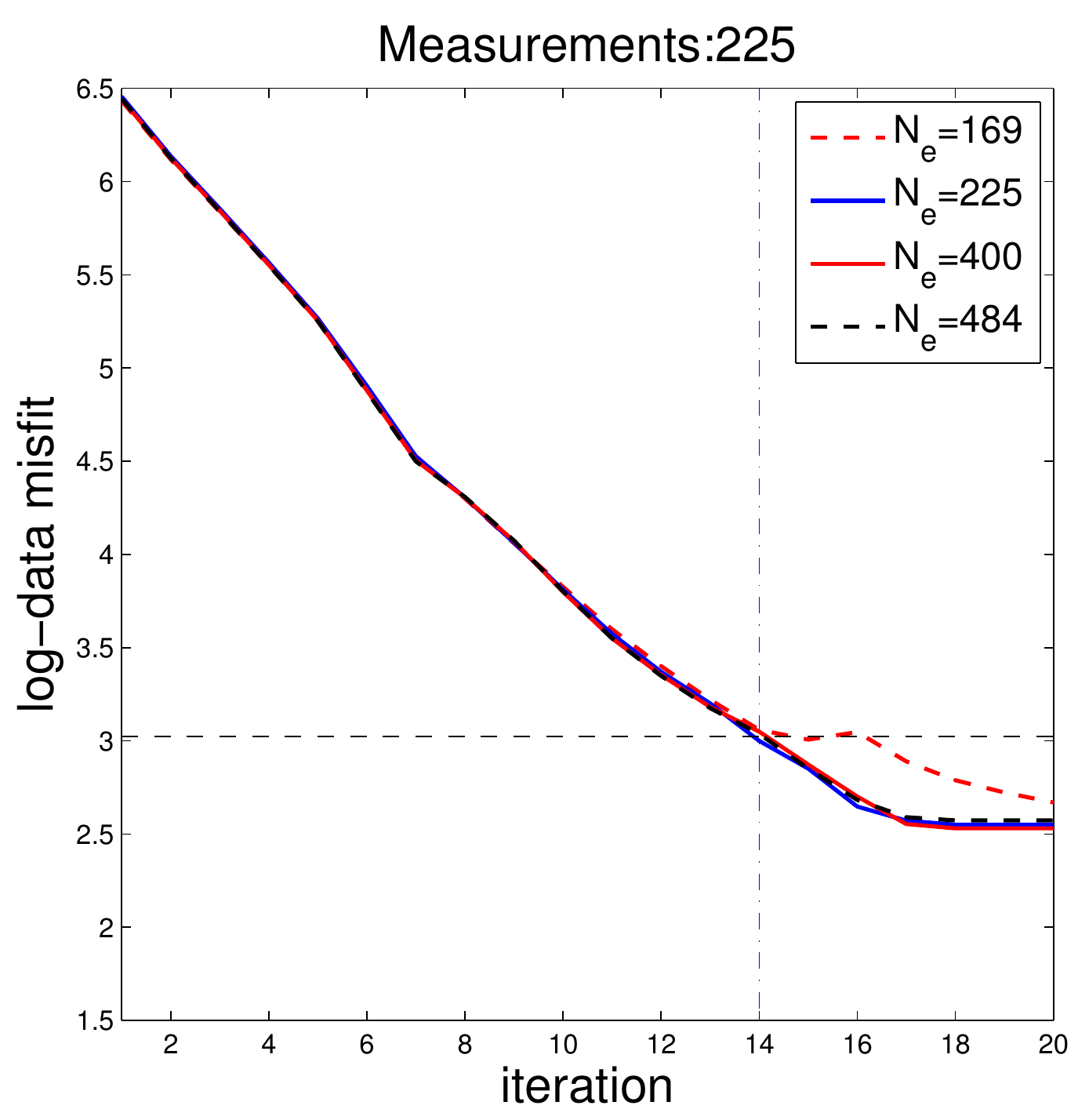}
\includegraphics[scale=0.3]{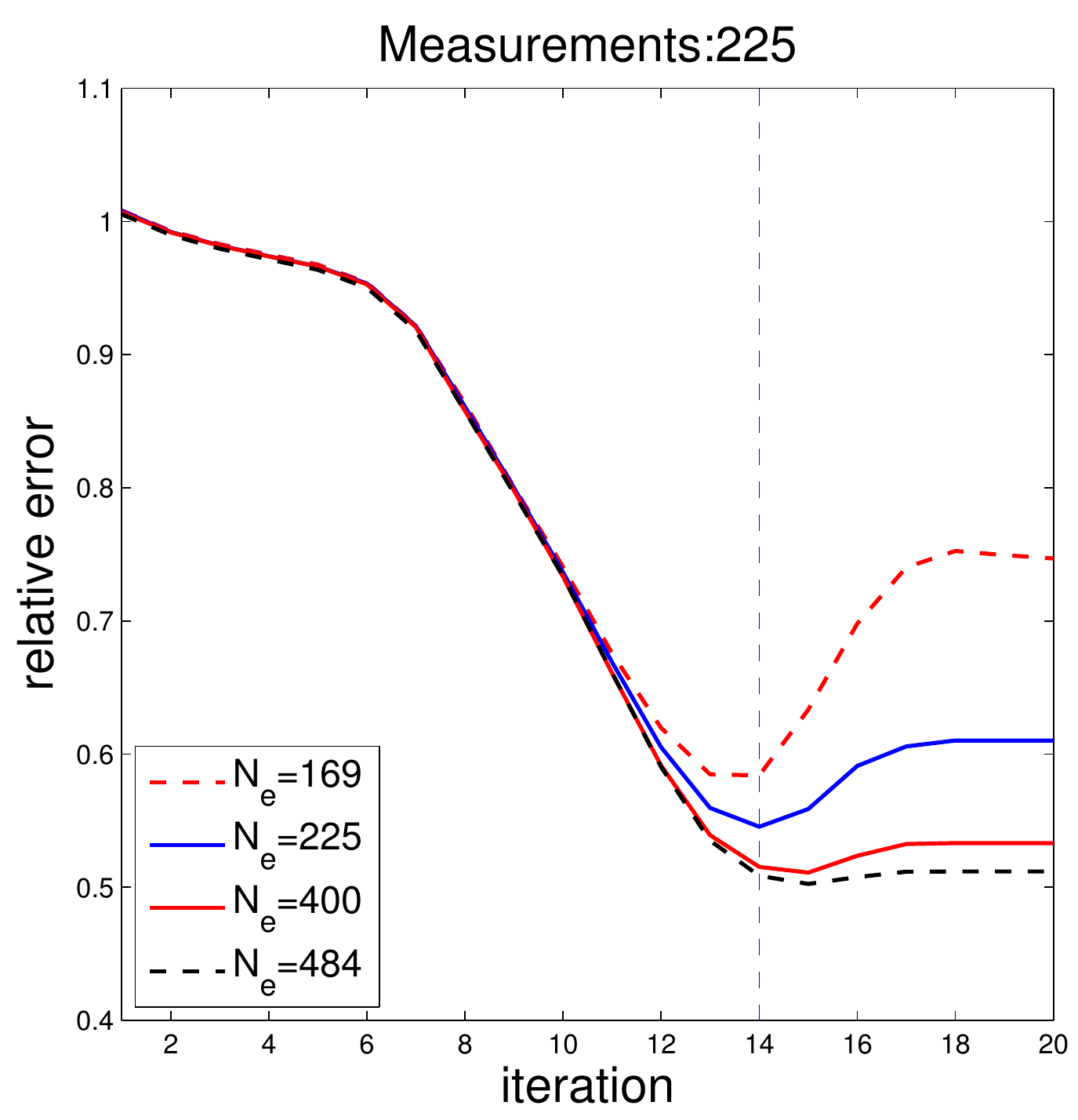}\\
\includegraphics[scale=0.32]{Well_Loc_12}
\includegraphics[scale=0.3]{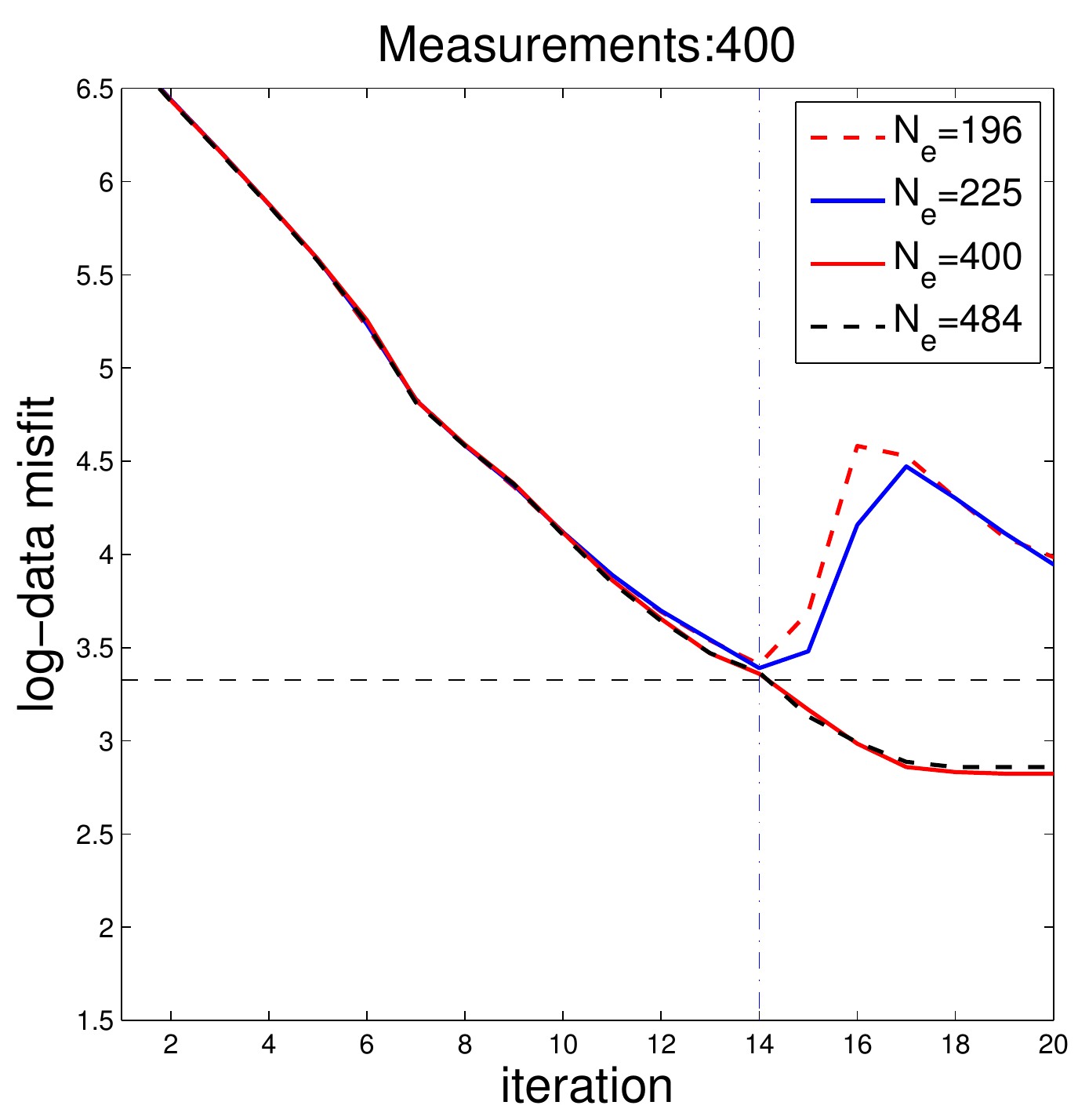}
\includegraphics[scale=0.3]{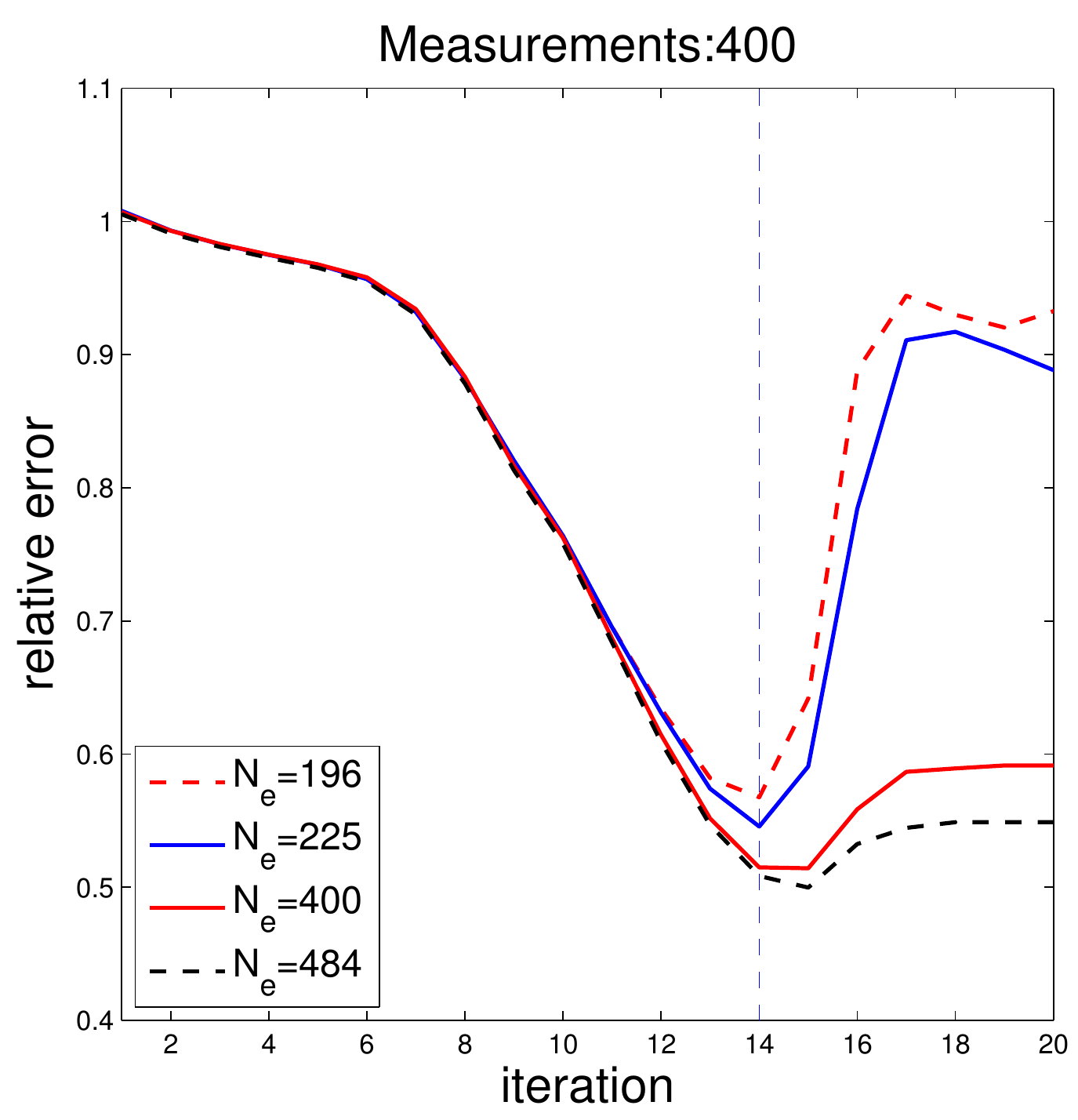}\\
\includegraphics[scale=0.32]{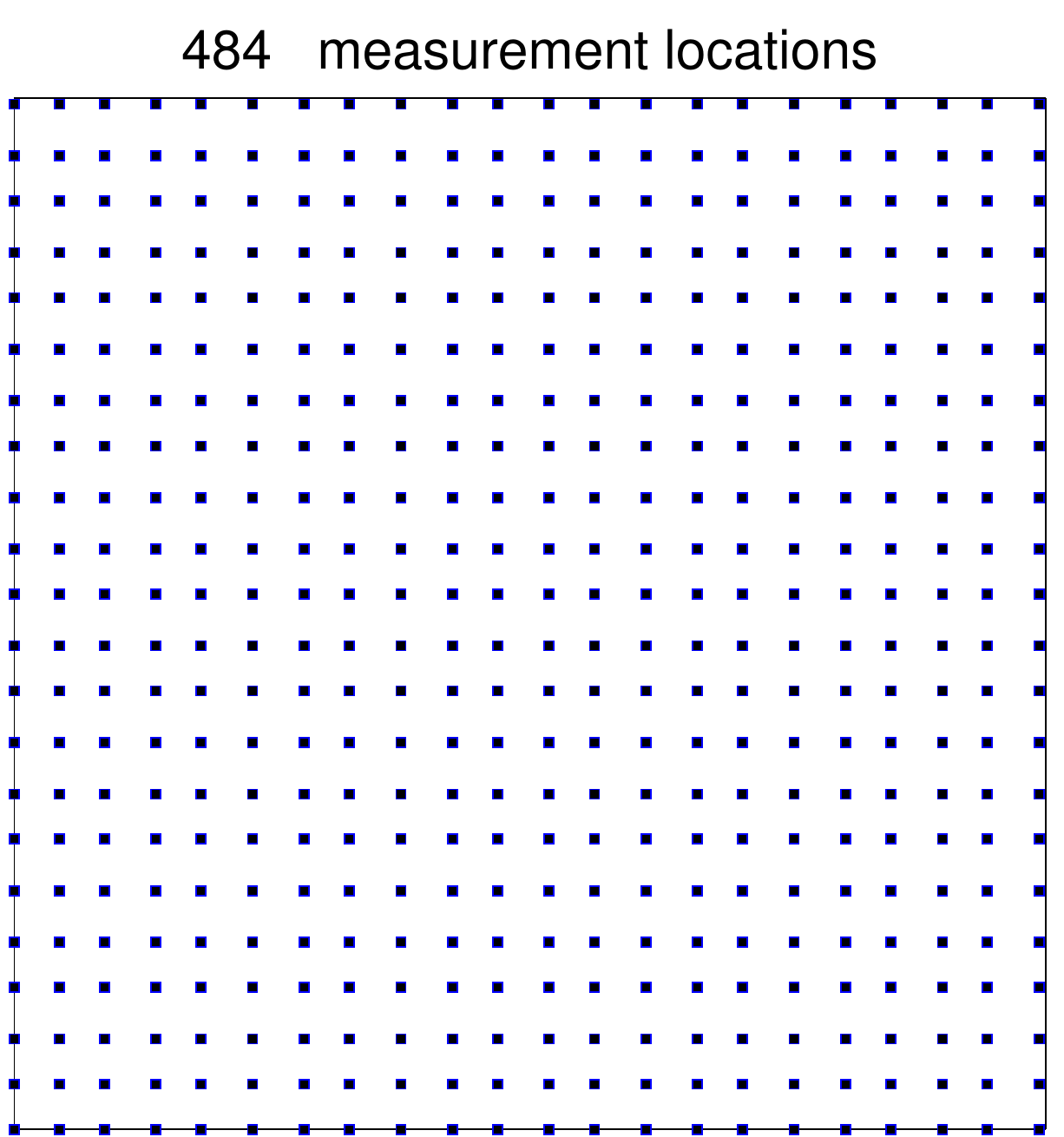}
\includegraphics[scale=0.3]{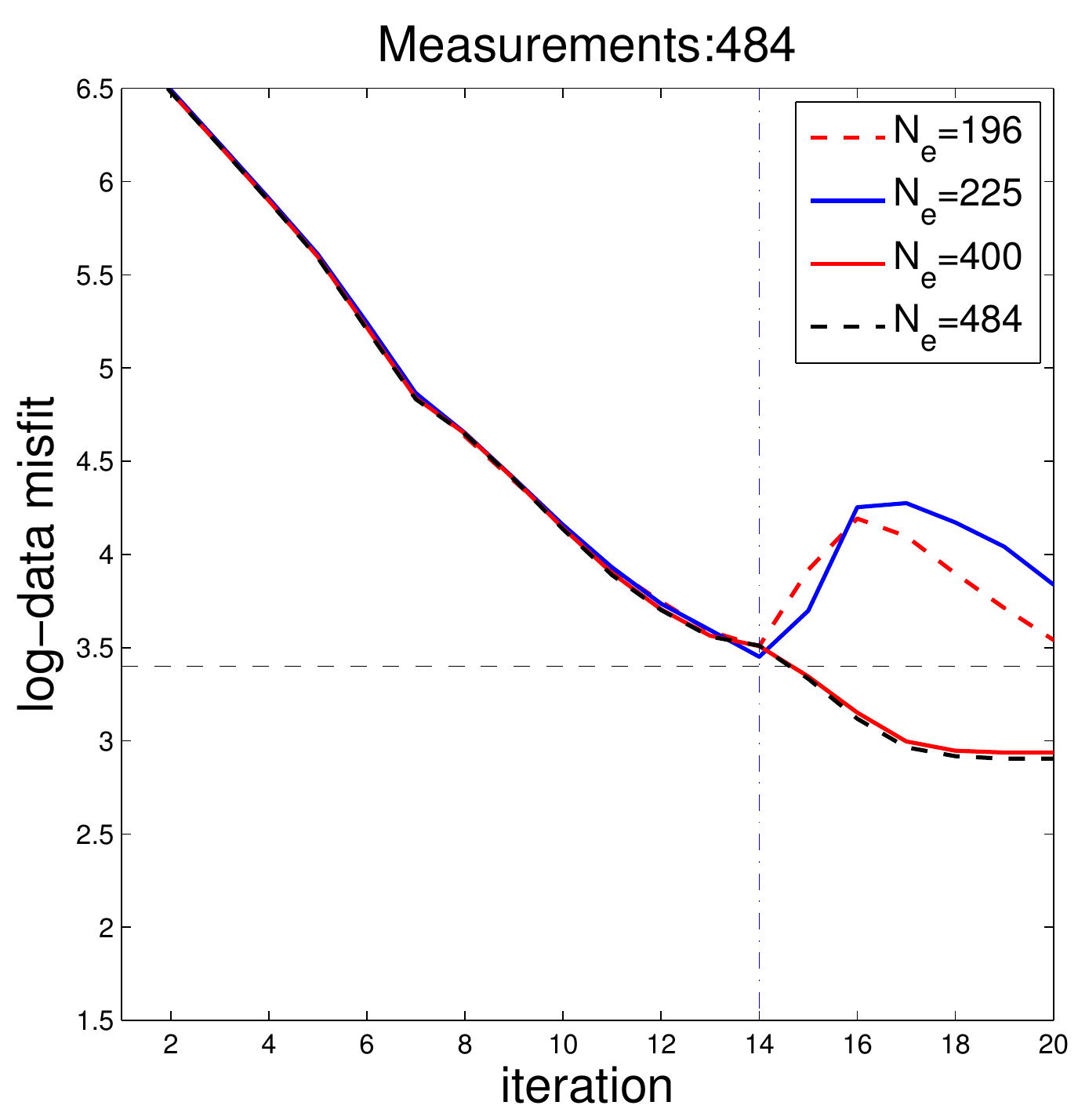}
\includegraphics[scale=0.3]{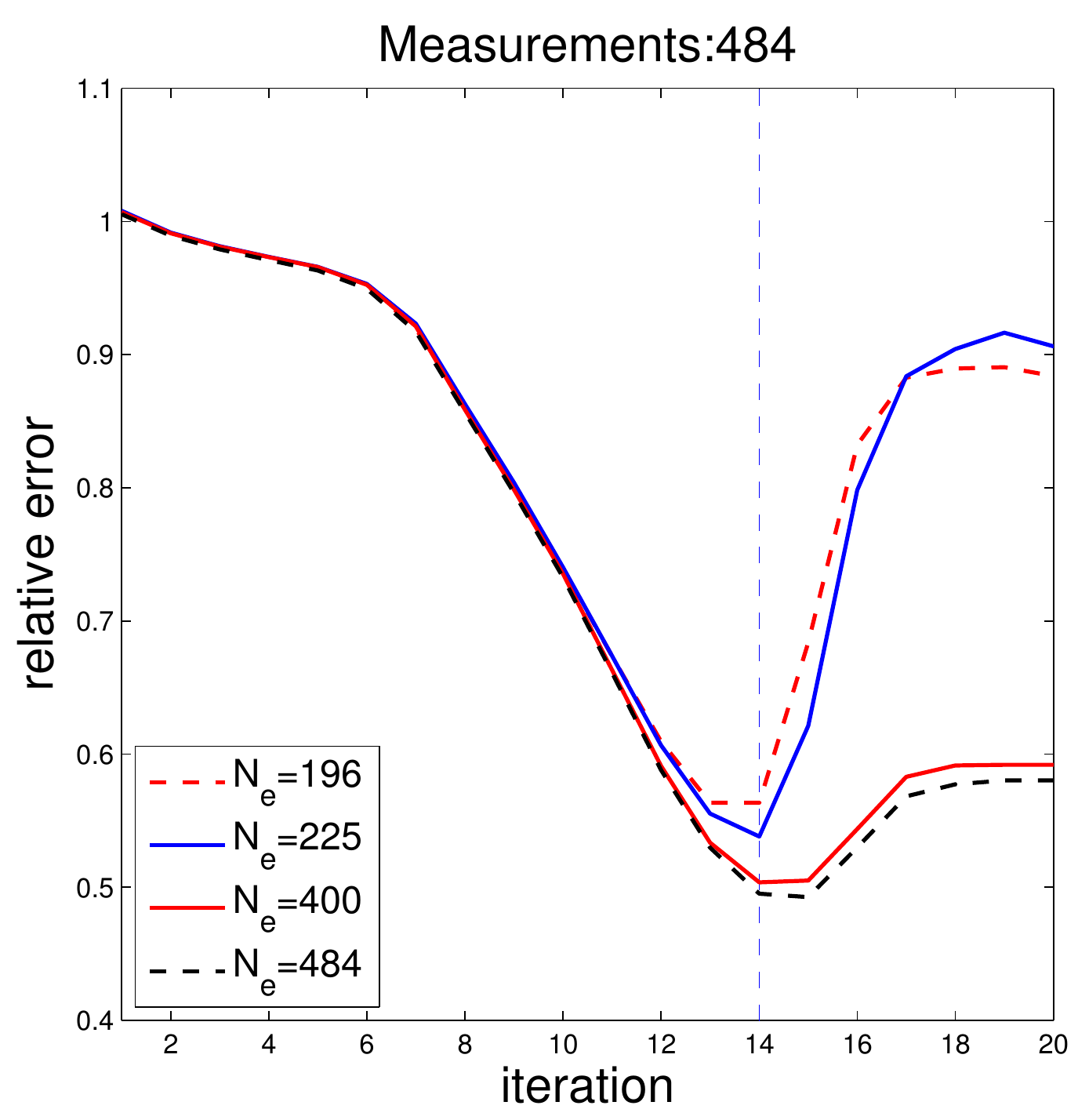}
 \caption{Numerical results with the proposed method for $\rho=0.7$ and different measurement configurations and ensemble sizes. Left column: measurement configuration. Middle column: relative error w.r.t. the truth. Right column: Log - data misfit.   Quantities are averaged at each iteration, over 40 experiments corresponding to different experiment. }
   \label{Fig6}
\end{center}
\end{figure}

We conduct a set of experiments, each of them similar to the one described in the preceding subsection  but with synthetic data generated from different measurement configurations. For each measurement configuration and ensemble size, we conduct a set of 40 experiments corresponding to different initial ensembles generated as described earlier. Then, log-data misfit and relative error with respect to the truth were computed and averaged (over the 40 experiments) at each iteration of the scheme. For clarity we only report these averaged quantities corresponding to each ensemble size and each measurement configurations. In the left column of \Fref{Fig5} and \Fref{Fig6} we display the measurement configuration for each experiment. The middle and right columns of Figures \ref{Fig5} and \ref{Fig6} show the (averaged over different experiments) log data misfit and error w.r.t truth for four different ensemble sizes (increasing from top to bottom in the legend of these Figures). The dotted horizontal line in the middle columns of  \Fref{Fig5} and \Fref{Fig6} denotes the value of the data misfit corresponding to $\eta /\rho$. The vertical line in the middle and right columns indicate the iteration after which the error w.r.t truth increases for the second choice (from to top to bottom) of ensemble size displayed on the legend of these figures. It is important to note that, irrespective of the ensemble size, the proposed algorithm does stabilize the estimates obtained in the early iterations. In other words, it avoids uncontrolled estimates that some standard unregularized Kalman-based methods may display in the first couple of iterations. However, the relation between ensemble size and number of measurements has severe effects on the ability of the stopping criteria (\ref{eq:m15}) to successfully terminate the algorithm under the guidelines that arise from the application of iterative regularization (i.e. with $\tau \approx 1/\rho$). Let us, for example, examine the case $M\leq 100$ reported in Figure \ref{Fig5}. For smaller measurement locations (i.e. $M=25,36$), an ensemble size equal to the number of measurements $N_{e}= M$, on average, will provide an estimate whose error with respect to the truth will decrease in the first iterations but then starts increasing long before the data-misfit has reached the value $\eta/\rho$. A similar but less drastic behavior is observed for larger number of measurements ($M=49, 64$). However, as we increase $M$, we find that a selection $N_{e}=M$ yields results that seemed to be stabilized whenever the scheme is stopped according to (\ref{eq:m15}) with  $\tau=1/\rho$. Note that for these measurement configurations with $M\leq100$, the ensemble size needs to be larger than the number of measurements so that the computational solution is properly stabilized with the proposed stopping criteria. 

Let us consider now the case $M>100$ reported in Figure \ref{Fig6}. In contrast to the previous case where an ensemble of a size larger than that number of measurements is needed for the stabilization of the computations, here certain choices of $N_{e}$ with $N_{e}<M$ may suffice. Note for example that for $M=225$, an ensemble with $N_{e}=169$ will result in stable computations by using (\ref{eq:m15}) with $\tau\approx 1/\rho$. Interestingly, for larger $M$ it is more evident that the selection $N_{e}>M$ is unnecessary for the algorithm to exhibit the regularization properties based on the stopping criteria mentioned above. Indeed, for $M=400$ and $N=484$ we note that $N_{e}=196$ and $N_{e}=225$ will clearly provide regularized estimates in the sense that their error will not increase before the data misfit reaches the value $\eta /\rho$. For these cases with $M> 100$, a selection of $N_{e}\ge M$ will clearly result in not only stable but an accurate identification. However, it is worth reiterating that for large-scale applications, the smallest ensemble size that produce a regularized solution is desirable due to the high computational cost of the forward simulations. In summary, there seems to be an apparent critical number of measurements ($M=100$) below of which a selection of $N_{e}>M$ is needed for the proposed method to fully stabilize the scheme. However, for larger number of measurements the regularizing properties of the scheme enable us to use a selection of $N_{e}<M$ which results in a reasonable computational cost.

\subsubsection{The tunable parameters $\rho$.}

\begin{figure}[htbp]
\begin{center}
\includegraphics[scale=0.2]{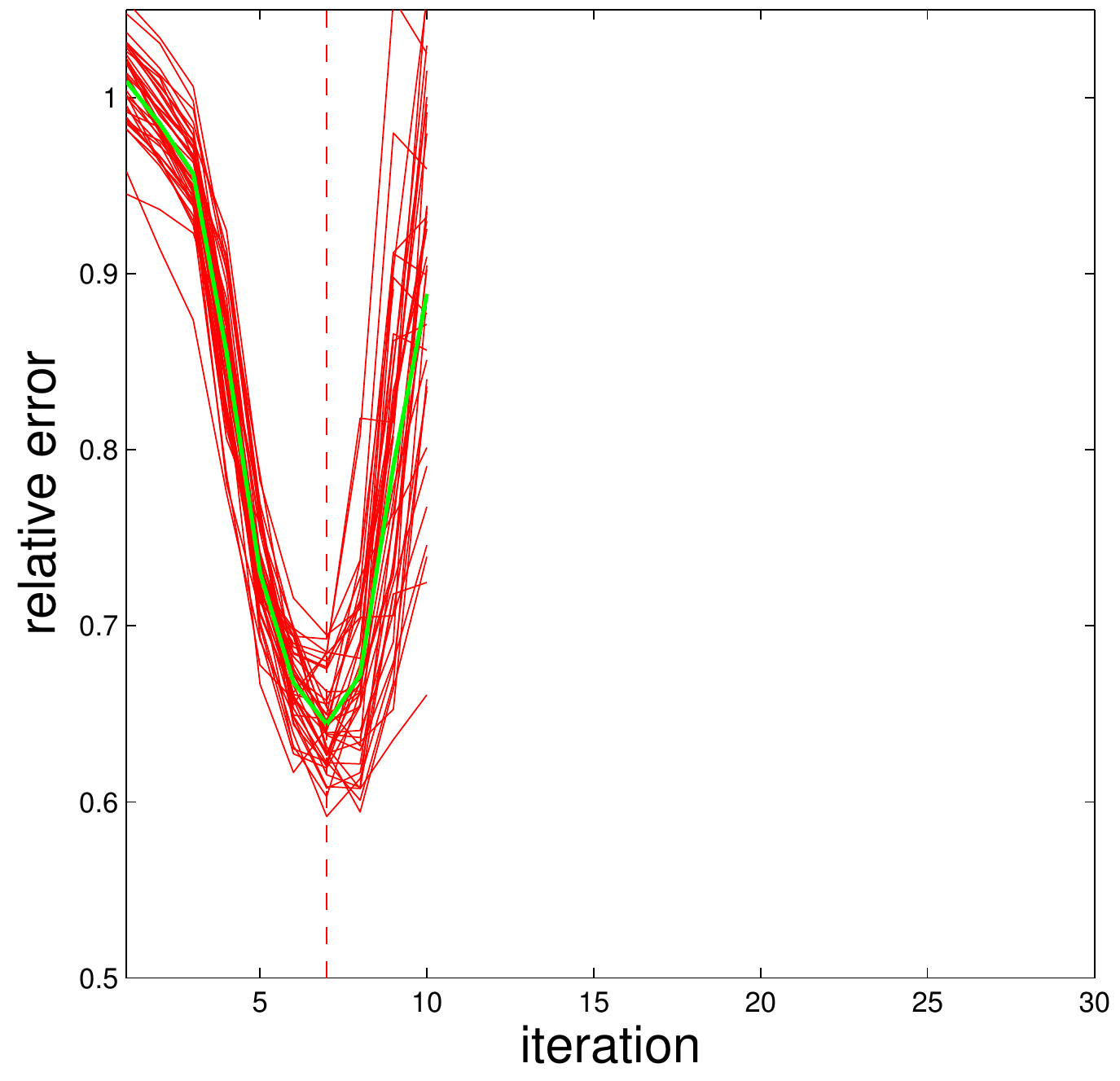}
\includegraphics[scale=0.2]{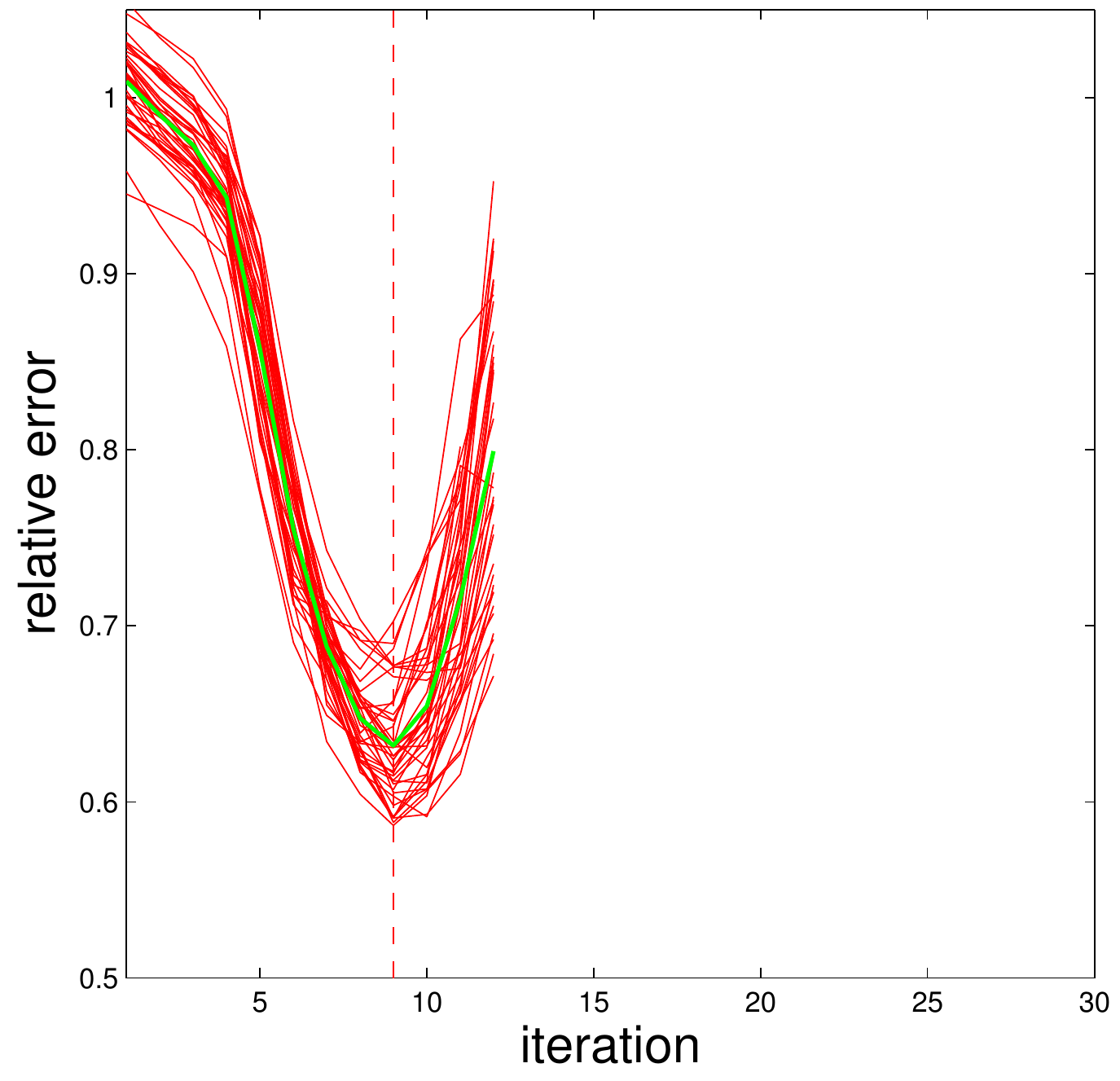}
\includegraphics[scale=0.2]{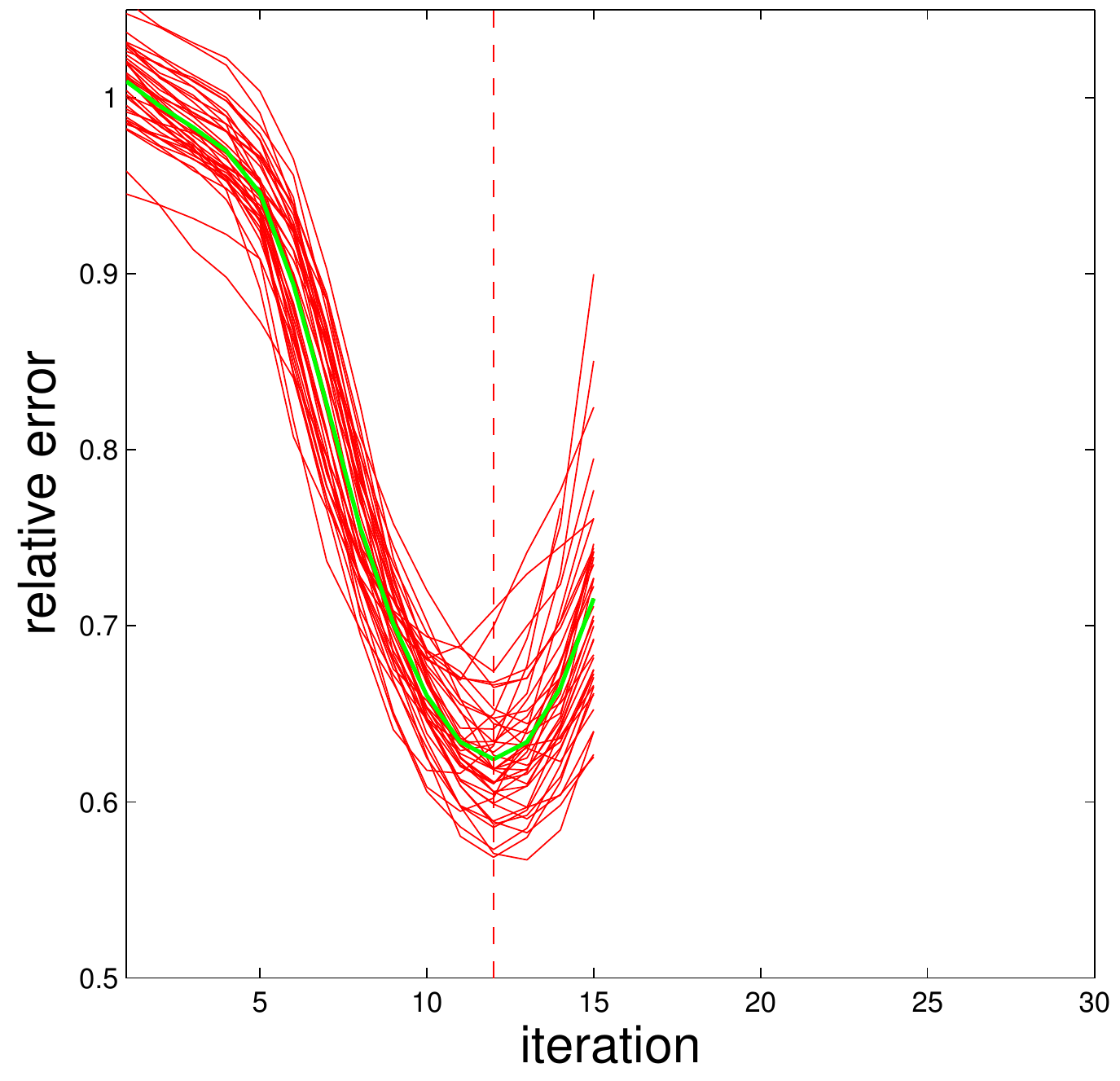}
\includegraphics[scale=0.2]{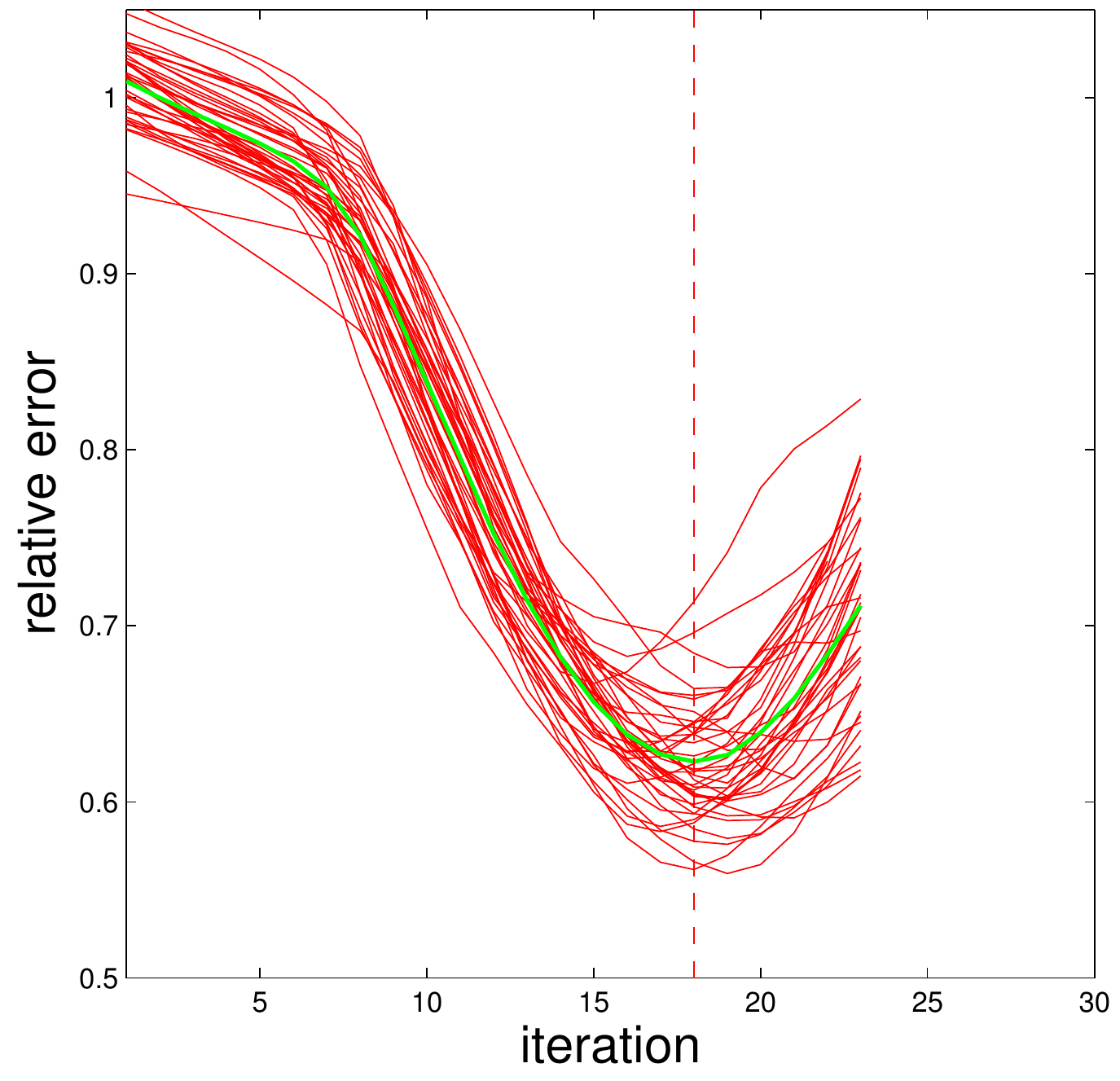}
\includegraphics[scale=0.2]{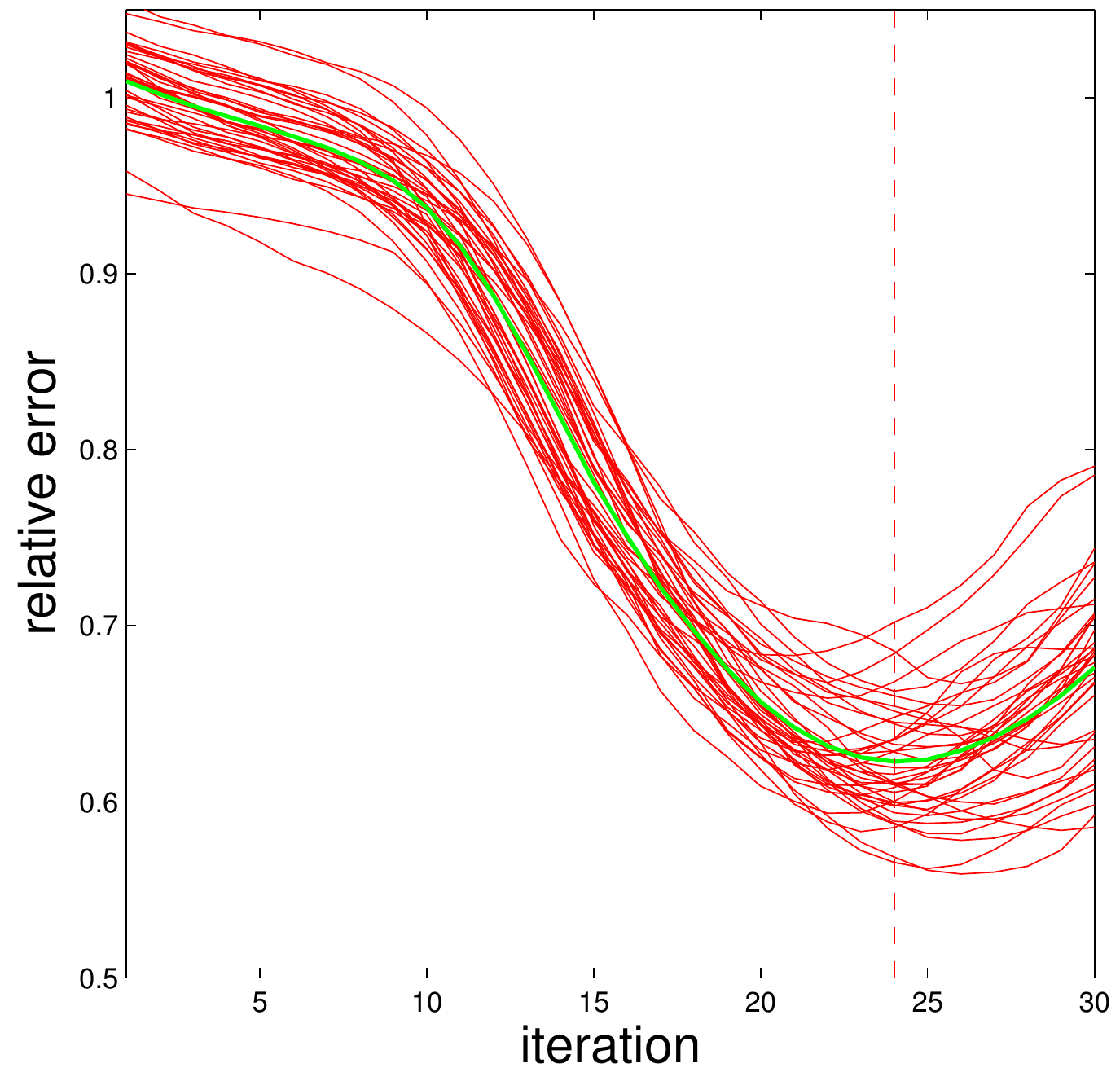}\\
\includegraphics[scale=0.2]{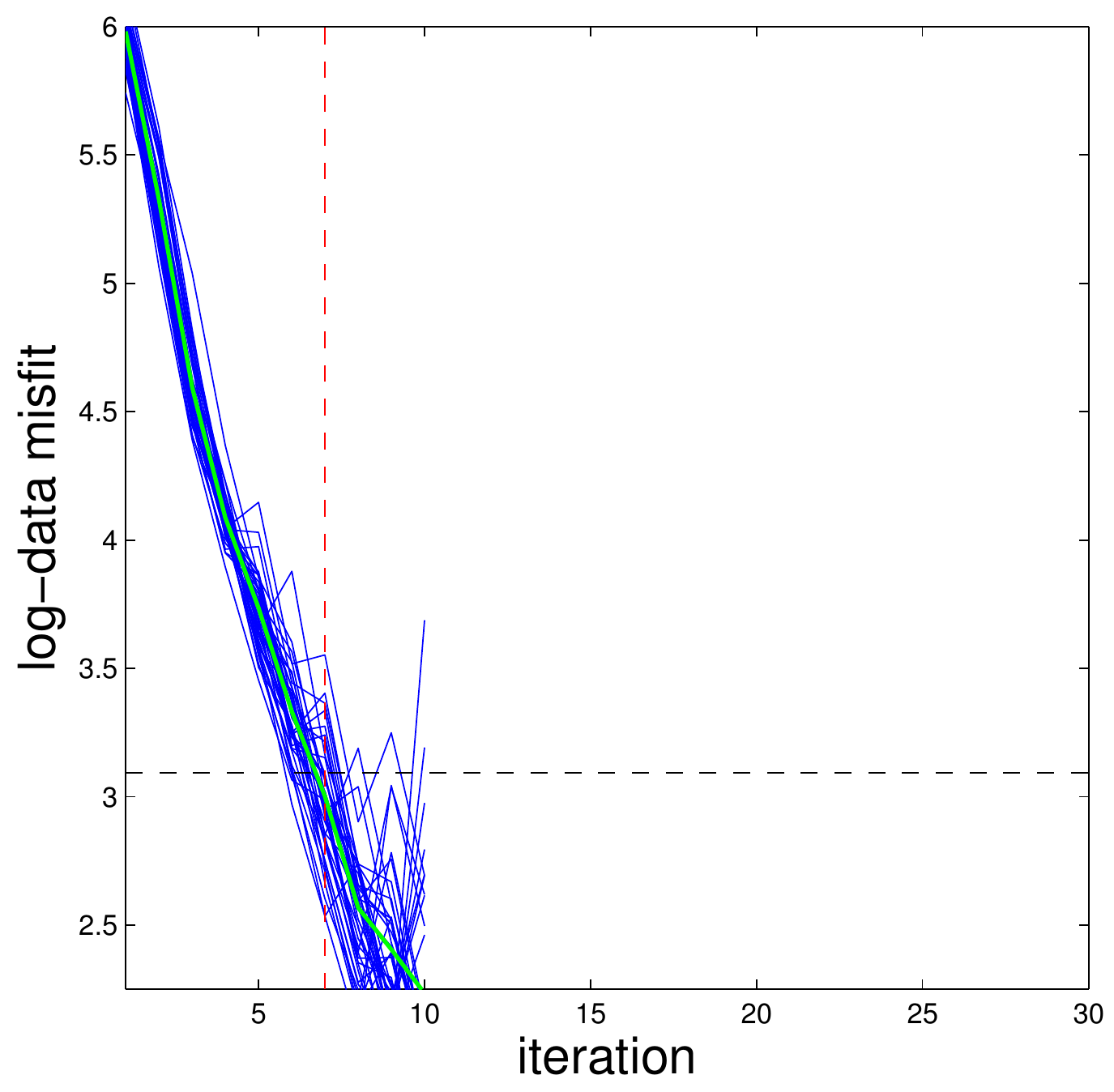}
\includegraphics[scale=0.2]{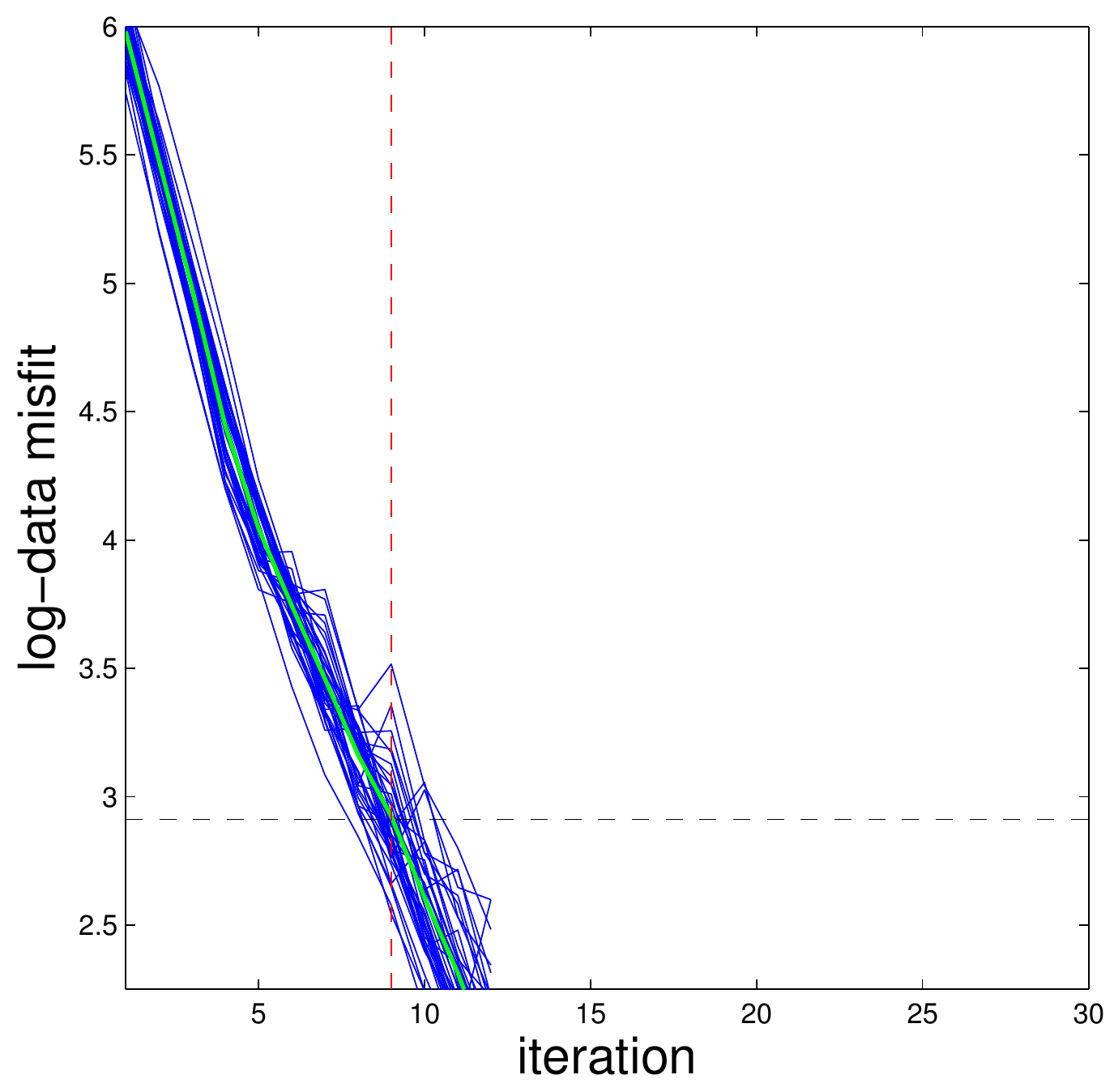}
\includegraphics[scale=0.2]{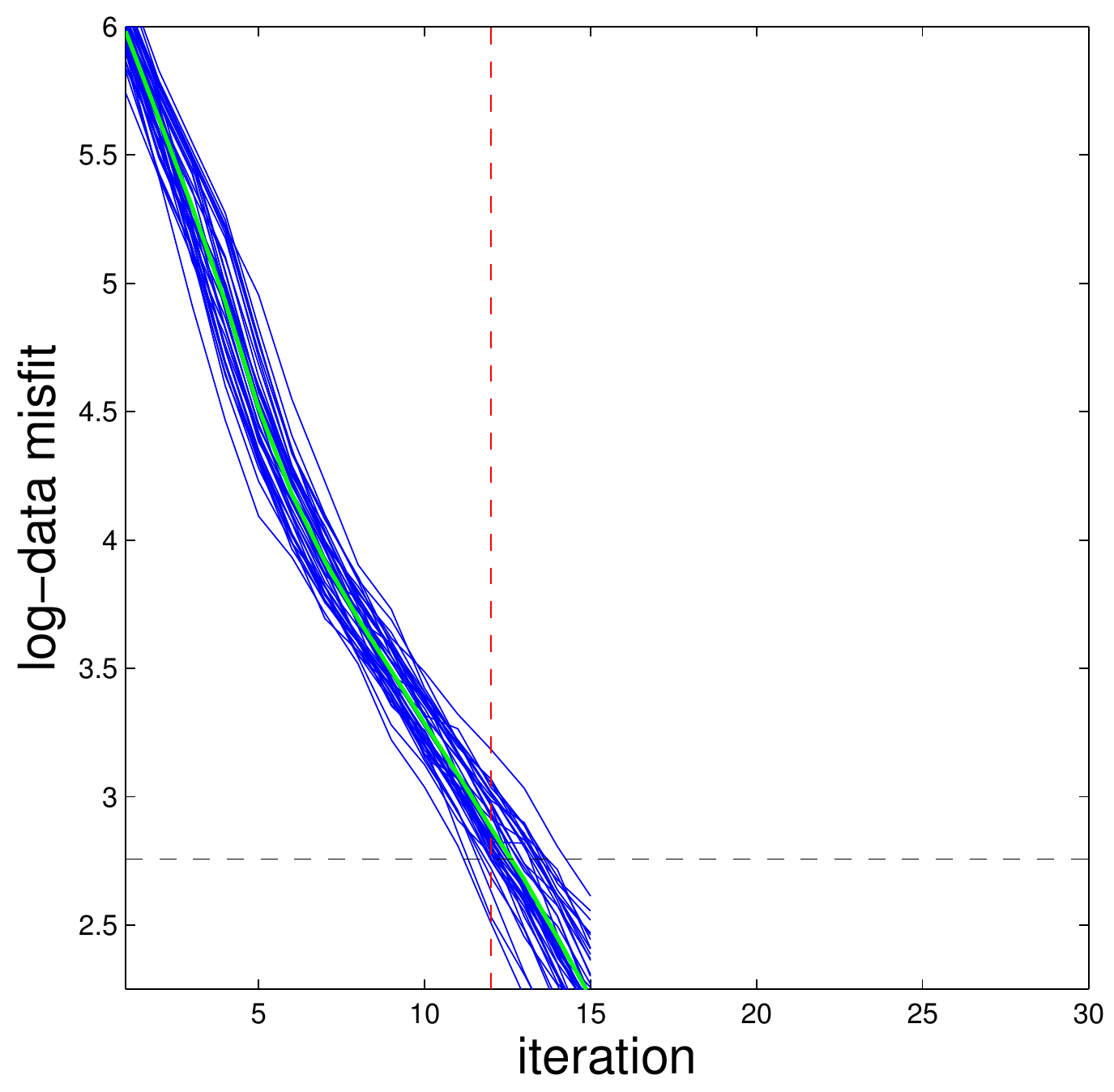}
\includegraphics[scale=0.2]{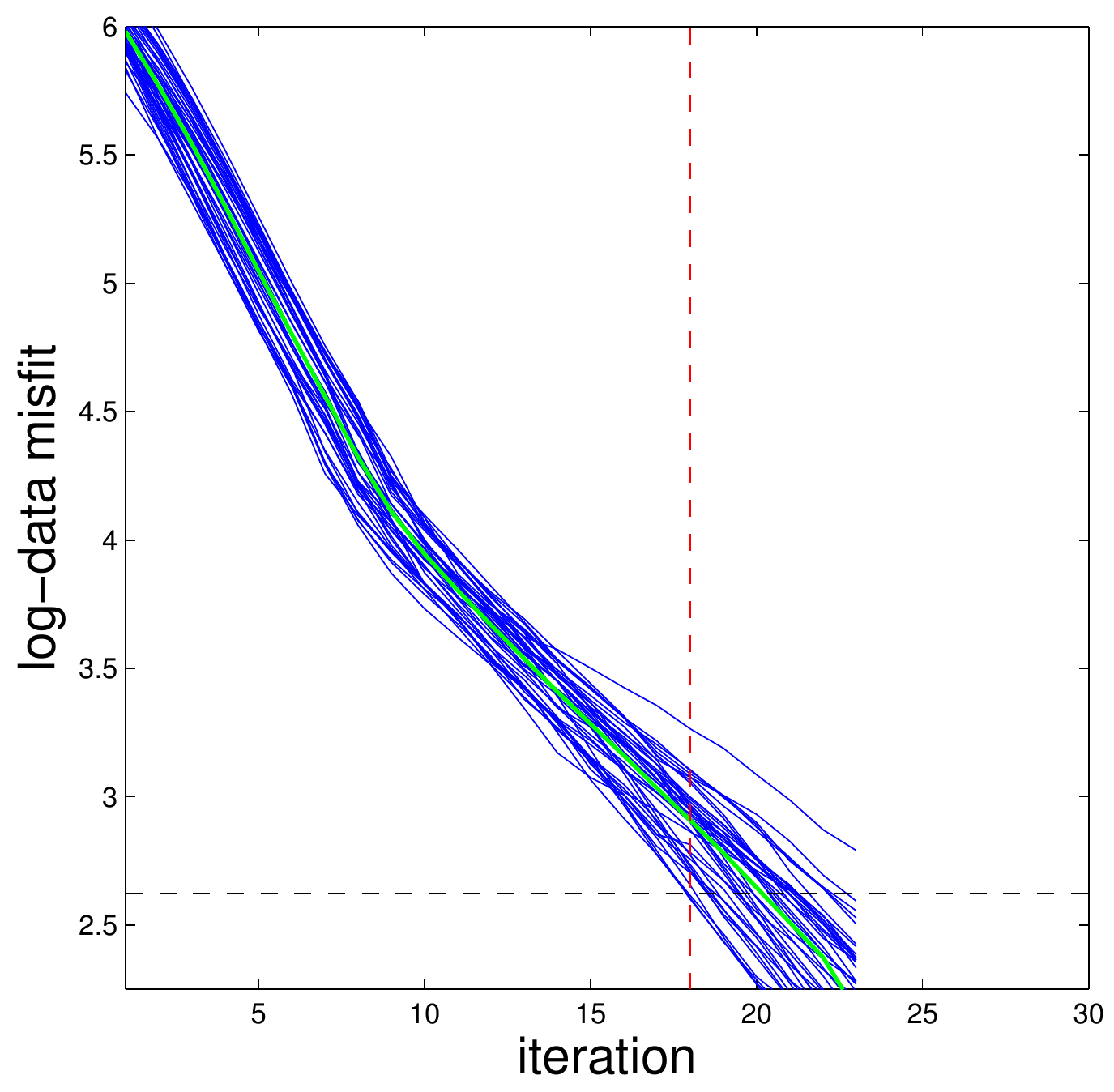}
\includegraphics[scale=0.2]{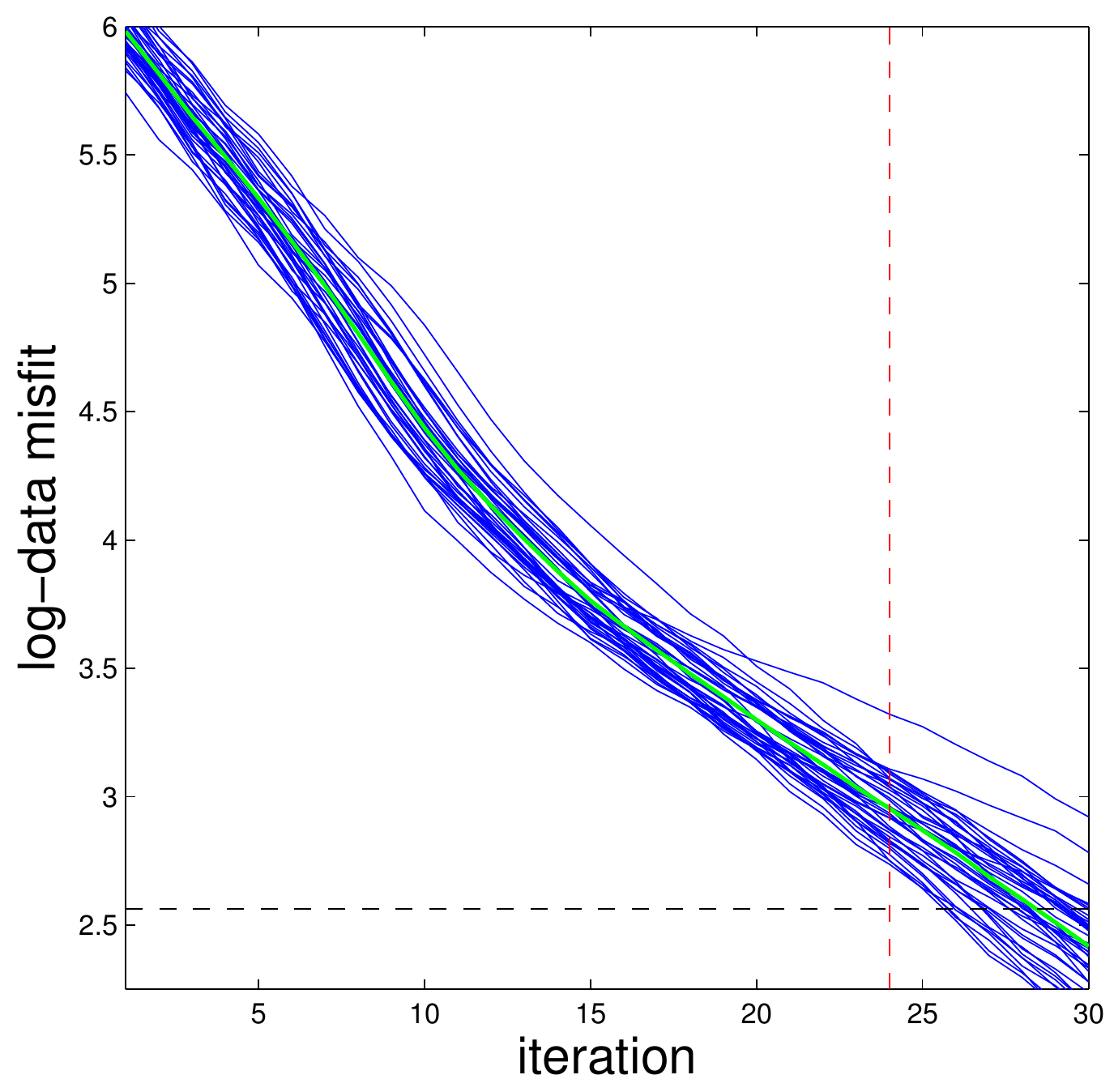}
 \caption{Regularizing ensemble method applied with $N_{e}=150$, $100$ measurements and $\rho$ (from left to right) $\rho=0.5,0.6,0.7,0.8,0.85$. Top row: relative error w.r.t. truth. Bottom row: Log - data misfit.  }
 \label{Fig7}
\end{center}
\end{figure}

As we discussed in subsection \ref{3_2}, the parameter $\rho$ controls the ensemble updates. The choice of this parameter has a significant effect on both accuracy and cost of the proposed scheme. We now investigate this effect with our Darcy flow model with the measurement configuration of $100$ locations displayed in Figure \ref{Fig1} (right). The results from the preceding section indicate that an ensemble of size $N_{e}>100$ will provide the stabilization needed with the choice $\tau\approx 1/\rho$ in the stopping criteria. We select $N_{e}=150$ for conducting a set of experiments where Algorithm \ref{Al1} is applied with several choices of $\rho$. For each choices of $\rho$, the experiment is repeated 40 times with a different selection of the initial ensemble generated as we previously described. In Figure \ref{Fig7} we display, for these 40 experiments, the log-data misfit (bottom) and the relative error with respect to the truth (top) of the ensemble mean obtained with the proposed scheme with different choices of $\rho$. The horizontal dotted line in Figure \ref{Fig7} (bottom) represent the log of $\eta/\rho$. Note that the stopping criteria (\ref{eq:m15}) with $\tau \approx 1/\rho$ provides a reasonable stabilization of algorithm. For $\rho<0.7$ we observe that the error with respect to the truth increases quite rapidly when the data misfit drops below the value $\eta/\rho$. We note that, on average, there is a slight decrease of the error with respect to the truth as we increase $\rho$. However, this slight increase in the accuracy of the identification has associated an increase in the computational cost. For the Darcy flow model, our experiments suggest that $\rho=0.7$ represents a reasonable choice in terms of accuracy and cost.

\subsubsection{Regularizing properties of the proposed method in the small-noise limit}

In this section we conduct experiments to study the regularization properties of the proposed scheme in the small noise limit. The previous experiments show that with a reasonable selection of the ensemble size with respect to the measurement configurations, the mean of the proposed ensemble method achieve stable approximations of the truth.  In other words, the early termination of the scheme, say at iteration $n^{\star}$,  yields an approximation $u^{\eta}\equiv \overline{u}_{n^{\star}}$ given by the mean of the ensemble at the stopped iteration. As discussed earlier, the aim of our iterative regularizing algorithm is to provide a stable approximation in the sense that, as $\eta\to 0$, $u^{\eta}\to u$ where $G(u)=G(u^{\dagger})$. Due to the lack of uniqueness of the identification problem we cannot claim that $u=u^{\dagger}$ (the estimate converges to the truth). Nevertheless, we expect that the elements $u\in X$ such that $G(u)=G(u^{\dagger})$ will possess similar spacial features of the truth. 

\begin{figure}[htbp]
\begin{center}
\includegraphics[scale=0.265]{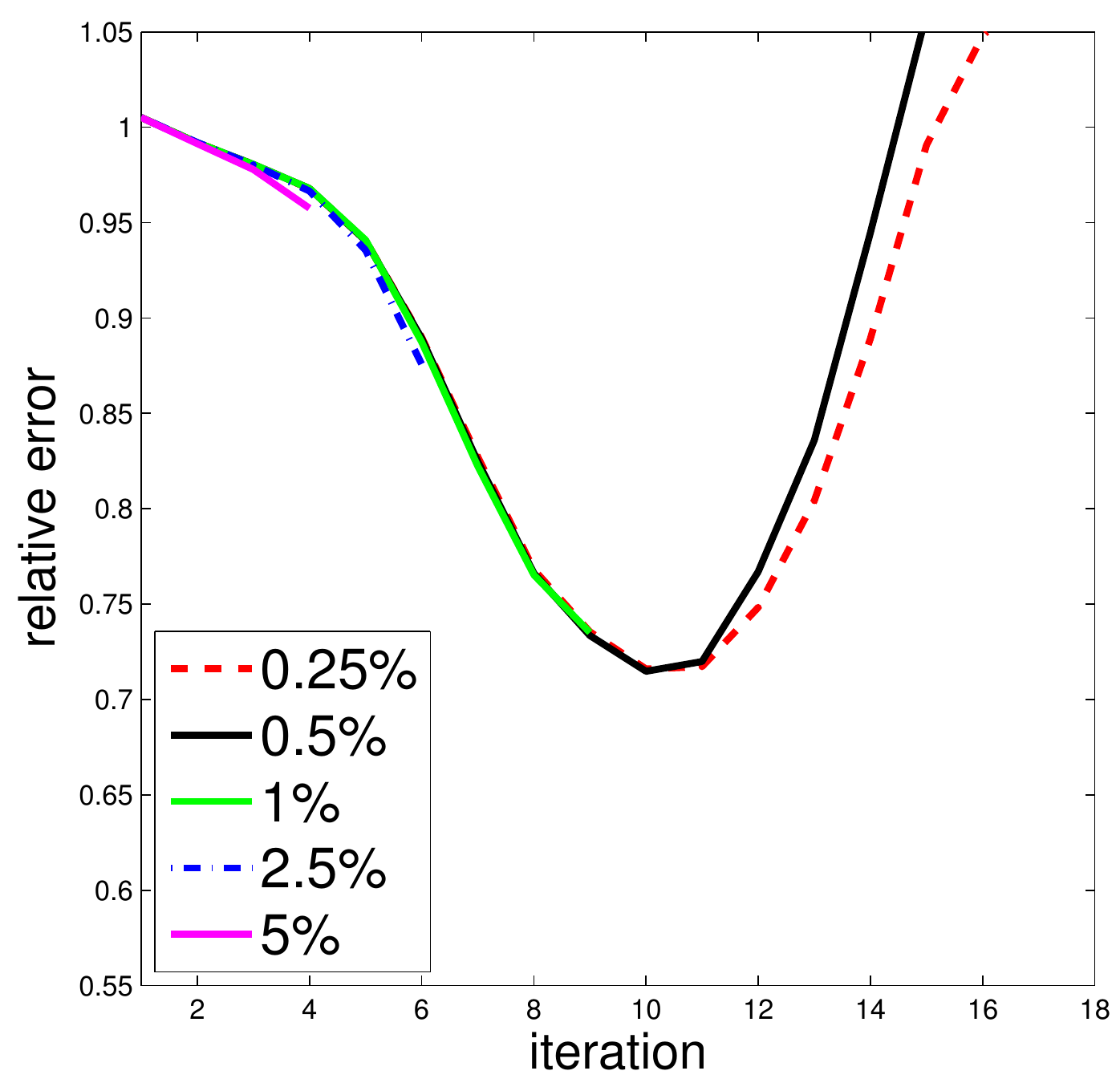}
\includegraphics[scale=0.265]{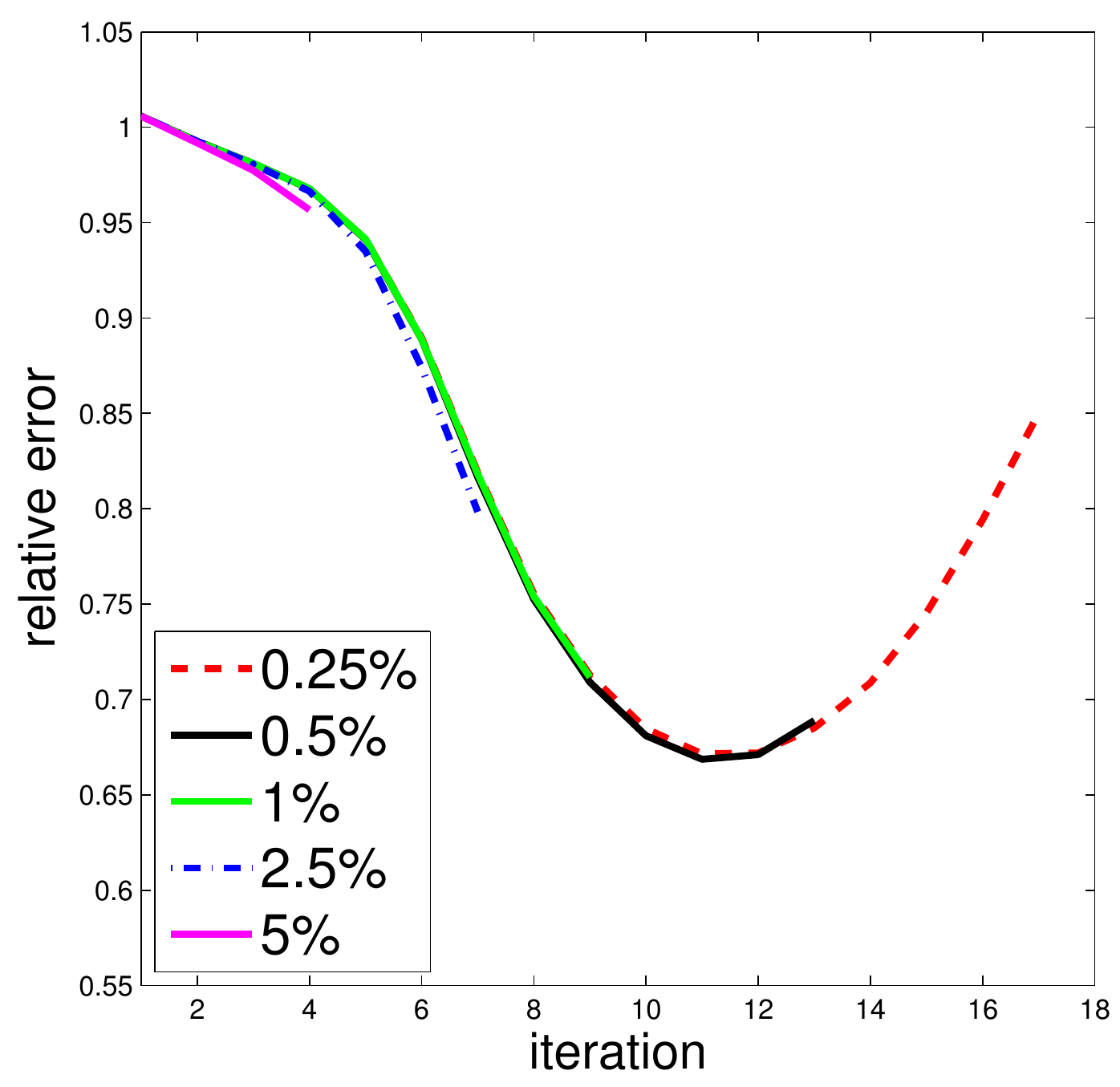}
\includegraphics[scale=0.265]{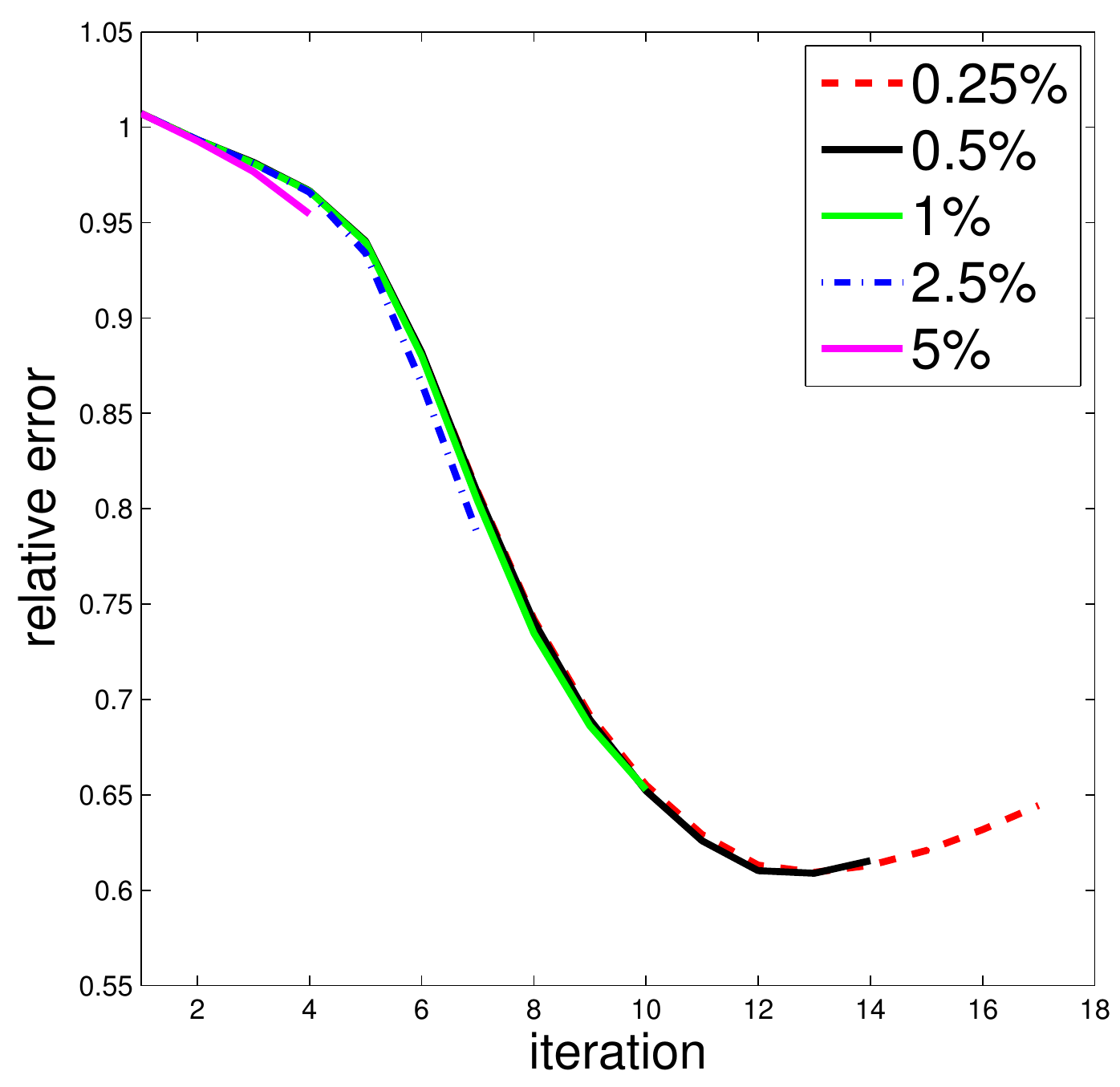}
\includegraphics[scale=0.265]{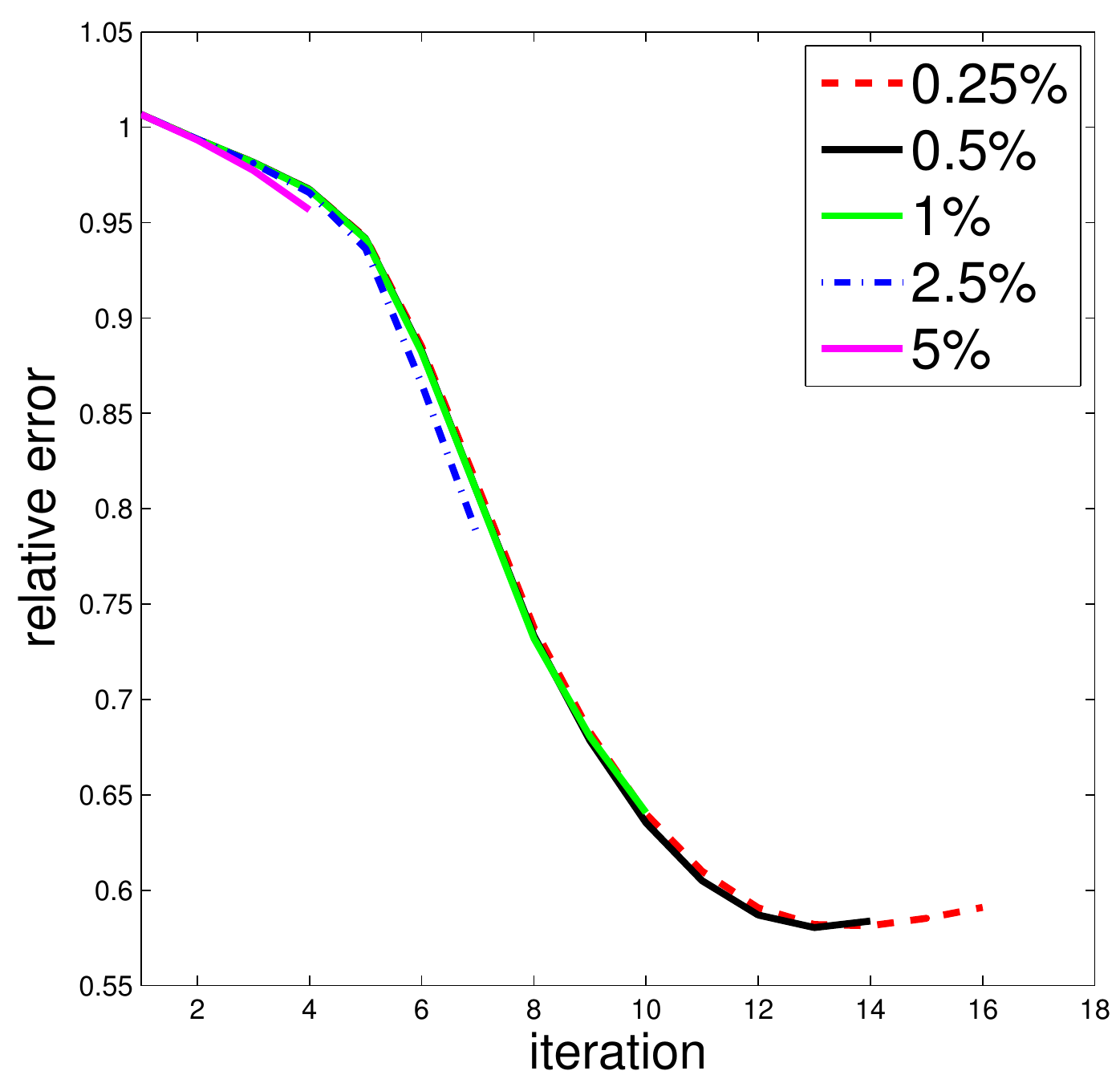}\\

\includegraphics[scale=0.265]{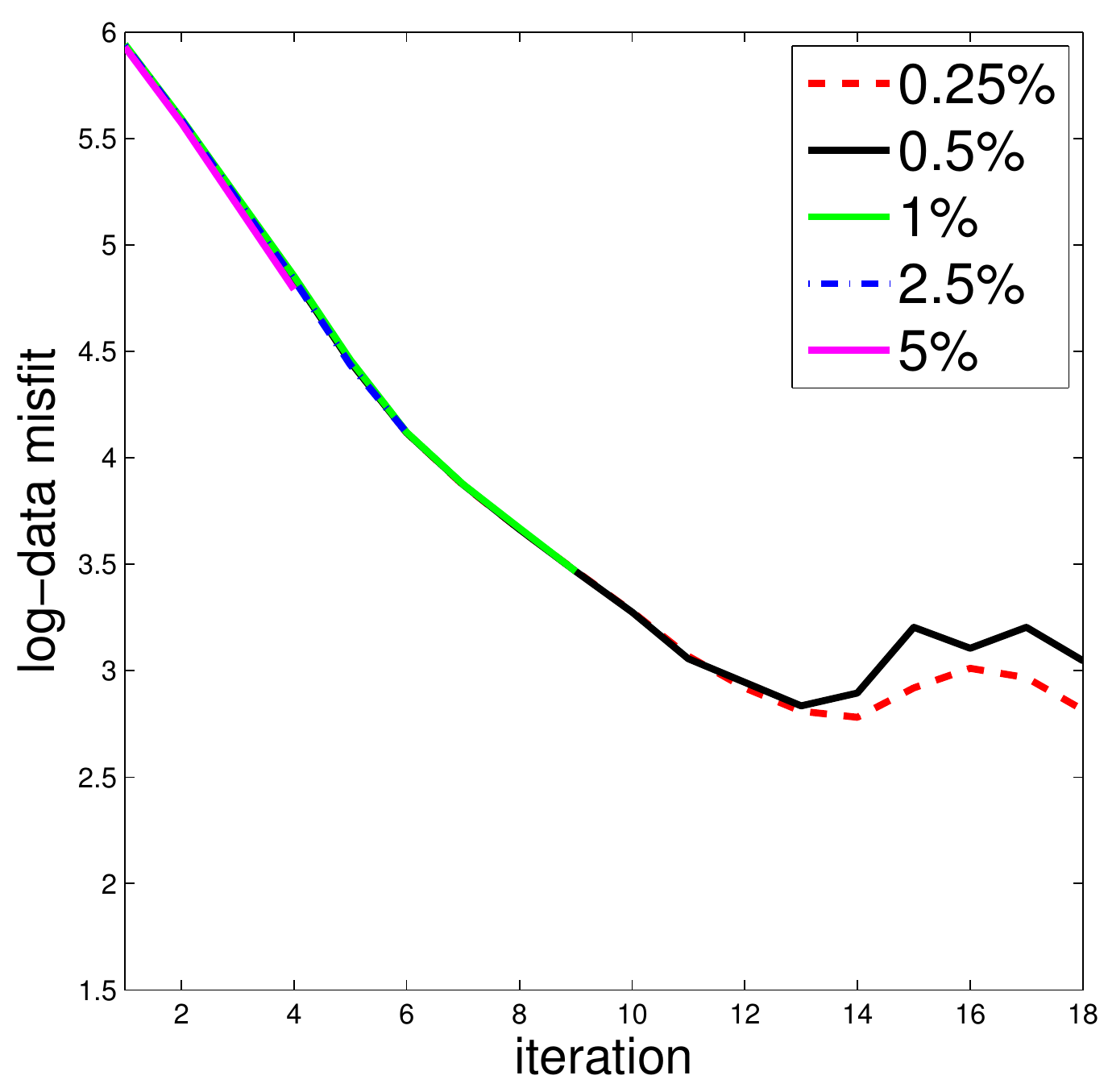}
\includegraphics[scale=0.265]{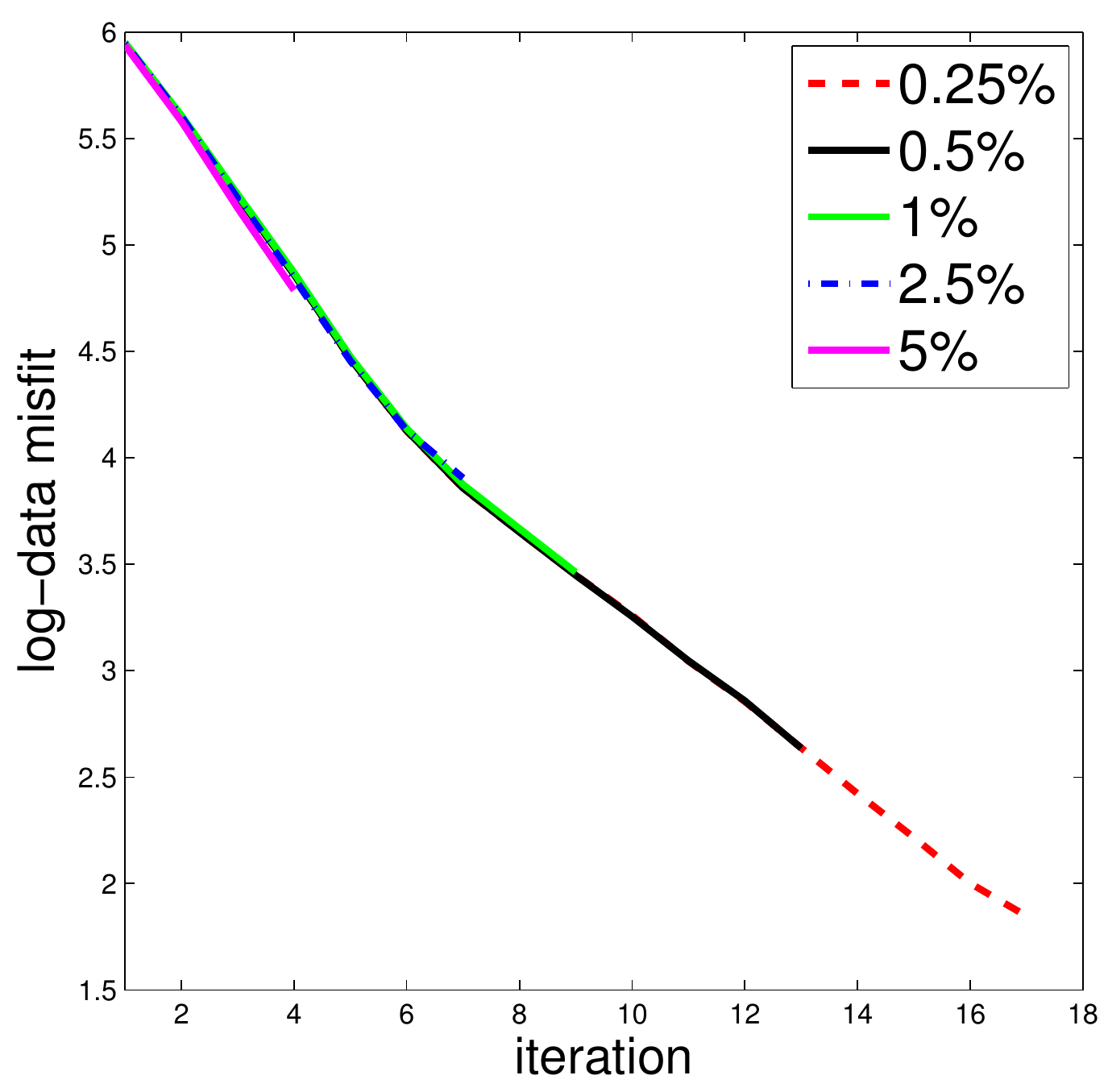}
\includegraphics[scale=0.265]{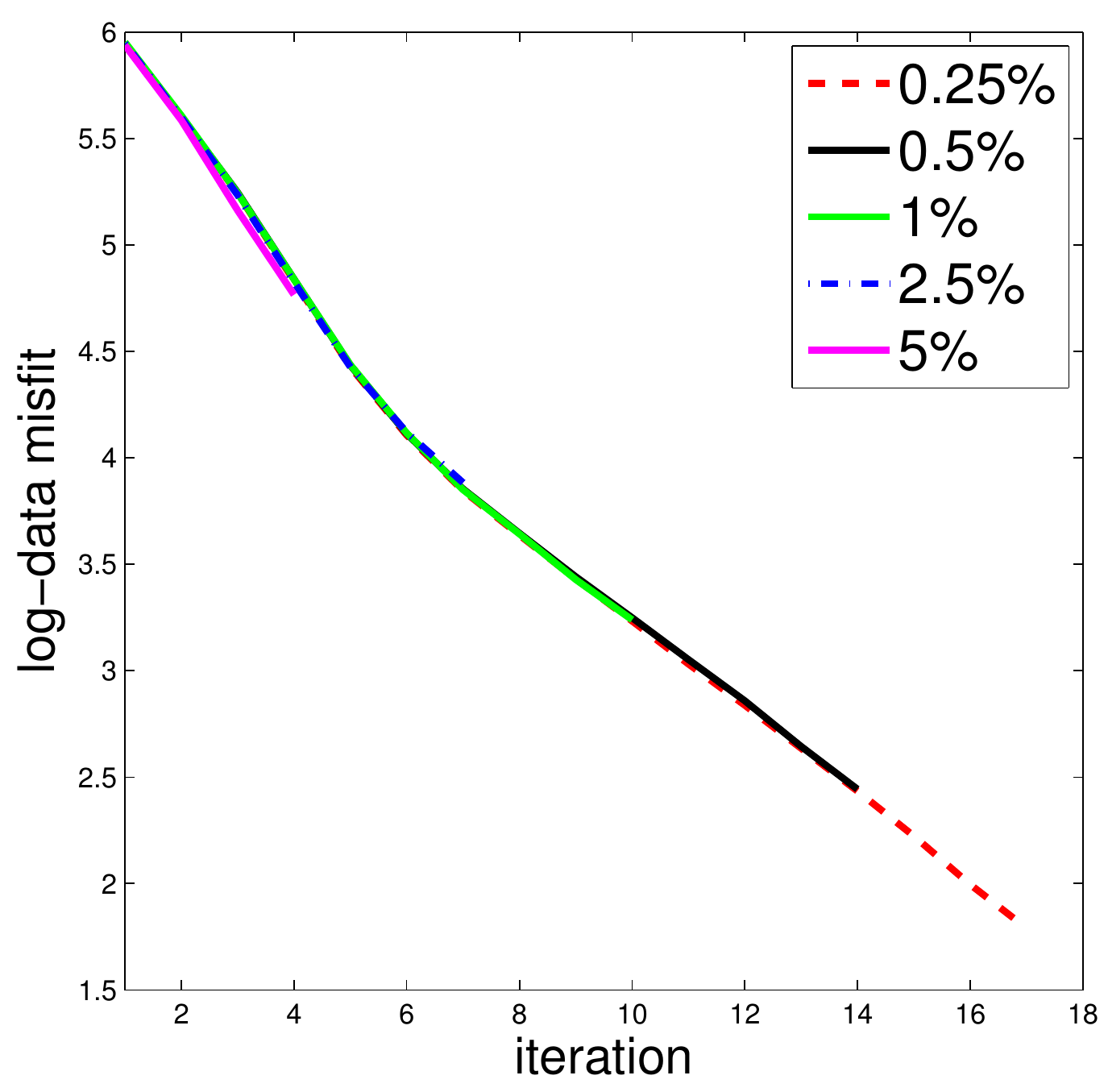}
\includegraphics[scale=0.265]{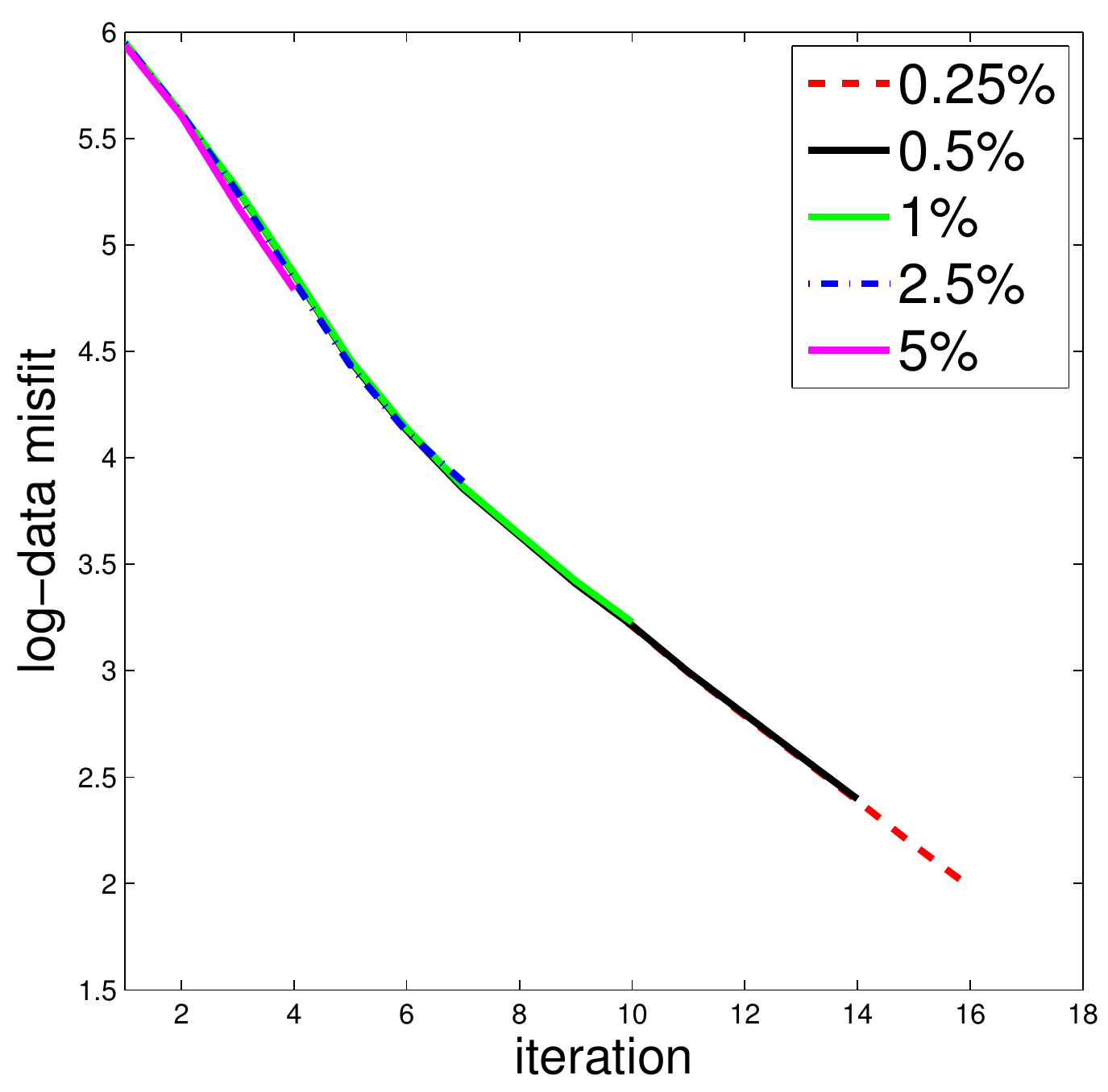}
 \caption{Log - data misfit (bottom) and relative error w.r.t. the truth (top) from 40 experiments with different initial ensembles of size of (from left to right) 75, 100, 150, 200 and different choices of noise level.}   \label{Fig8}
\end{center}
\end{figure}

We consider again the configuration with $M=100$ as before and we vary the ensemble size. In Figure \ref{Fig8} we show the relative error with respect to the truth (top) and the log-data misfit (bottom) obtained for each ensemble size with a different selection of noise levels. The level of noise is selected so that the norm of the noise relative to the data is of the percentage indicated in the plots. The matrix $\Gamma$ as selected as described earlier but kept fixed for these experiments associated to different noise levels. For this particular experiment the algorithm is stopped according to (\ref{eq:m15}) with $\tau=1/\rho$. We note that for the smallest noise considered here $0.5\%$ an ensemble of size $N_{e}=150$ was needed to properly stabilize the computations. Again, this shows that the proposed method, for sufficient large ensemble,  inherits the regularization properties of the regularization LM scheme of \cite{Hanke}.

\subsubsection{Comparison with the regularizing LM scheme.}
\begin{figure}[htbp]
\begin{center}
\includegraphics[scale=0.45]{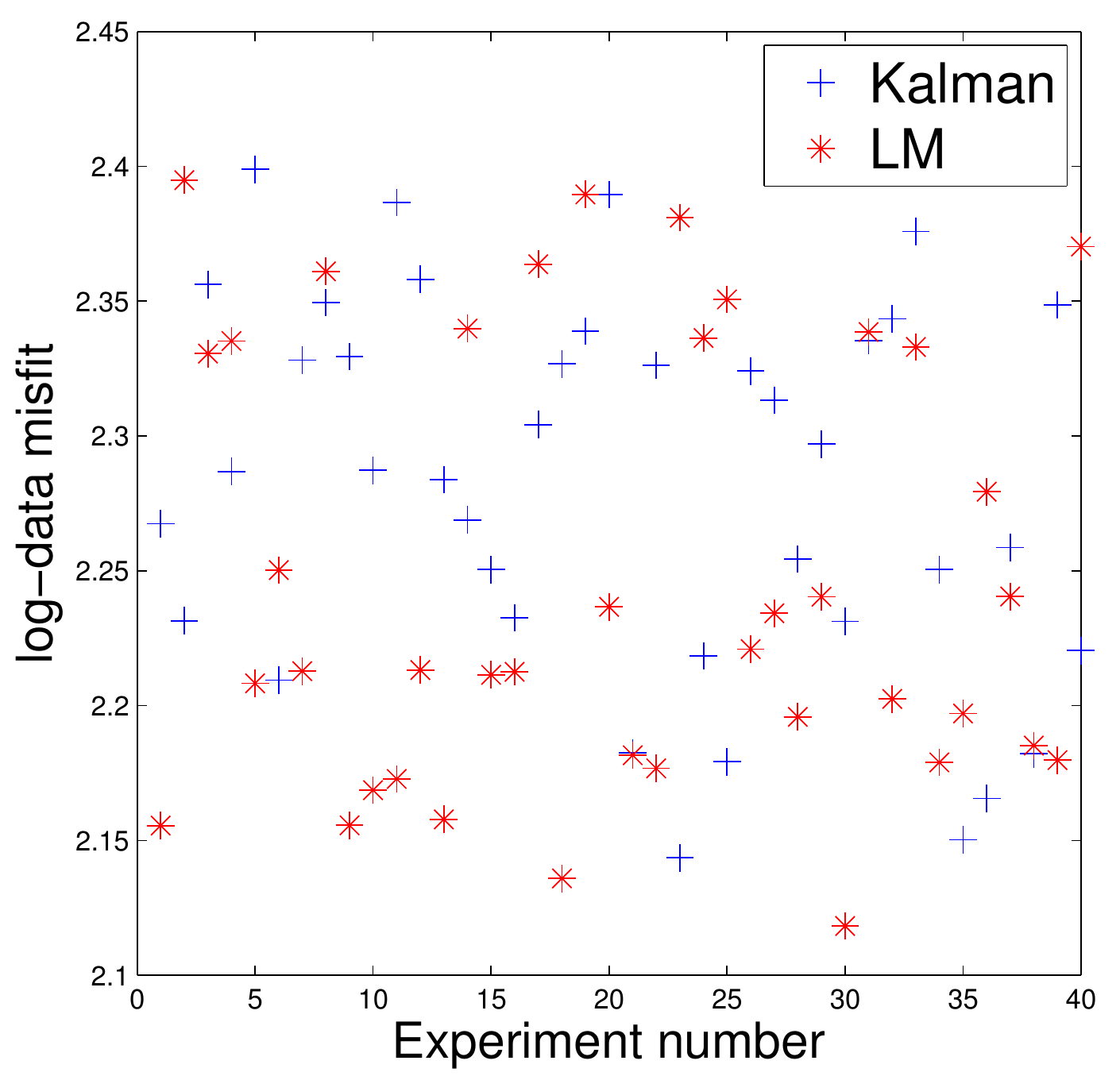}
\includegraphics[scale=0.45]{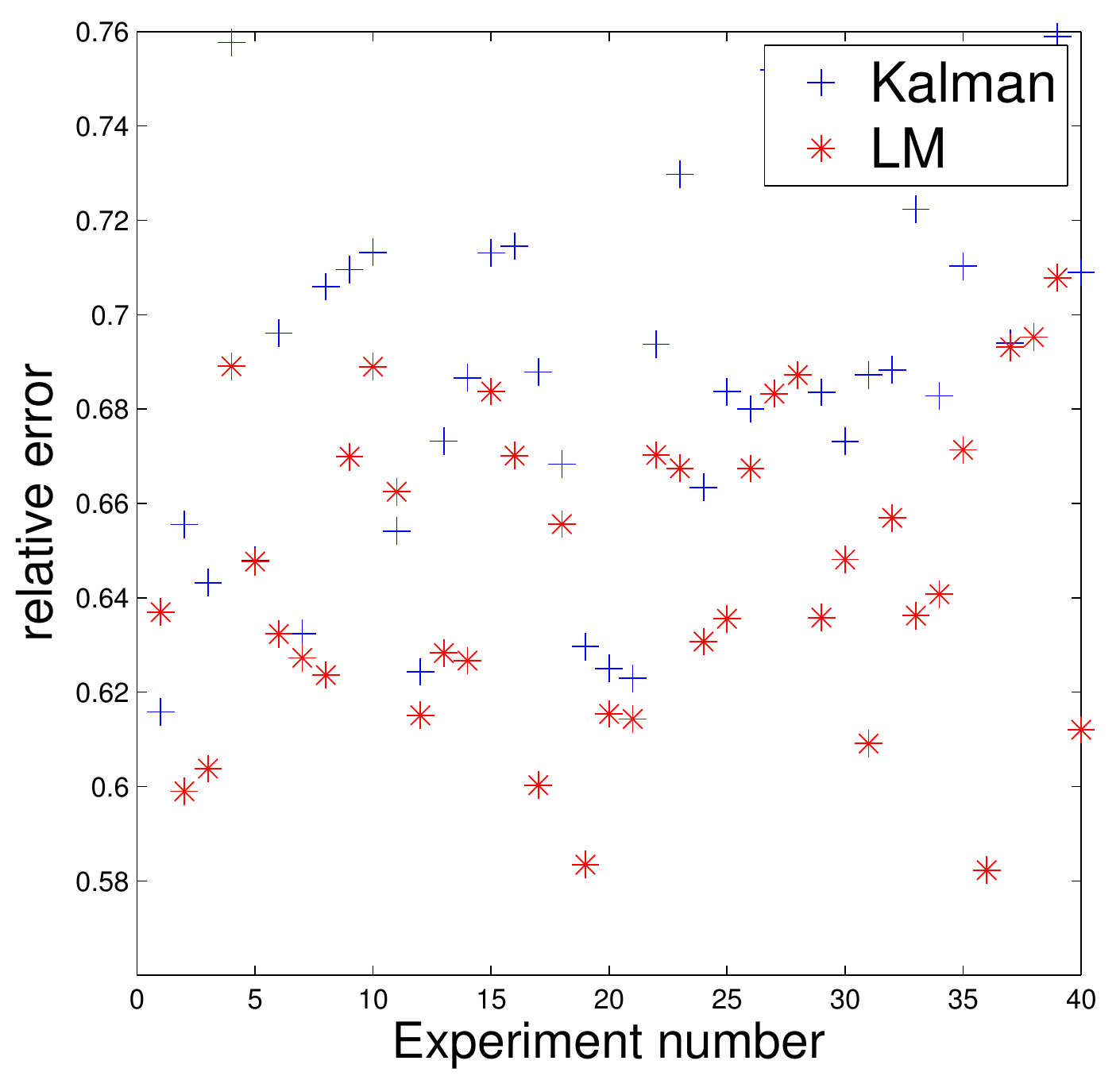}
 \caption{Log-data misfit (left) and rel. error w.r.t truth (right) from estimates obtained with the regularizing LM scheme and the proposed ensemble method (Results from 50 experiments with different initial ensembles of size $N_{e}=150$).}
    \label{Fig9}
\end{center}
\end{figure}

In Figure \ref{Fig9} we display the numerical results from the application of the proposed scheme to identify the log-conductivity described in the preceding subsections. As before, we used 40 different experiments corresponding to different initial ensembles of size $N_{e}=150$. In this case, however, we compare with the regularizing LM scheme applied to the solution of the parameter identification but constrained to the subspace generated by the initial ensemble. The numerical results displayed in Figure \ref{Fig9} suggests that the proposed ensemble Kalman method produces a derivative-free approximation of the regularizing LM scheme where the Fr\'echet derivative is approximated by the ensemble covariance. Moreover, these results show that the proposed ensemble Kalman algorithm is minimizing the data misfit is a stable fashion. 

\subsection{Test Model II: EIT}\label{num_EIT}

\begin{figure}[htbp]
\begin{center}
\includegraphics[scale=0.3]{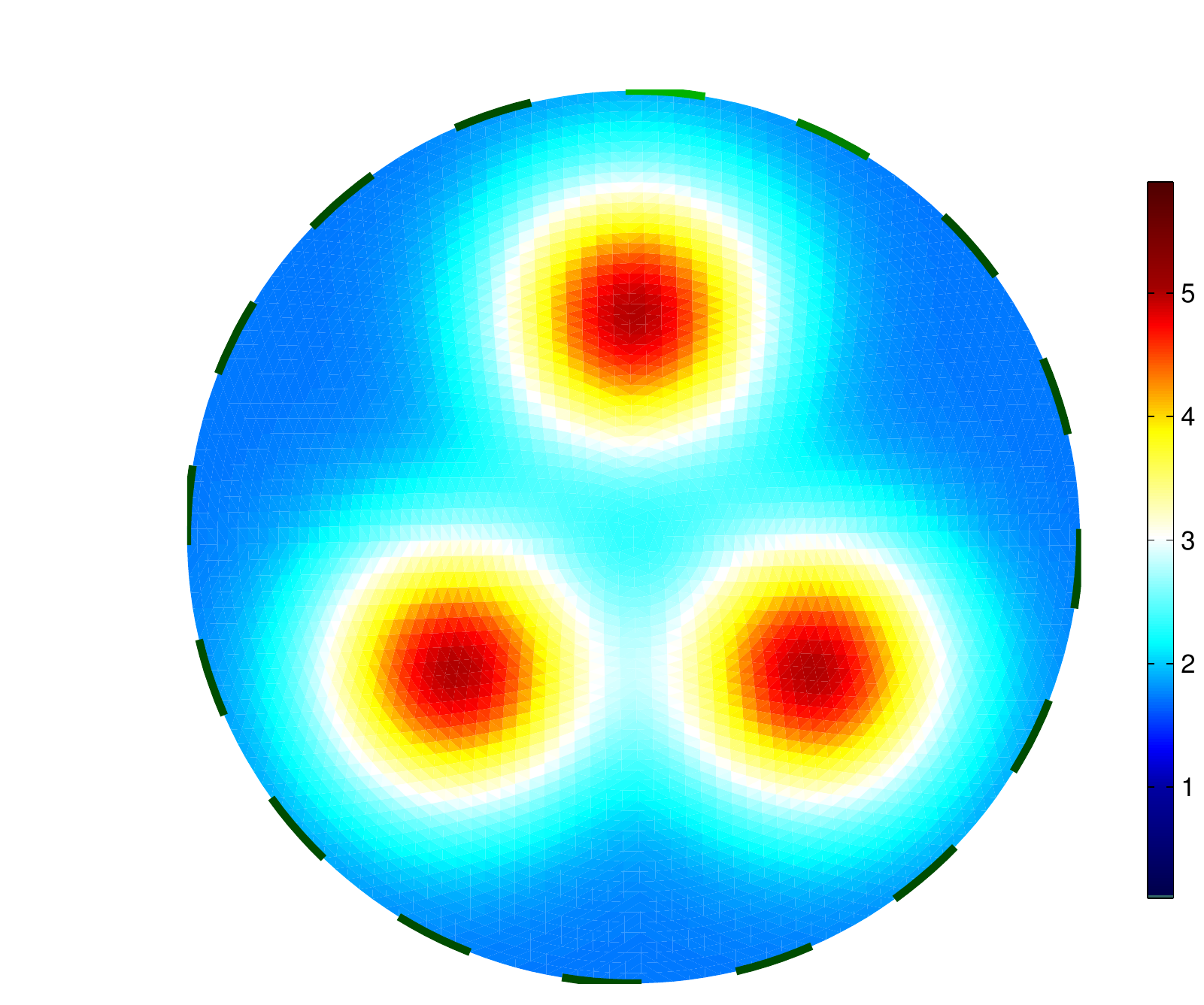}
\includegraphics[scale=0.3]{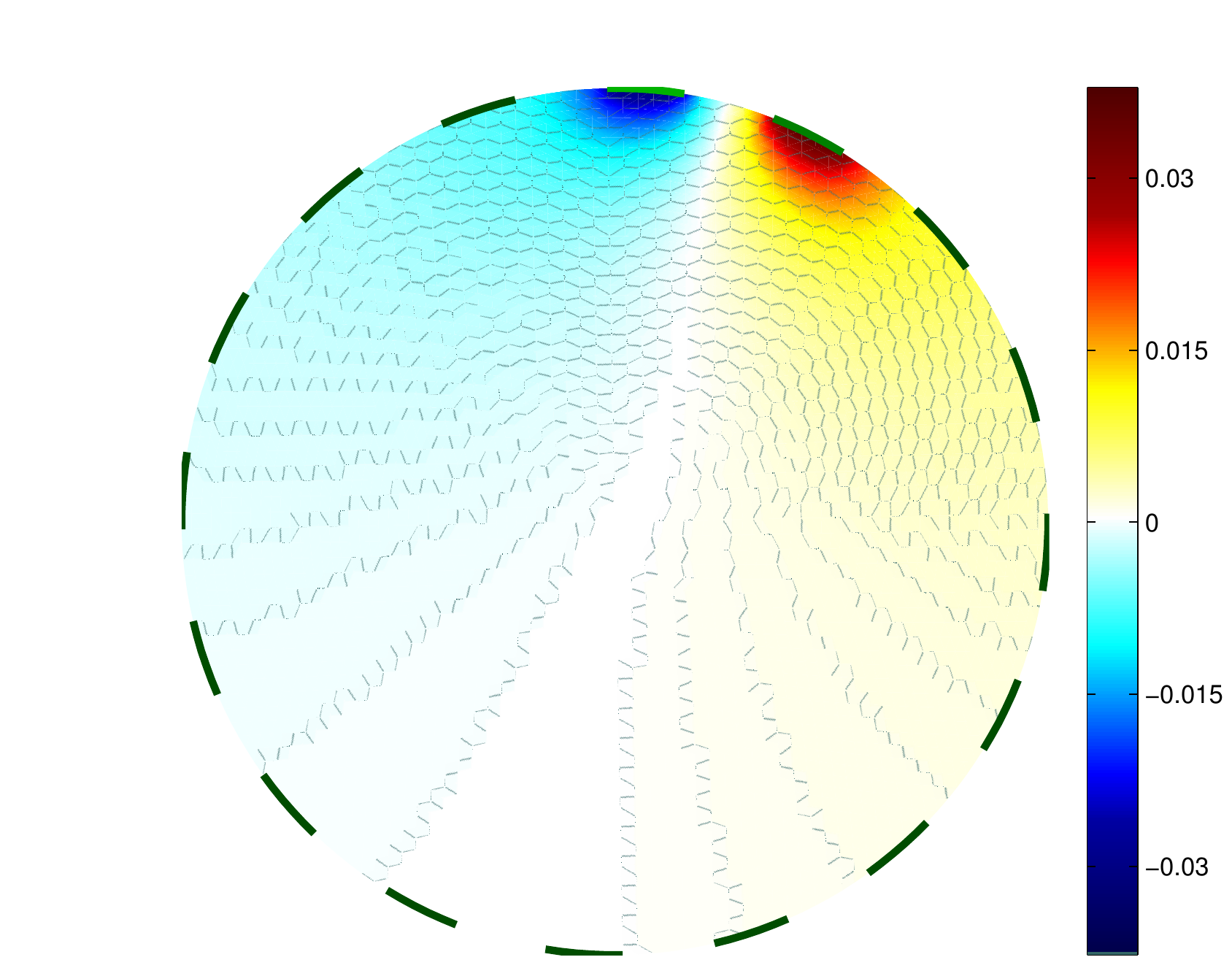}
\includegraphics[scale=0.3]{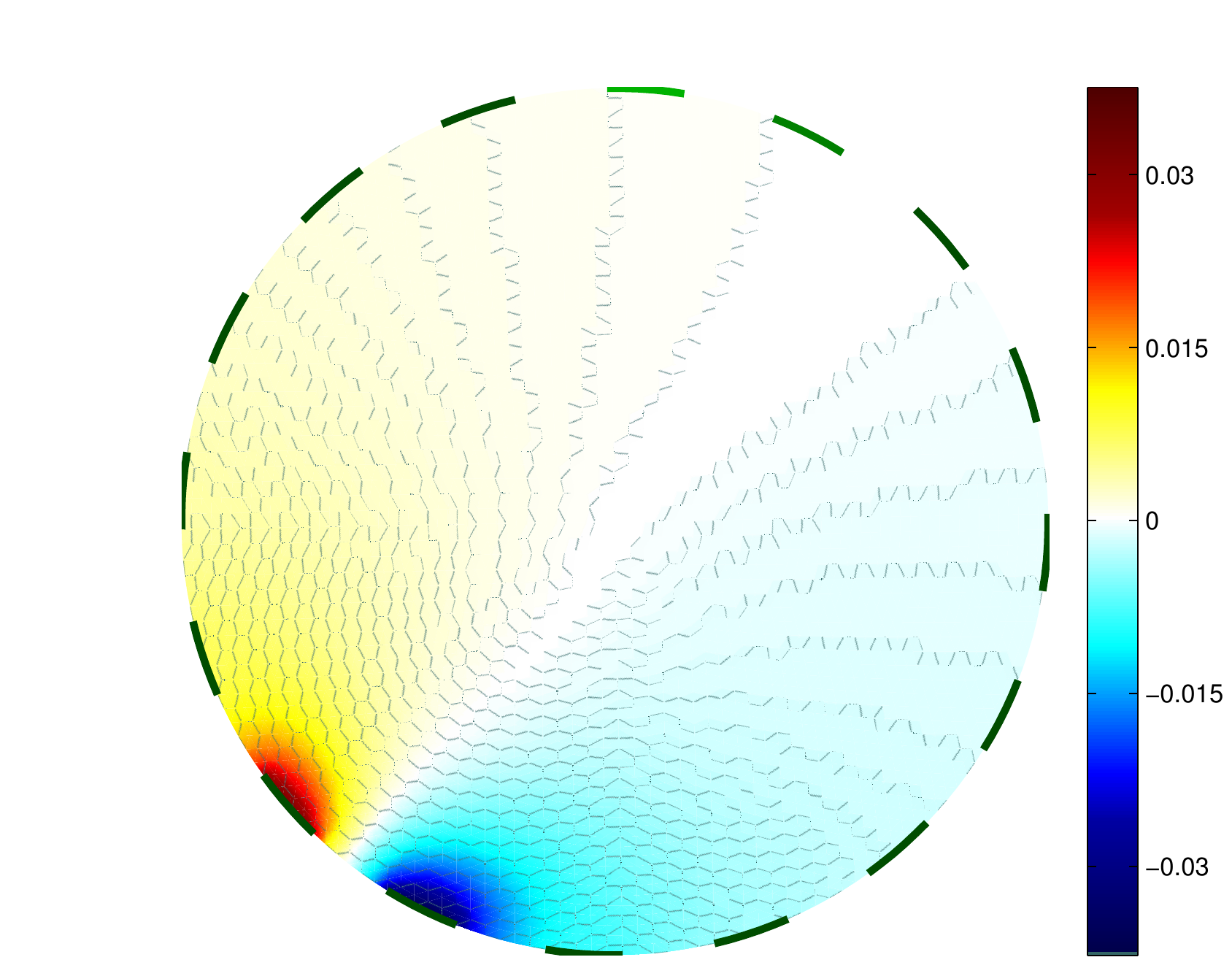}

 \caption{Left: true log-conductivity. Middle and Right: Two of the 15 voltage simulations patterns computed with the true conductivity. }    \label{FigEIT1}
\end{center}
\end{figure}

In this section we investigate the performance of the regularizing ensemble Kalman method for the solution of the EIT problem described in subsection \ref{EIT}. Similar to the previous subsection, our aim is to understand fundamental aspects of the regularizing properties of the scheme. In particular, in this subsection our aim is to observe the effect of the selection of the initial ensemble on the accuracy and cost of the proposed technique. In contrast to the previous subsection where the true unknown field was a random draw from a distribution that we used as well for the generation of the initial ensemble, here we prescribe a conductivity $\kappa$ that consist of three cicular inclusions on a circular domain $D$ of unit radius similar to the one used in \cite{WM}. This true conductivity is displayed in \Fref{FigEIT1} along with the configuration of 16 electrodes used in the CEM described in subsection \ref{EIT}. We consider 15 current patterns where current is injected between a pair of adjacent electrodes. For each current pattern we collect measurements on all the electrodes thereby having $M=240$ observations. We choose a value $z_{k}=0.01$ for the contact impedances of all the electrodes.

In the middle and middle-right of \Fref{FigEIT1} we show some of the true voltages obtained with the FEM solver from EIDORS MATLAB framework \cite{EIDORS}. Synthetic data were generated from the aforementioned voltages by adding Gaussian random noise of standard deviation of $2\%$ of the signal. Inverse crimes are avoided by using a finer mesh (with 7744 elements) than the one used for the application of Algorithm \ref{Al1} (mesh of 6400 elements). The experiments in this subsection are focused on the electrode and measurement configuration described earlier. 

\begin{figure}[htbp]
\begin{center}
\includegraphics[scale=0.18]{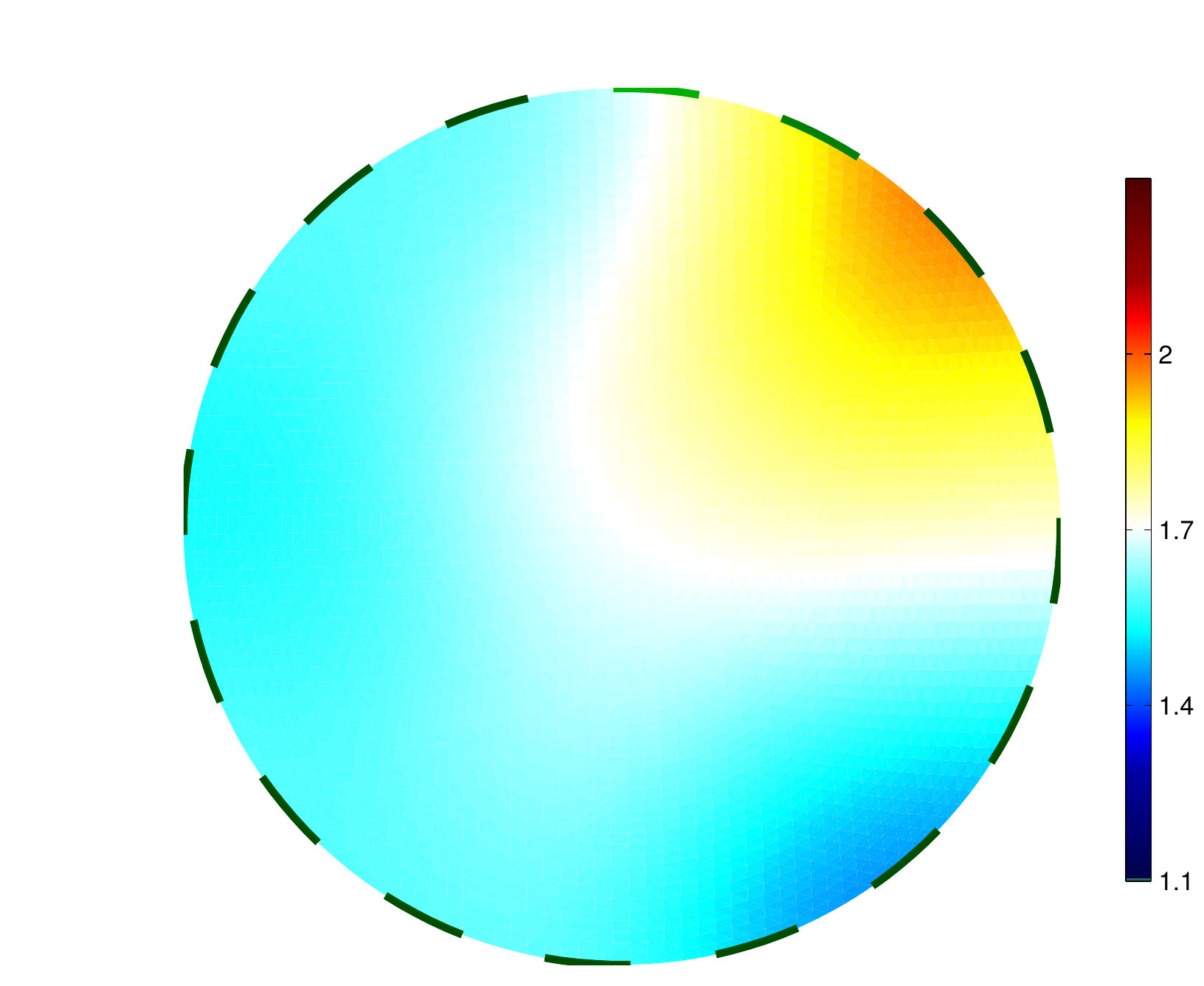}
\includegraphics[scale=0.18]{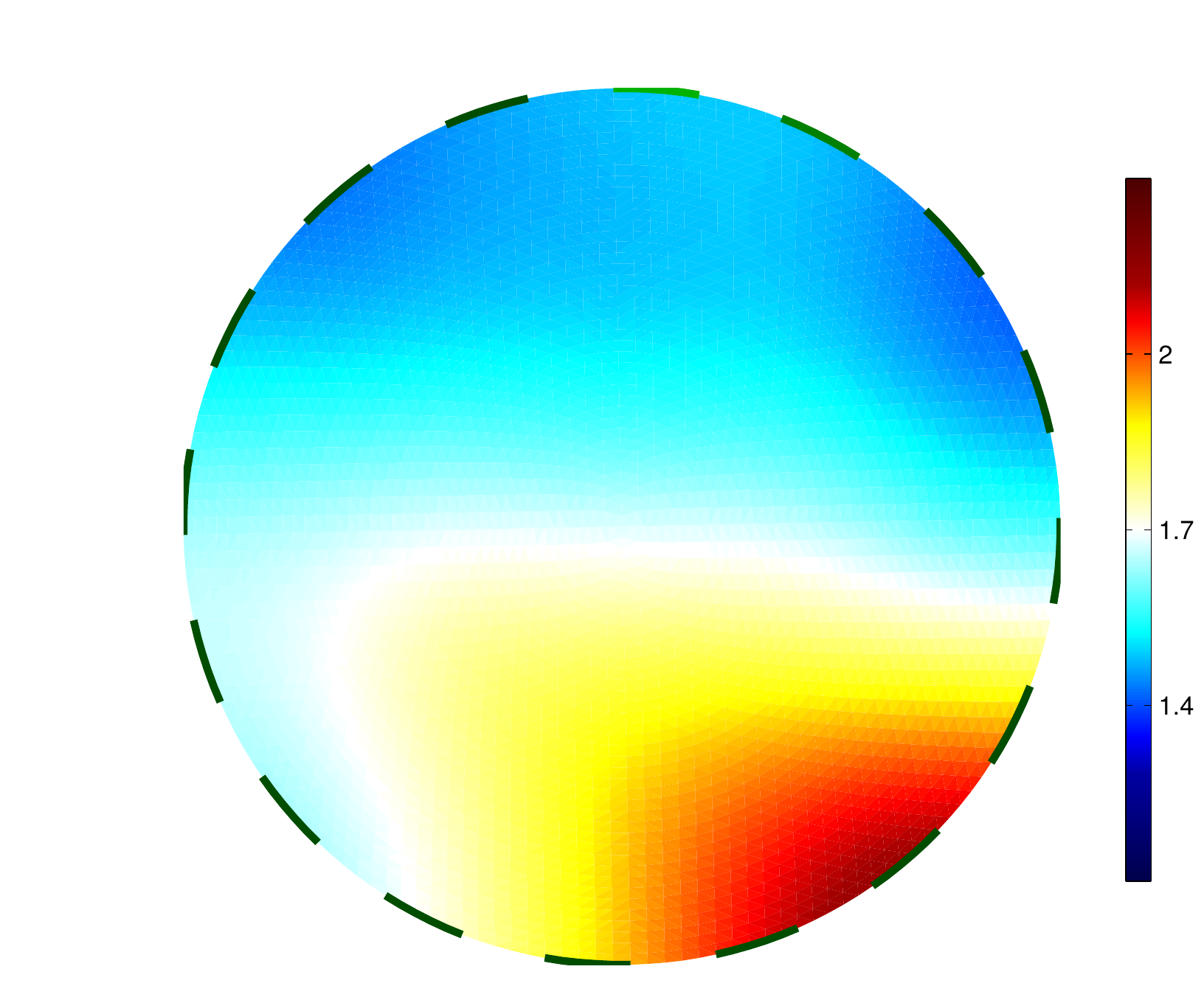}
\includegraphics[scale=0.18]{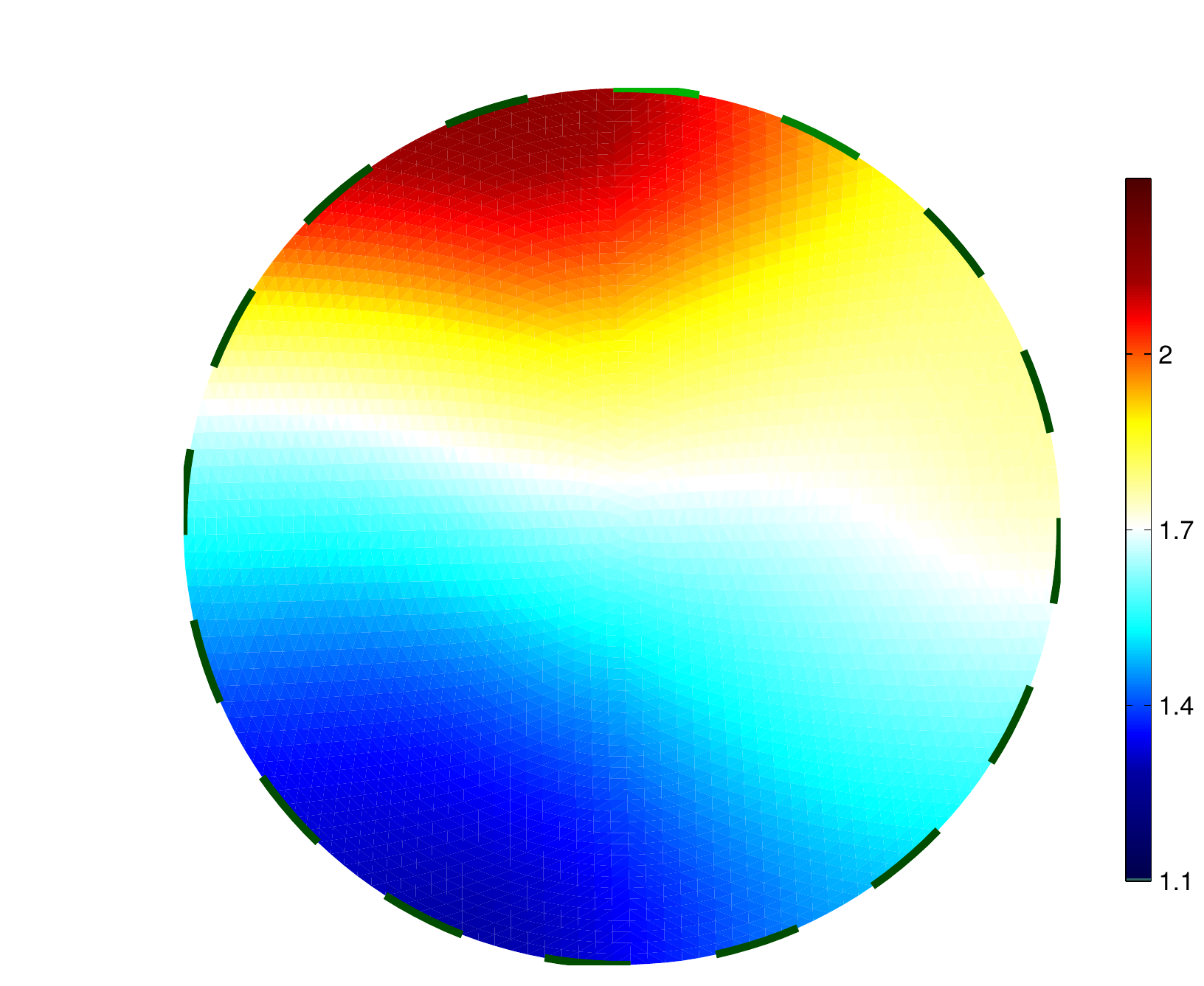}
\includegraphics[scale=0.18]{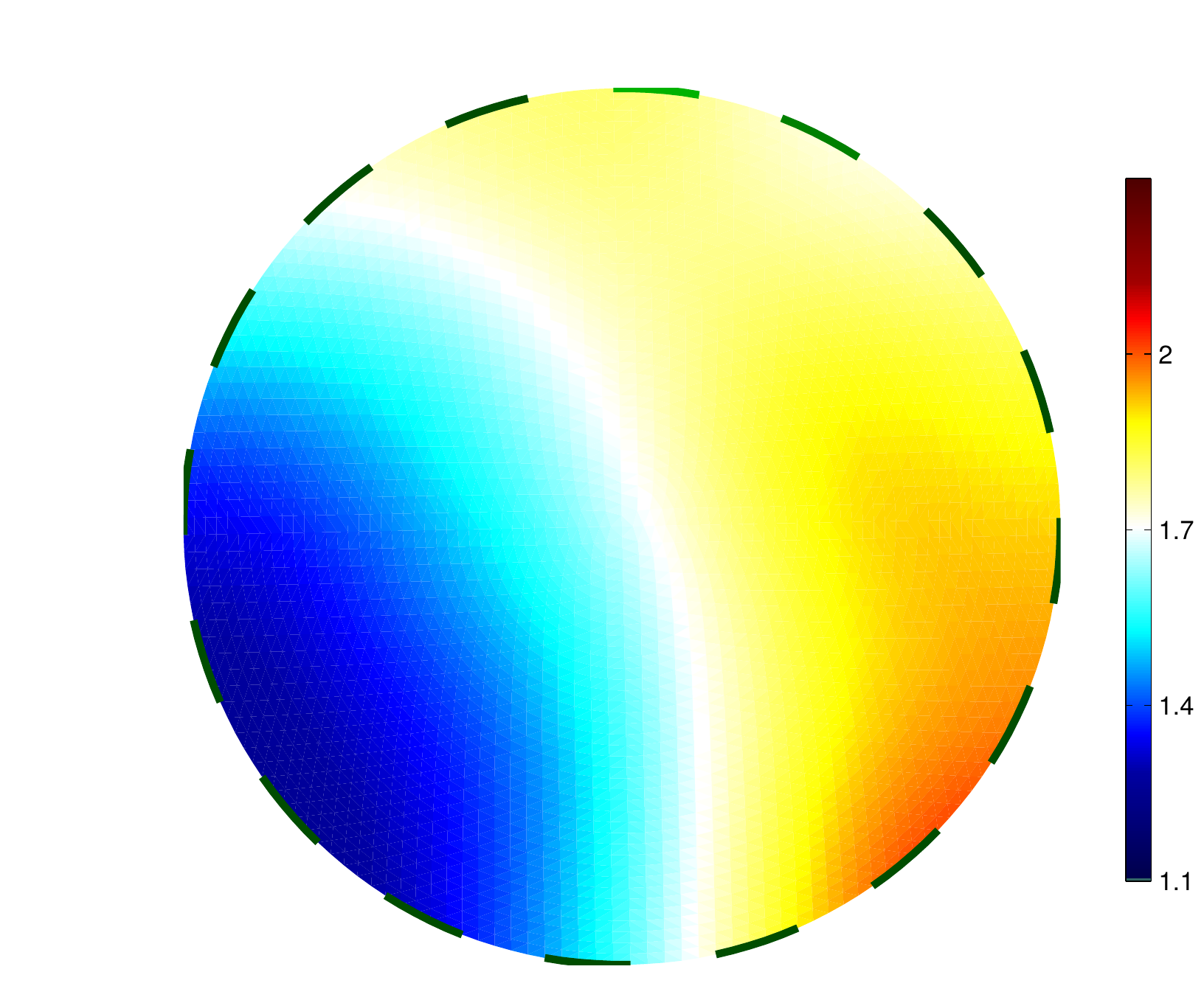}
\includegraphics[scale=0.18]{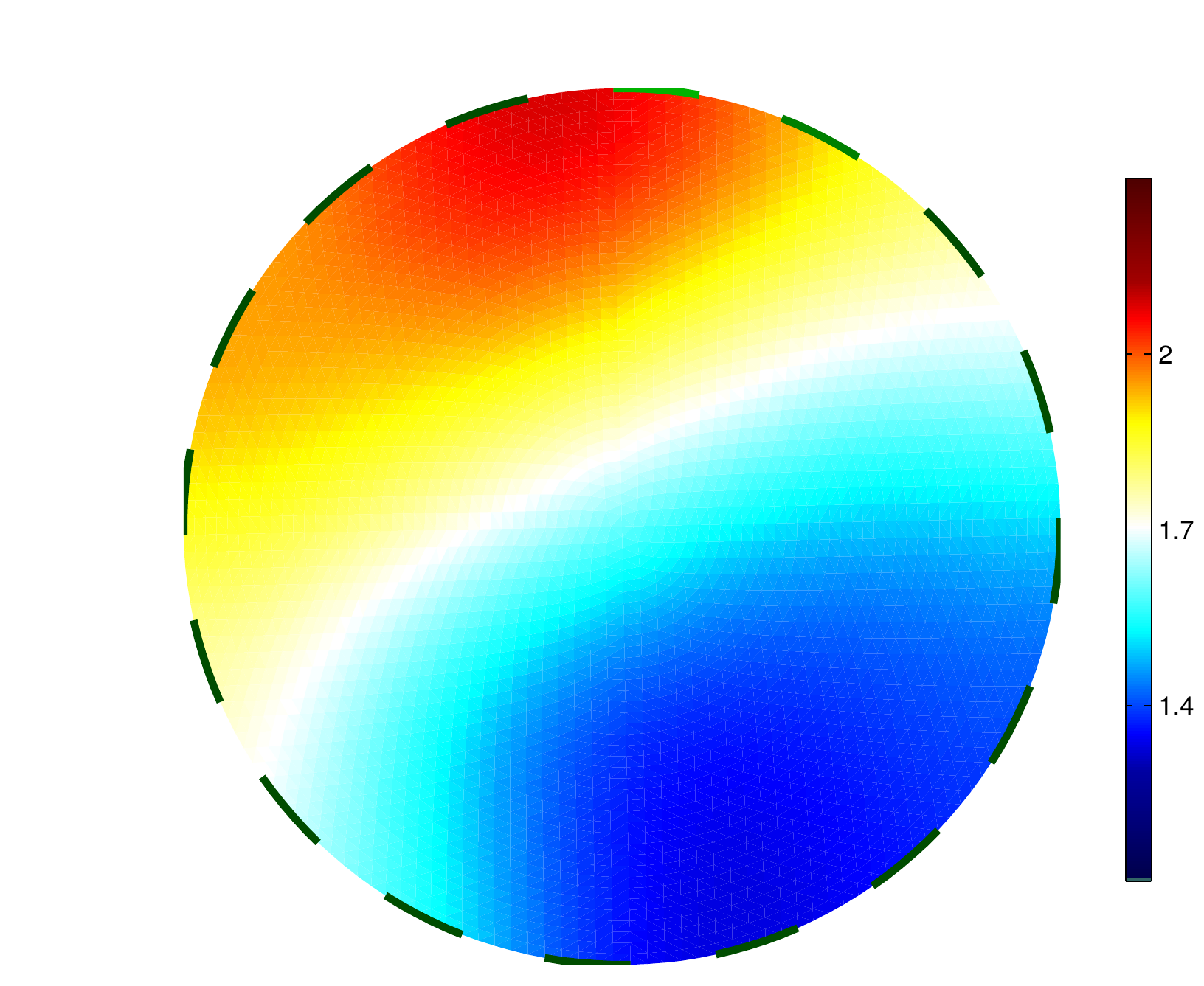}
 \caption{Draws from a Gaussian prior of the form (\ref{eq:cova1}) with $L=0.2$ and $\theta=5$. }
 \label{FigEIT2}
\end{center}
\end{figure}

Note that the experiments from this subsection comprise the more general case where there is no characterization of the truth in terms of a probability distribution. In this case, we still consider an initial ensemble generated from samples from a Gaussian measure $N(0,C)$. However, this probability distribution is completely artificial and defined only for the purpose of generating an initial ensemble. In other words, we assume no link between the truth and our choice of initial ensemble. We consider the covariance  operator that arises from the Wittle-Matter correlation function \cite{WM}. In 2D, such covariance operator takes the form
  \begin{equation}\label{eq:cova1}
C =\omega L^2( 1-L^2 \Delta)^{-\theta}
\end{equation}
where $\Delta$ is the Laplacian operator, $L$ is a correlation length. For the purpose of this work we regard $\theta$ as a parameter that controls the regularity of the samples and $\omega$ as a  scaling constant. Some samples from such distribution with parameter $L=0.2$, $\theta=5$ and $\omega=0.1$ are displayed in \Fref{FigEIT2}. As before, samples from $N(0,C)$ are obtained by means of the KL expansion of random fields $N(0,C)$. Samples are generated on a regular rectangular domain that contains the computational domain shown in \Fref{FigEIT1} and projected on the FEM domain used for the computation of the CEM.

\subsubsection{Ensemble size $N_{e}$, tunable parameter $\rho$  and small noise limit.}

\begin{figure}[htbp]
\begin{center}
\includegraphics[scale=0.21]{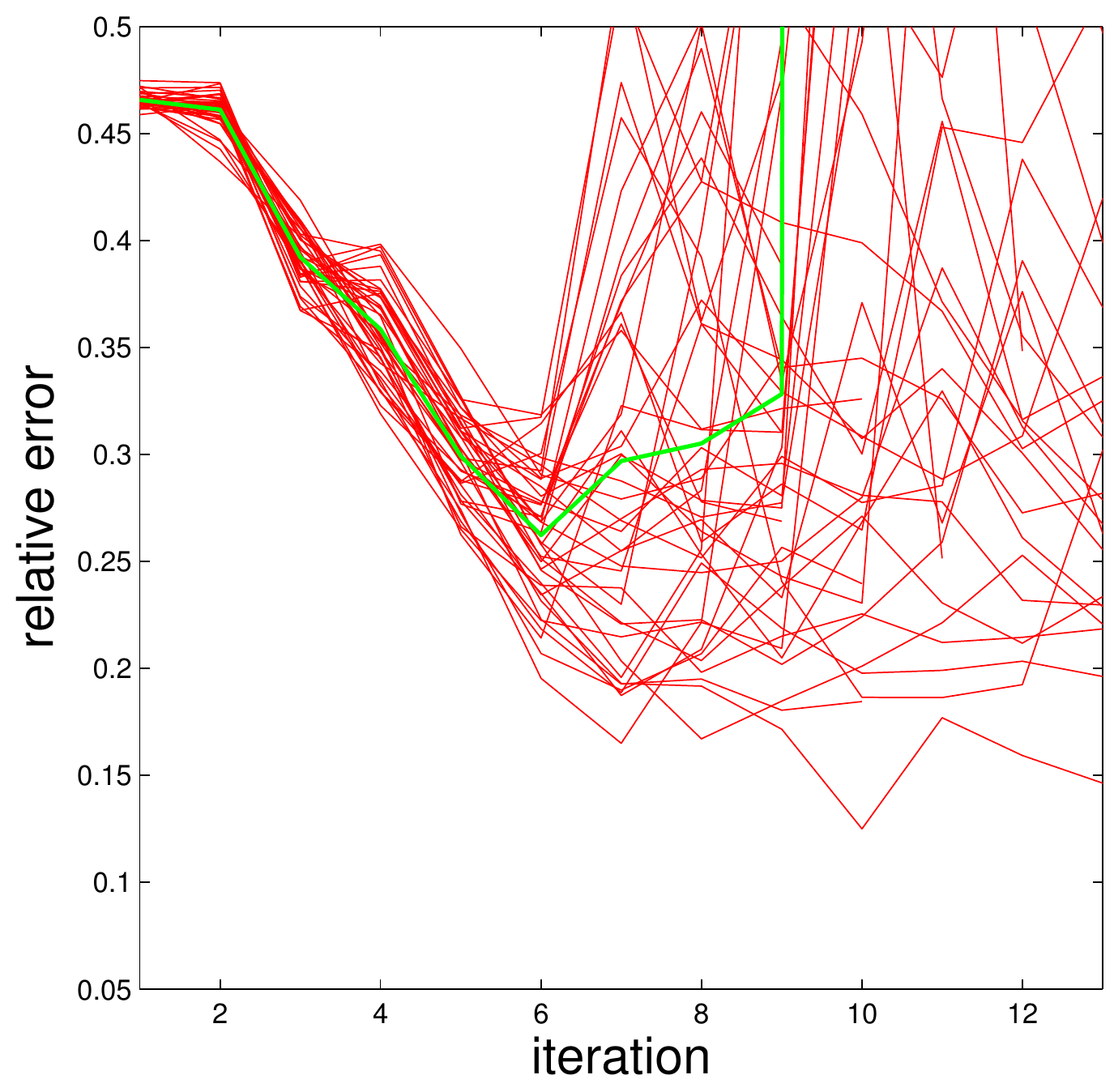}
\includegraphics[scale=0.21]{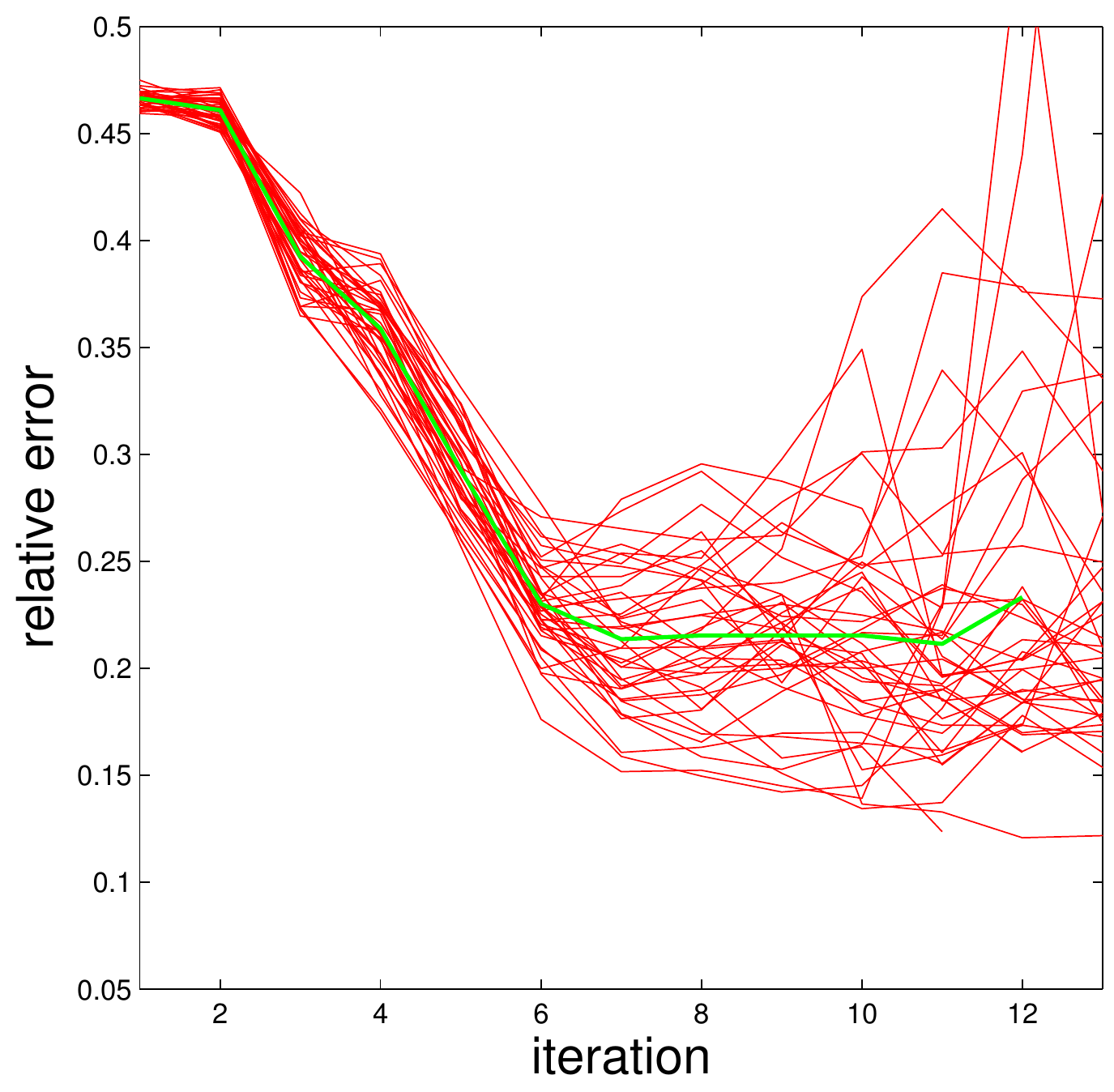}
\includegraphics[scale=0.21]{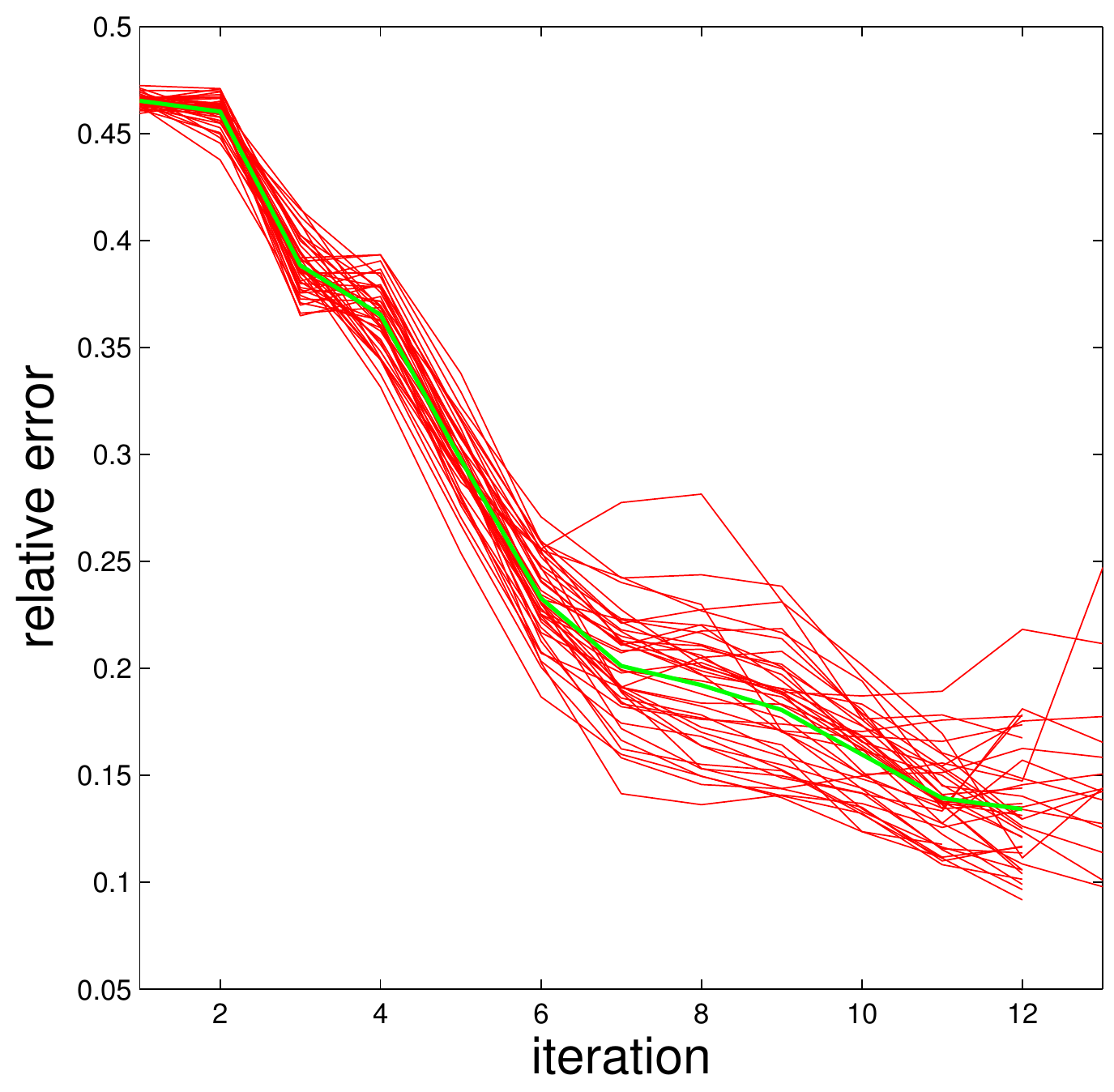}
\includegraphics[scale=0.21]{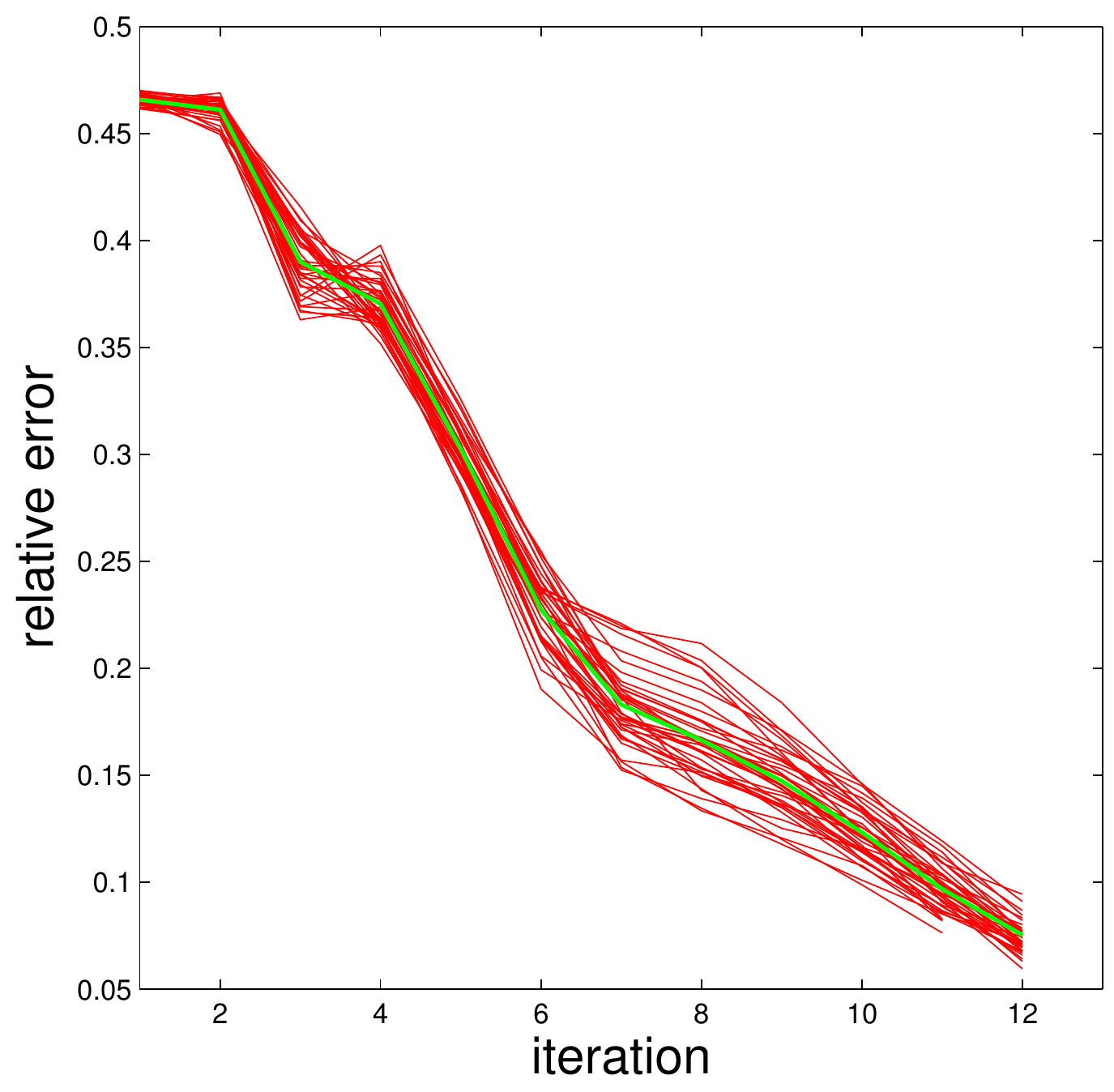}
\includegraphics[scale=0.21]{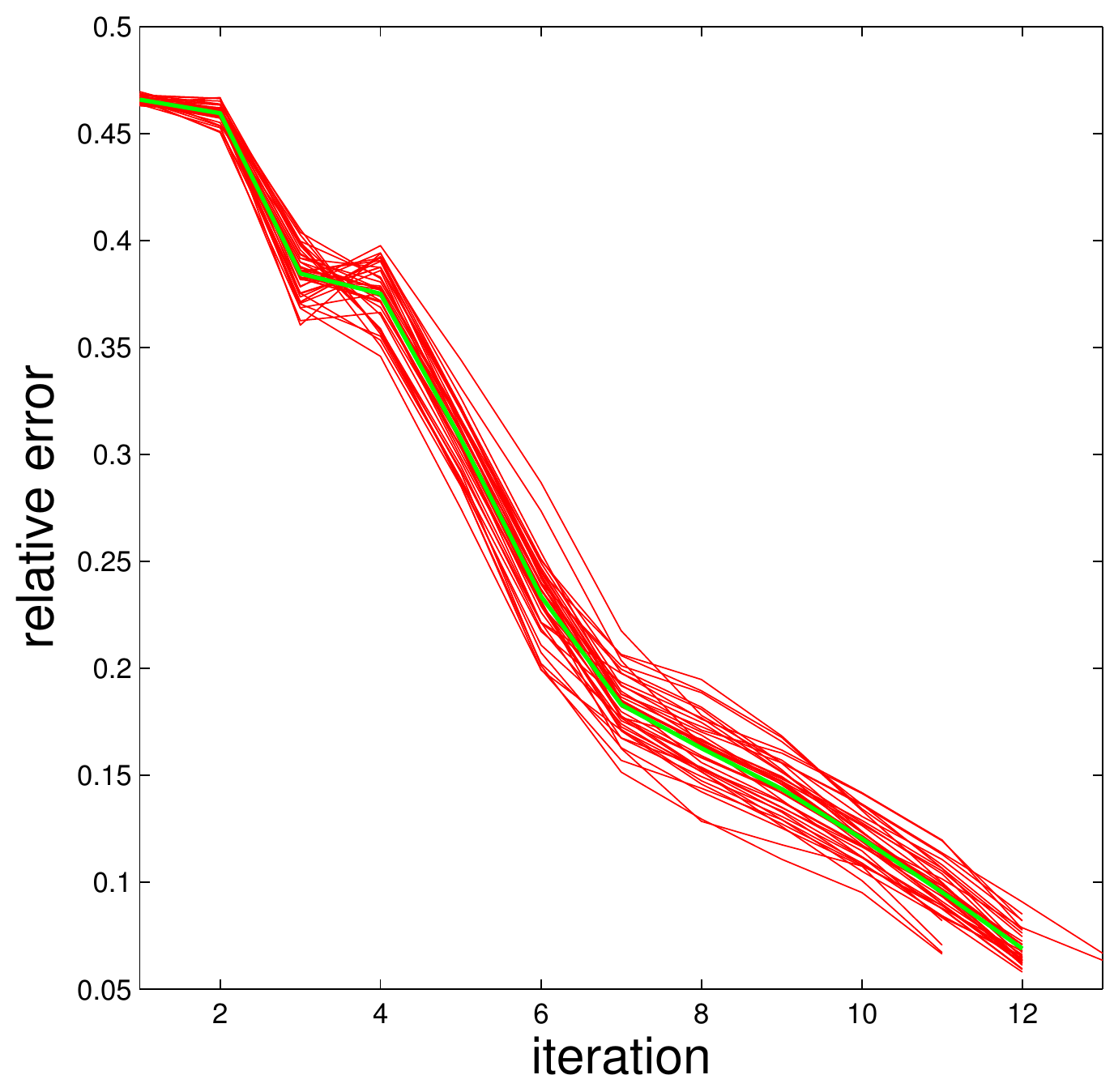}\\
\includegraphics[scale=0.21]{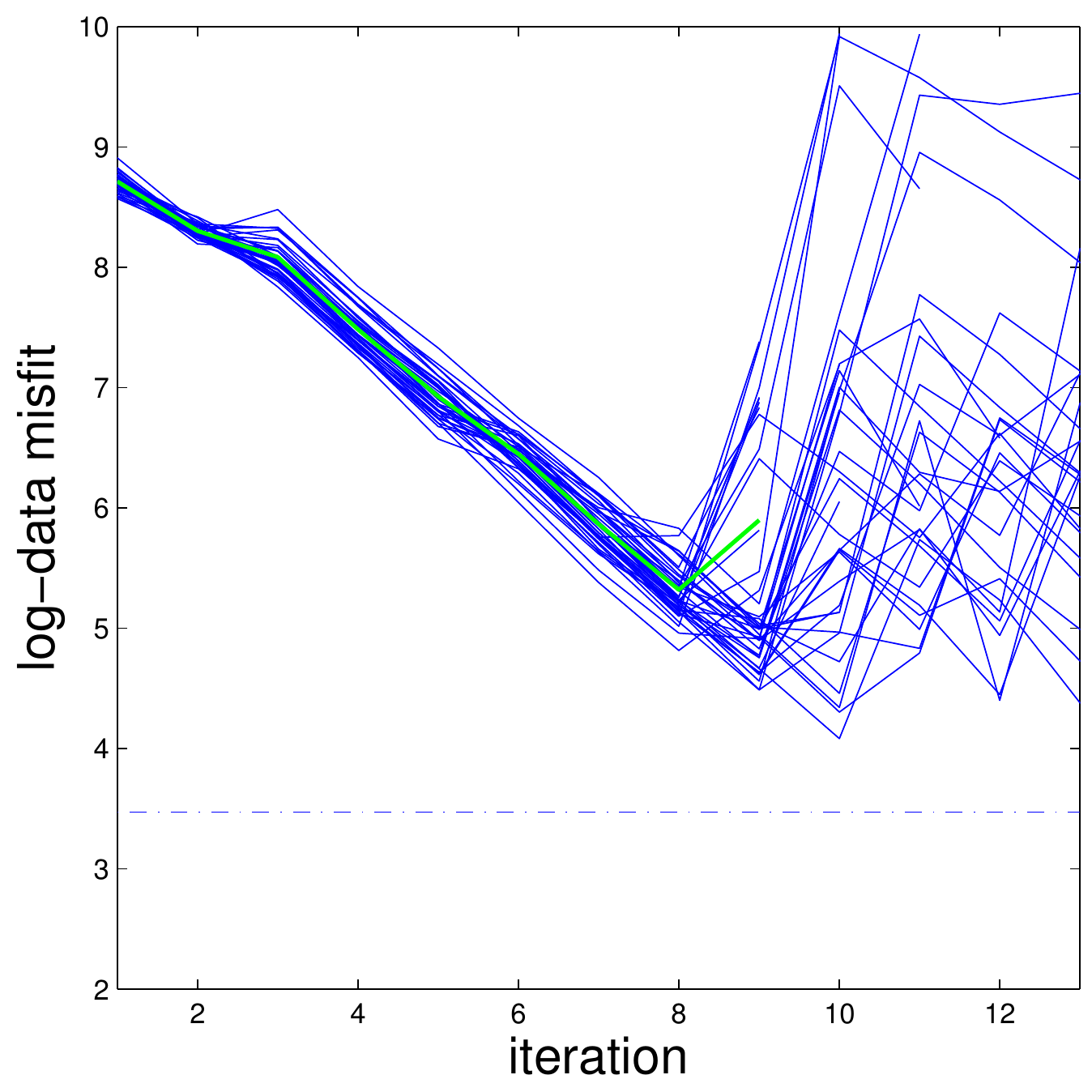}
\includegraphics[scale=0.21]{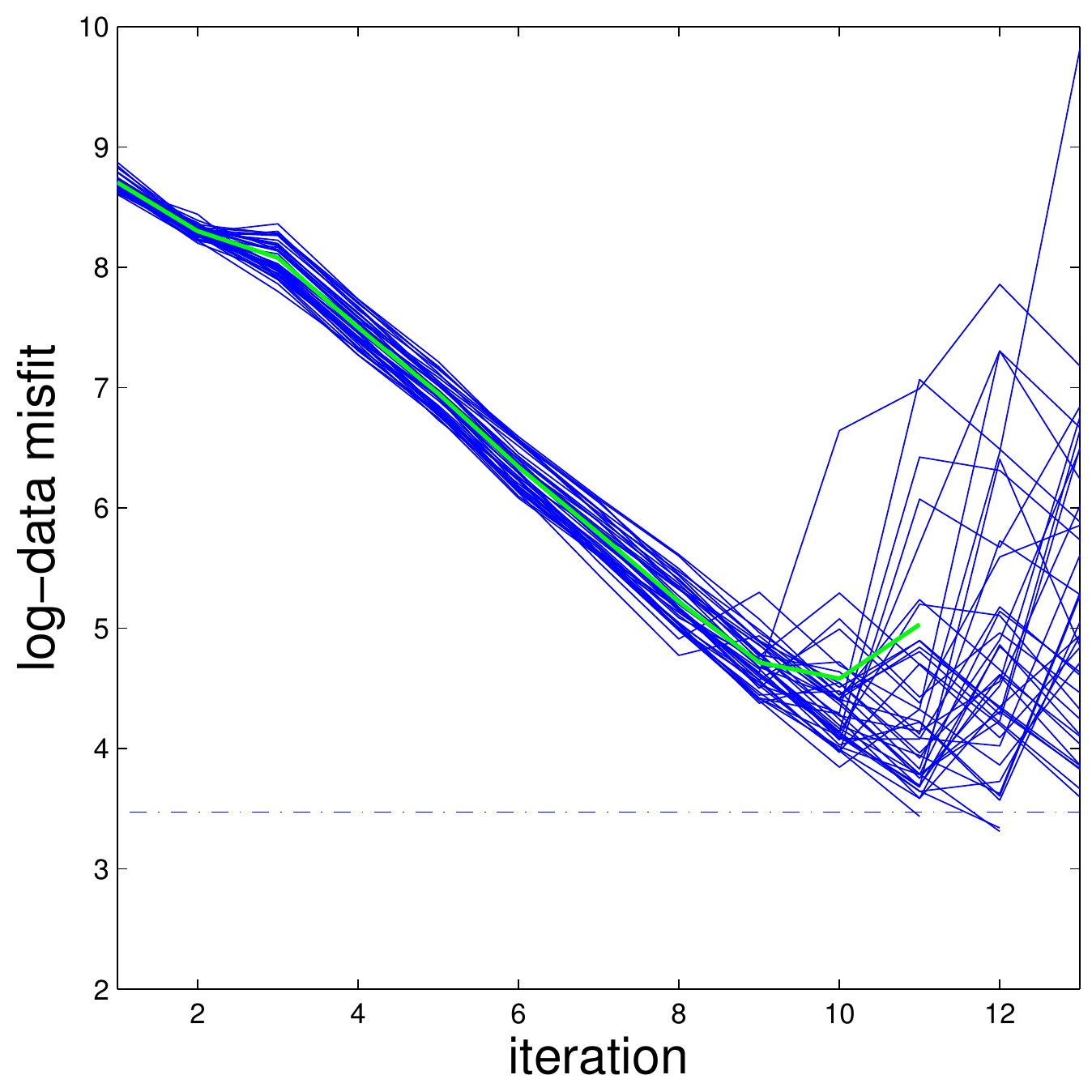}
\includegraphics[scale=0.21]{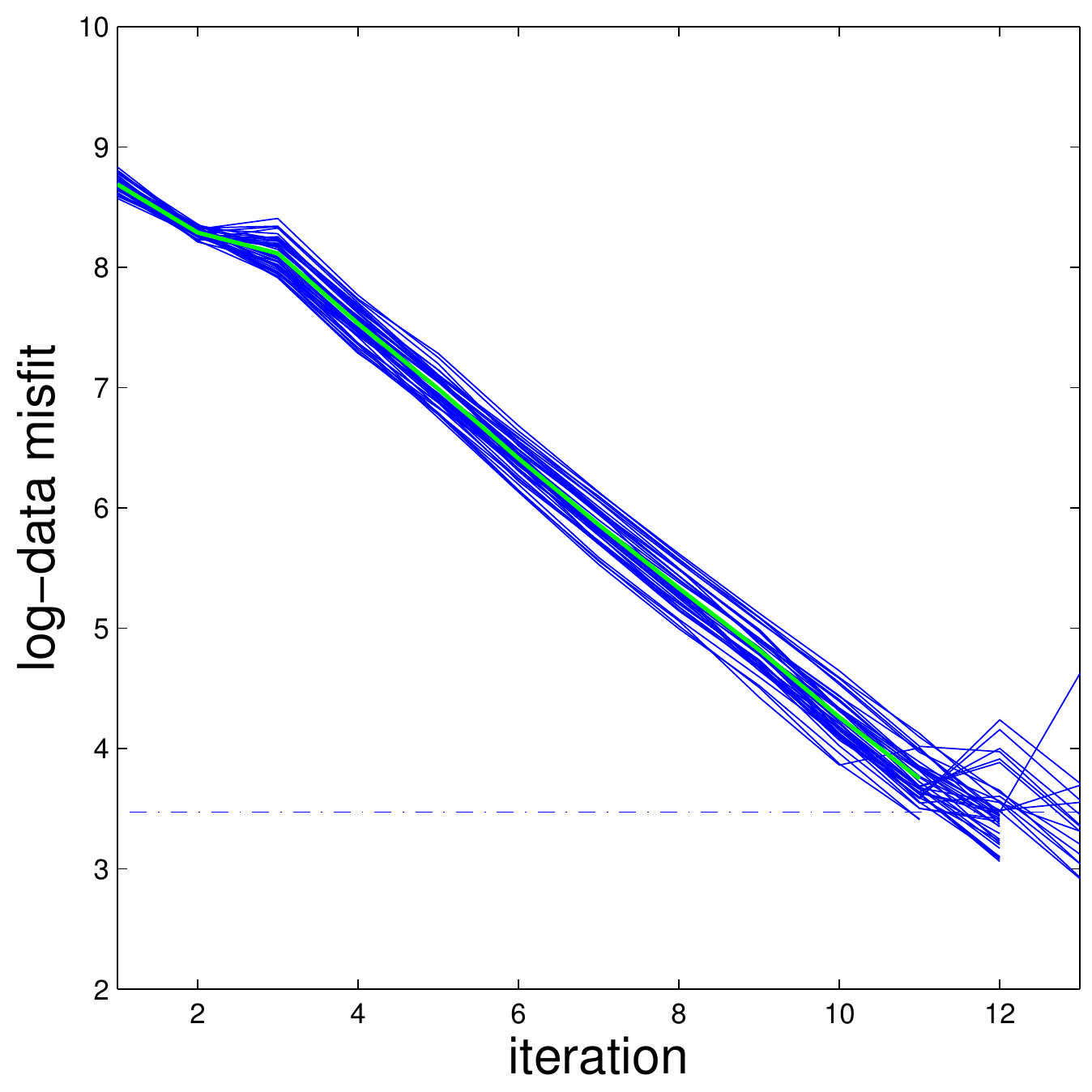}
\includegraphics[scale=0.21]{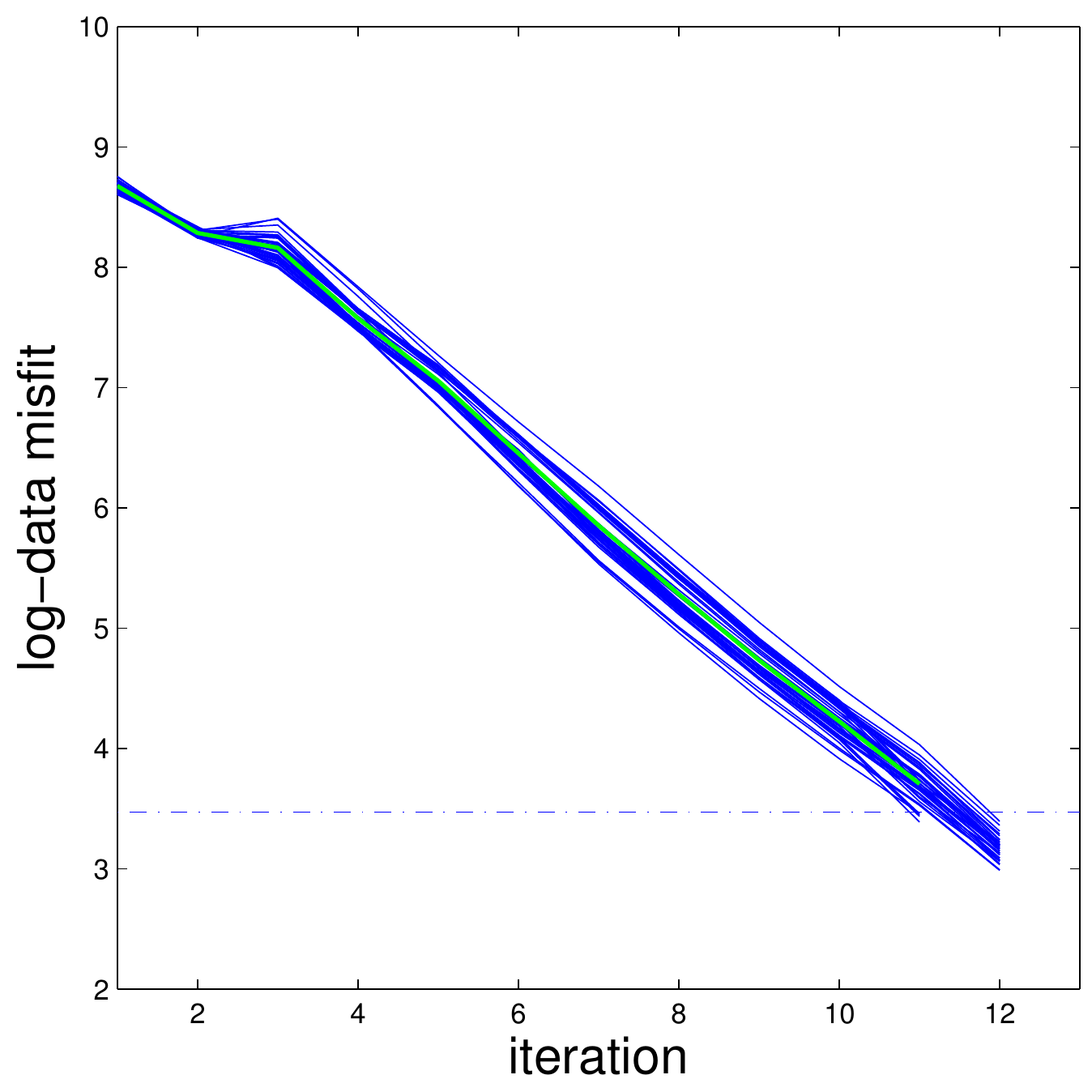}
\includegraphics[scale=0.21]{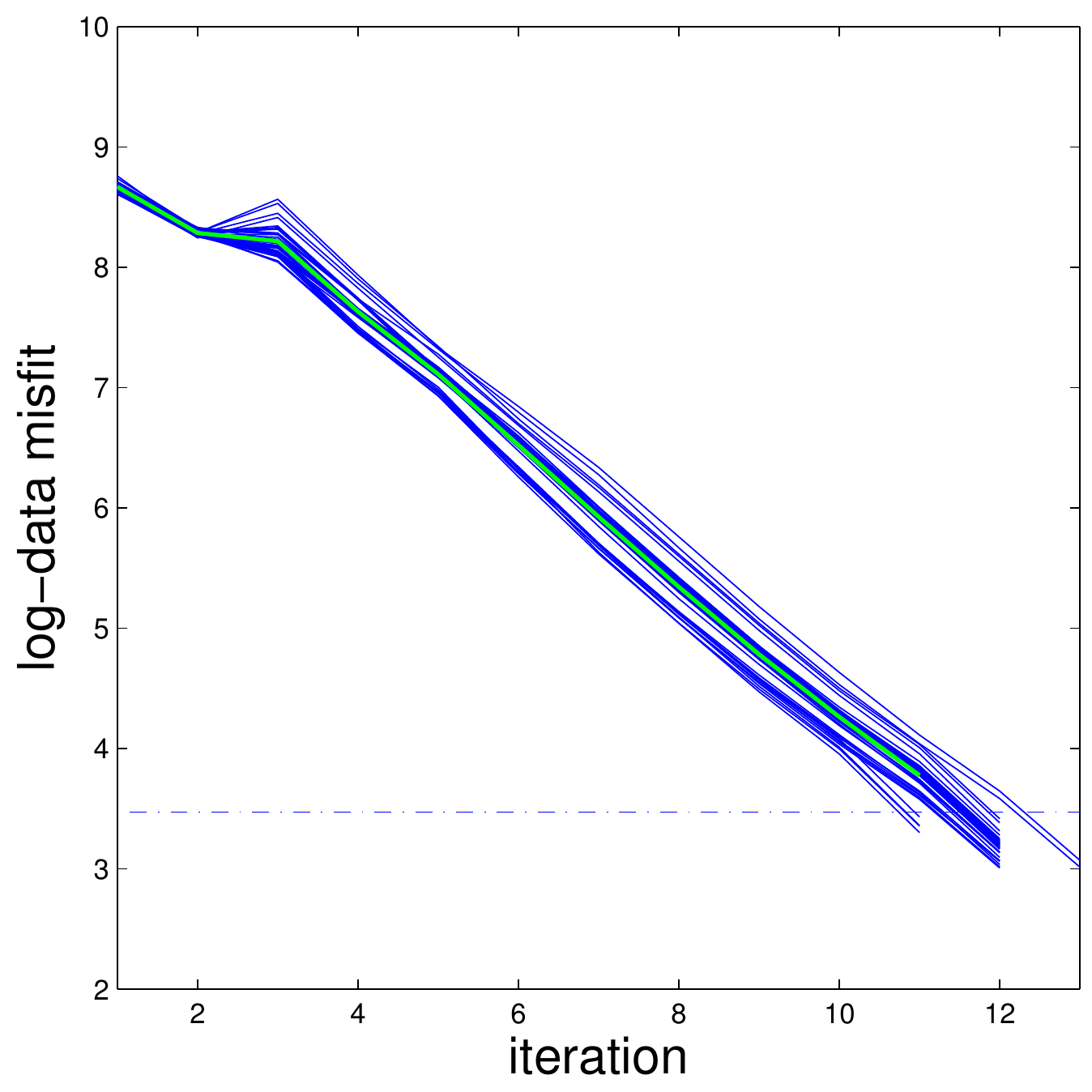}
 \caption{Relative error w.r.t. truth (top) and log - data misfit (bottom) obtained from Algorithm \ref{Al1} with (from left to right) $N_{e}=50, 60, 75, 150,250 $. These plots display the results from 40 experiments with different initial ensembles.}   \label{FigEIT3}

\end{center}
\end{figure}

\begin{figure}[htbp]
\begin{center}
\includegraphics[scale=0.18]{True_EIT}\\
\includegraphics[scale=0.18]{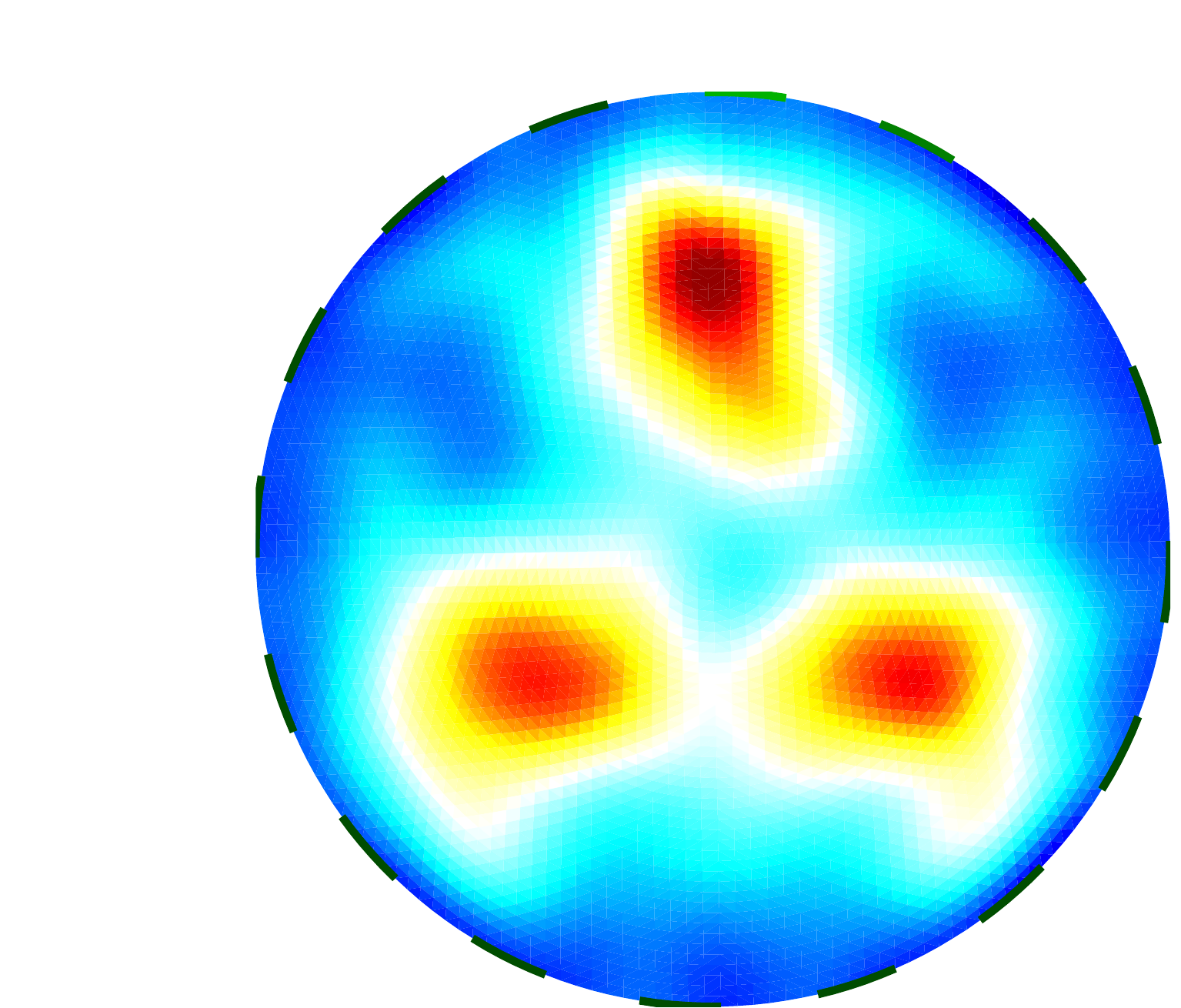}
\includegraphics[scale=0.18]{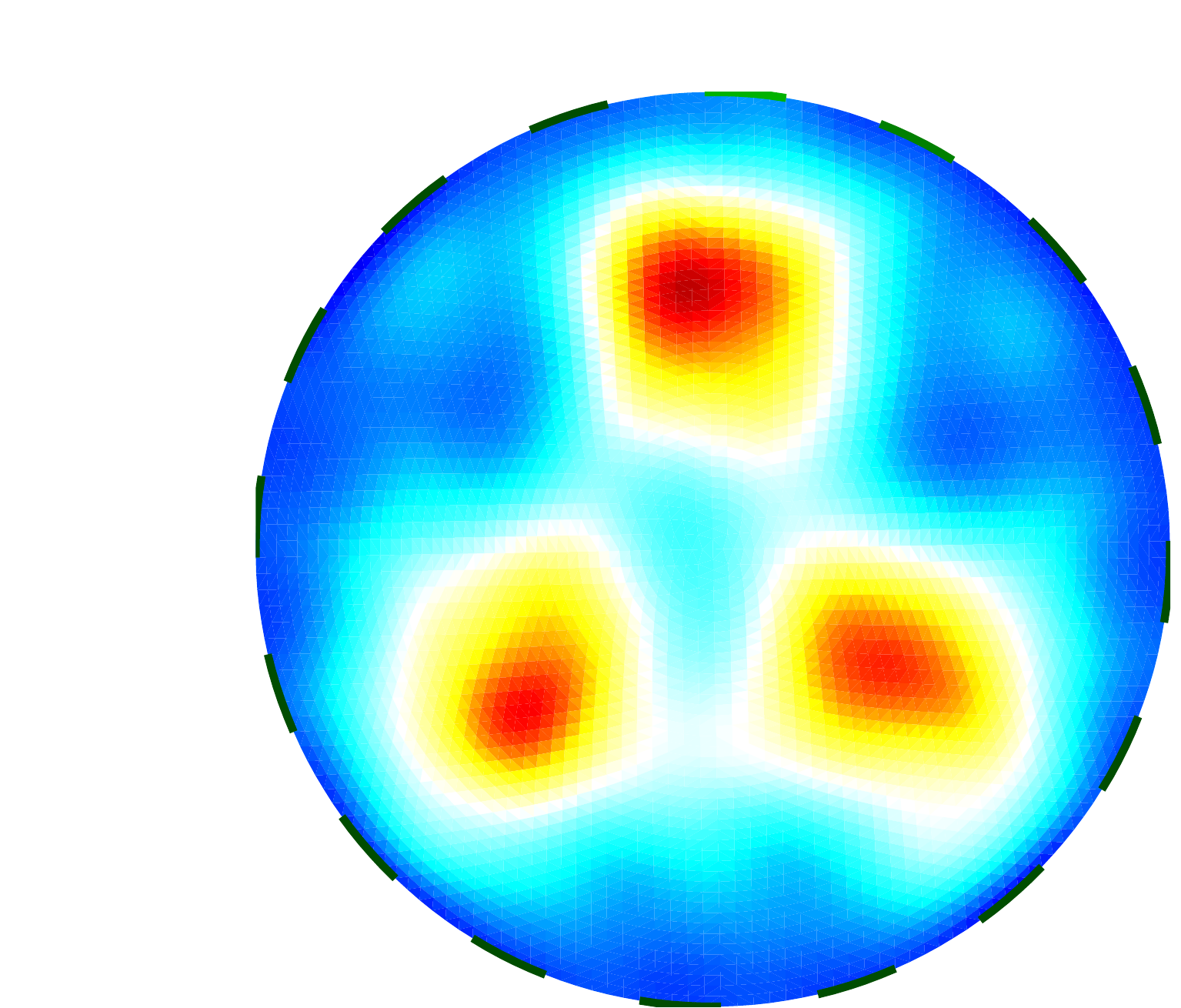}
\includegraphics[scale=0.18]{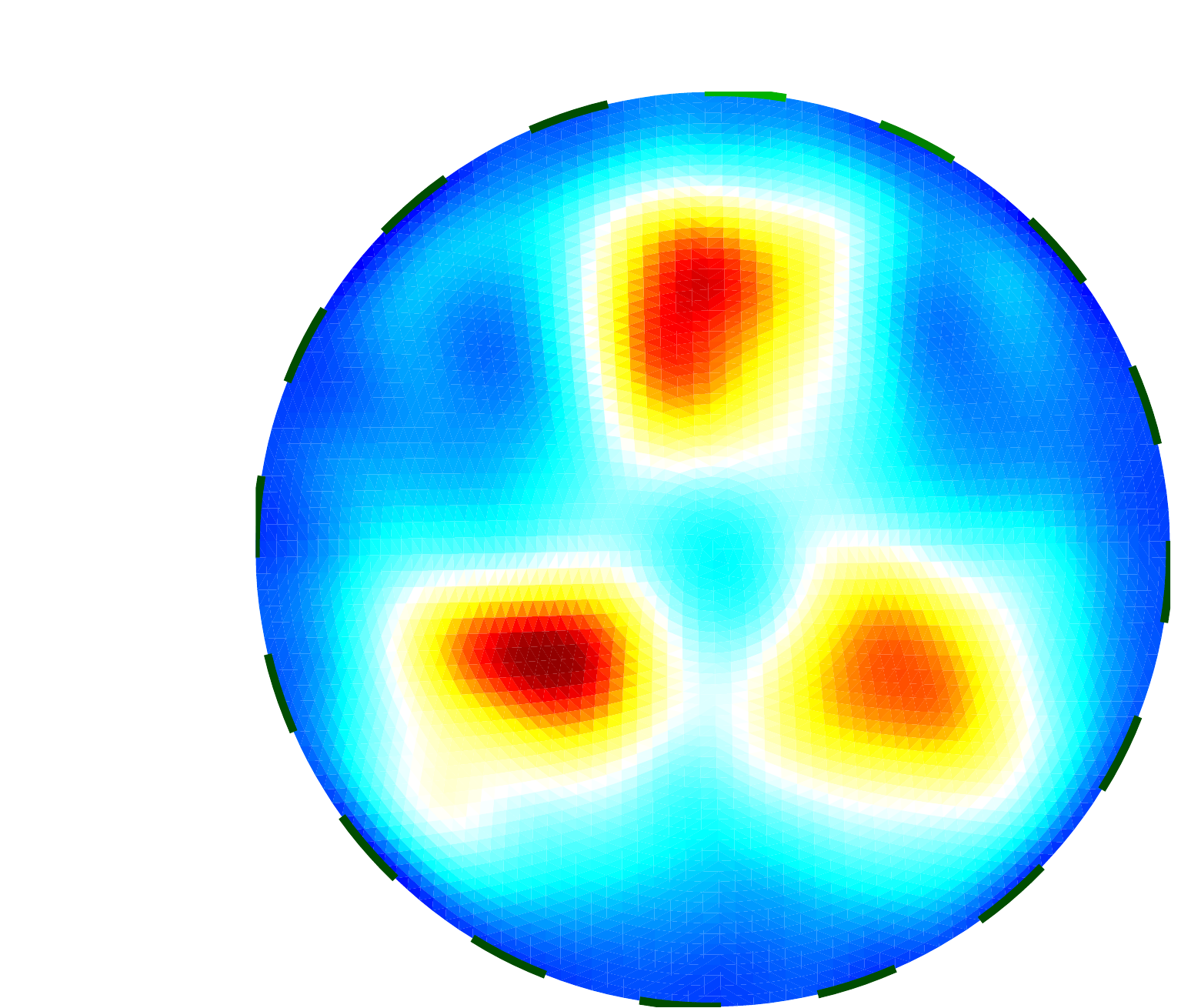}
\includegraphics[scale=0.18]{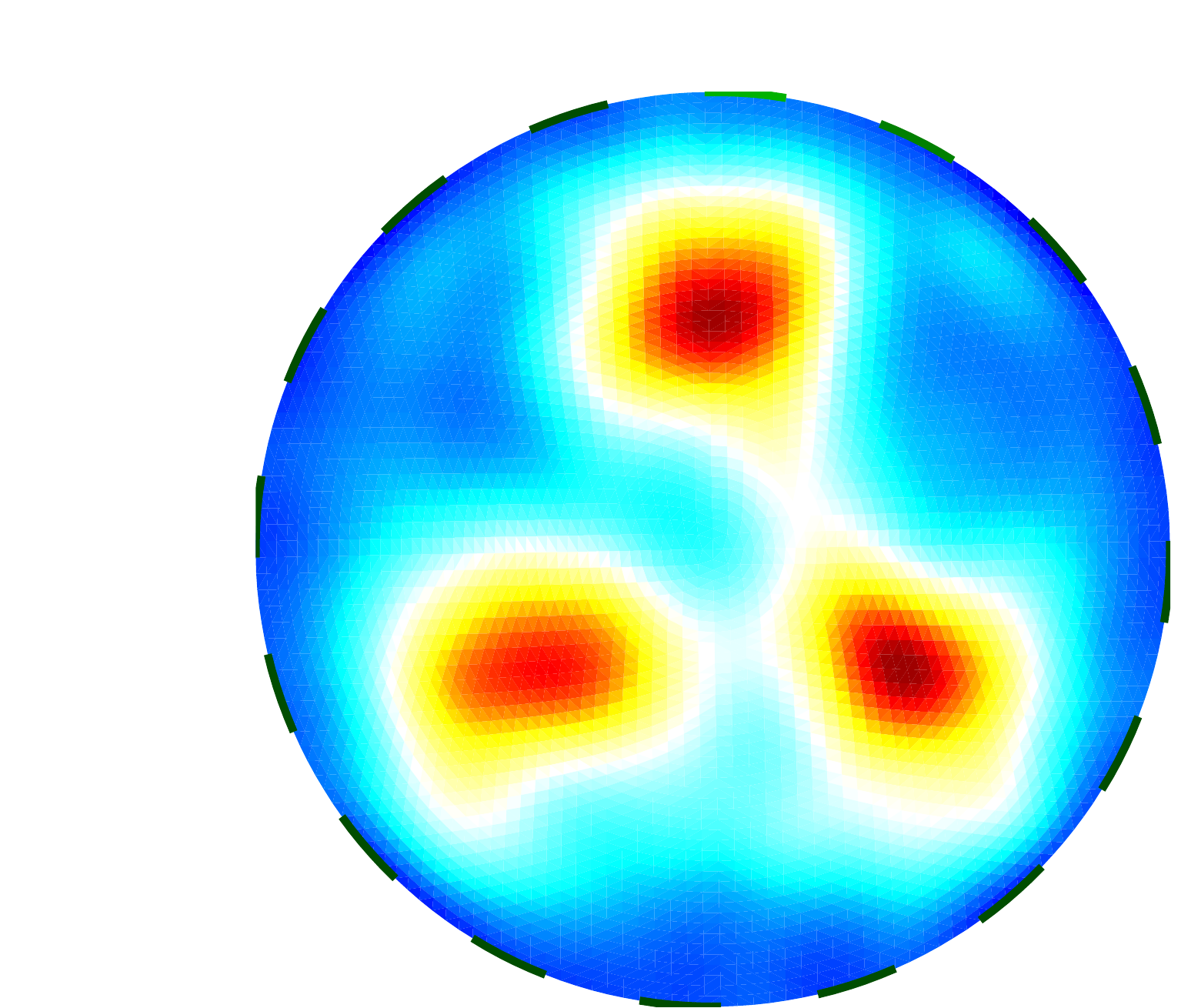}
\includegraphics[scale=0.18]{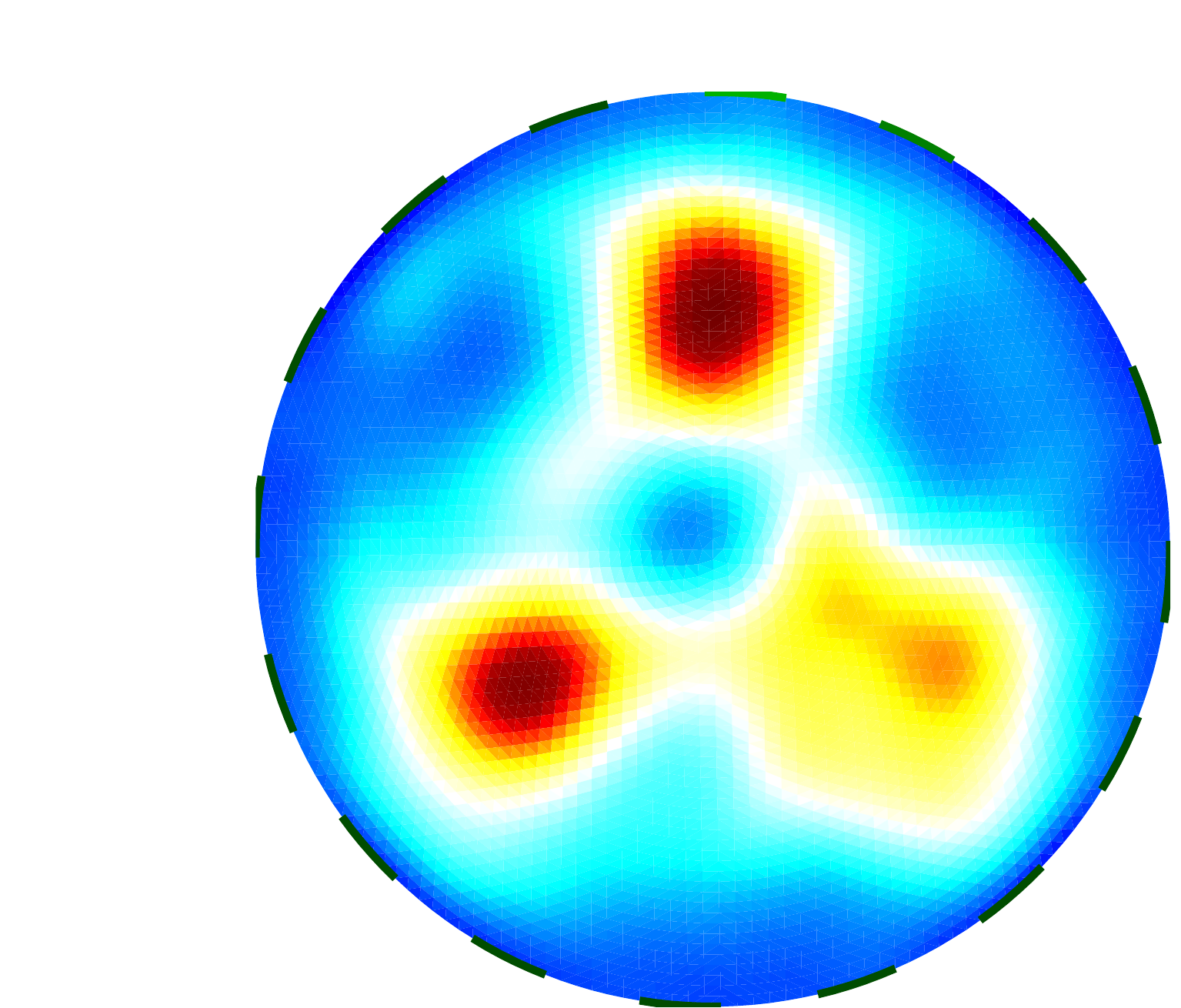}
 \caption{Top: True log-conductivity. Middle and bottom: Estimates of log-conductivity obtained from 5 experiments with different initial ensembles. The proposed method is used with  $\rho=0.5$ and $N_{e}=100$. } \label{FigEIT4}
\end{center}
\end{figure}

In subsection \ref{Darcy_Num} we have extensively studied the role of ensemble size and tunable parameters. The aim here is to briefly verify that similar results are obtained when Algorithm \ref{Al1} is applied to the EIT problem. For a fixed parameter $\rho=0.5$ we conduct a set of experiments where the scheme is applied to 40 different initial ensembles generated from the Gaussian measure described above. The corresponding log-data misfit (bottom) and relative error (top) from these experiments are displayed in  \Fref{FigEIT3}. As with the previous test model, the ensemble size is crucial to the performance of the proposed iterative scheme. We observe a critical ensemble size of $N_{e}=75$ above which the proposed scheme is properly stabilized when the early termination is carried out according to (\ref{eq:m15}) with the selection of $\tau\approx 1/\rho$ provided by the theory of the regularizing LM scheme. In this case we note again that the critical size above which the method exhibits regularizing properties was smaller than the number of measurements; this confirms that the proposed scheme addresses the small ensemble effect described in subsection \ref{unreg}. The ensemble mean $\overline{u}_{n}$ from 5 experiments (from different initial ensembles) with $N_{e}=100$ are displayed in \Fref{FigEIT4}.

\begin{figure}[htbp]
\begin{center}
\includegraphics[scale=0.25]{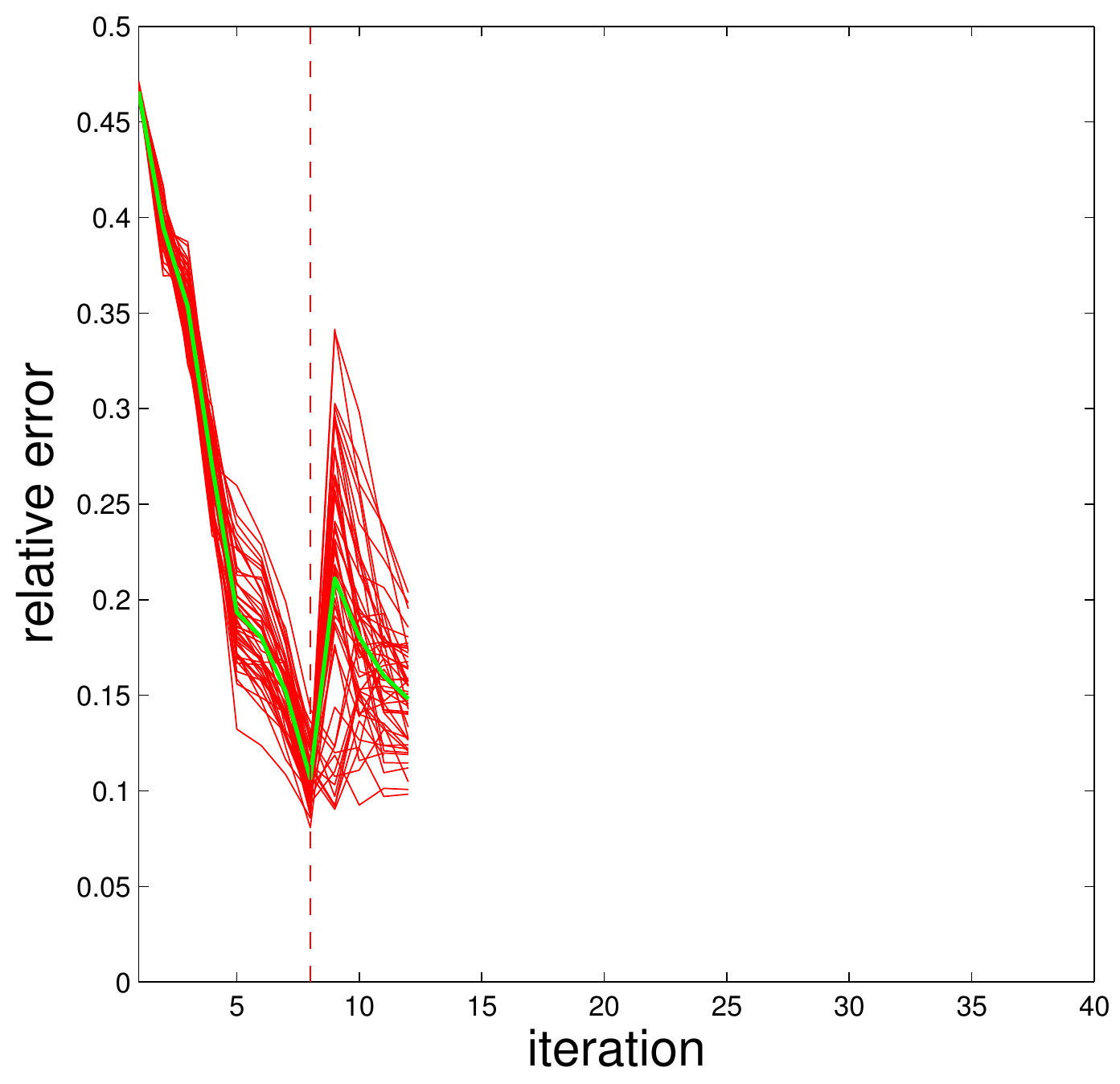}
\includegraphics[scale=0.25]{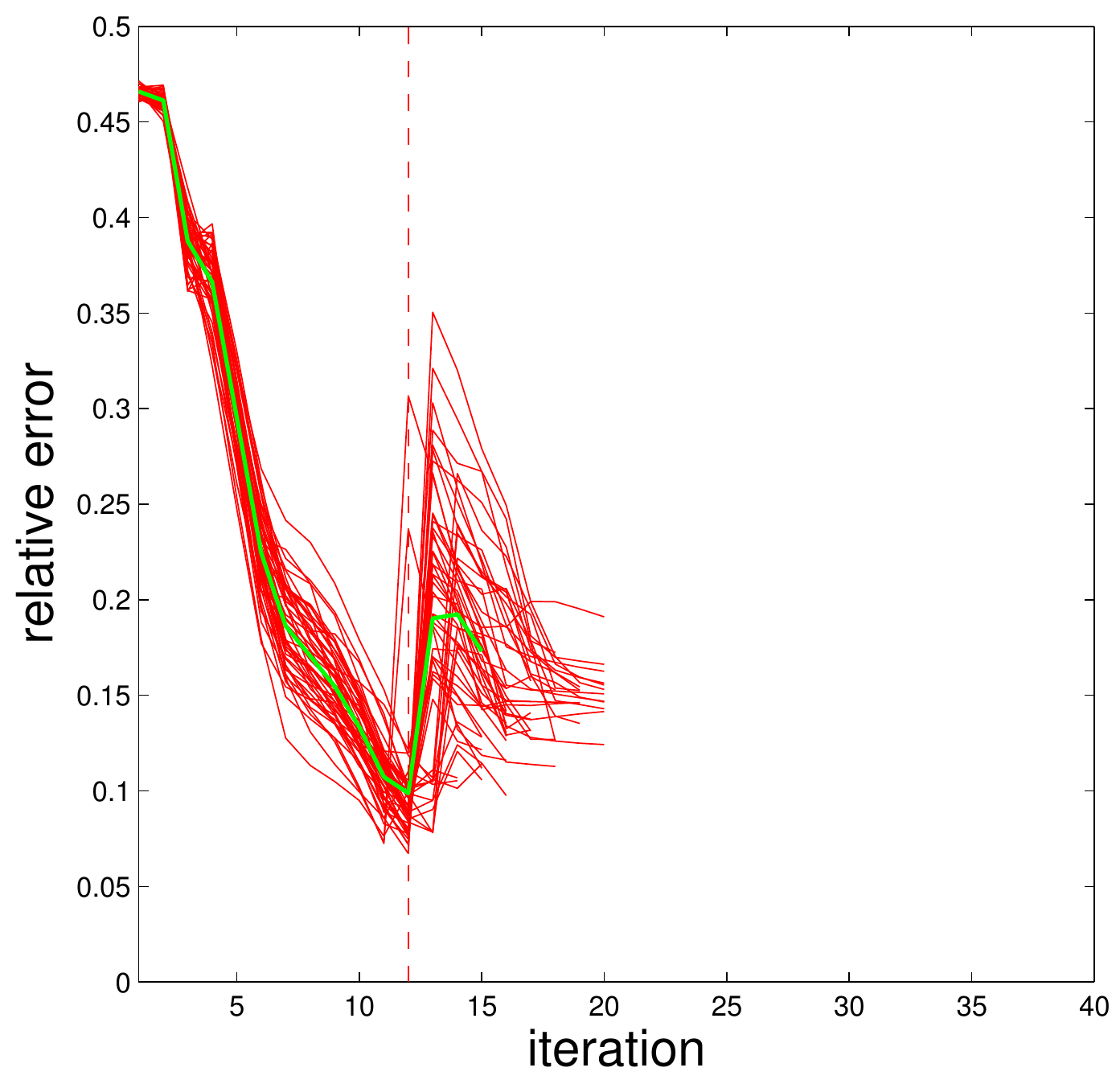}
\includegraphics[scale=0.25]{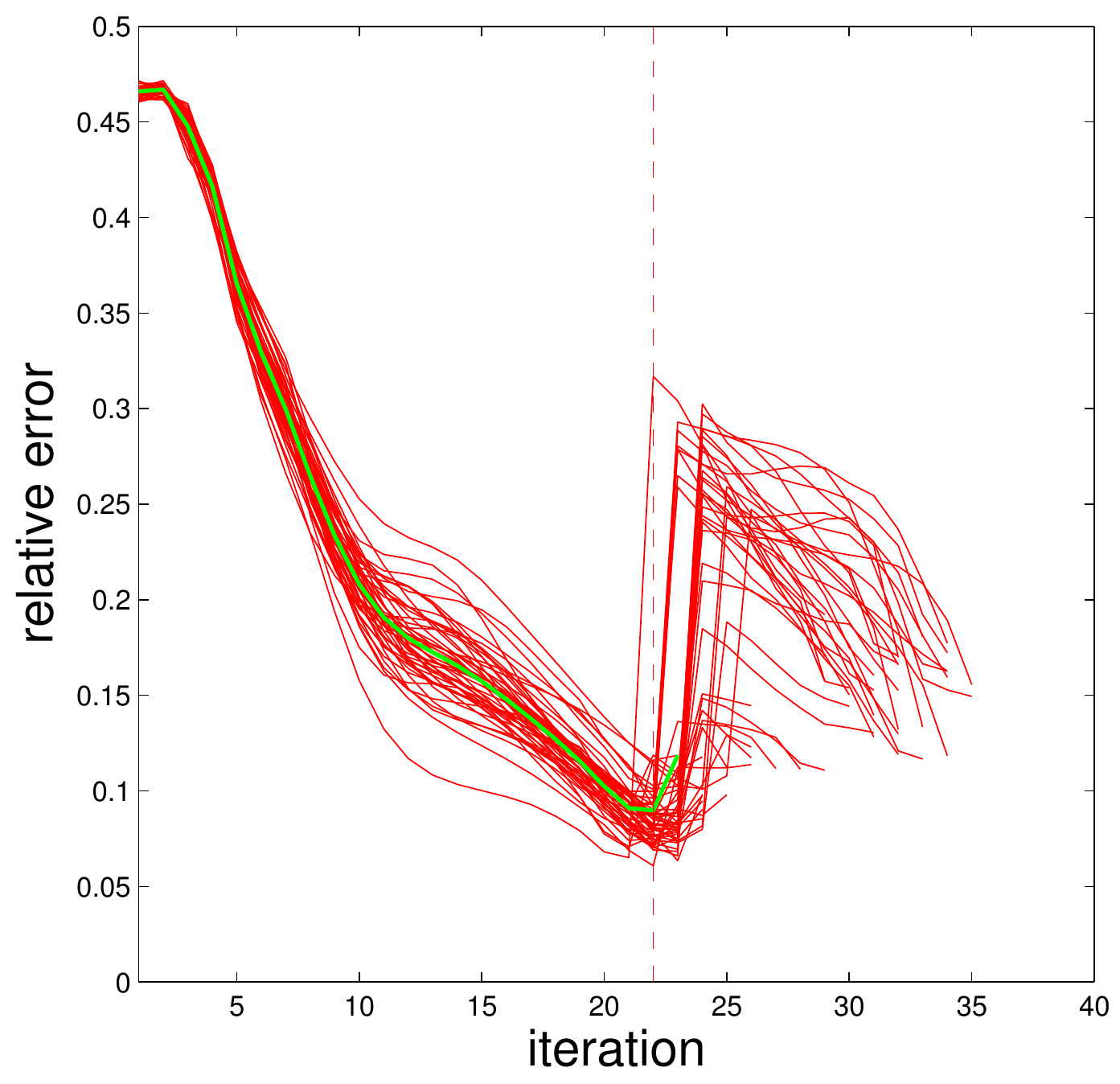}
\includegraphics[scale=0.25]{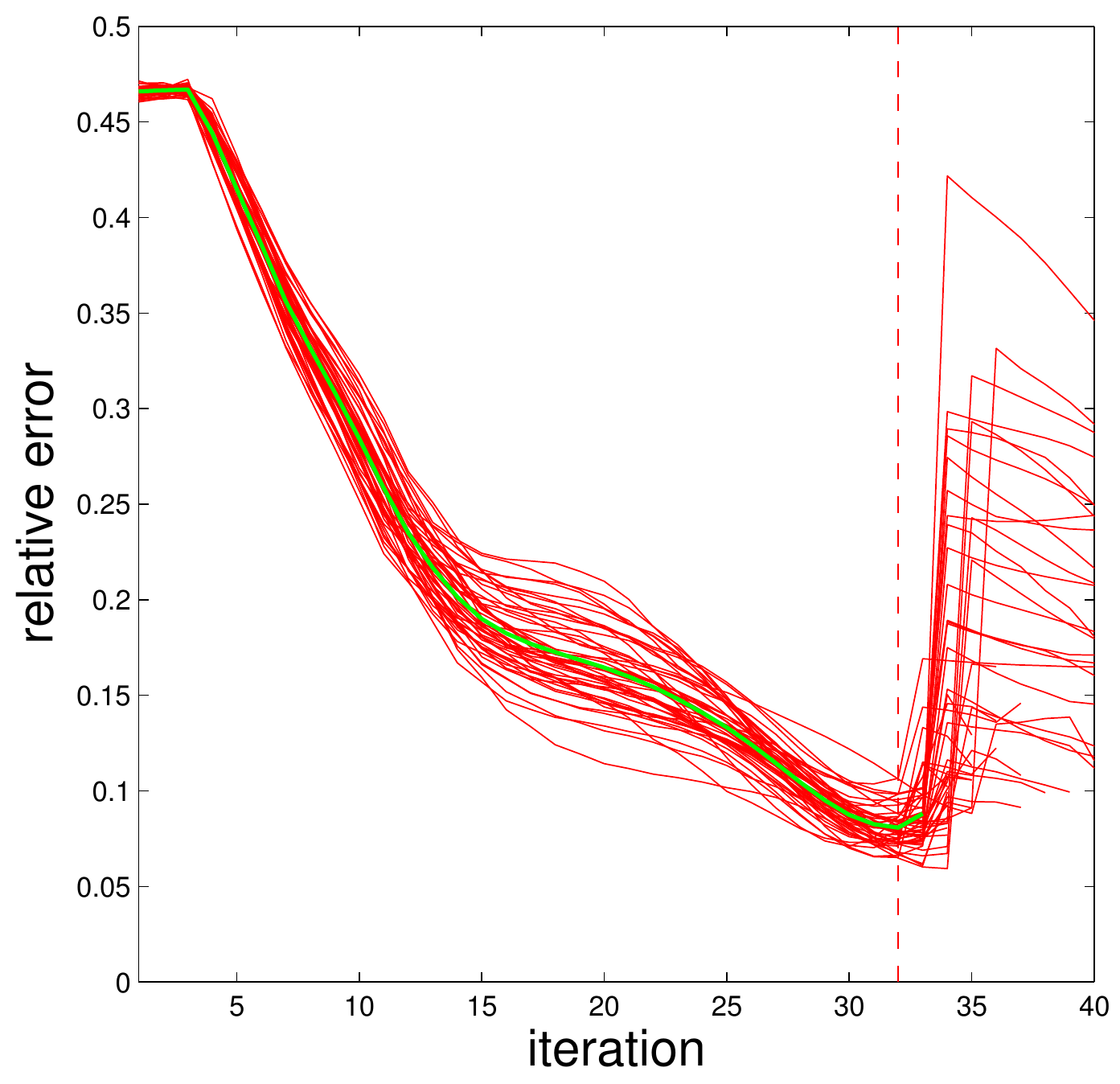}\\
\includegraphics[scale=0.25]{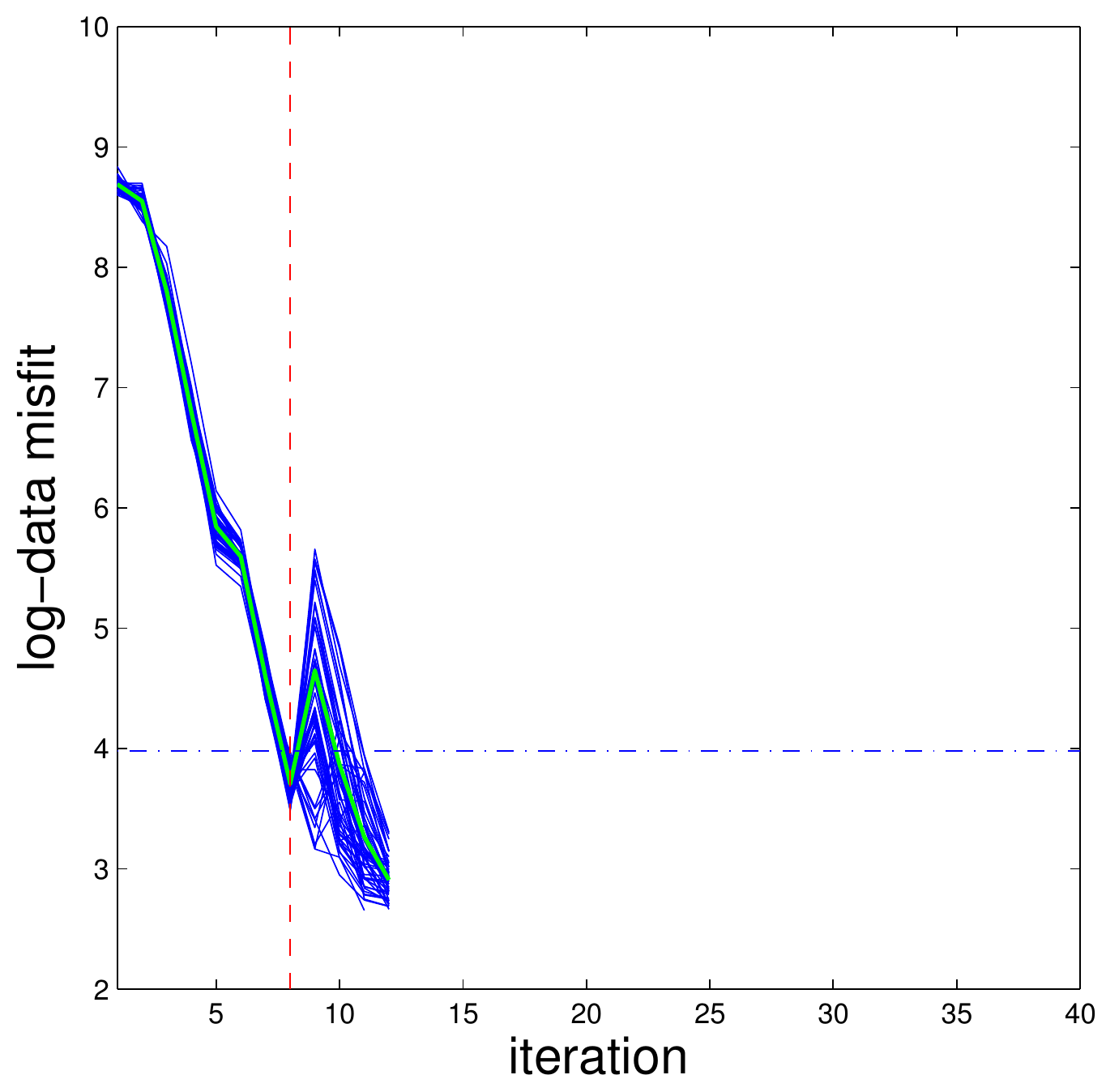}
\includegraphics[scale=0.25]{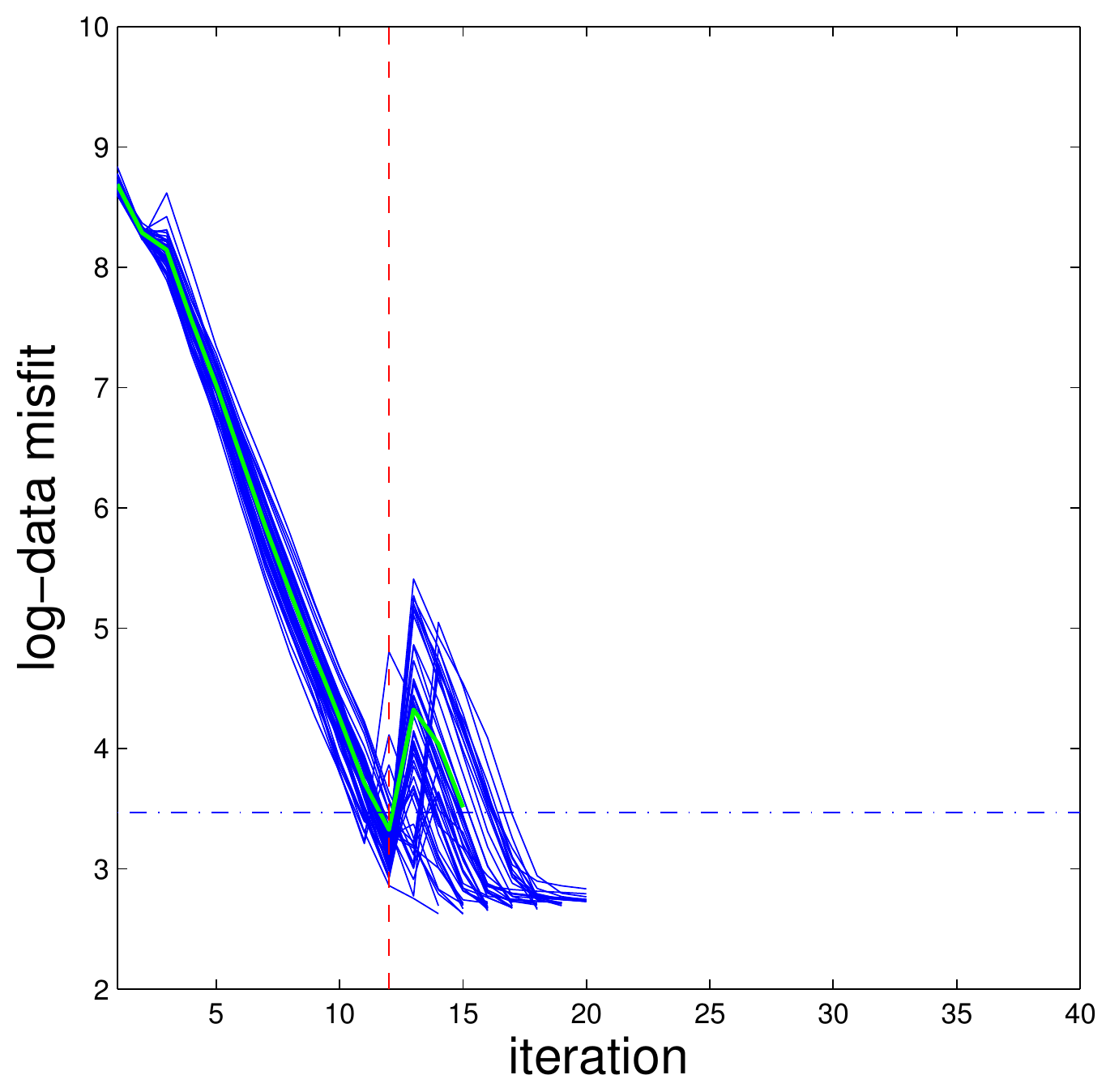}
\includegraphics[scale=0.25]{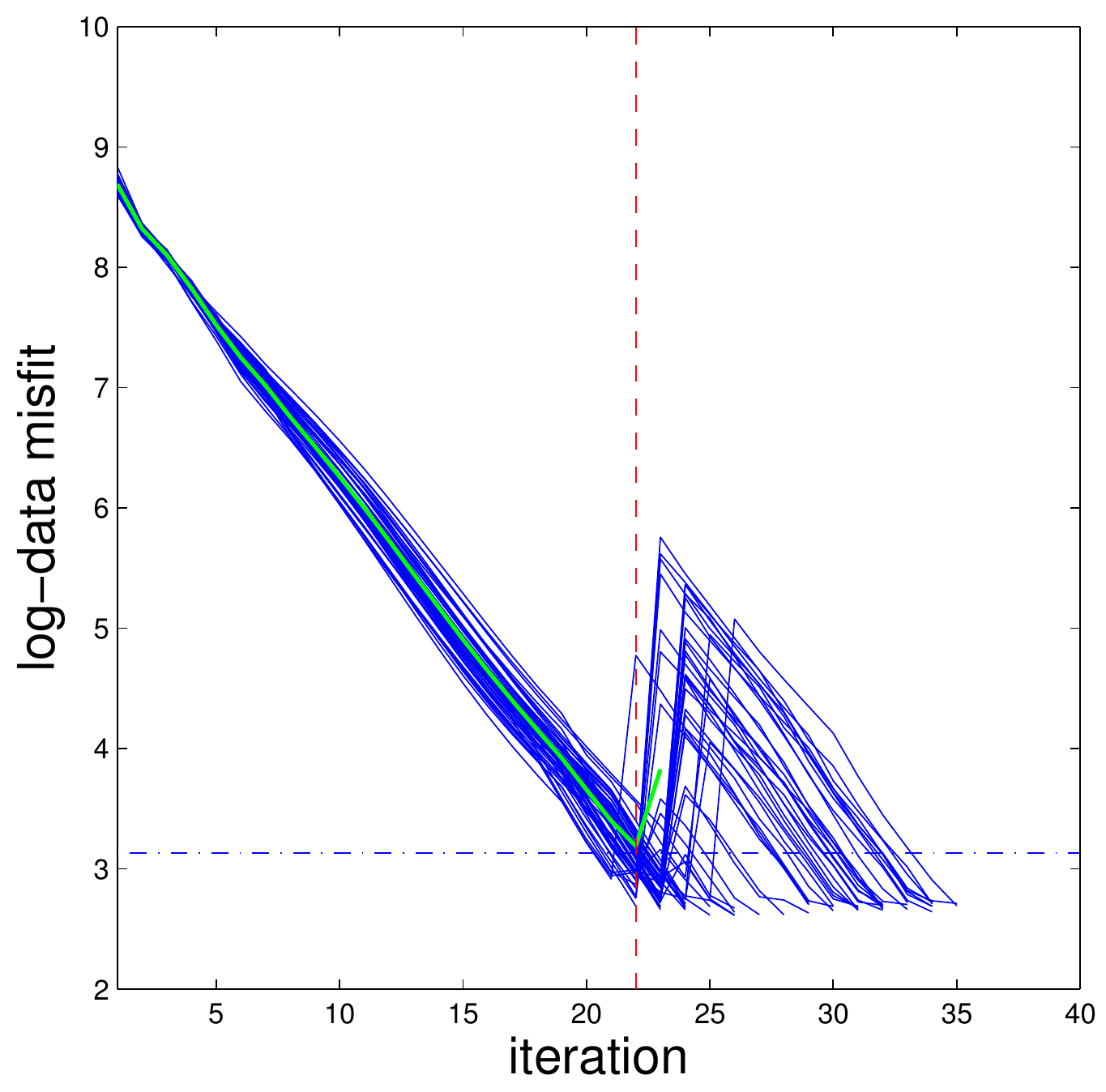}
\includegraphics[scale=0.25]{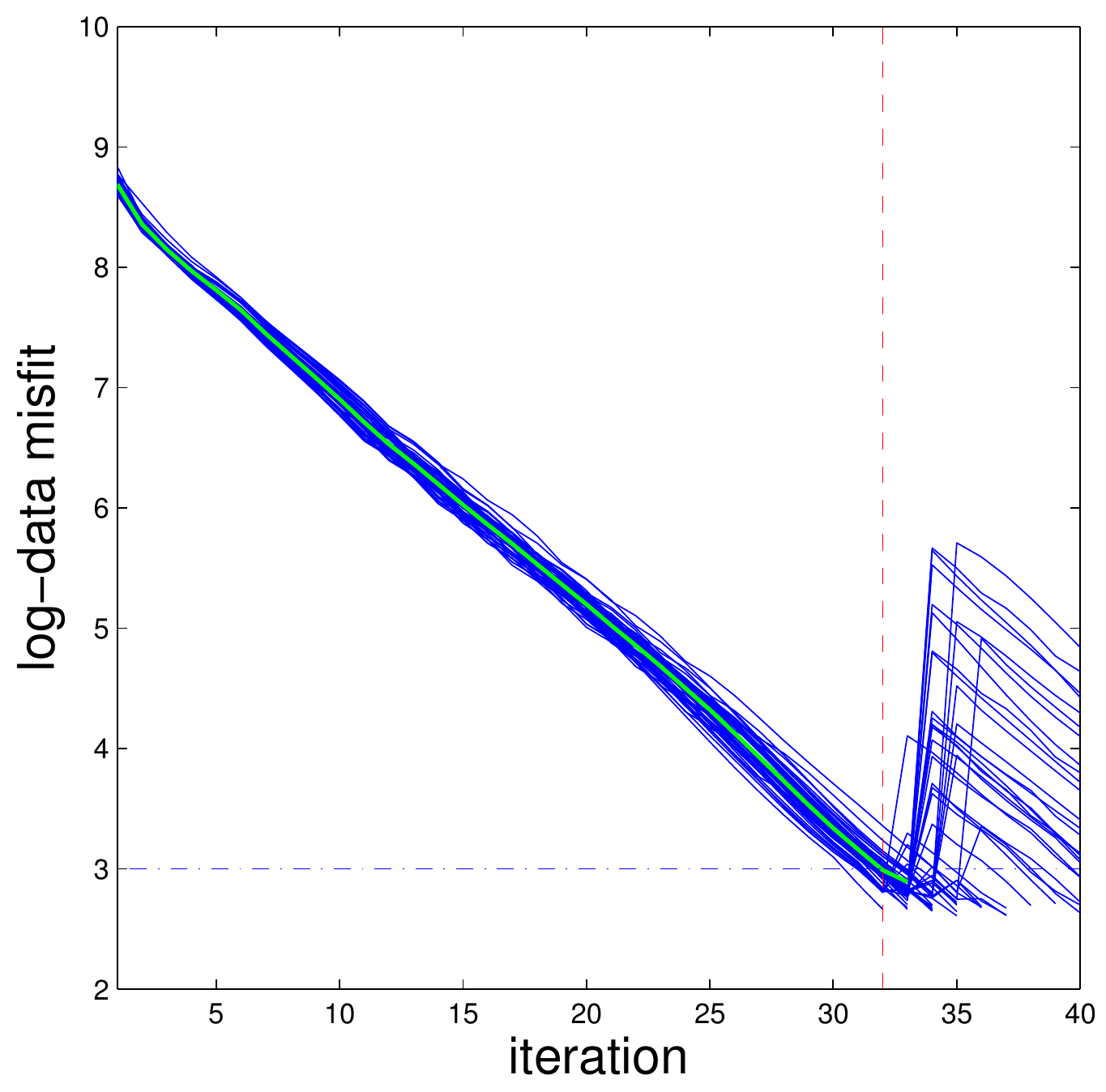}\\
 \caption{Error w.r.t. truth (top) and log - data misfit (bottom) obtained with Algorithm \ref{Al1} with (from left to right) $\rho=0.3,0.5, 0.7, 0.8 $ from 40 experiments with different initial ensembles of size $N_{e}=100$}   \label{FigEIT5}
\end{center}
\end{figure}

In \Fref{FigEIT5} we fix $N_{e}=100$ and now perform a set of experiments with different selections of $\rho$ in (\ref{eq:m12}). As before, we confirm that for this sufficiently large ensemble size, the selection of $\tau>1/\rho$ in the termination of the scheme ensures stabilized estimates. It is worth noticing that the aforementioned selection of $\tau$ is essential as the relative error increases abruptly once the data misfit drops below the value $\eta/\rho$.  While there is a slight increase in the accuracy (in terms of the error w.r.t truth) as we increase $\rho$, we clearly observe that the number of iterations and so the computational cost of the scheme increases substantially. For this case we observe that $\rho=0.5$ provides a reasonable balance between computational cost and accuracy.

\begin{figure}[htbp]
\begin{center}
\includegraphics[scale=0.4]{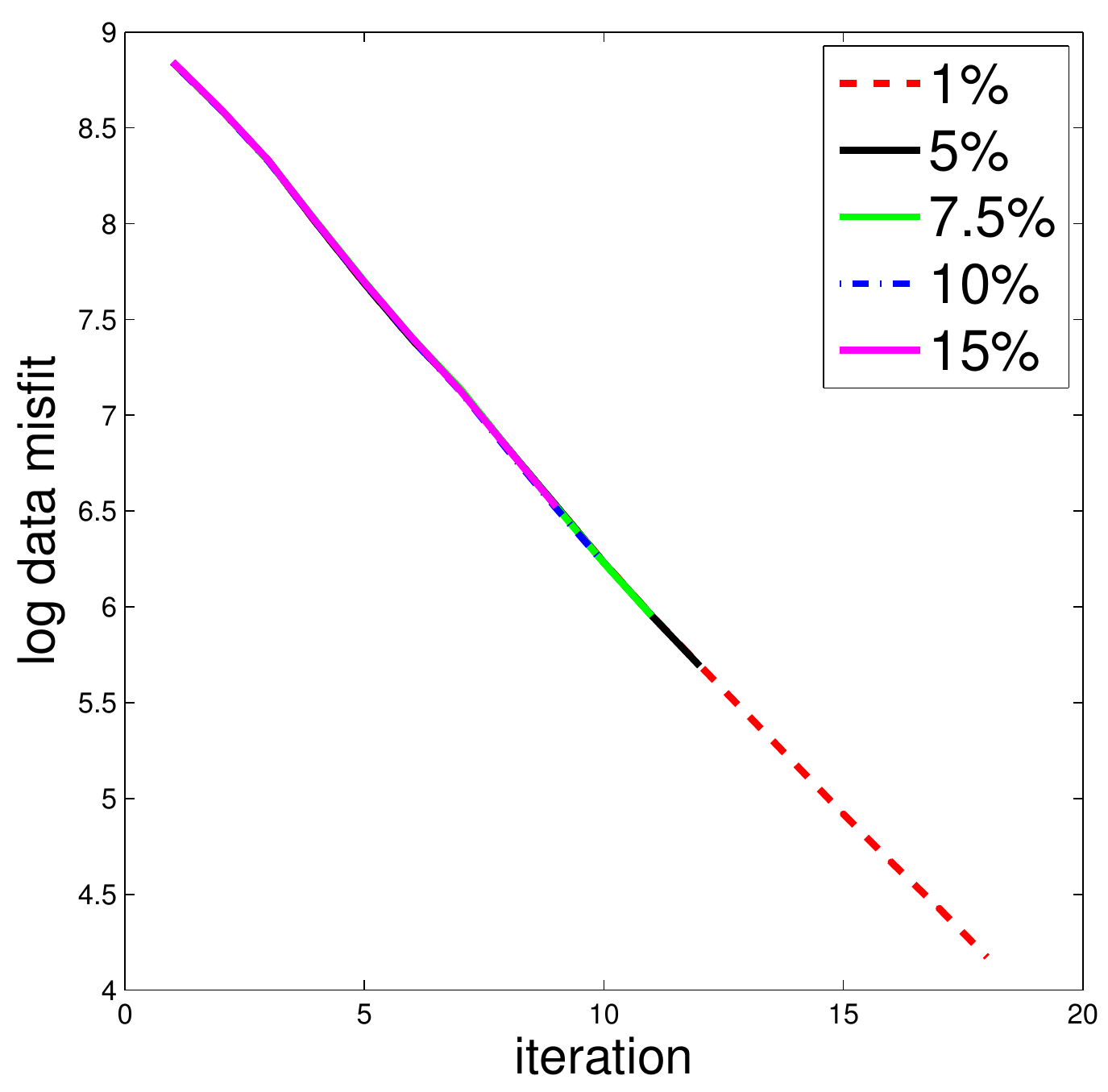}
\includegraphics[scale=0.4]{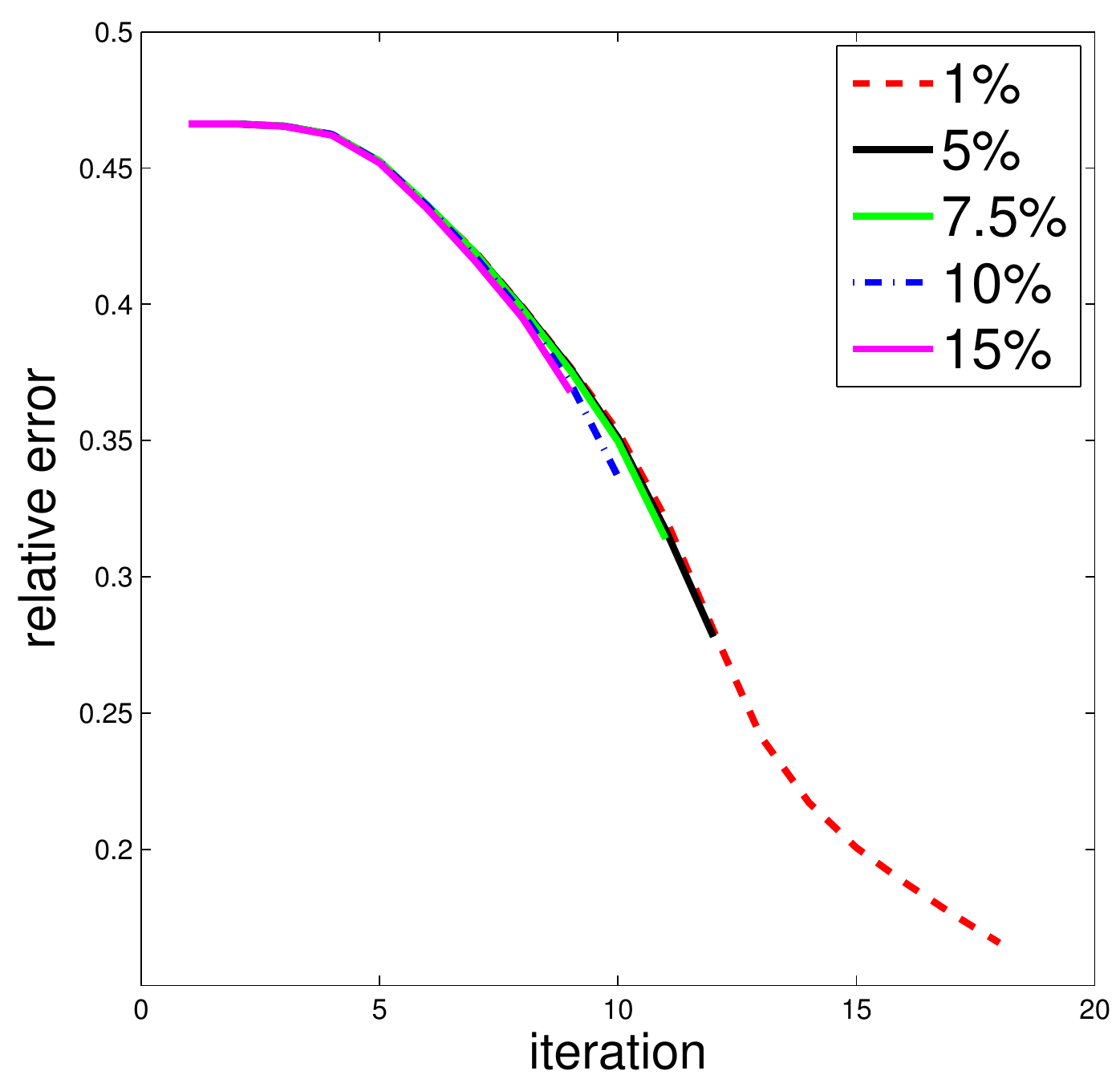}
 \caption{Log data misfit (right) and error w.r.t truth (left) from the average over 40 experiments where Algorithm \ref{Al1} was applied  with different initial ensembles to synthetic data of noise levels of $1\%, 5\%, 7.5\%, 10\% 15\%$. } \label{FigEIT6}
\end{center}
\end{figure}

Finally, in \Fref{FigEIT6} we show the effect of the stabilization of the scheme as the noise decreases. These results corresponds to the application of the scheme with synthetic data generated from the electrode configuration described before but with different noise levels. The results from \Fref{FigEIT6} corresponds to averages, at each iteration, from 40 experiments obtained from 40 different selections of initial ensembles. As we expected, the proposed ensemble scheme inherits the regularizing properties of the LM method and produces stable computations that converge, in a stable fashion as $\eta\to 0$.

\subsubsection{The selection of the initial ensemble}

\begin{figure}[htbp]
\begin{center}
\includegraphics[scale=0.18]{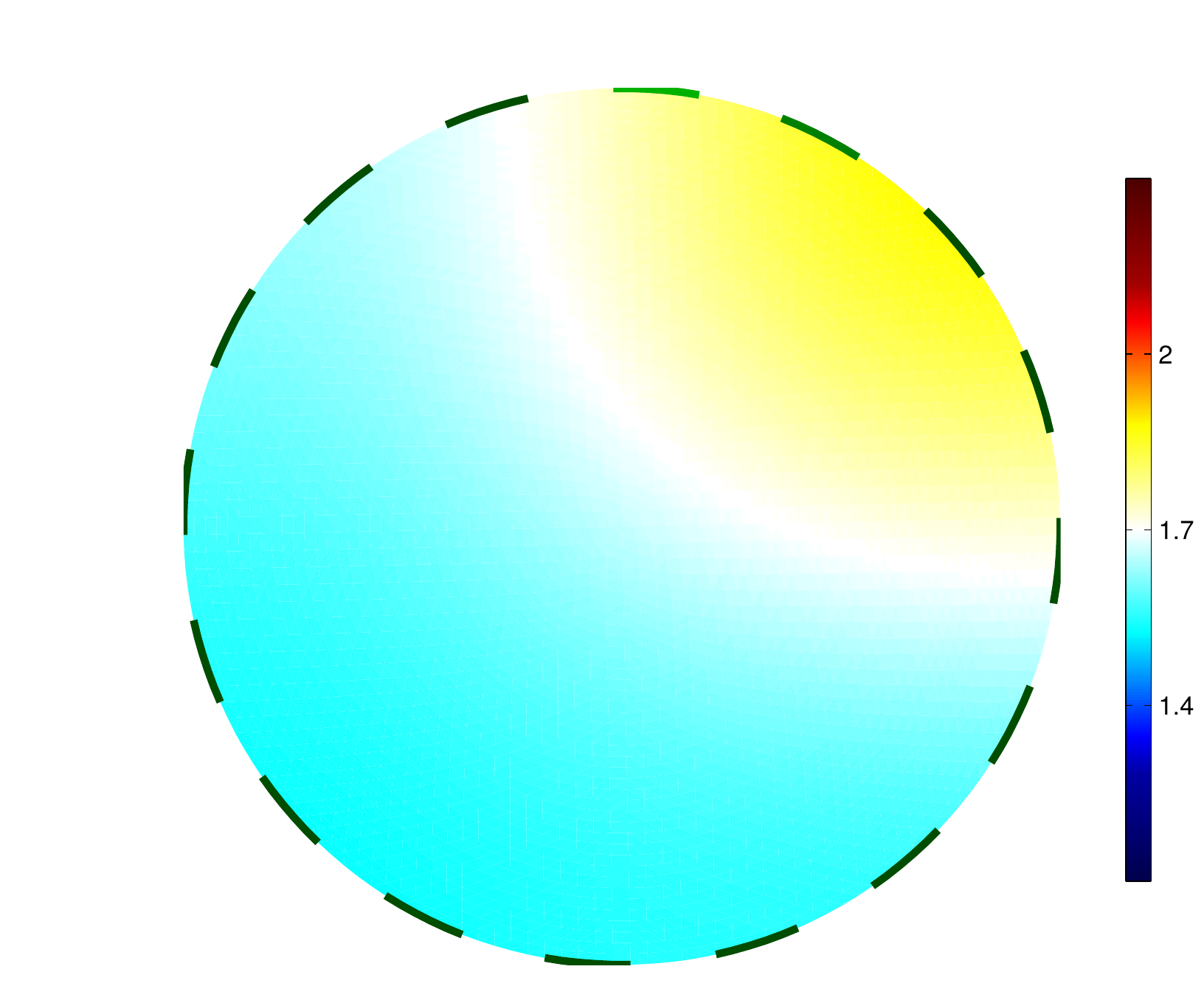}
\includegraphics[scale=0.18]{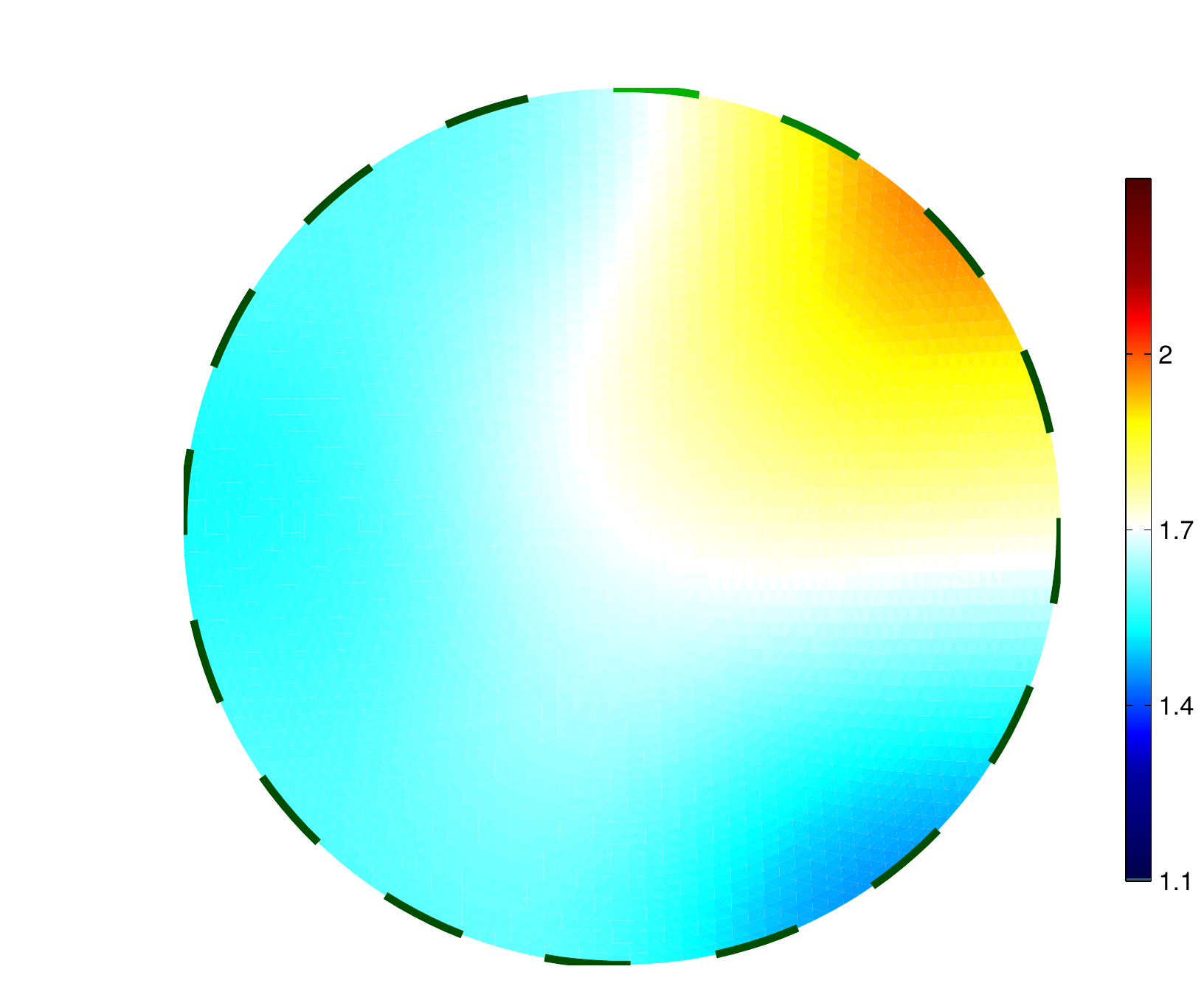}
\includegraphics[scale=0.18]{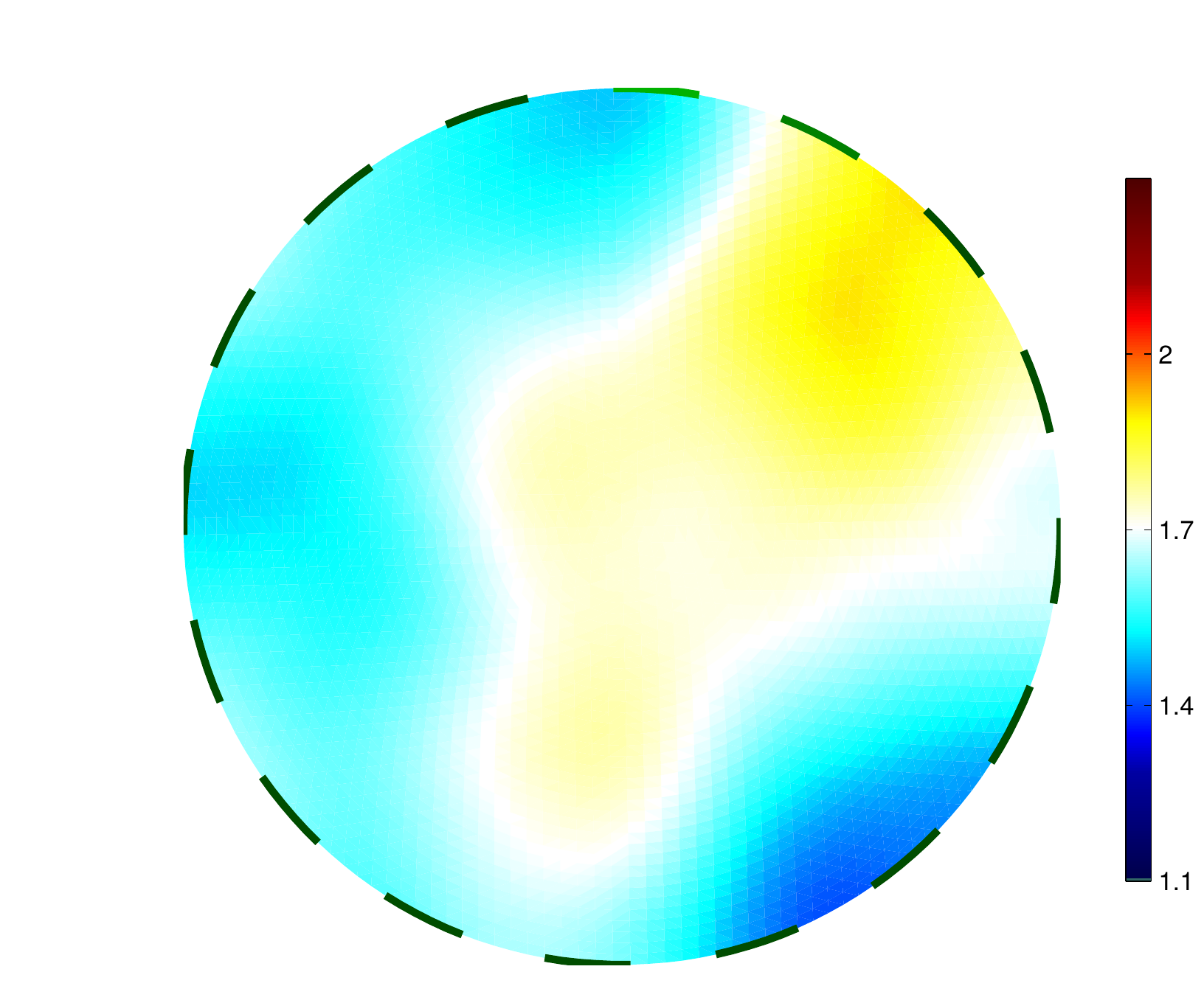}
\includegraphics[scale=0.18]{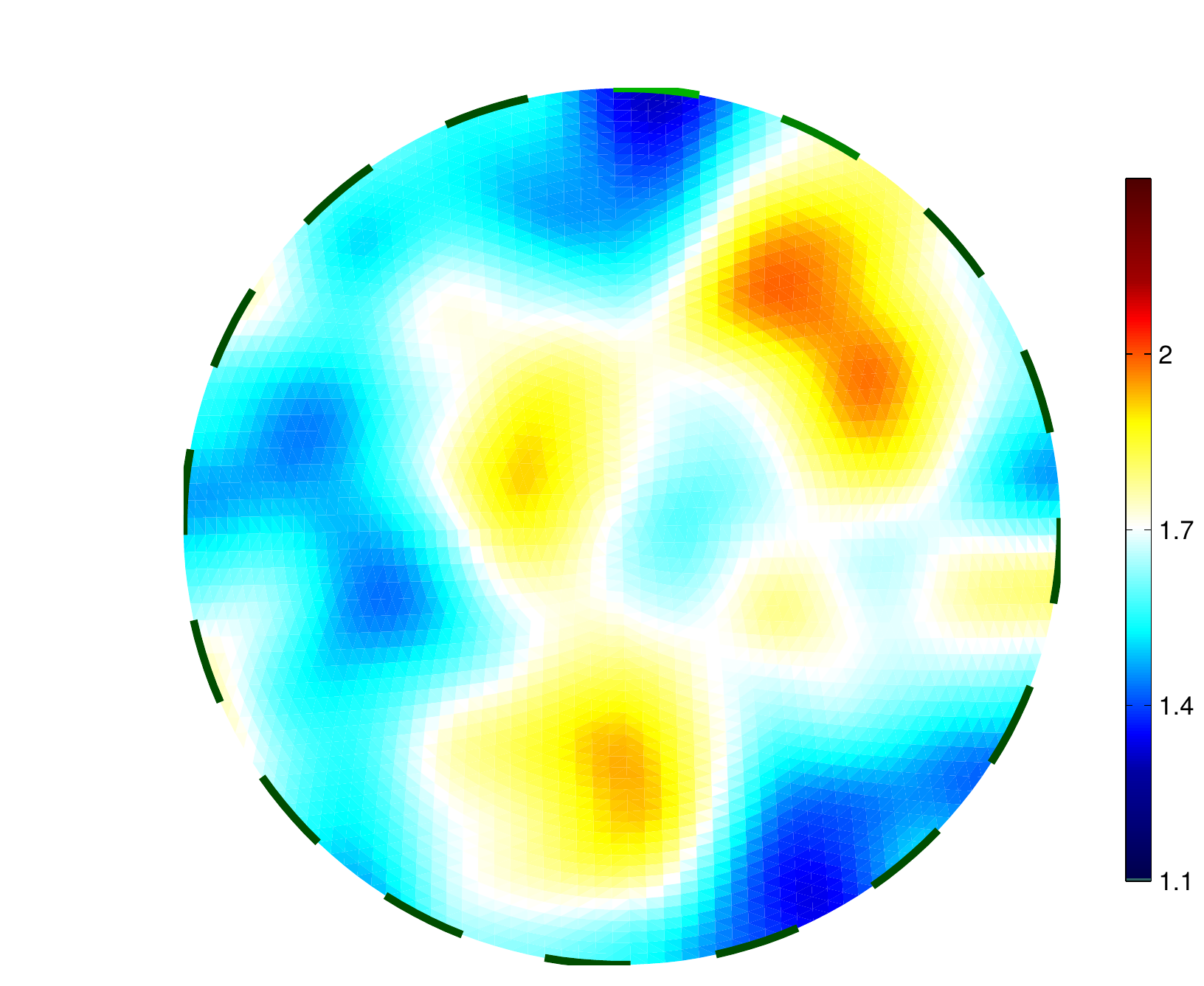}
\includegraphics[scale=0.18]{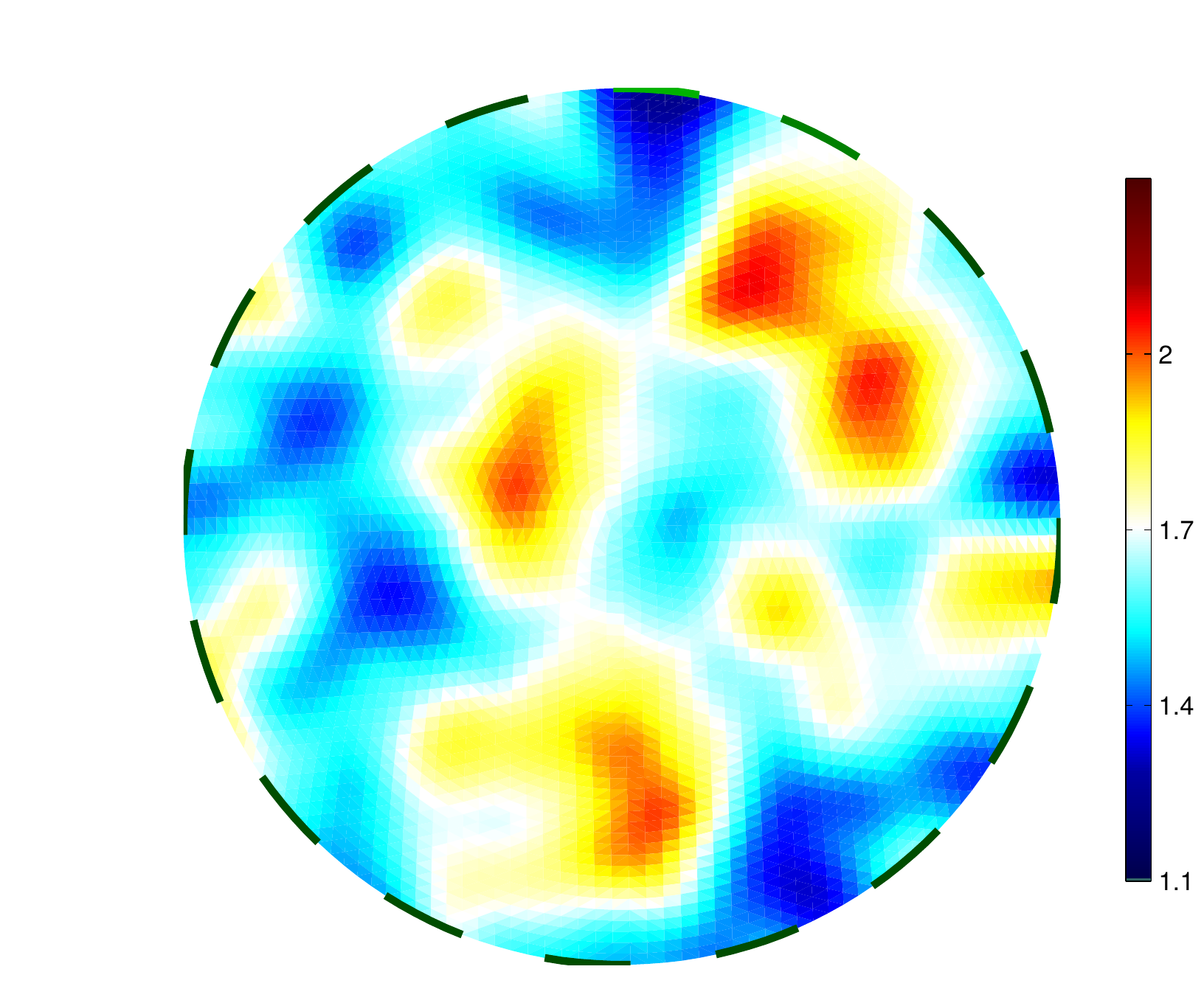}\\
\includegraphics[scale=0.18]{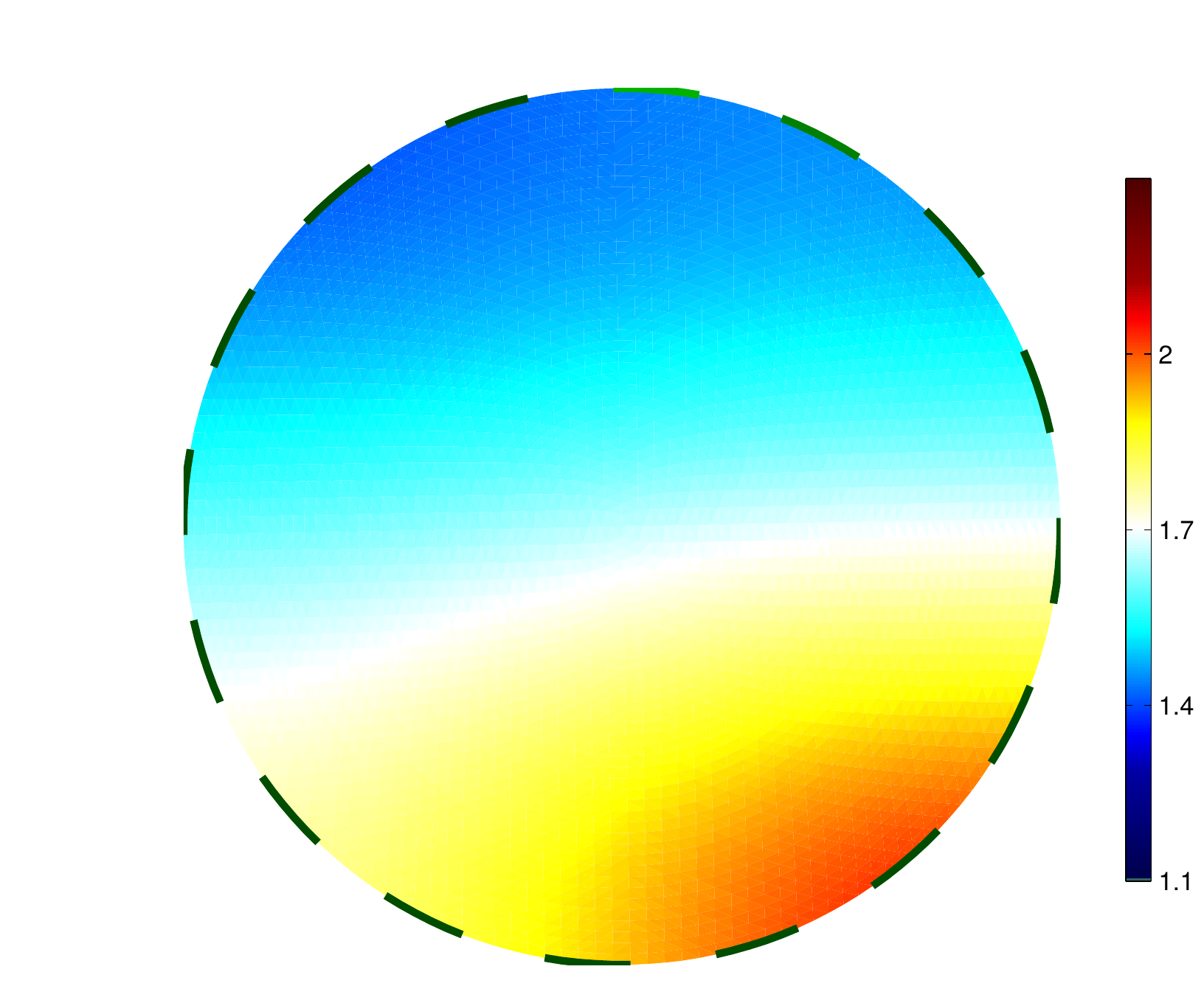}
\includegraphics[scale=0.18]{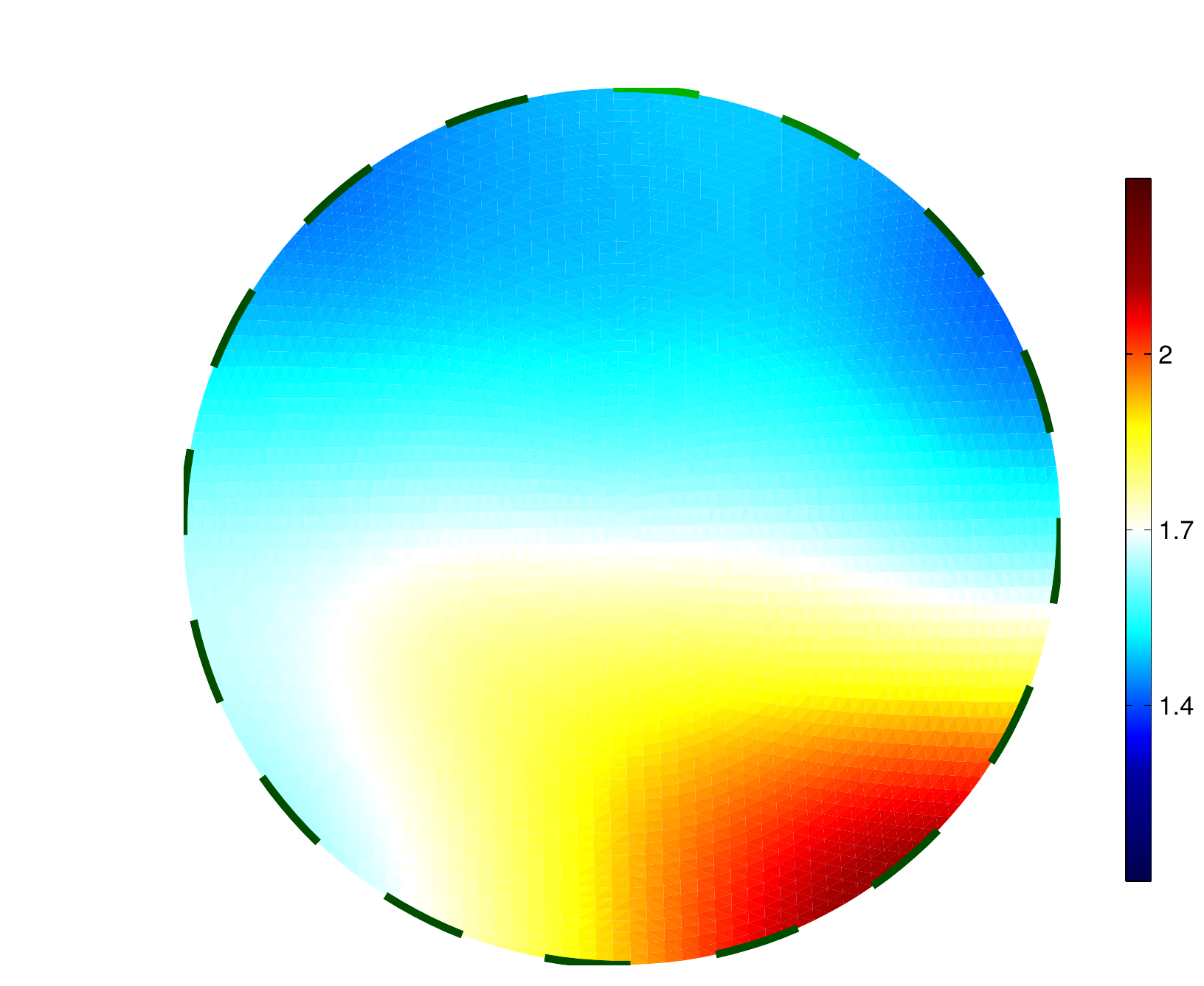}
\includegraphics[scale=0.18]{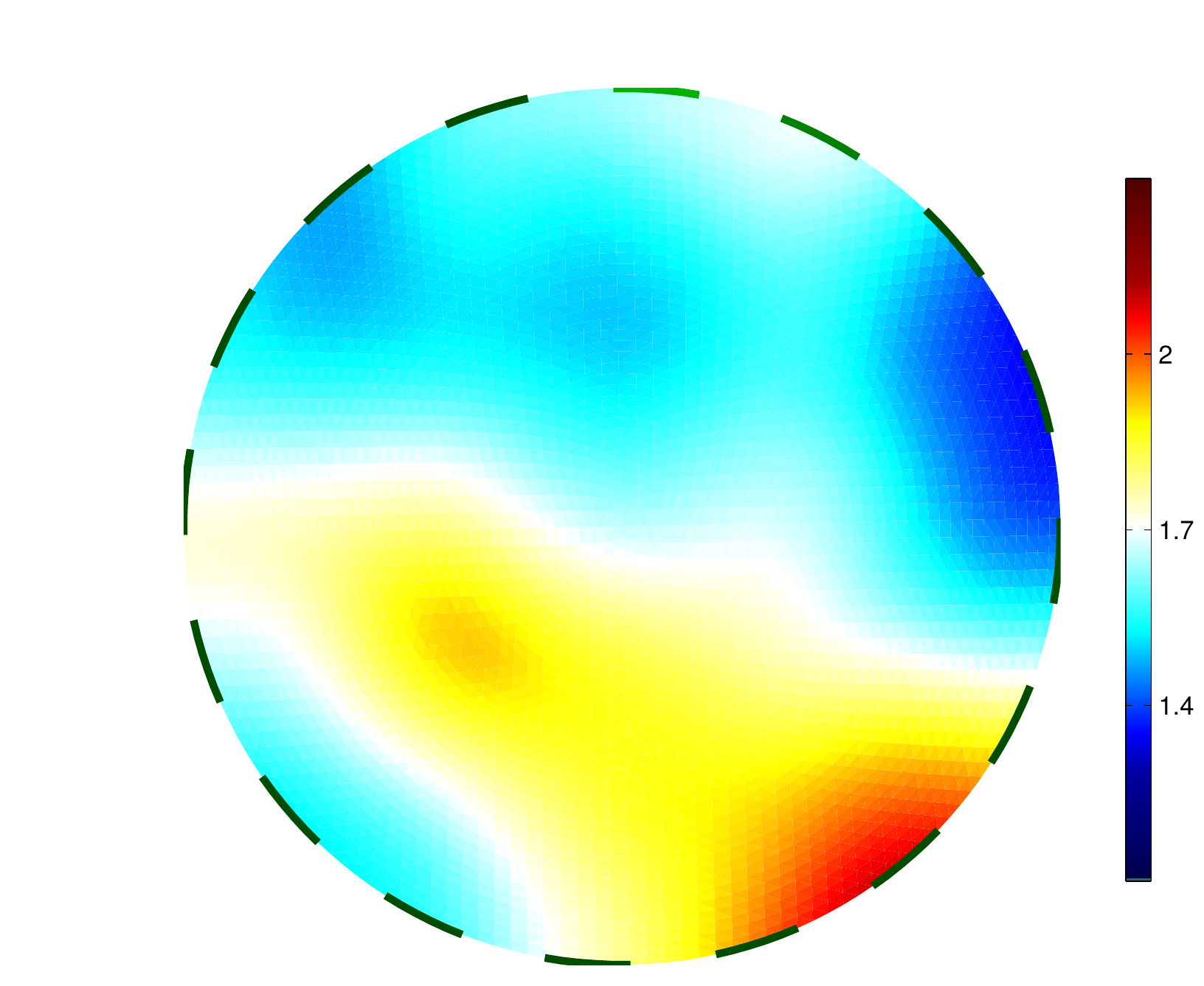}
\includegraphics[scale=0.18]{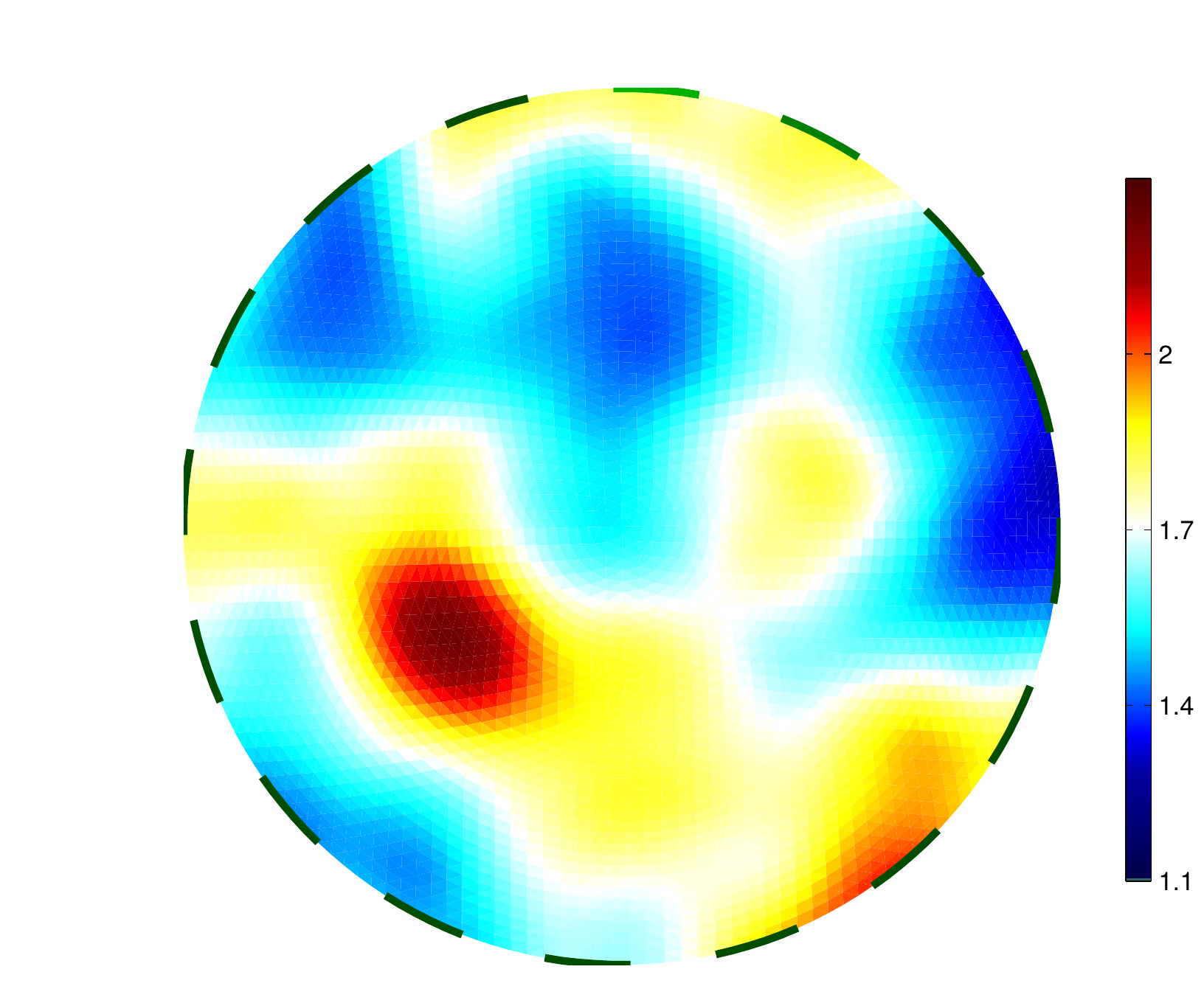}
\includegraphics[scale=0.18]{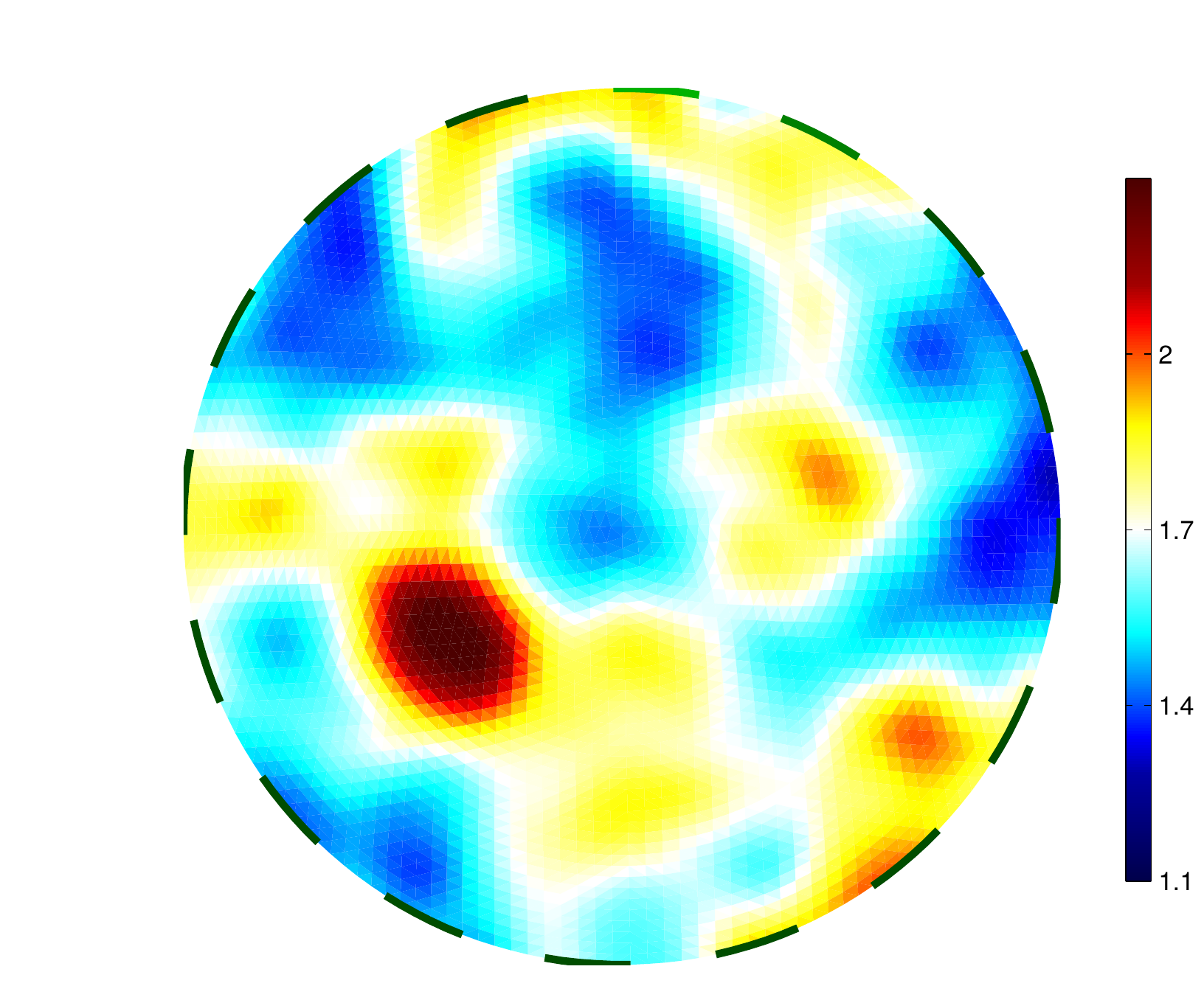}\\
\includegraphics[scale=0.18]{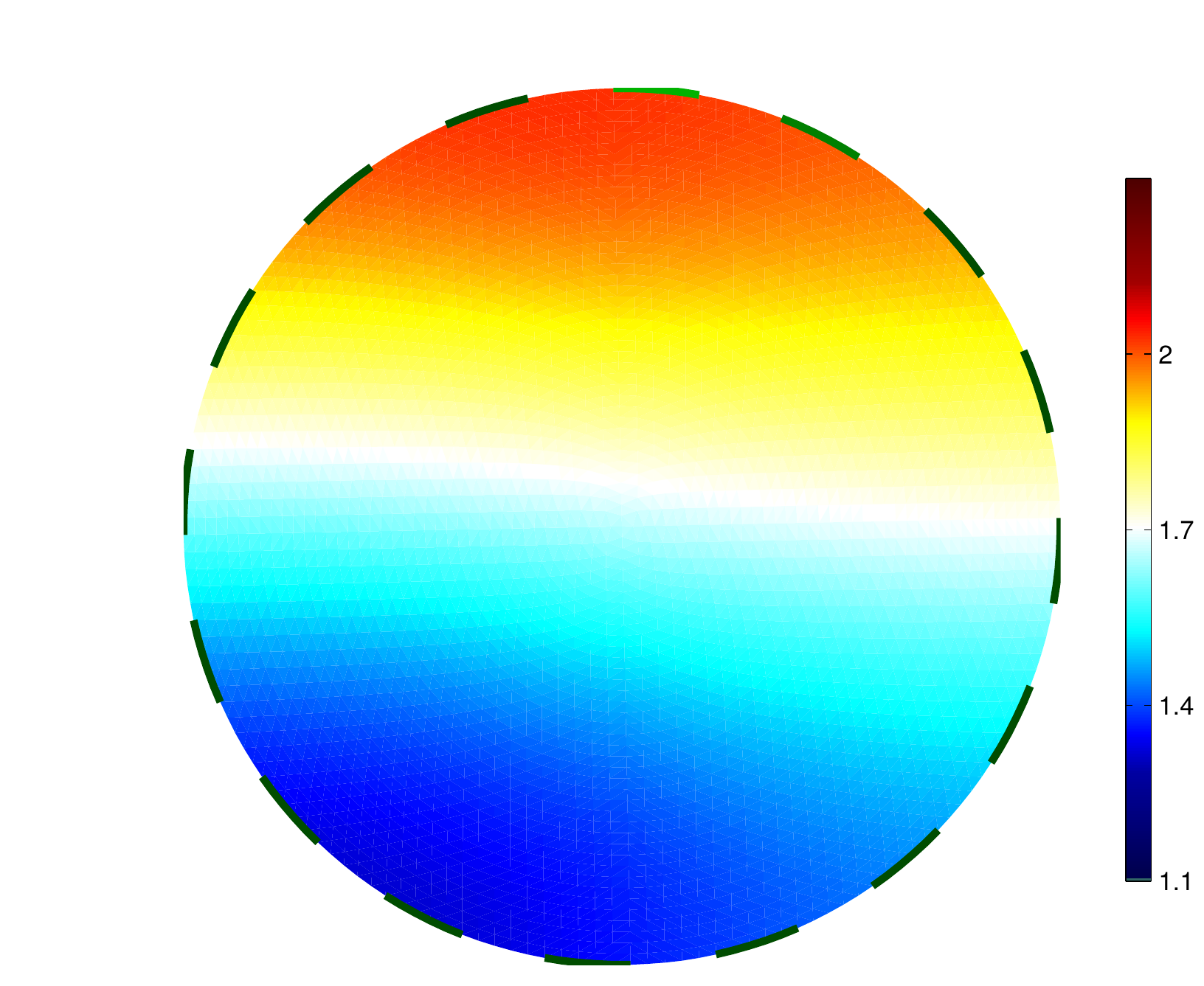}
\includegraphics[scale=0.18]{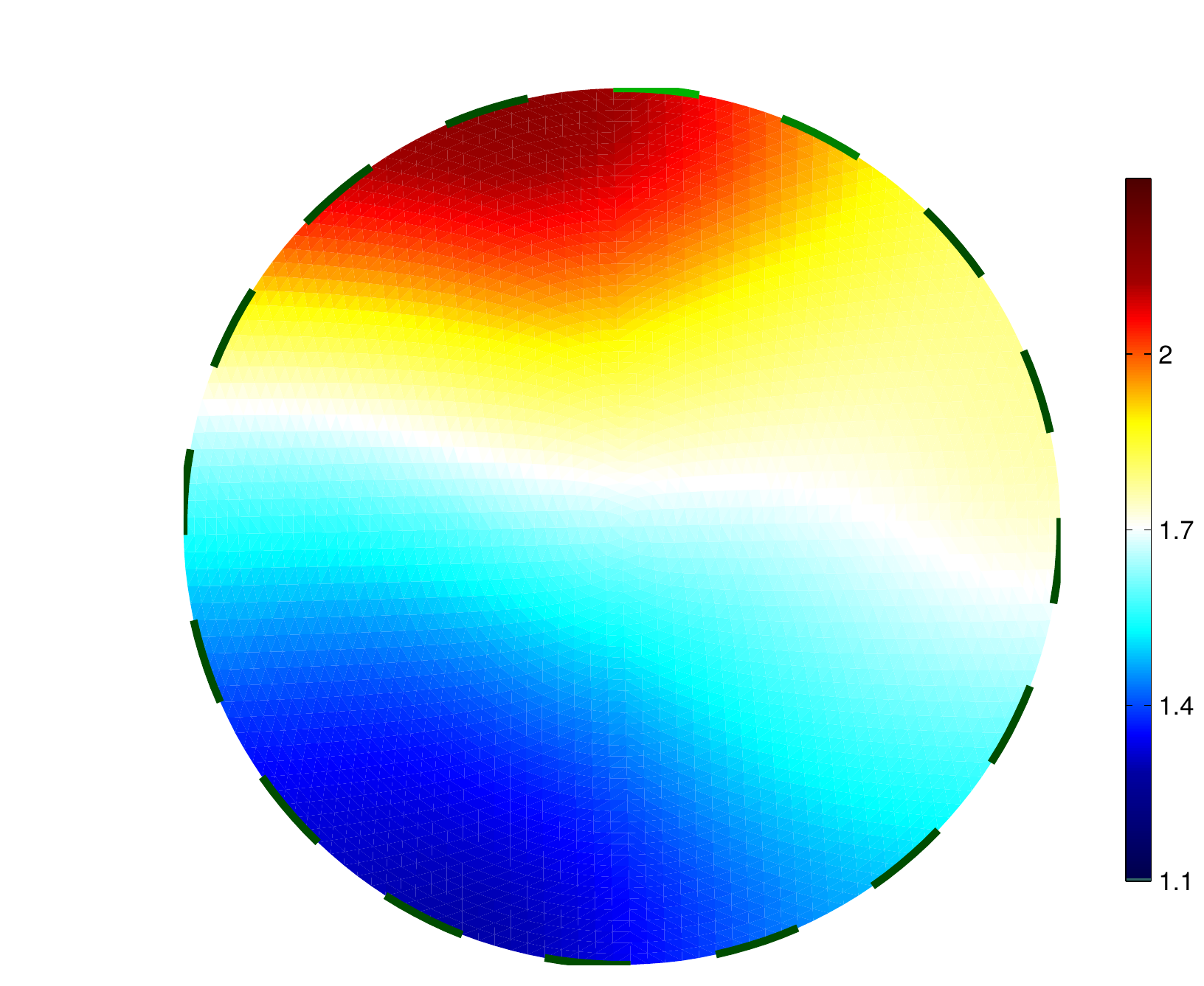}
\includegraphics[scale=0.18]{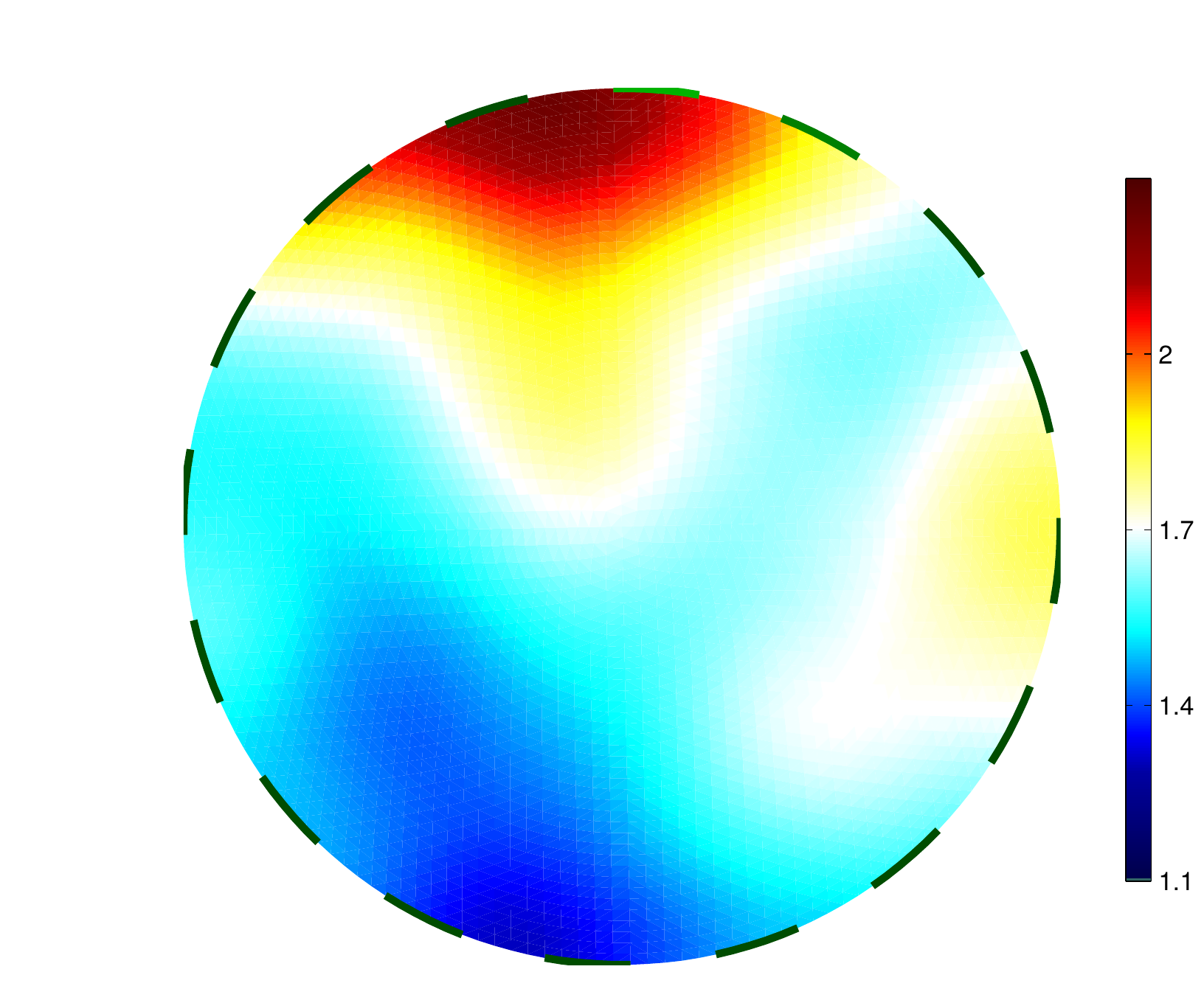}
\includegraphics[scale=0.18]{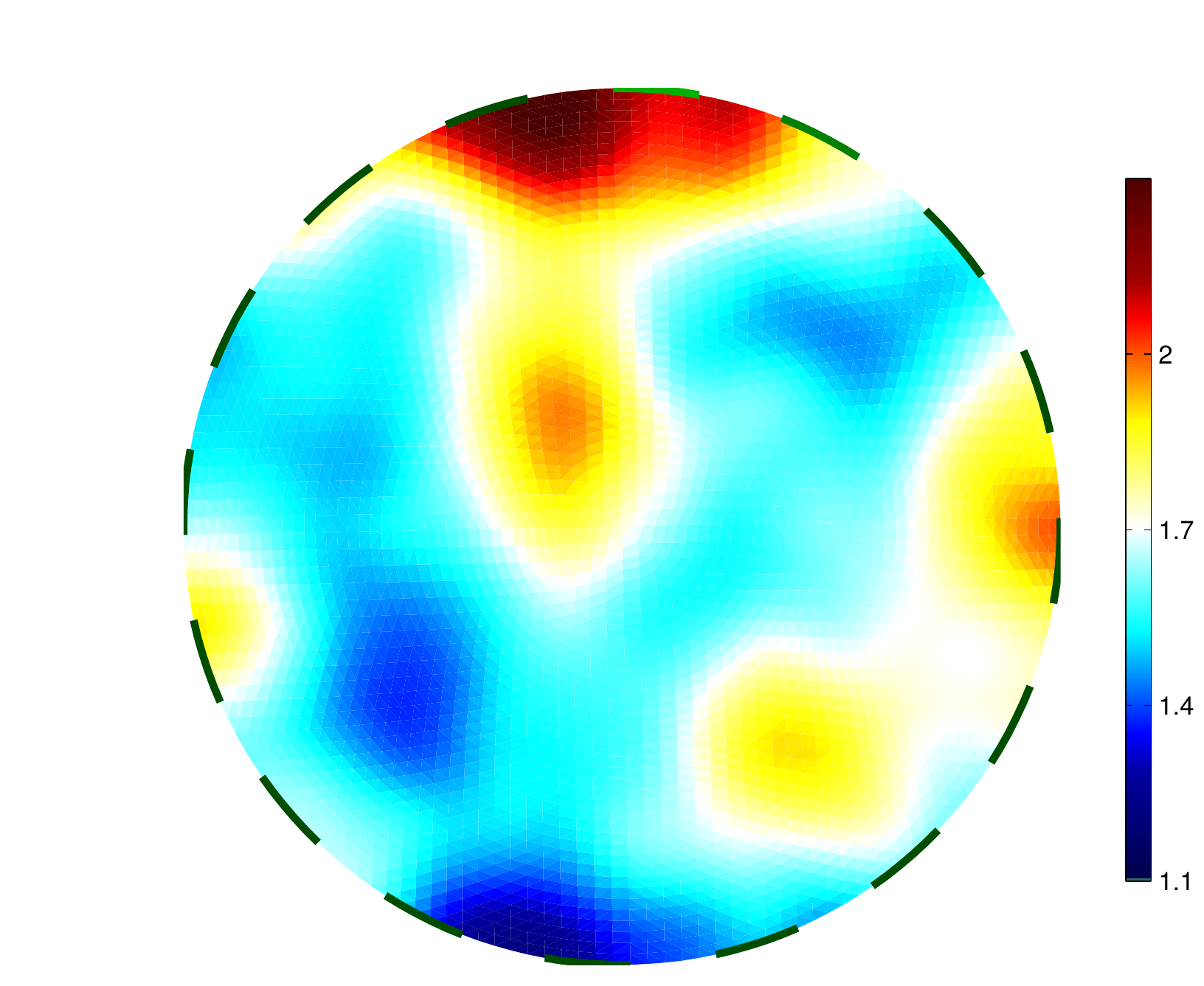}
\includegraphics[scale=0.18]{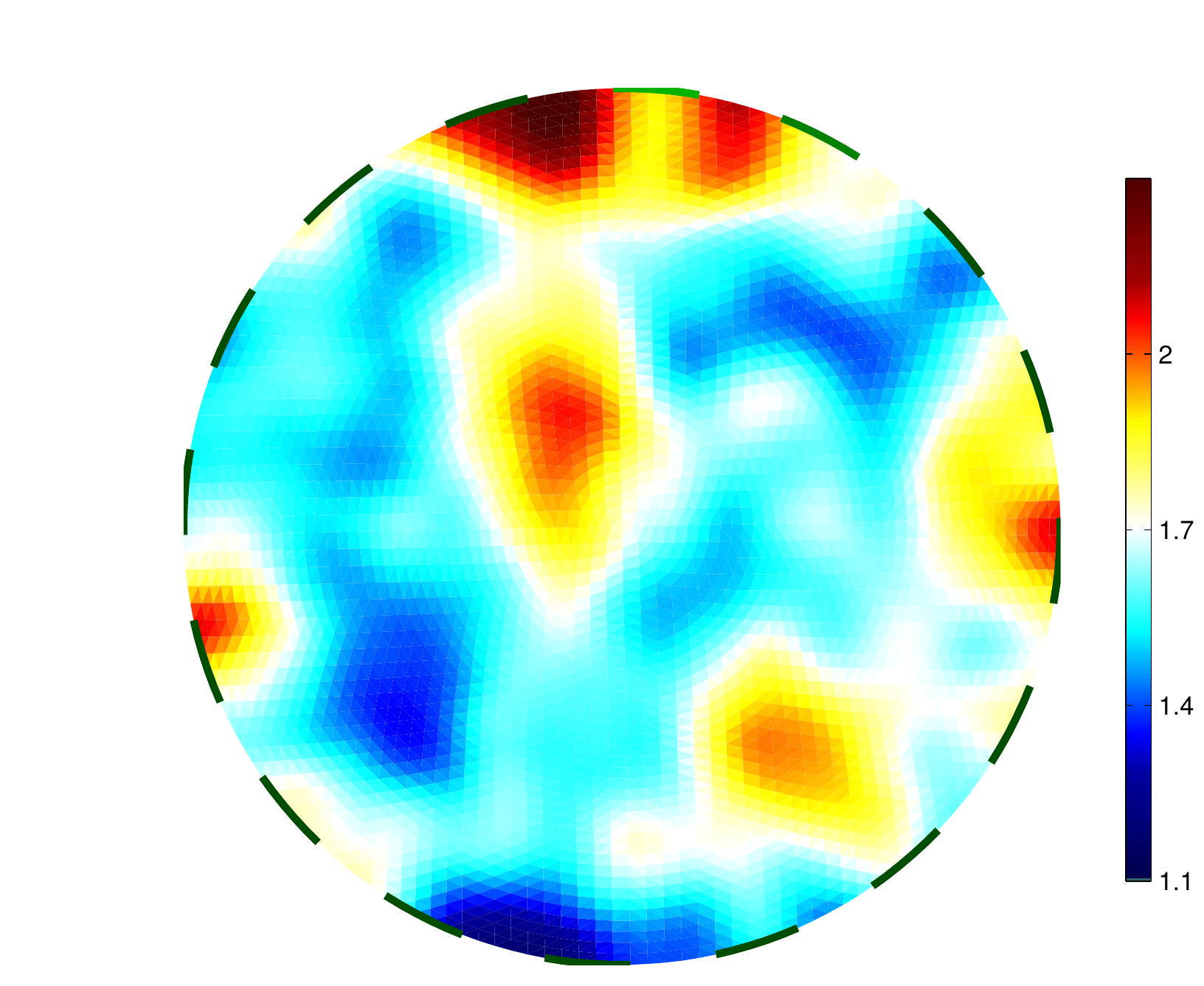}
 \caption{Draws from a Gaussian prior of the form (\ref{eq:cova1}) with (from left to right) $L=1, 0.2, 0.1, 0.06, 0.05$. Each row corresponds to a different realization of the KL coefficients in the parameterization of the Gaussian field.}
 \label{FigEIT7}

\end{center}
\end{figure}

As we have indicated earlier, the proposed scheme is highly dependent on the selection of the initial ensemble. In the experiments from the previous subsection, we attempt to reduce such dependence by selecting the initial ensemble from the same distribution that we use to generate the truth. This selection on the initial ensemble, from Proposition \ref{Prp2} ensures that the estimate from Algorithm \ref{Al1} has the same regularity/spatial structure of the truth. We wish now to investigate the effect of the selection of the initial ensemble on the estimates produced by our method. We therefore require a systematic way to generate substantially different initial ensembles. This can be achieved, for example, by modifying the parameters $L$ and $\theta$ in (\ref{eq:cova1}) that, in turn, control the spatial correlation and regularity of the samples that we use as members of such initial ensemble. For simplicity we focus only on the selection of $L$  and consider initial ensembles with the covariance (\ref{eq:cova1}) for different choices of  $L$. In \Fref{FigEIT7} we present some samples generated with several choices of $L$ (from left to right: $L=1, 0.2, 0.1, 0.06, 0.05$). Each row corresponds to a fixed set of realizations of coefficients in the KL expansion that we use to generate those Gaussian fields.   We can observe visually how the  spatial correlation of the samples decreases with $L$.

\begin{figure}[htbp]
\begin{center}
\includegraphics[scale=0.214]{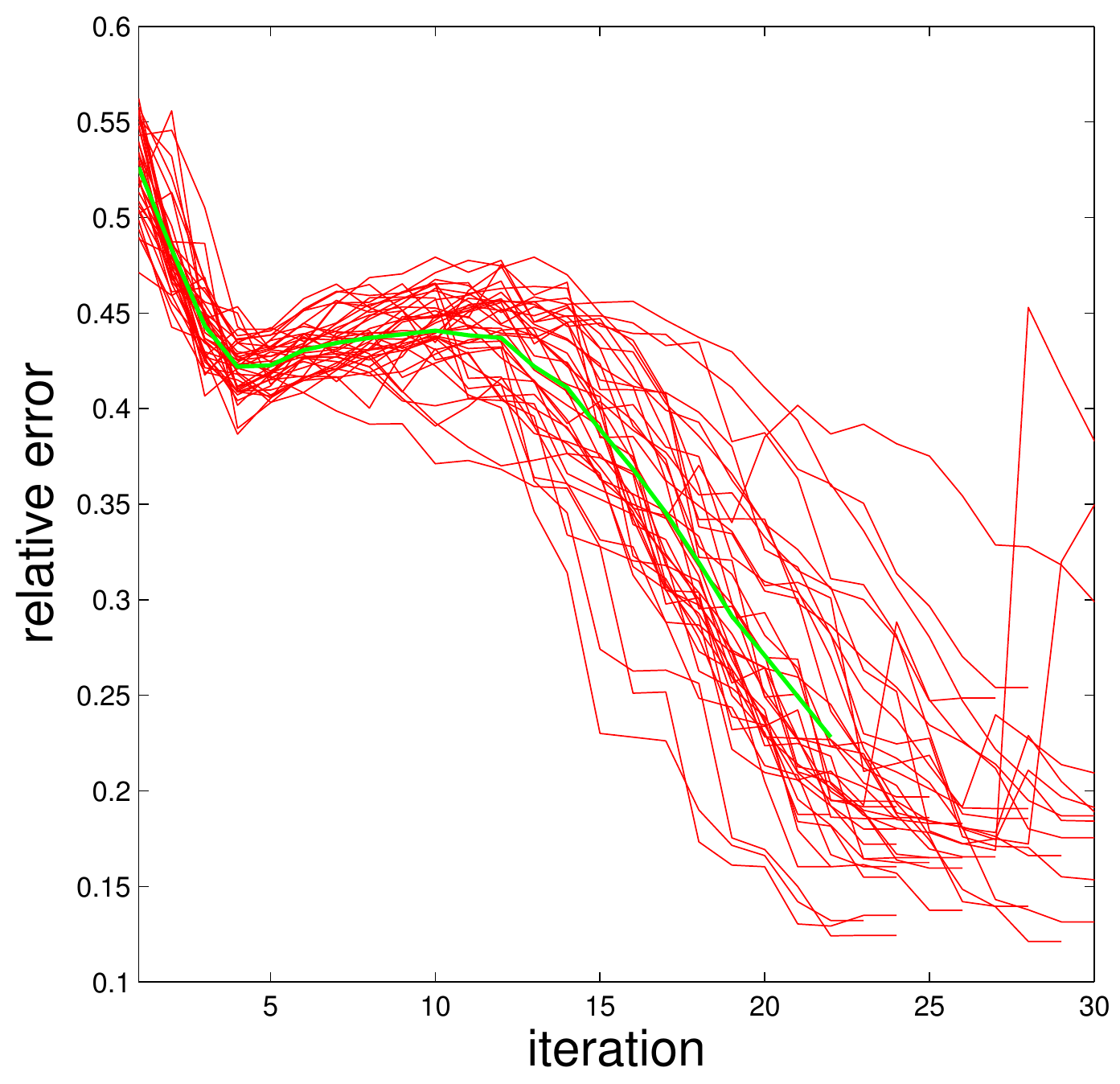}
\includegraphics[scale=0.214]{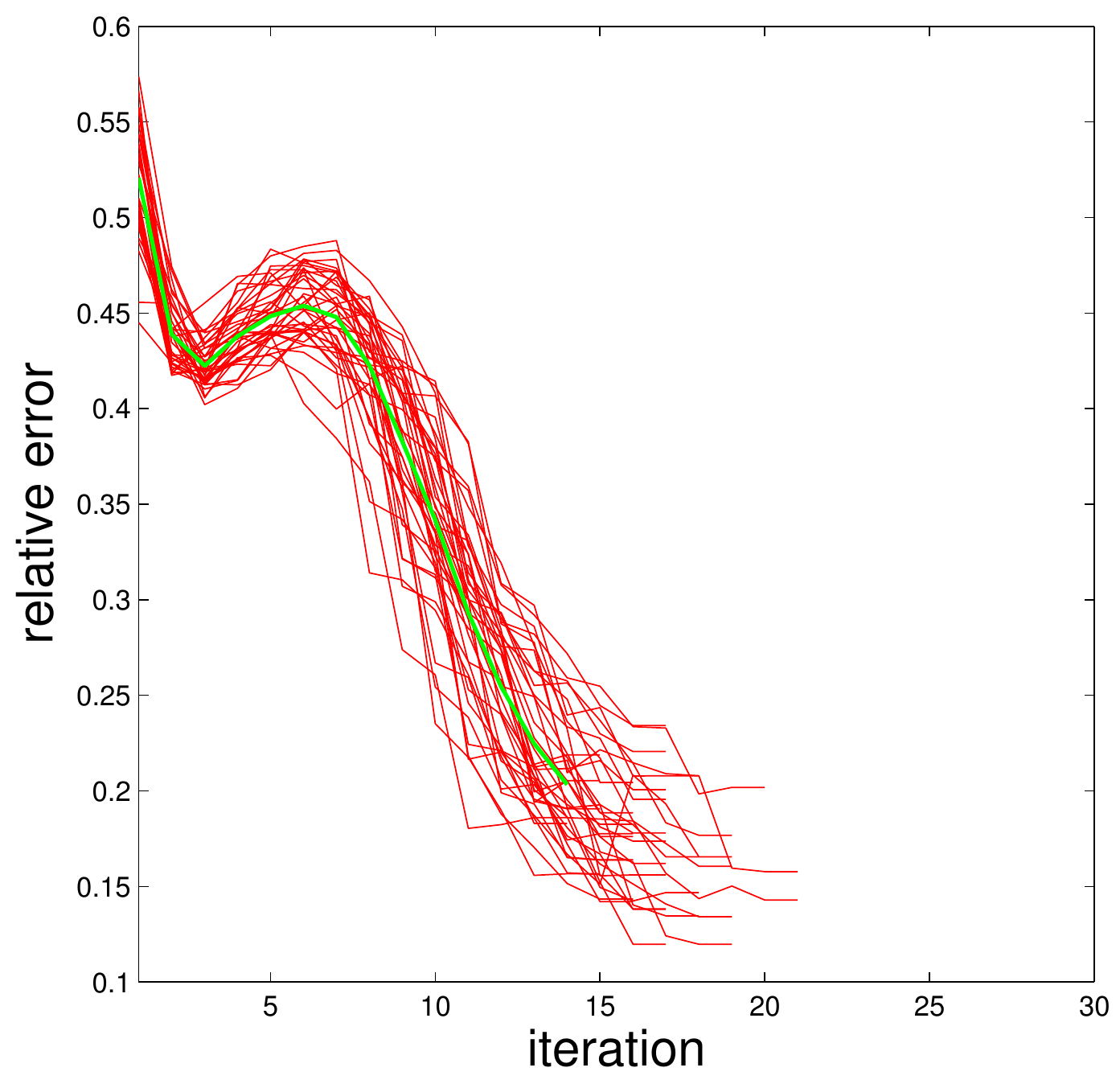}
\includegraphics[scale=0.214]{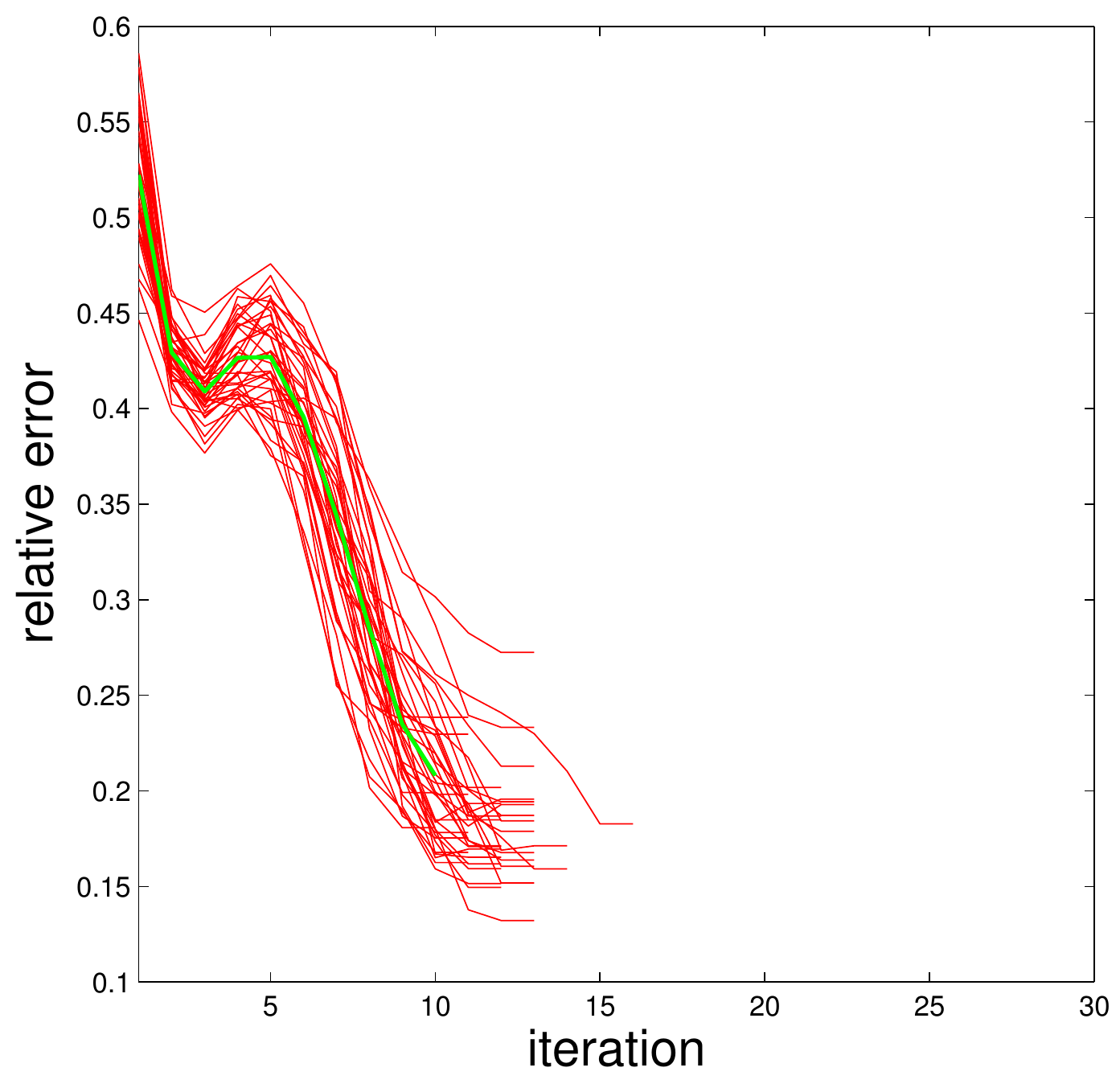}
\includegraphics[scale=0.214]{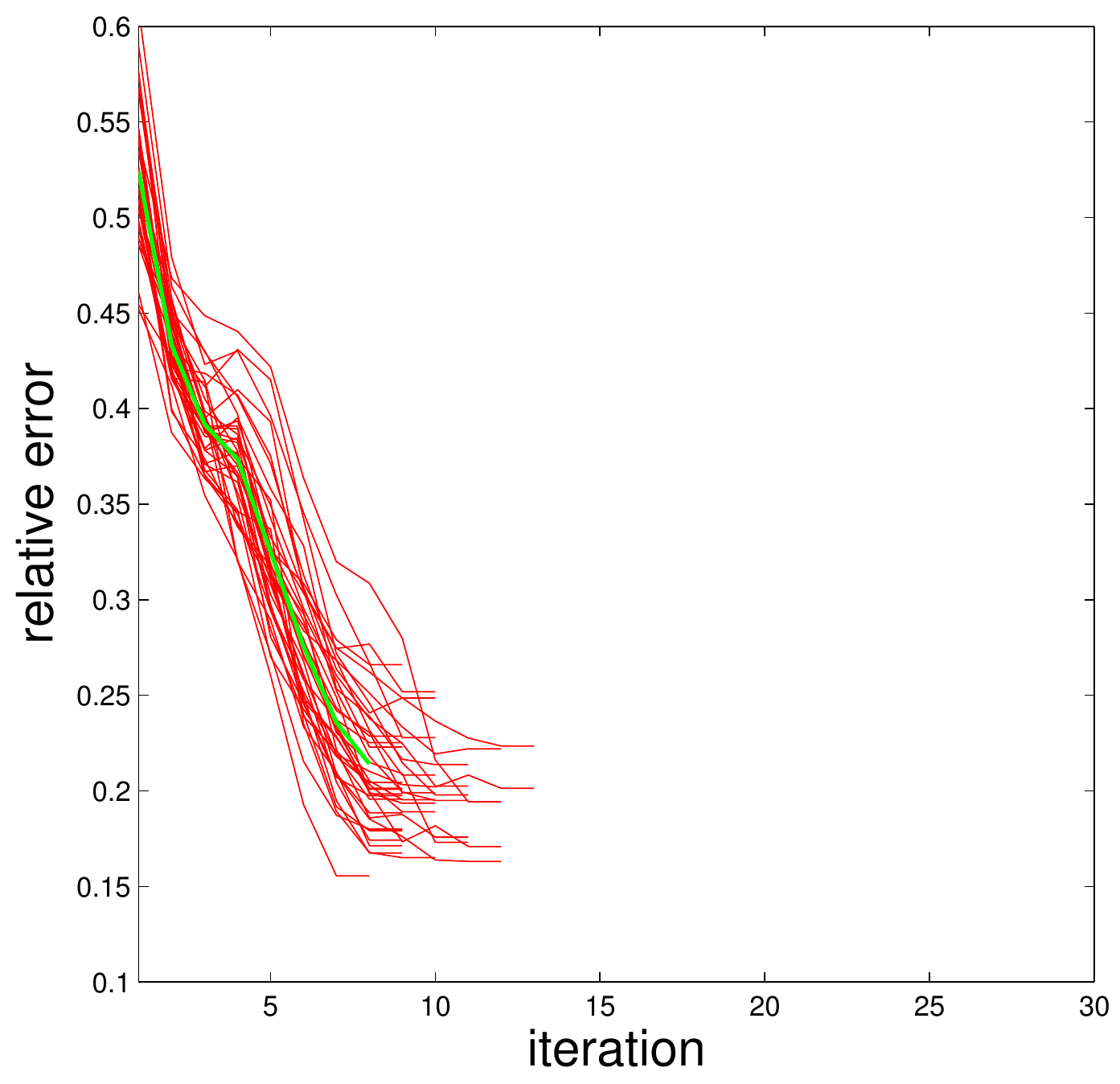}
\includegraphics[scale=0.214]{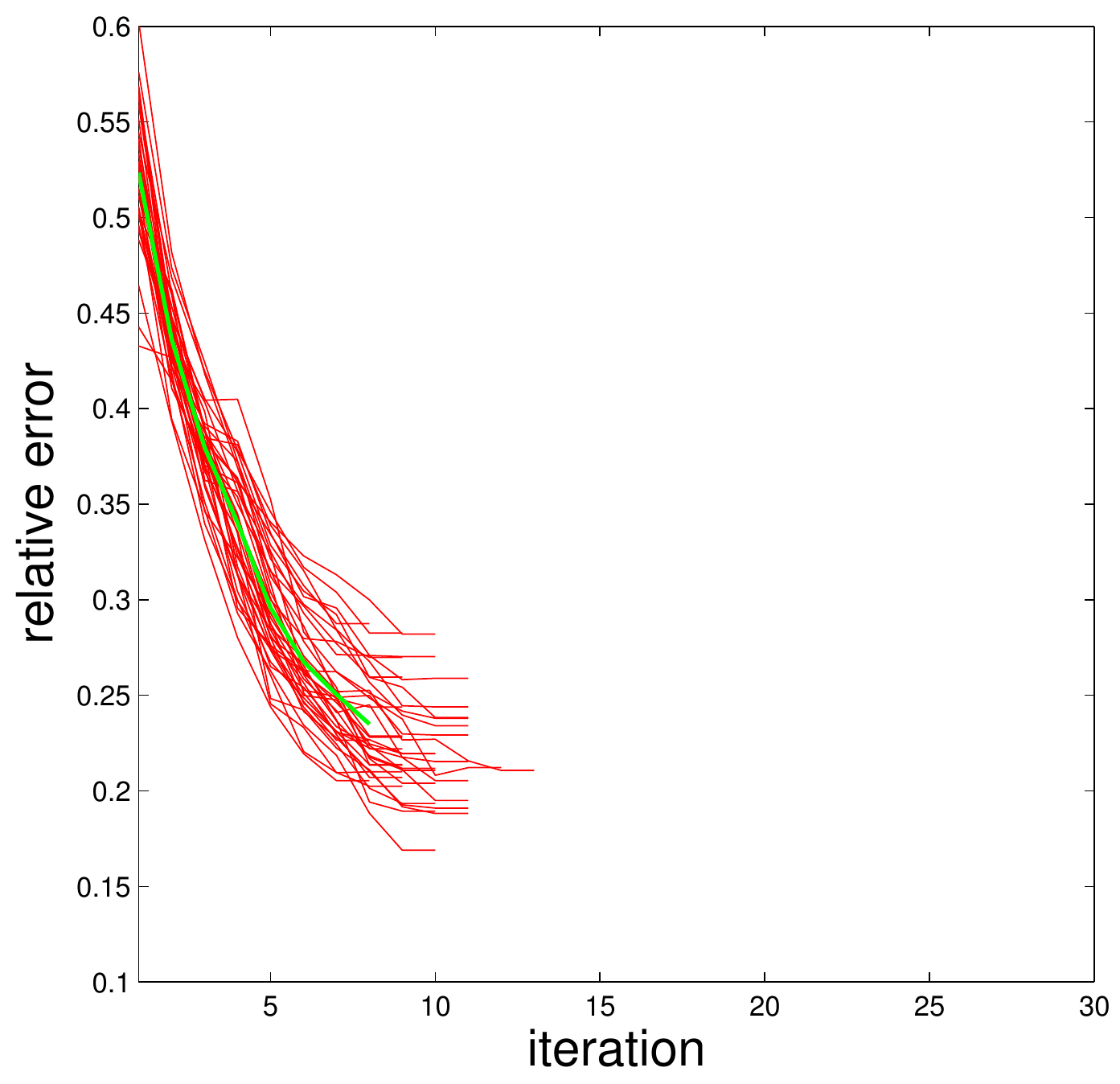}\\
\includegraphics[scale=0.214]{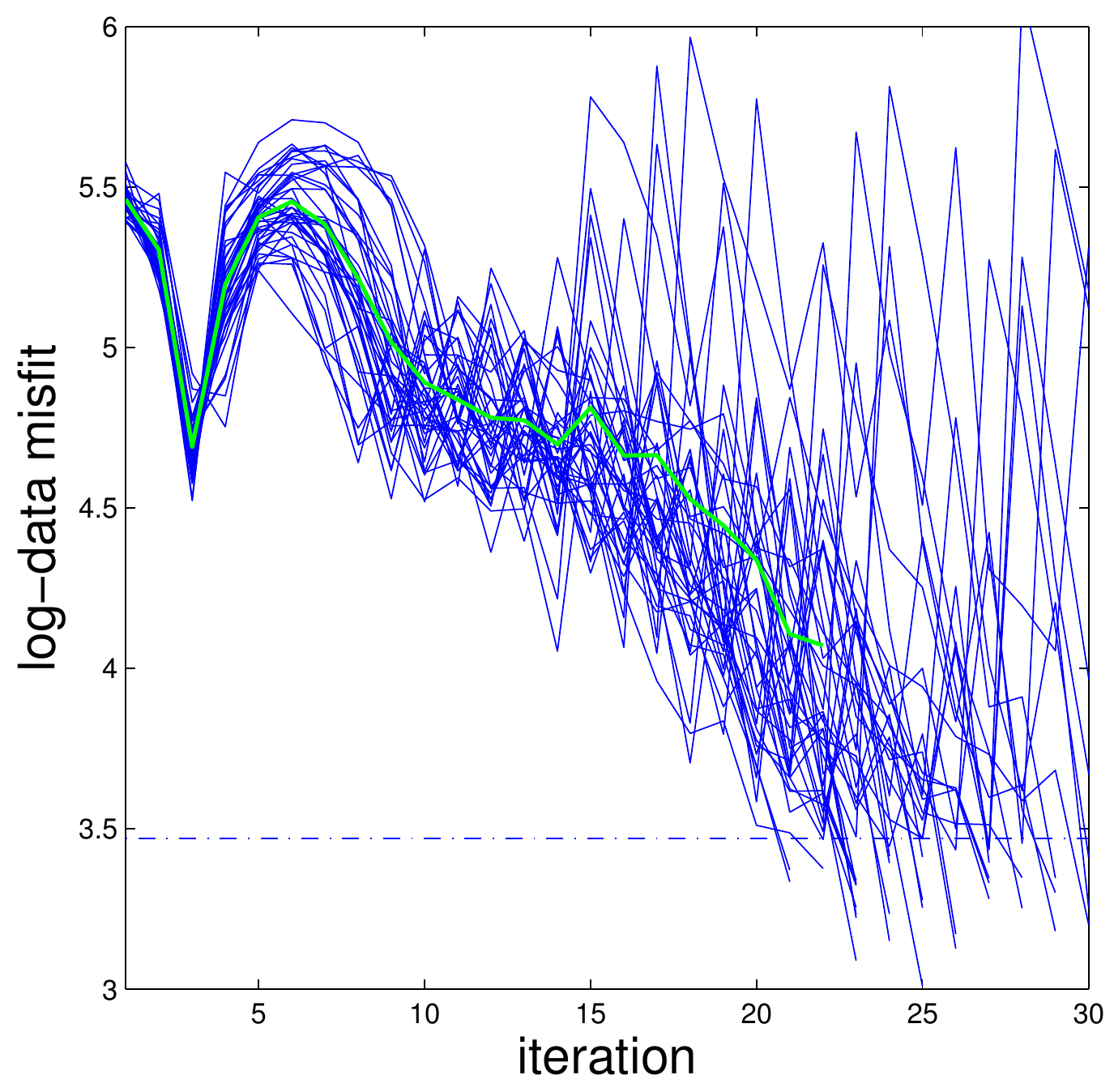}
\includegraphics[scale=0.214]{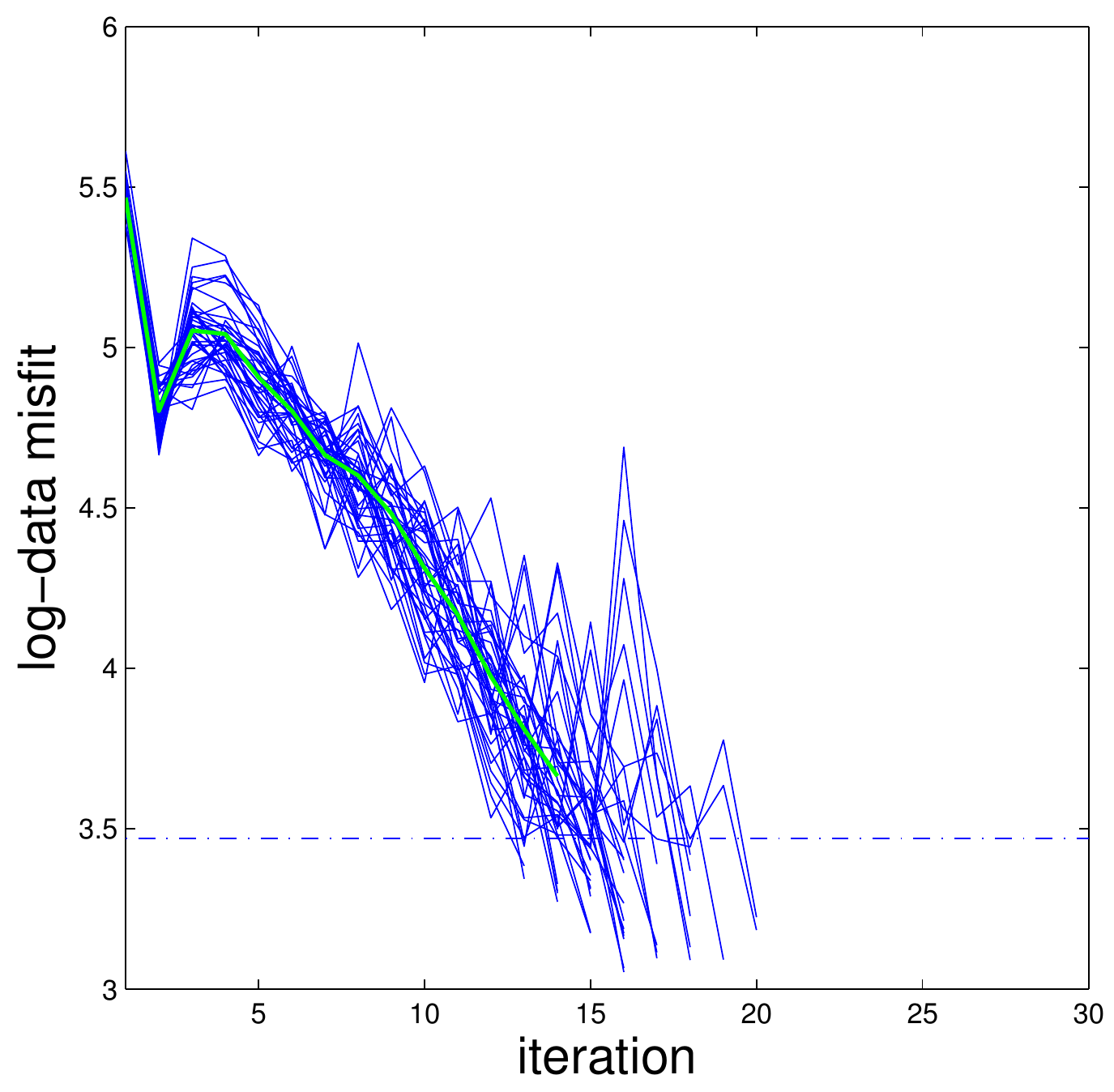}
\includegraphics[scale=0.214]{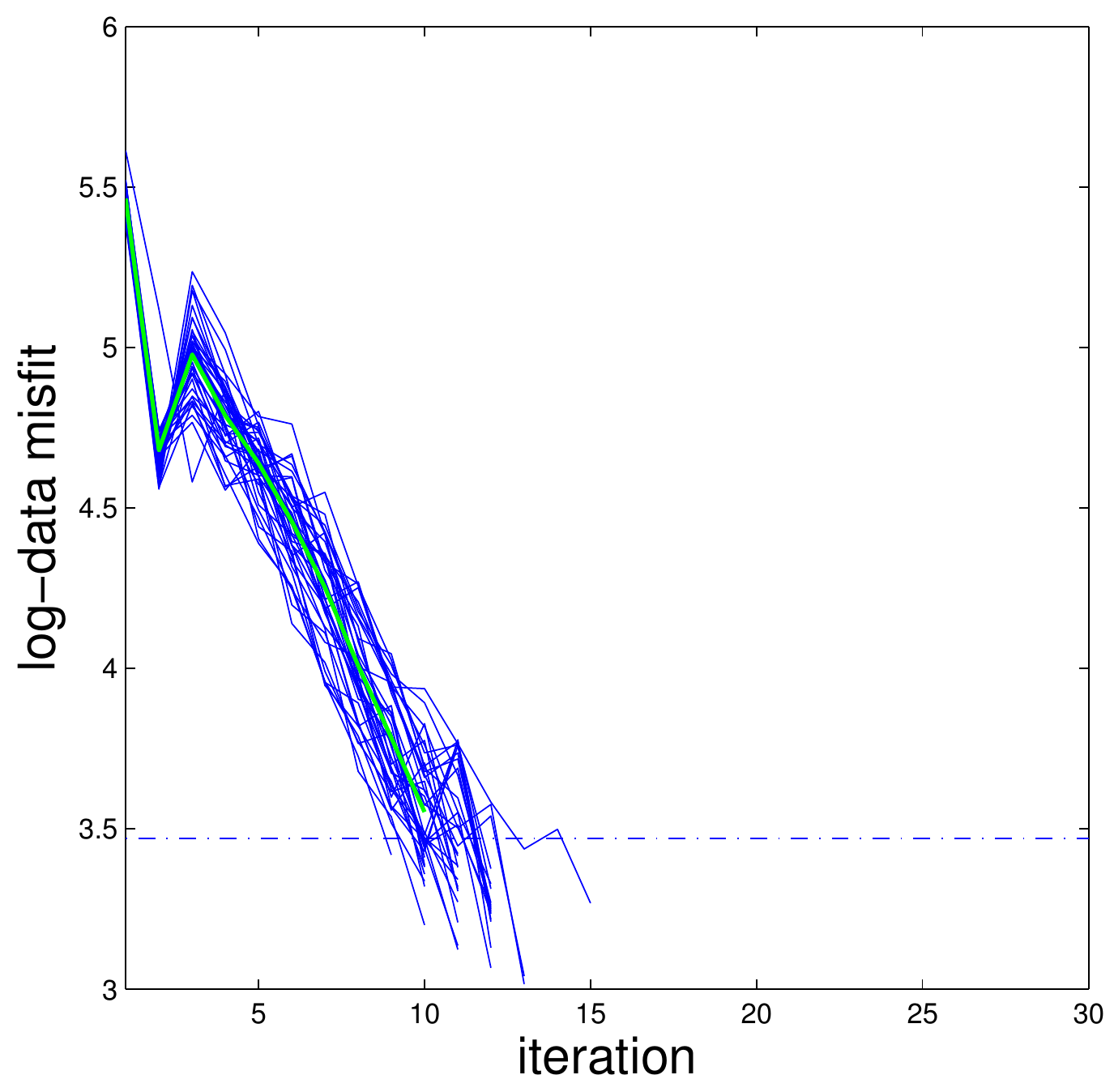}
\includegraphics[scale=0.214]{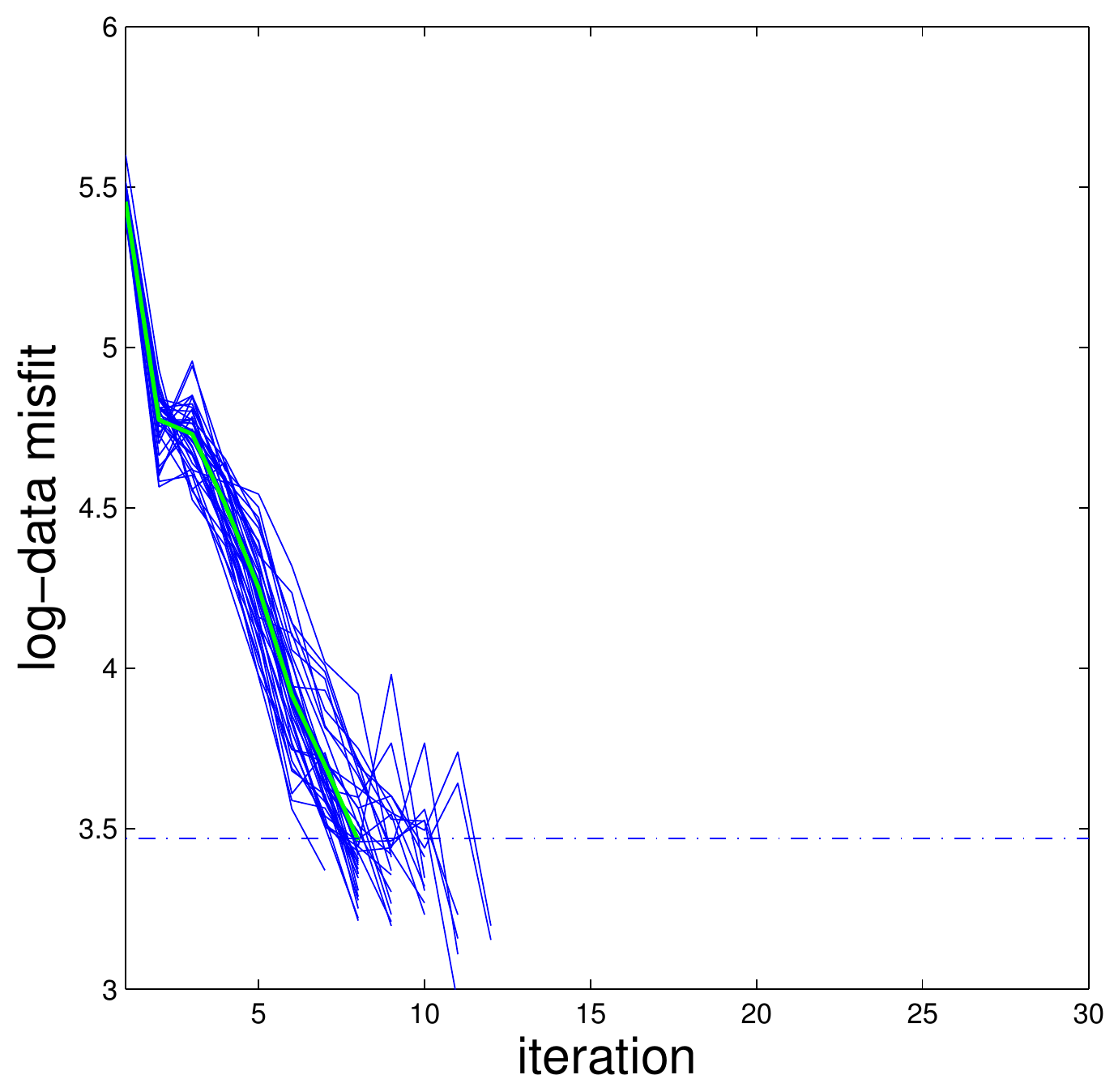}
\includegraphics[scale=0.214]{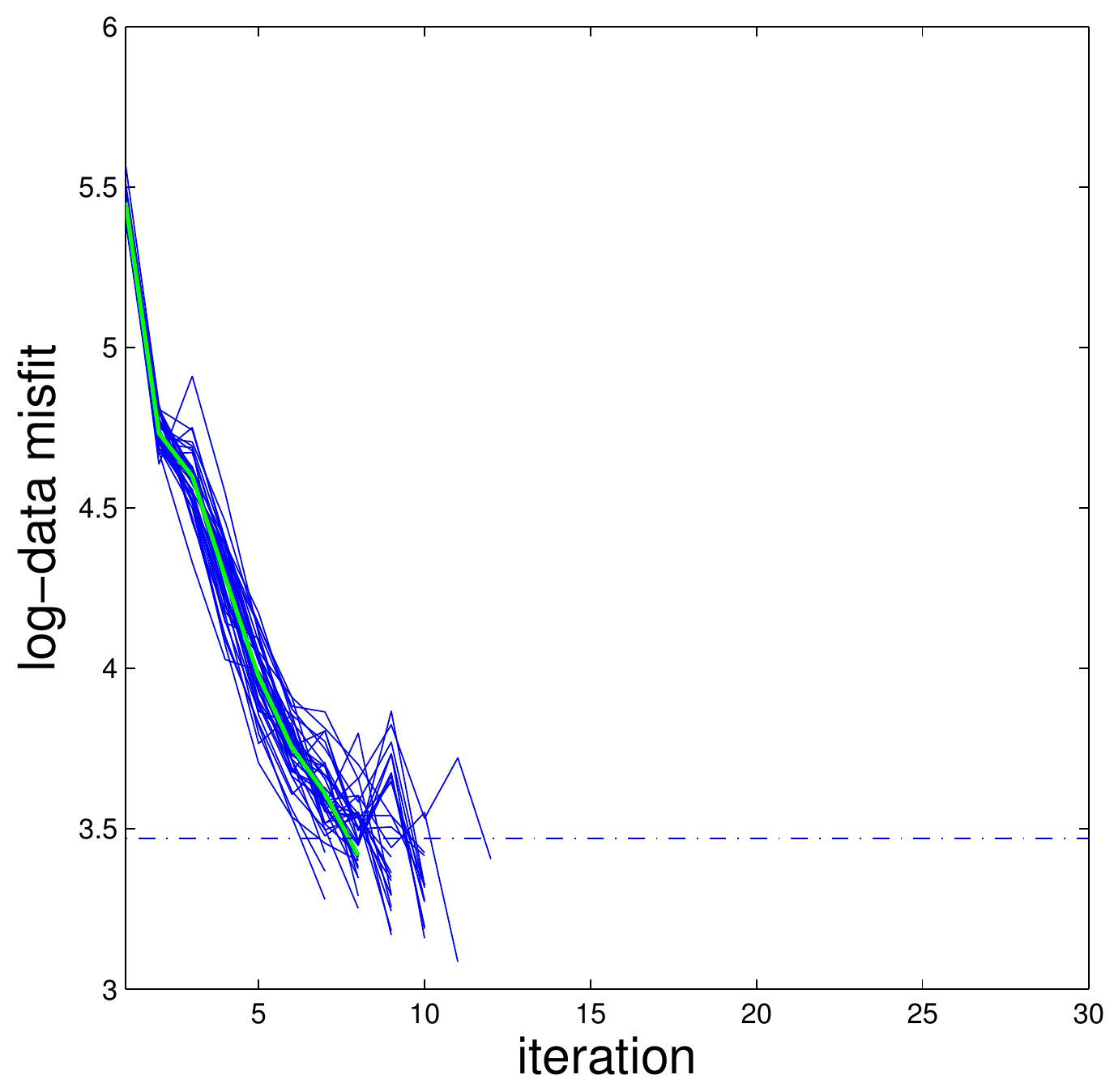}

 \caption{Error with respect to the truth (top) and log data misfit (bottom) from 40 experiments with different initial ensembles generated from samples of $N(0,C)$ with values of $L$ in (\ref{eq:cova1}) (from left to right) $L=0.2, 0.1, 0.06,0.04, 0.03$. }\ \label{FigEIT8}

\end{center}
\end{figure}

We apply Algorithm \ref{Al1} for the solution of the EIT problem with an ensemble of size $N_{e}=100$ generated from the Gaussian measure described earlier with different choices of parameter $L$. For each value of $L$, we report the results from 40 experiments corresponding to different initial ensembles (i.e. different draws from the Gaussian measure). \Fref{FigEIT8} displays the resulting log-data misfit (bottom) and error w.r.t truth (top) . The selection of the initial ensemble based on $L$ has a significant effect on the performance of the method. We note that the accuracy degrades as we decrease the correlation length. For the largest correlation length considered here ($L=1$) we note that after a few iterations, before the data misfit reaches the value $\eta/\rho$, the relative error shows a slight increase which results in large fluctuations of the data misfit. As we decrease the correlation length the relative error is stabilized and the data misfit reduced. However, we note that the accuracy decreases as we use initial ensembles generated with members that have smaller correlation lengths ($L<0.1$). Indeed, for $L=0.05$, the relative error with respect to the truth increases before the data misfit drops below  $\eta/\rho$,. It is thus clear that the selection of the initial ensemble (in this case based on $L$) affects the regularizing properties of the scheme. 

\begin{figure}[htbp]
\begin{center}

\includegraphics[scale=0.18]{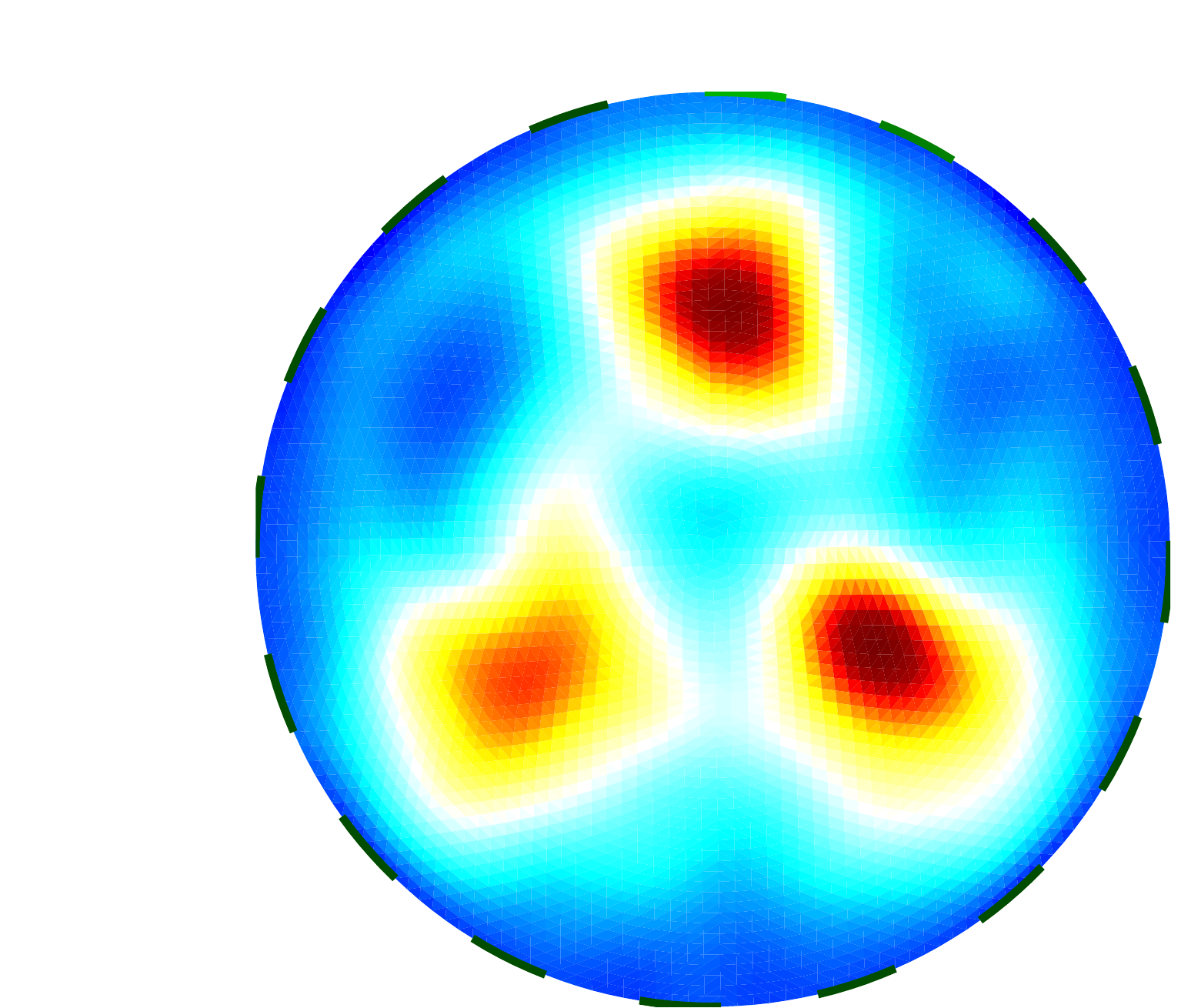}
\includegraphics[scale=0.18]{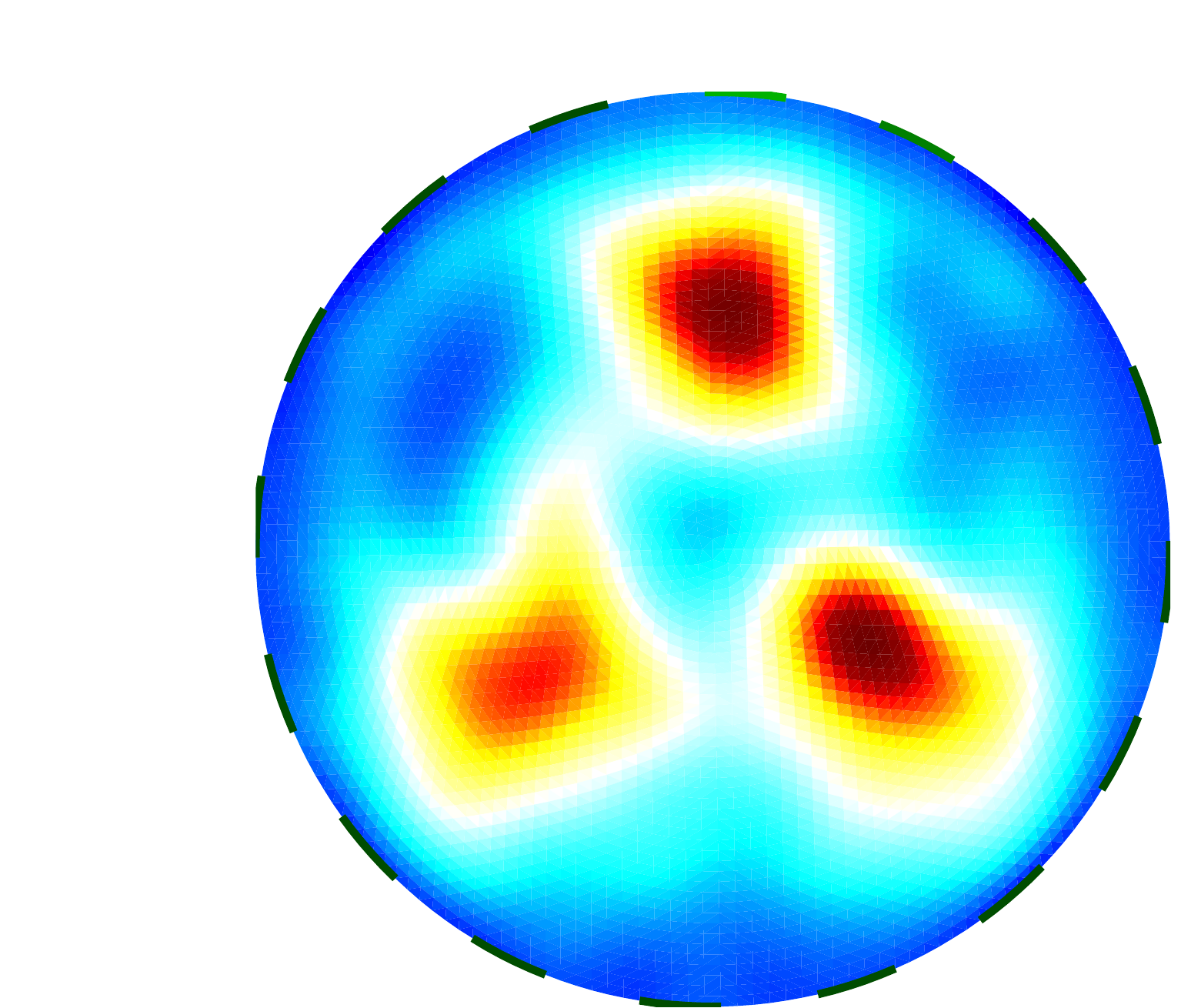}
\includegraphics[scale=0.18]{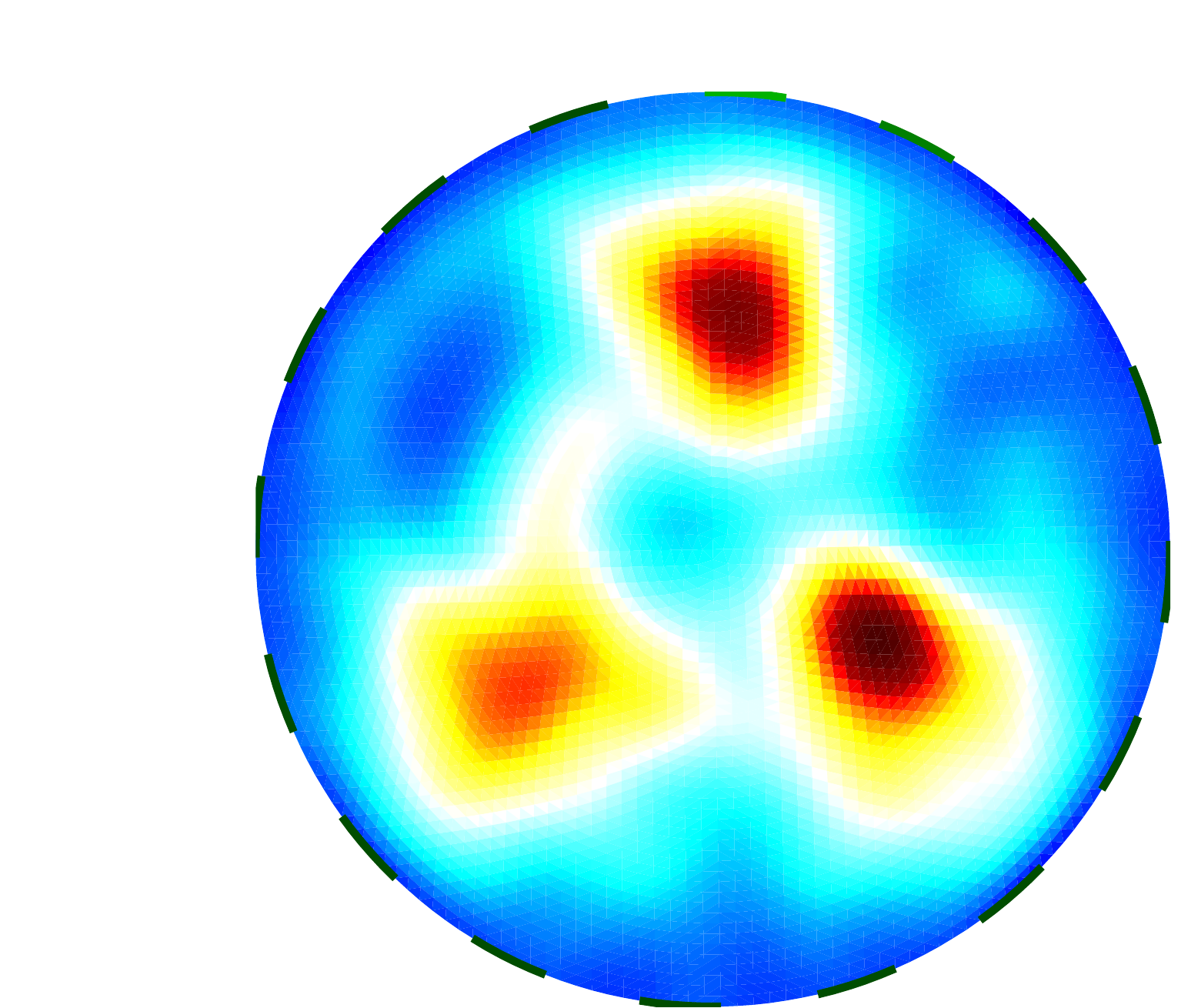}
\includegraphics[scale=0.18]{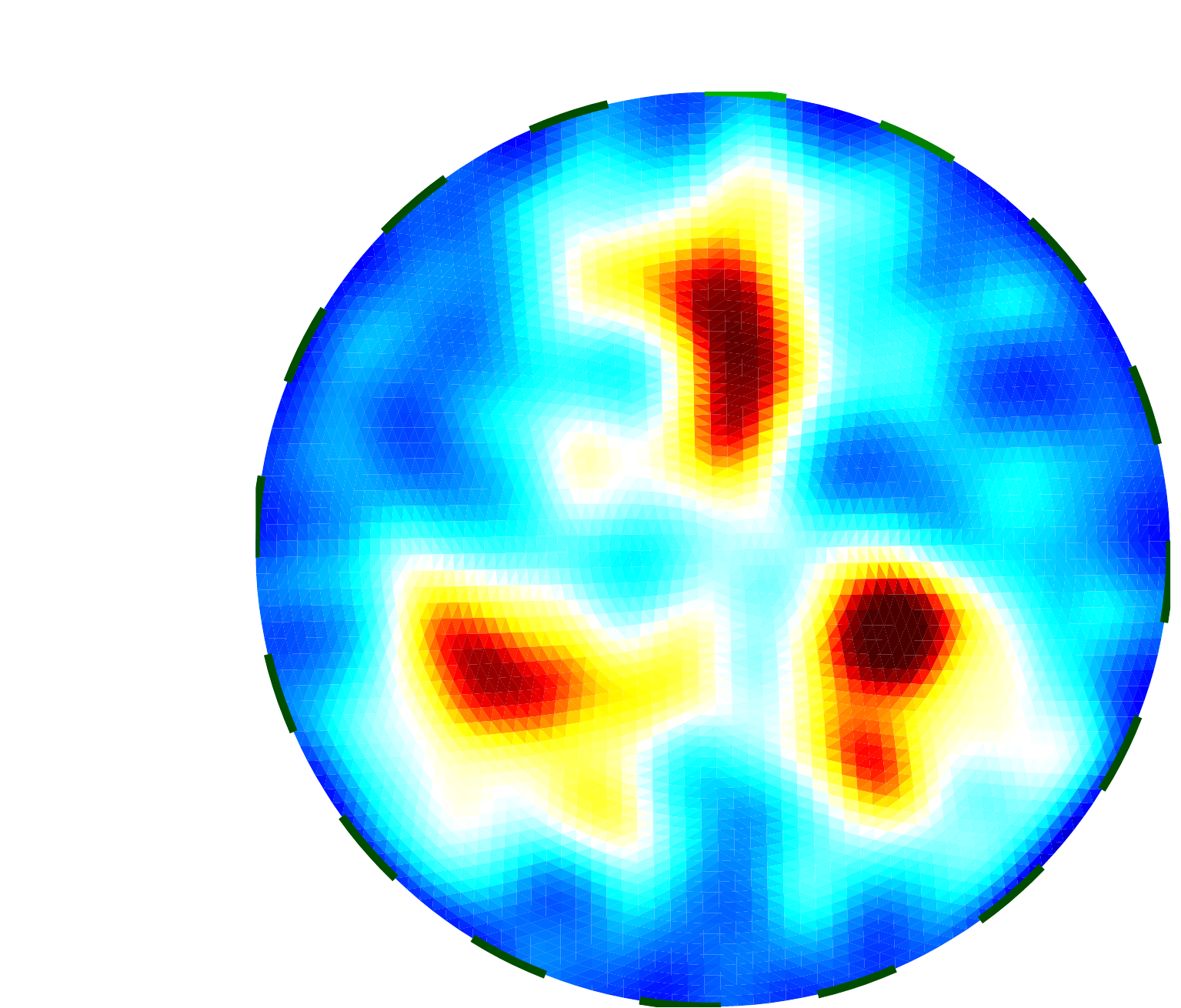}
\includegraphics[scale=0.18]{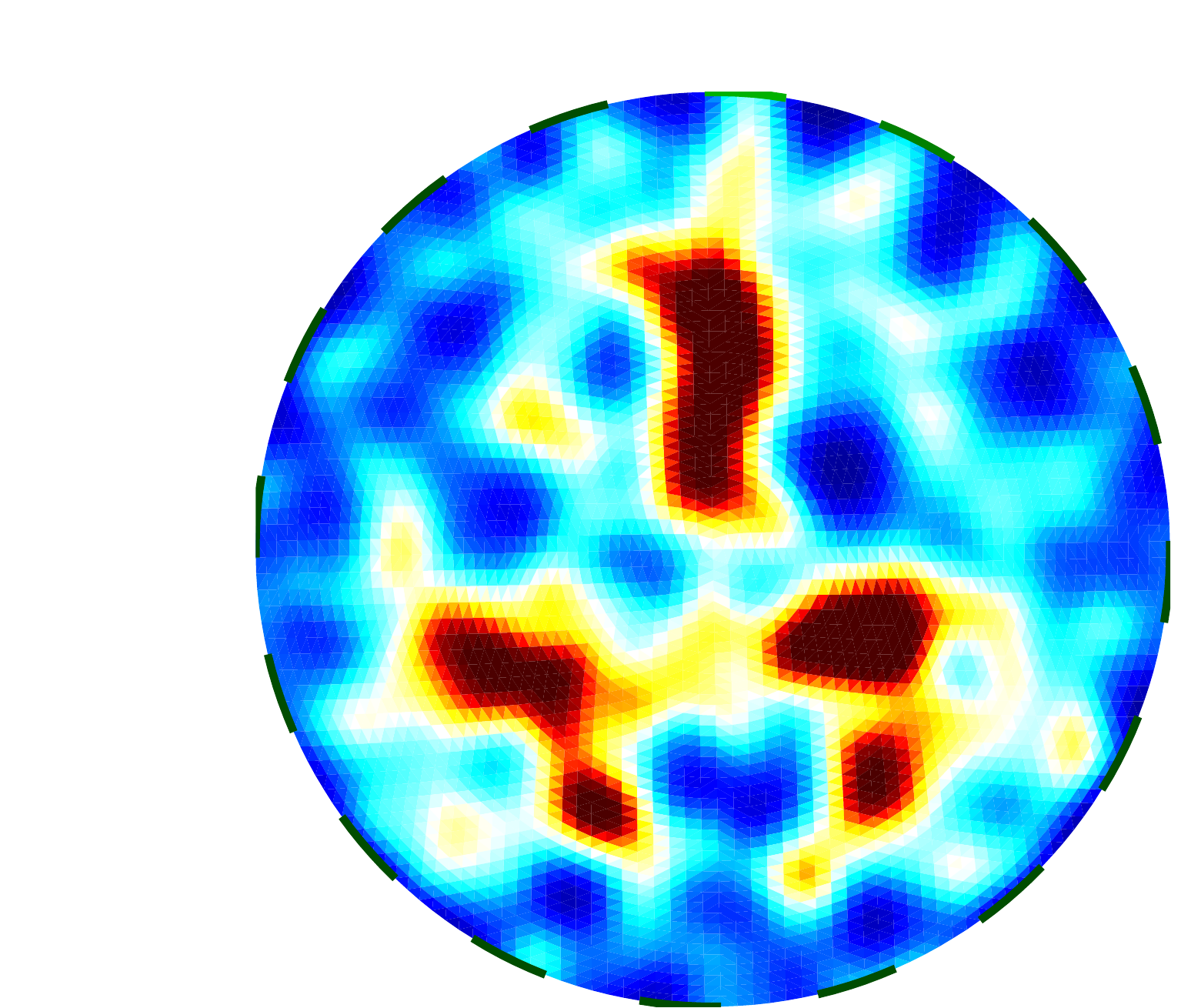}\\
\includegraphics[scale=0.18]{True_EIT}
\caption{Top: Estimates of level-set obtained from 5 experiments with different initial ensembles with samples of $N(0,C)$ with values of $L$ in (\ref{eq:cova1}) (from left to right) $L=0.2, 0.1, 0.06,0.05$. Bottom: True conductivity. }  \label{FigEIT9}

\end{center}
\end{figure}

In \Fref{FigEIT9} we display estimates (i.e. ensemble mean $\overline{u}_{n}$) of the log-conductivity obtained from one of the the previous experiments for the different choices of $L$. More concretely, these are the estimates obtained with different initial ensemble corresponding to different choices of $L$ but generated with the same realization of KL coefficients. From the invariance subspace property of the proposed scheme, it comes as no surprise that an initial ensemble of elements that have very small spatial correlation cannot possibly produce an estimate that characterizes the truth. Indeed, the inclusions of high conductivity from the truth have a diameter of approximately 0.3 which can be better identified when the initial ensemble is generated from smooth fields with a similar spatial correlation. Similar dependence on correlation lengths were obtained in the Bayesian level-set approach of \cite{Level_set_US}. While prior knowledge of the truth can be certainly used for the construction of the initial ensemble, further research should address the potential identification of parameters that determine the regularity of the initial ensemble that better characterize the truth. 



\section{Applications to geometric inverse problems}\label{Sec:Applications}

The numerical investigations of the preceding sections establishes that the proposed ensemble Kalman method, with a sufficiently large ensemble size, inherits the regularizing properties of the regularizing LM approach of \cite{Hanke}. Thus, our scheme can be potentially applied to produce stable and accurate estimates of PDE-constrained inverse problems. Moreover, those estimates can be computed at a reasonable computational cost without the need of the derivative of the forward map. In this section, we show the potential application of the proposed method to solve identification problems where the computation of such derivative can be very cumbersome. In concrete, in this section we investigate the application of the proposed regularizing ensemble Kalman method for the solution of shape identification problems.

\
We are interested in the identification of an unknown region $\Omega$  within $D$, the spatial domain of definition of the PDE under consideration, whose boundary $\partial \Omega$ defines a sharp discontinuity of an unknown parameter in the PDE model. That is the case, for example, in subsurface flow applications where the conductivity of the aquifer/reservoir may take substantially different (by several orders of magnitude) values on facies characterized with different geologies. While the nominal values on such geologic facies may be known a priori, the interface between geologic facies is usually unknown.

A common approach to parameterize an unknown interface (or shape) is to use an implicitly parameterization in terms of level-sets \cite{Osher}. In other words, we may assume that
\begin{eqnarray}\label{level}
\Omega = \{x\in D|~u(x)\leq 0\}
\end{eqnarray}
where here $u$ denotes the level-set function. The zero level-set provides the unknown interface which is now parameterized by a function $u$ and thus the identification can be posed within a functional-analytical framework. Level-set-based methods for shape identification problems have been extensively applied in the literature \cite{BurgerSurvey, BurgerGB,Iglesias3,Santosa}. However, to our best of knowledge, most of these approaches use a variational framework for solving shape-based least-squares problems which, in turn, requires computation of the shape derivatives of the forward map. In addition, these standard variational formulations require the solution of the level-set equation to evolve the shape so that it minimizes the, possibly regularized, squared data misfit. This approach often leads to computations of flat level-set functions that need to be redefined. In this section we show that the proposed ensemble Kalman method can be used to approximate the solution to identification problems without using neither the shape derivative of forward maps nor the level-set equation. More precisely, we apply the proposed iterative ensemble Kalman method for the solution of identification problems where the interface is parameterized by (\ref{level}). The update formulate (\ref{eq:m16}) with the controlled/regularized selection of $\alpha$ induces the motion of the shape. Moreover, we propose to apply the Kalman method with a selection of the initial ensemble of level-sets that consist of samples from Gaussian measures $N(0,C)$ which a covariance that enforces some smoothness of the level-set function and could potentially incorporate some intrinsic correlation length. Nonetheless, this probabilistic approach for the generation of the initial ensemble is conducted artificially for the sake of the implementation. However, it has been recently shown  \cite{Level_set_US} that there is zero probability of generating samples (and thus initial ensembles) corresponding to flat level-set functions. This, combined with the invariance subspace property of the proposed method, ensures that our estimates do not become flat thus overcoming the common issue of the flattening of the level-set function typical of standard level-set approaches.

\subsection{Estimation of geologic facies}\label{geo_facies}

\begin{figure}[htbp]
\begin{center}
\includegraphics[scale=0.28]{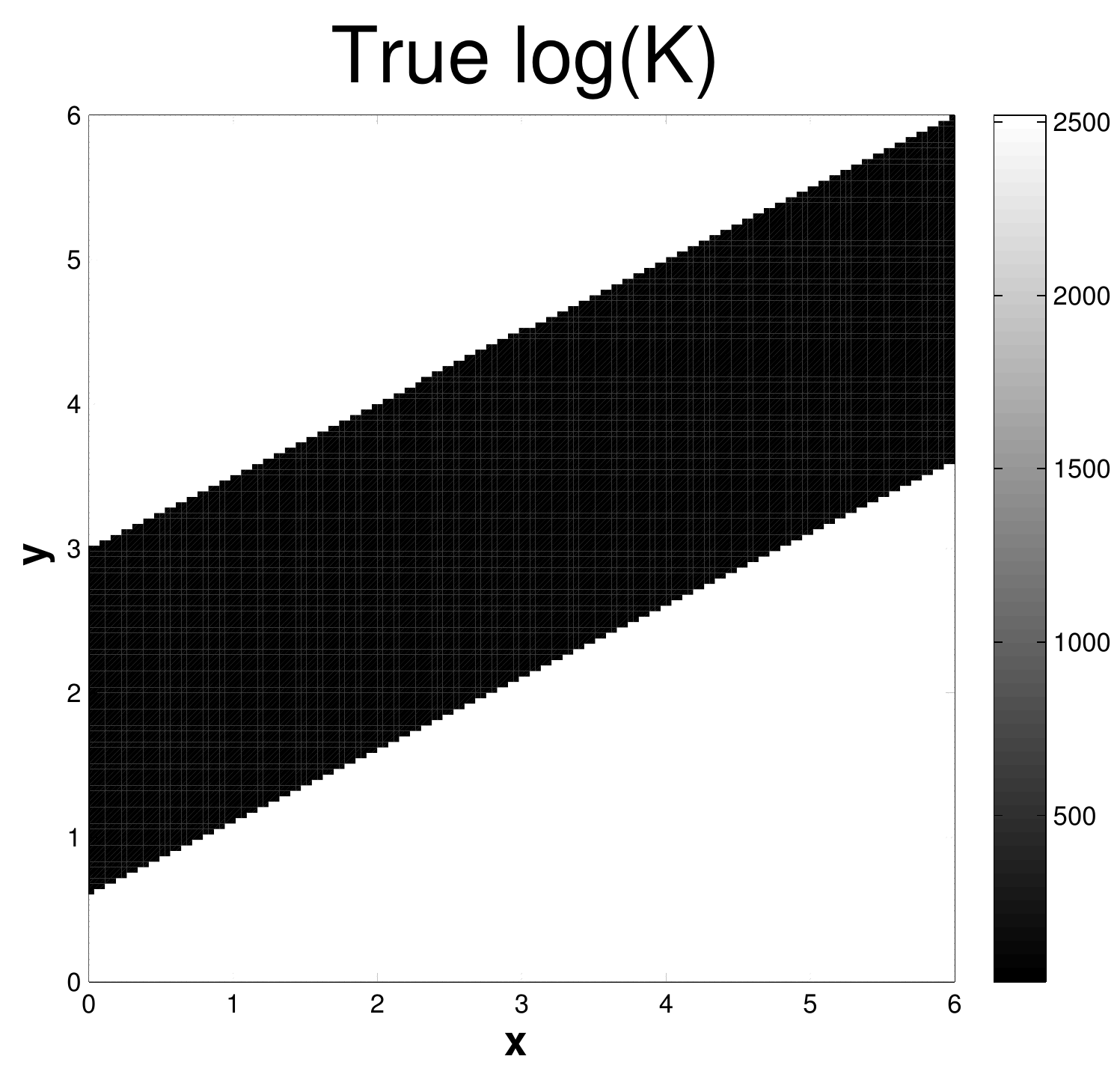}
\includegraphics[scale=0.28]{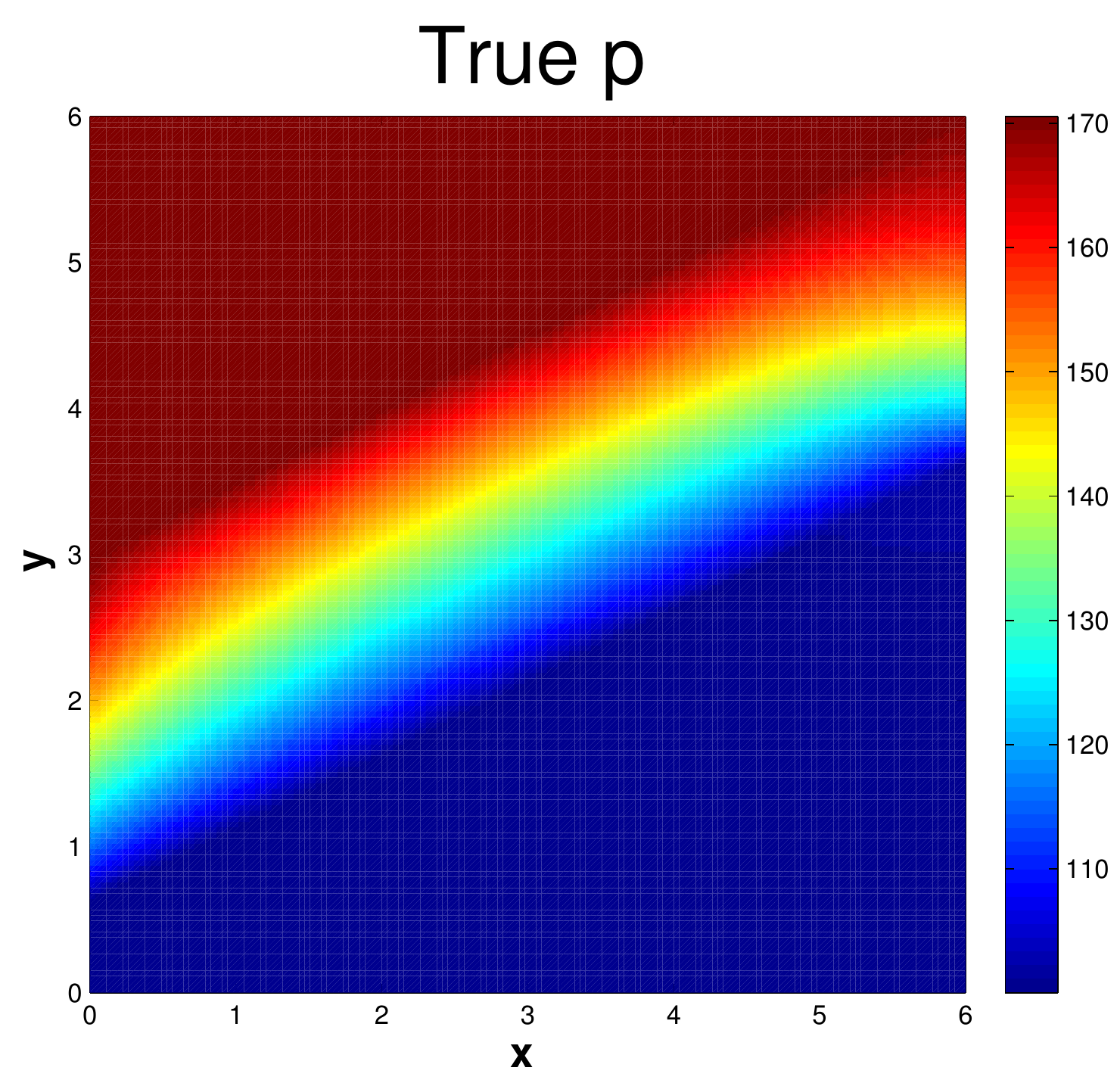}
\includegraphics[scale=0.28]{Well_Loc}
 \caption{Left: True hydraulic conductivity.  Middle: true hydraulic head. Right: Measurement locations. }
    \label{Fig_gw1}
\end{center}
\end{figure}

We are interested in the identification of geologic facies in the test model described in section \ref{test}. More precisely we want to estimate the region of high conductivity $\Omega$ (represented by (\ref{level})) in an aquifer $D$. For this example we prescribe a true conductivity $K^{\dagger}$ displayed in Figure \ref{Fig_gw1} (middle). This true conductivity simulates a typical layered structure in the geologic properties of an aquifer. We use $K^{\dagger}$ in (\ref{eq:a1})-(\ref{eq:a3}) to generate synthetic data that we wish to invert in order to identify the true conductivity. The measurement configuration is displayed in \Fref{Fig_gw1} (right). As before, inverse crimes are avoided by using a finer grid for the generation of the truth than the one used for the initial ensemble and so for the inversion. 

Our goal is now to apply the proposed ensemble Kalman method on the variable $u$ corresponding to the level-set function that parametrizes $\kappa$ as follows
\begin{eqnarray}\label{103}
\kappa(u)= \kappa_{i}\chi_{u\leq 0 }+\kappa_{e}\chi_{u>0},
\end{eqnarray}
where $\chi_{A}$ denotes the characteristic function of region $A$ and $\kappa_{i}$ and $\kappa_{e}$ are the positive constants that define the high and low conductivity values of the true conductivity \Fref{Fig_gw1} (left). We assume these two constants are known; the unknown is the interface between the regions of different conductivity.

\begin{figure}[htbp]
\begin{center}%
\includegraphics[scale=0.7]{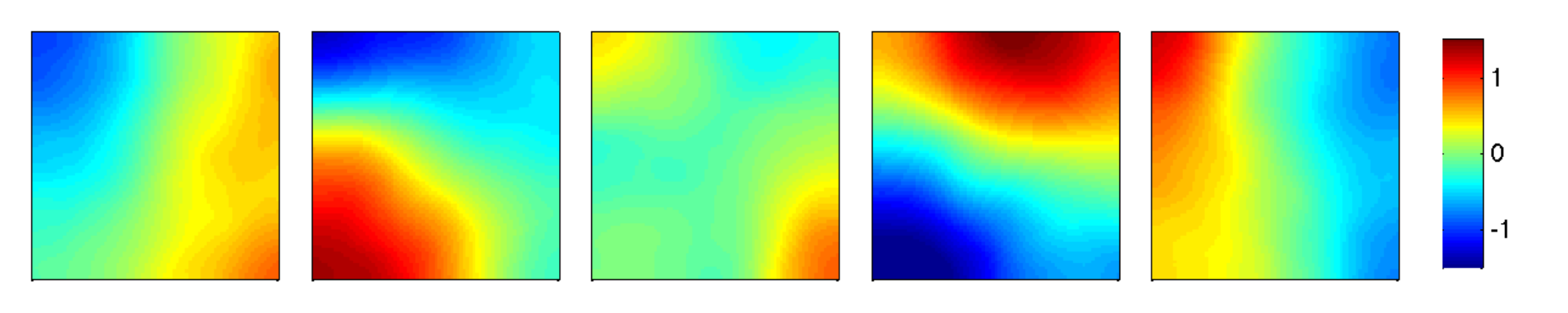}
\includegraphics[scale=0.7]{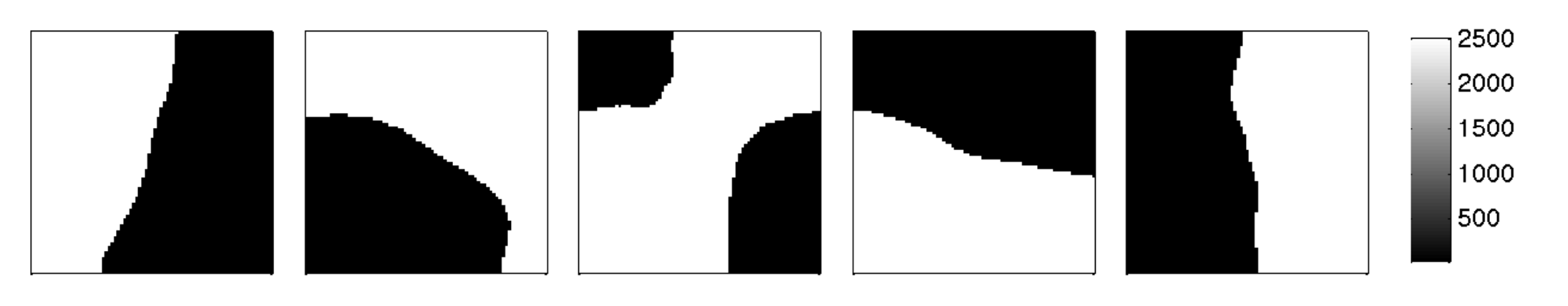}\\
 \caption{Top: Different realizations from level-set functions sampled of $N(0,C)$ with $C$ from (\ref{eq:cova3}). Bottom: Conductivities computed from these sample level-set functions (by means of (\ref{103})).} \label{Fig_gw2}

\end{center}
\end{figure}

As described earlier we prescribe an artificial Gaussian measure from which we generate members of an initial ensemble of level-set functions. More concretely we consider samples from $N(0,C)$ with 
  \begin{equation}\label{eq:cova3}
C =\omega (-\Delta)^{-\theta}
\end{equation}
where $\Delta$ is the Laplacian operator, $\omega$ is a  scaling constant and $\theta$ controls the regularity. In contrast to (\ref{eq:cova1}), the previous expressions defines a covariance with an intrinsic fixed correlation length that we are not able to vary. In Figure \ref{Fig_gw2} (top row) we show some samples from the prior distribution and in Figure \ref{Fig_gw2} (bottom row) we show the corresponding conductivities obtained from (\ref{103}).

\begin{figure}[htbp]
\begin{center}
\includegraphics[scale=0.2175]{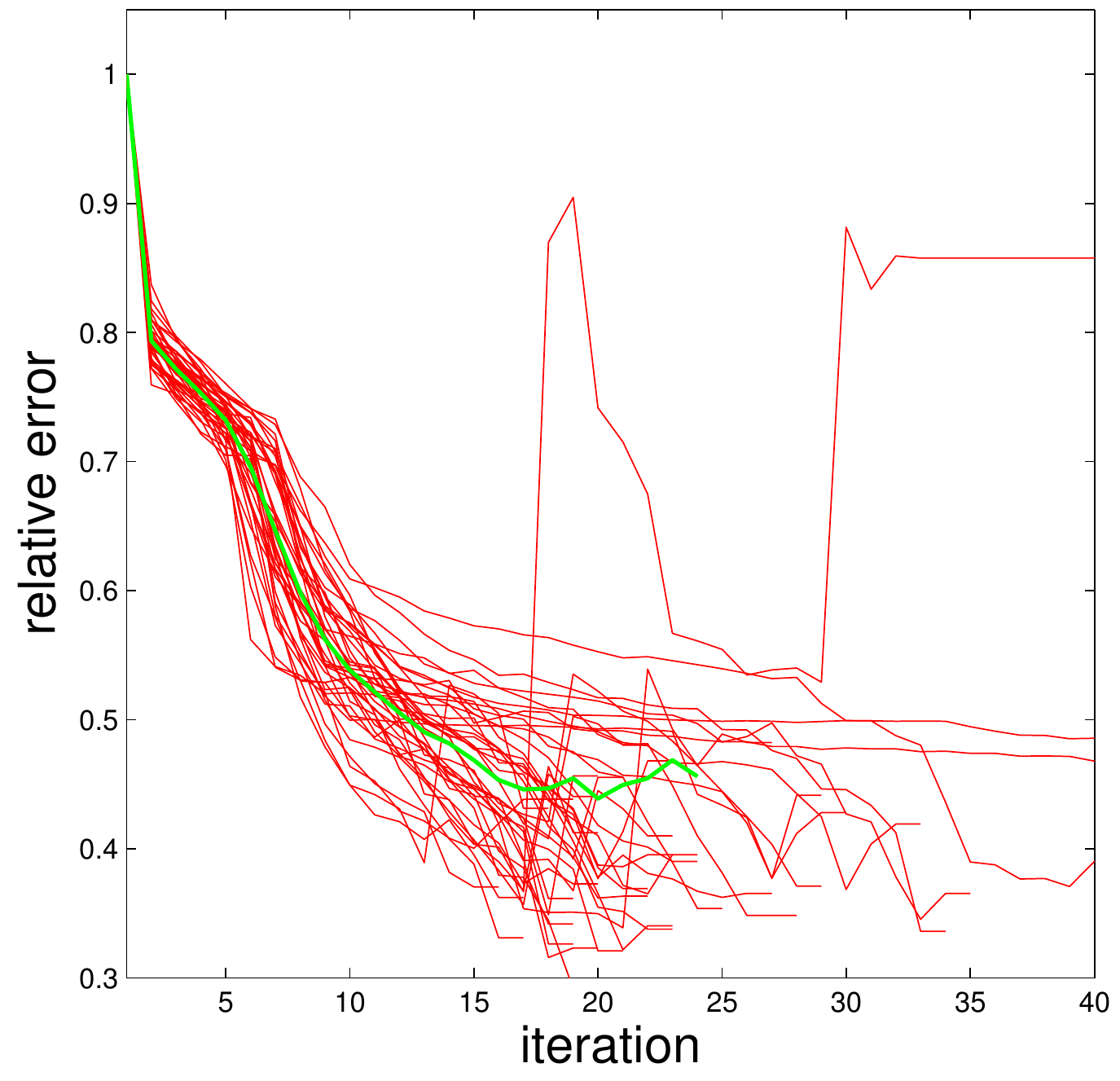}
\includegraphics[scale=0.2175]{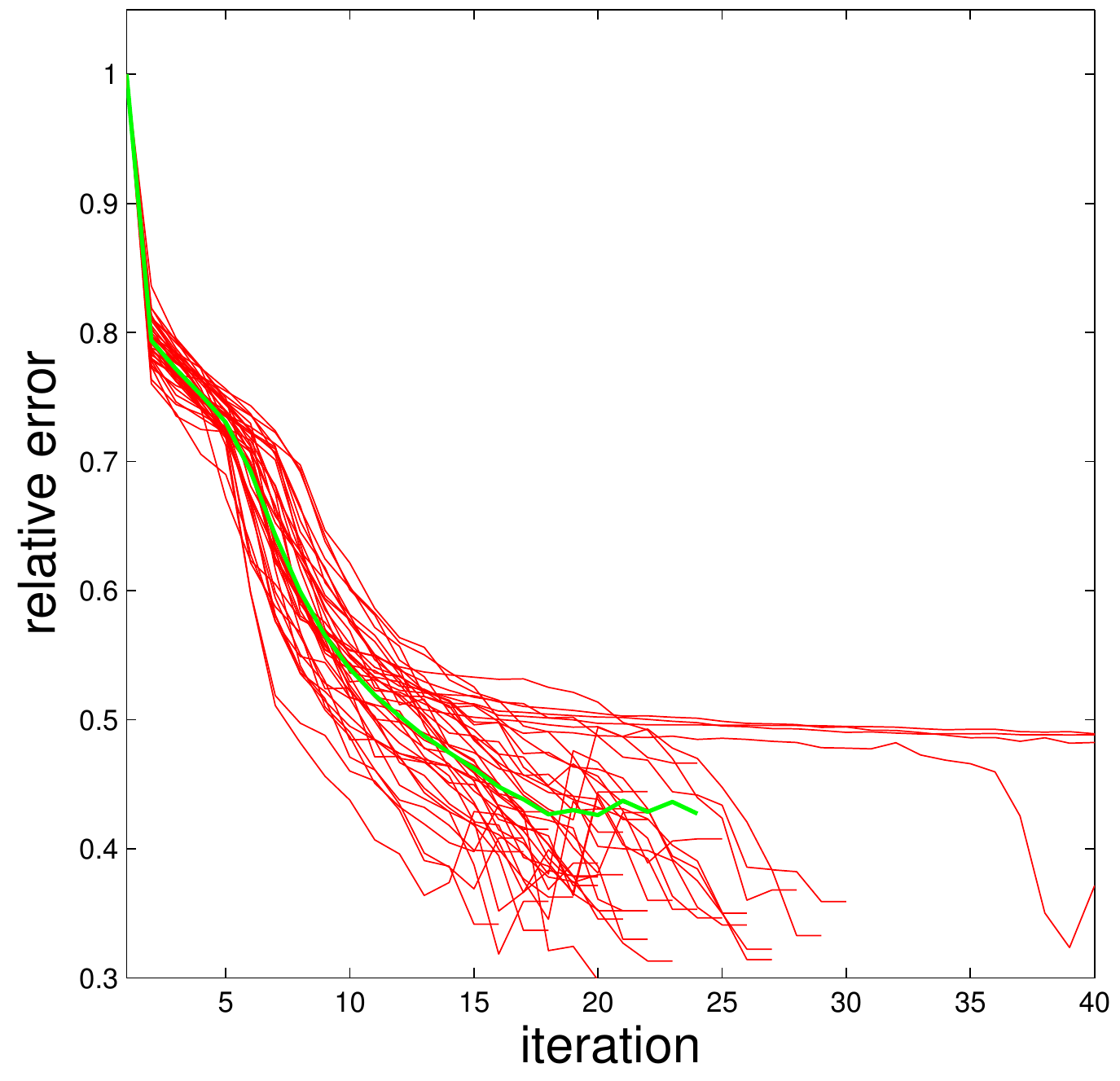}
\includegraphics[scale=0.2175]{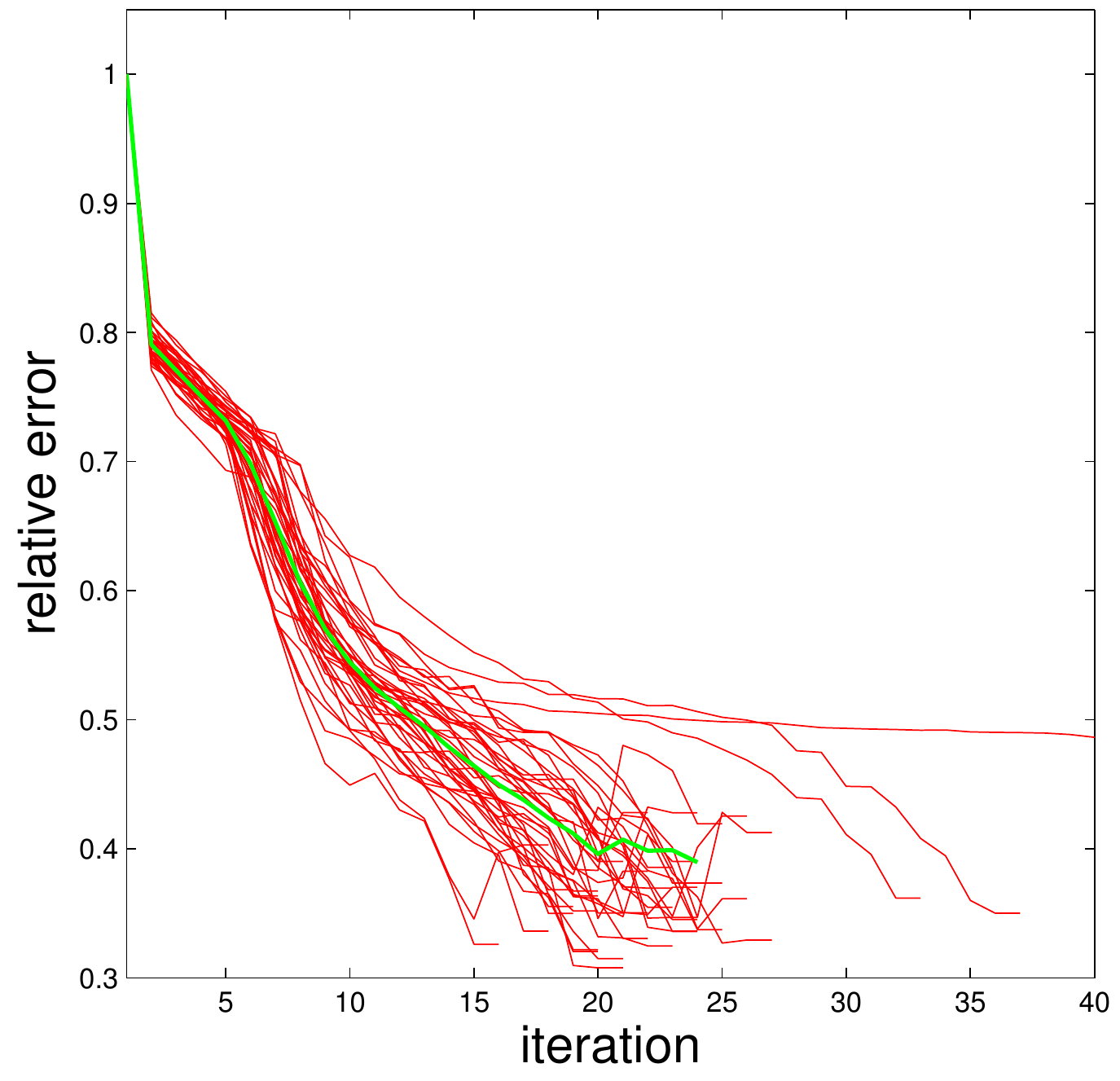}
\includegraphics[scale=0.2175]{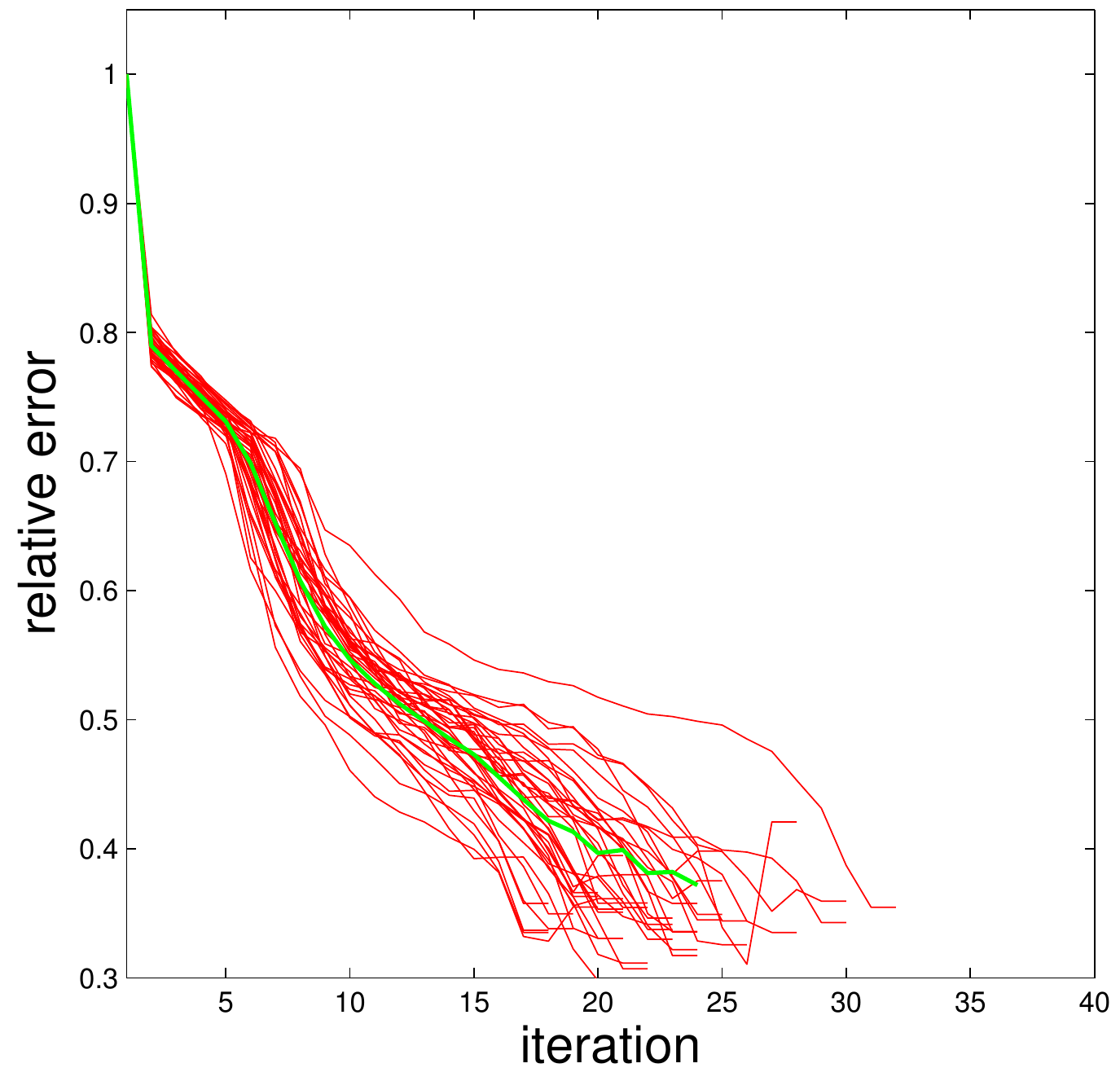}
\includegraphics[scale=0.2175]{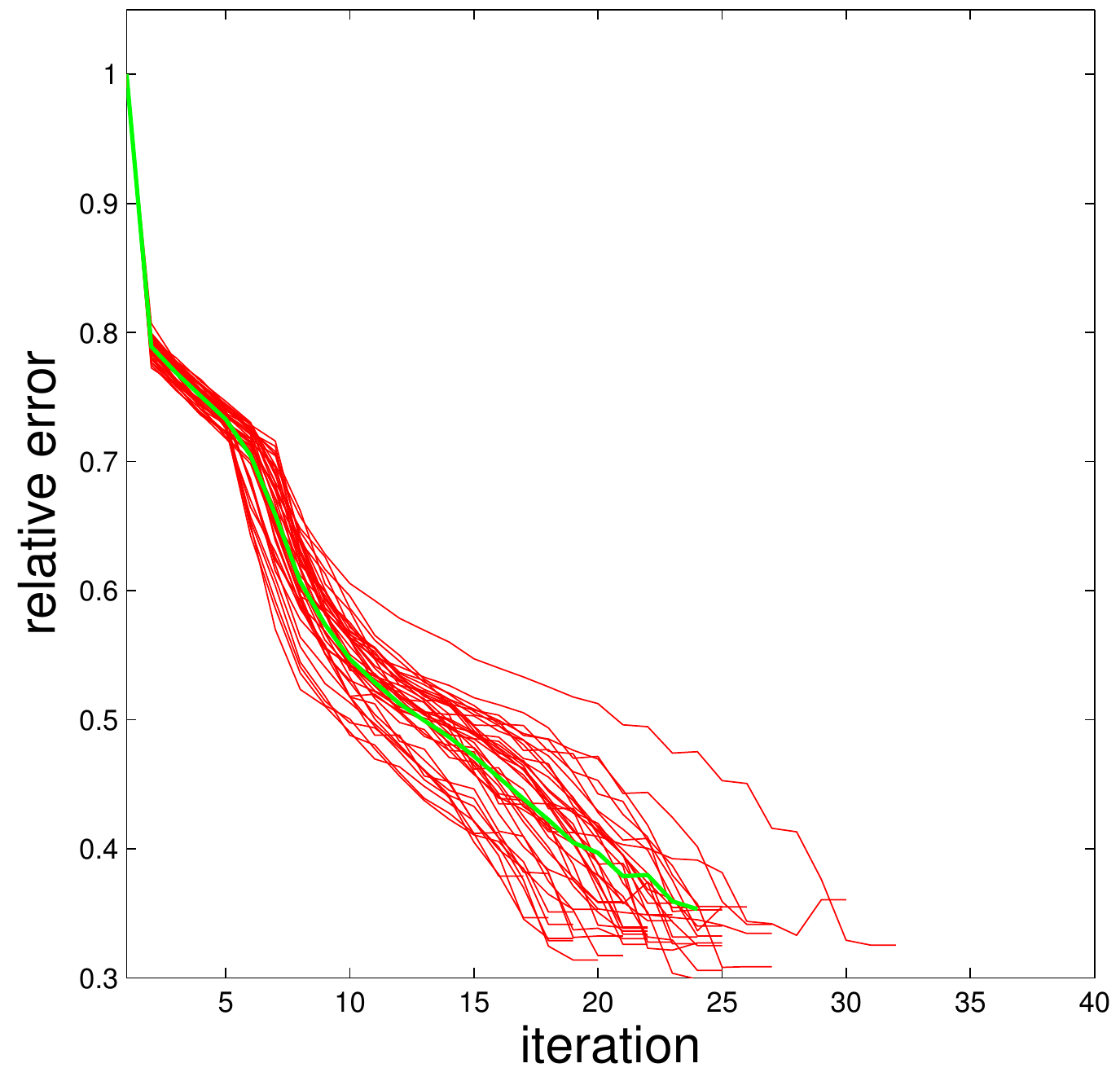}\\
\includegraphics[scale=0.2175]{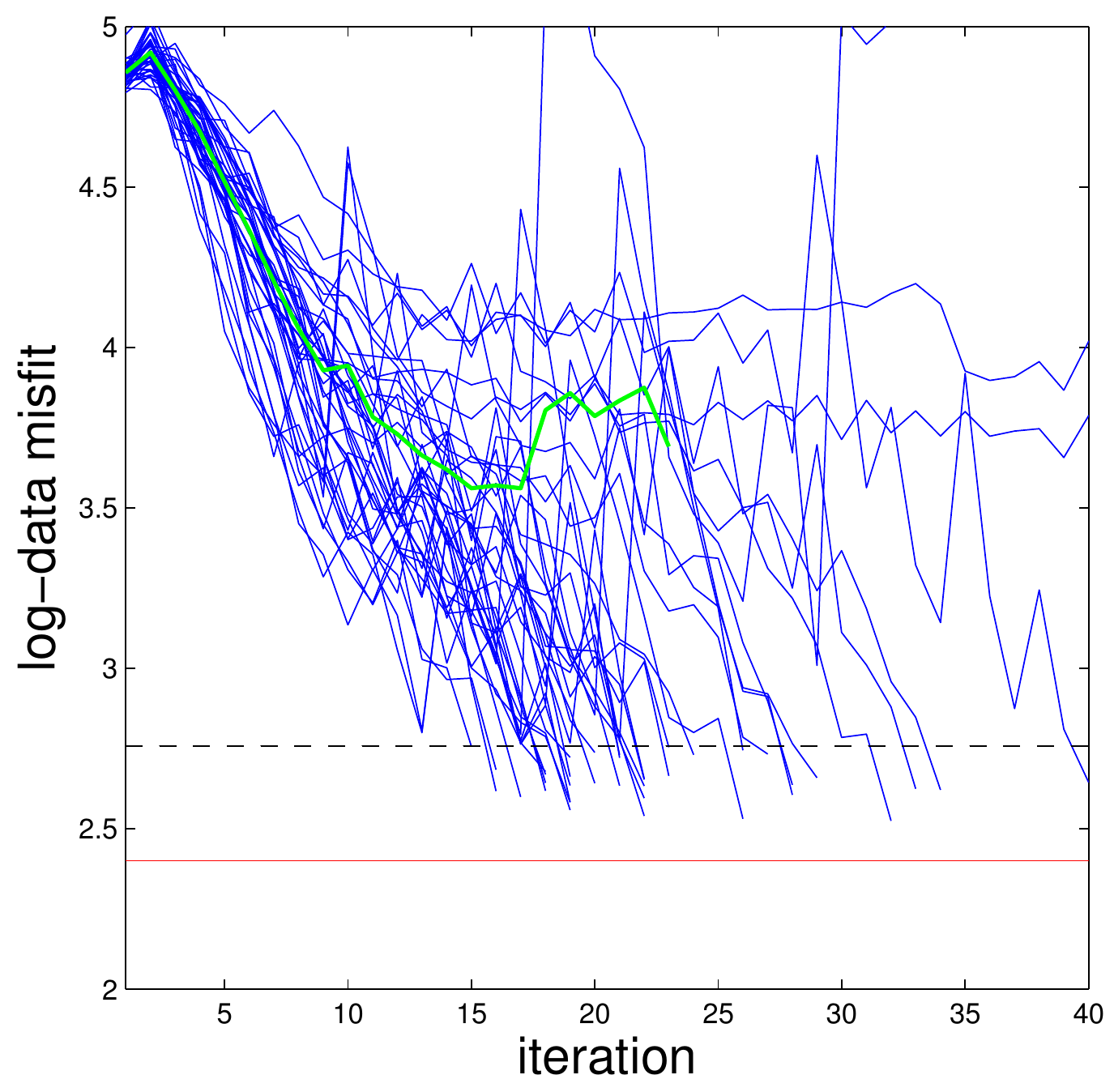}
\includegraphics[scale=0.2175]{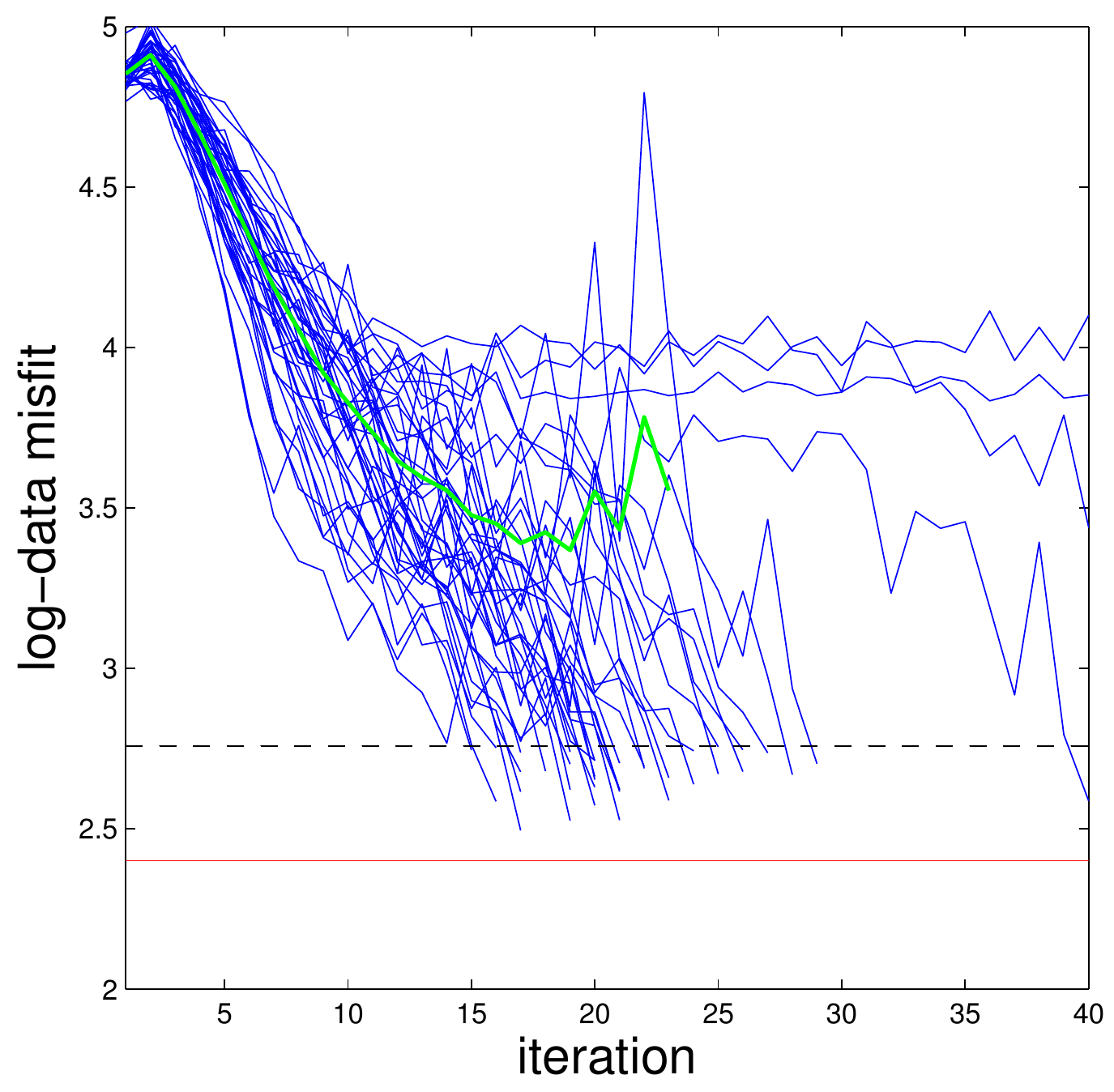}
\includegraphics[scale=0.2175]{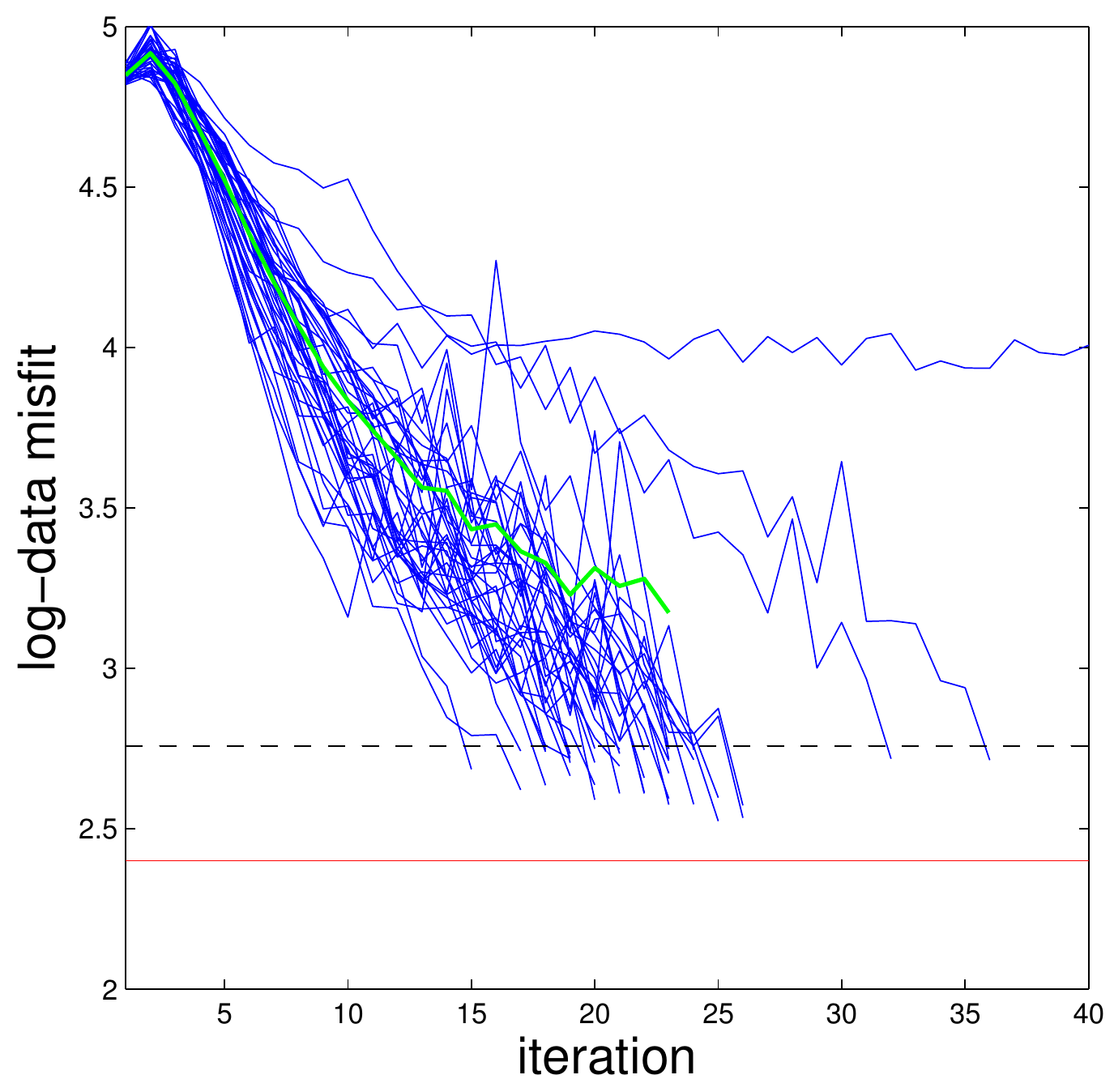}
\includegraphics[scale=0.2175]{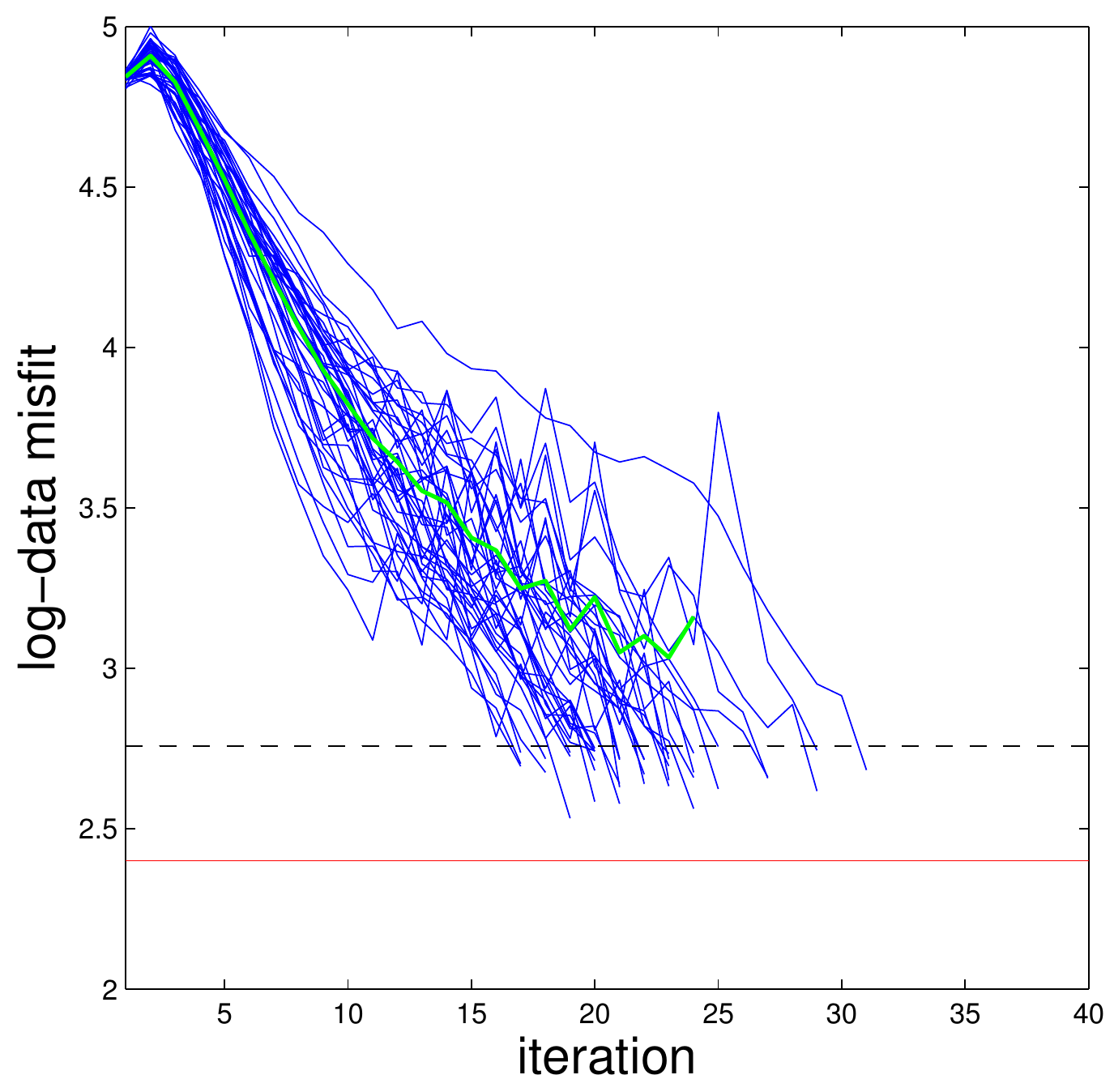}
\includegraphics[scale=0.2175]{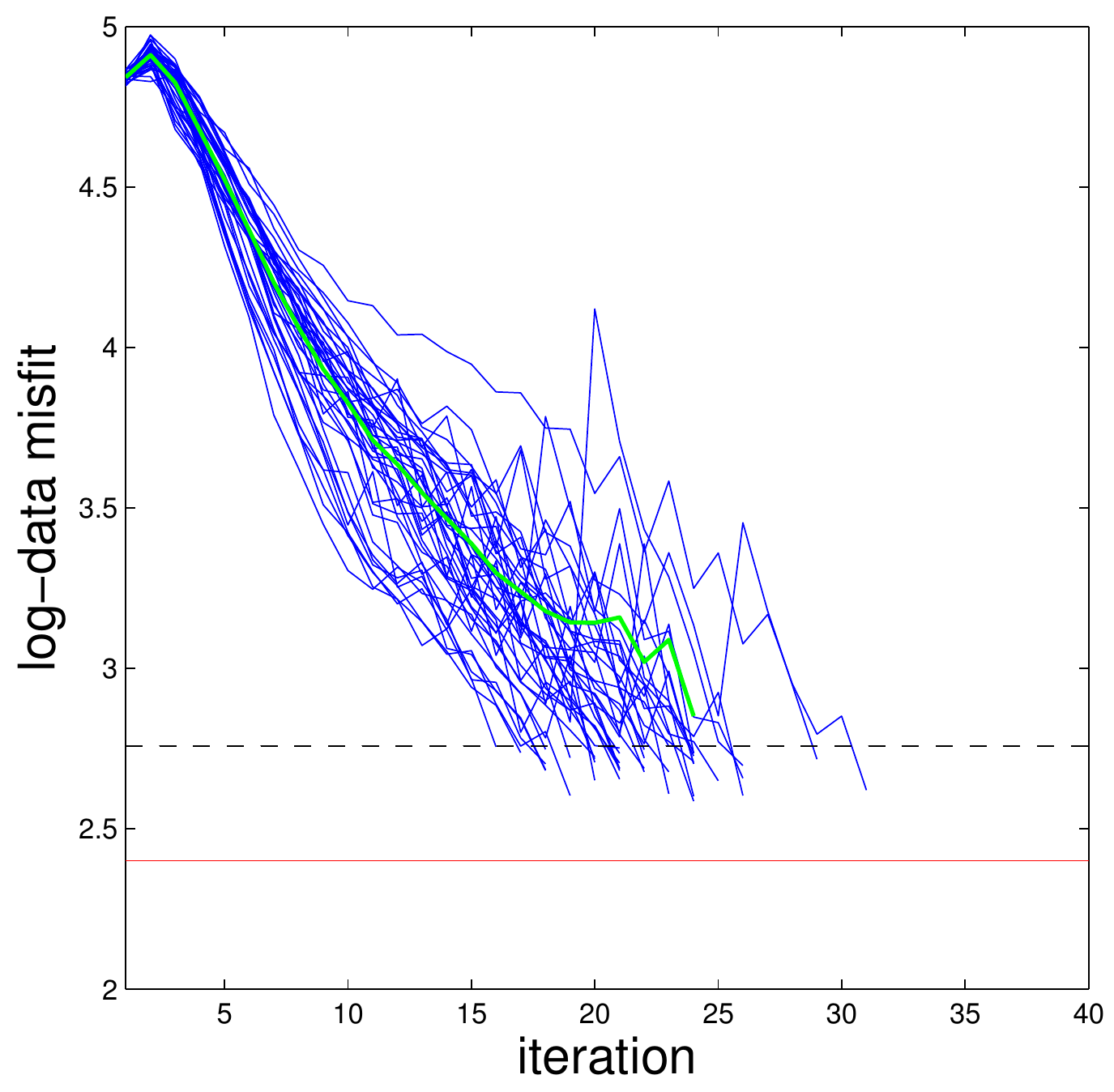} \caption{Error with respect to the truth (top) and log data misfit (bottom) from 40 experiments with different initial ensembles of size (from left to right) $N){e}=75, 100,150, 250, 350$. }    \label{Fig_gw3}

\end{center}
\end{figure}

\begin{figure}[htbp]
\includegraphics[scale=0.7]{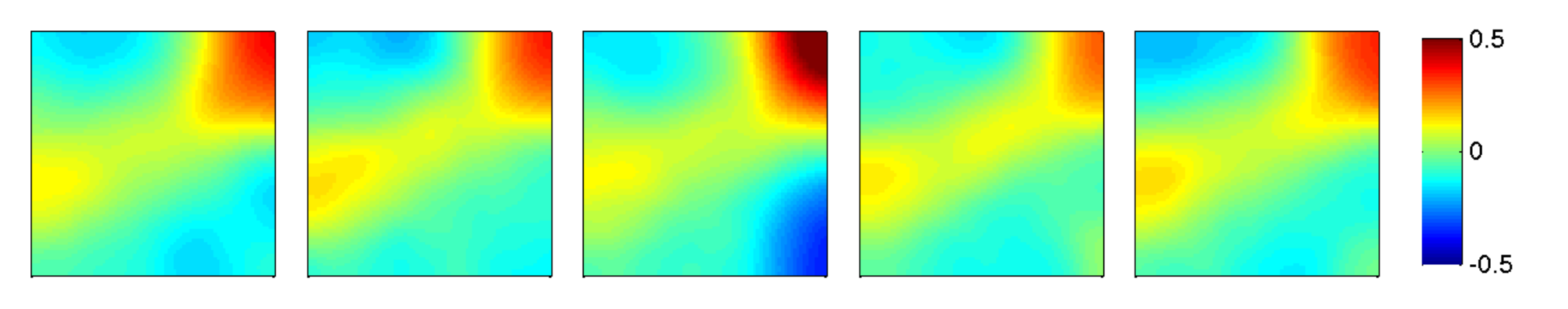}\\
\includegraphics[scale=0.7]{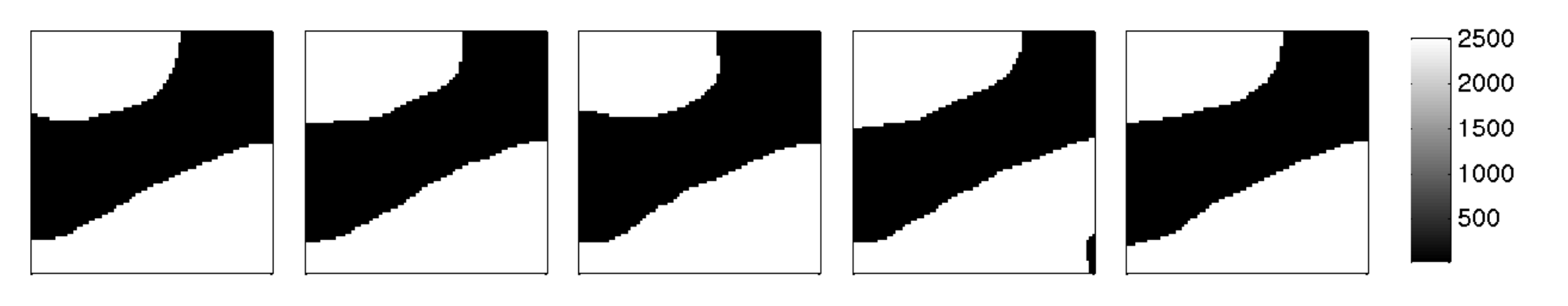}\\
\begin{center}
\includegraphics[scale=0.175]{True_LS_GW}
\caption{Top: Estimates of level-set obtained from 5 experiments with different initial ensembles with samples of $N(0,C)$ with $C$ from (\ref{eq:cova2}).  Middle: Conductivities computed from these estimated level-set functions (by means of (\ref{103})). Bottom: True hydraulic conductivity. }     \label{Fig_gw4}

\end{center}
\end{figure}

We apply the regularizing ensemble Kalman method proposed with different ensemble sizes $(N_{e}=75, 100,150,250,350)$ and with a fixed parameter $\rho=0.7$. As before, we consider 40 experiments with different initial ensembles and the corresponding log-data misfit and error with respect to the truth are displayed in Figure \ref{Fig_gw3}. Although we are applying the ensemble algorithm to the estimation of the level-set function $u$, we consider the relative error with respect to the truth of the corresponding conductivity given by (\ref{103}). Note that a large ensemble is needed in order to obtain stable computations in the scheme until the data misfit reaches the value $\eta/\rho$ that we use to stop the algorithm. This value is indicated with the dotted horizontal line in the bottom of Figure \ref{Fig_gw3}. For sufficiently large ensemble size ($N_{e}>150$), on average both data misfit and error decrease. Stable computations are obtained when the algorithm is stopped via (\ref{eq:m15}) with $\tau \approx 1/\rho$. In the top and middle rows of \Fref{Fig_gw4} we display estimates of the level-set function and the corresponding conductivity obtained with Algorithm \ref{Al1} (for $N_{e}=150$) from five experiments with different initial ensemble sizes (with initial ensemble generated from the same Gaussian measure). The visual agreement between the estimates of conductivity and truth (\Fref{Fig_gw4} bottom) is remarkable.

\subsection{EIT. Sharp interfaces.}

\begin{figure}[htbp]
\begin{center}

\includegraphics[scale=0.27]{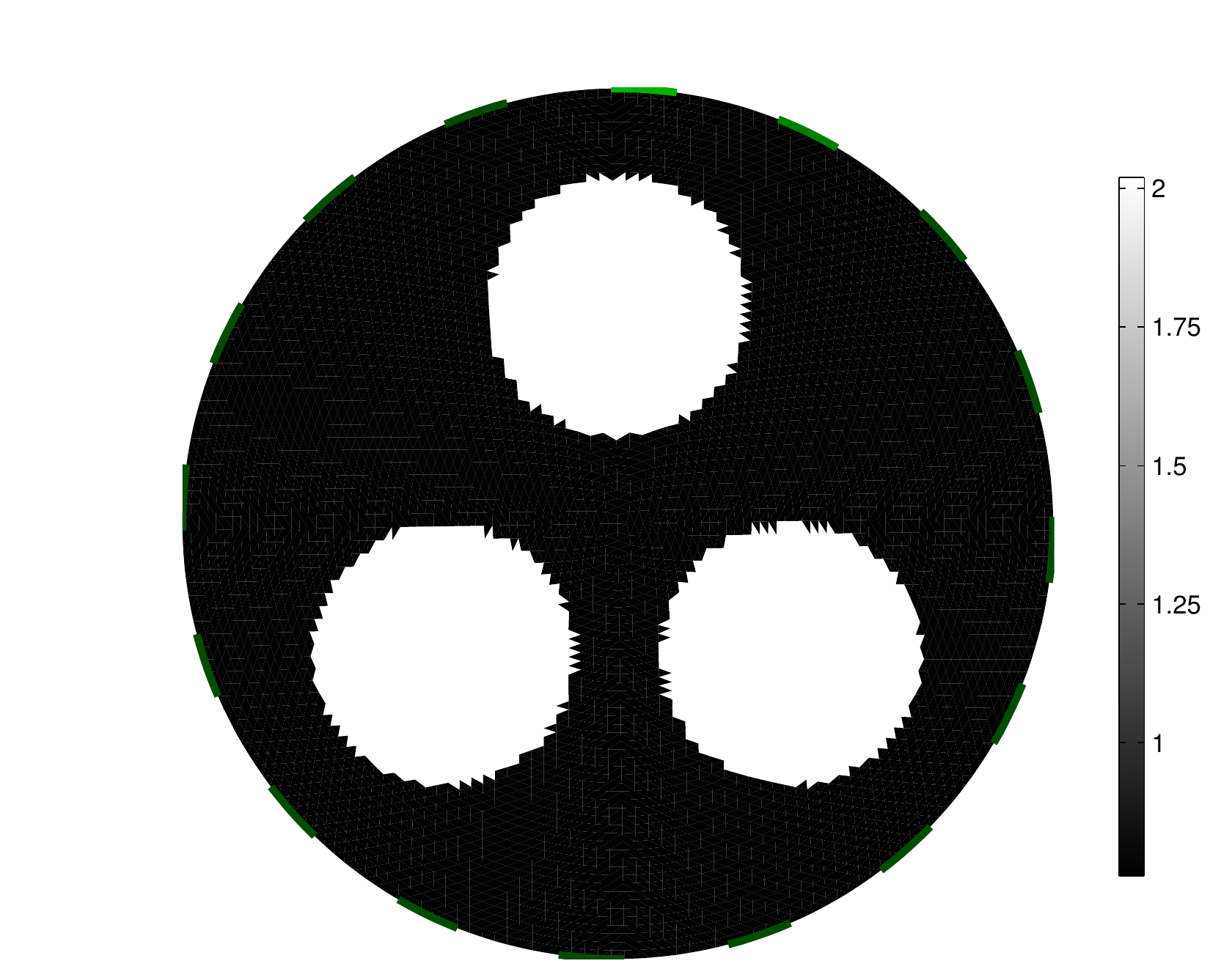}
\includegraphics[scale=0.27]{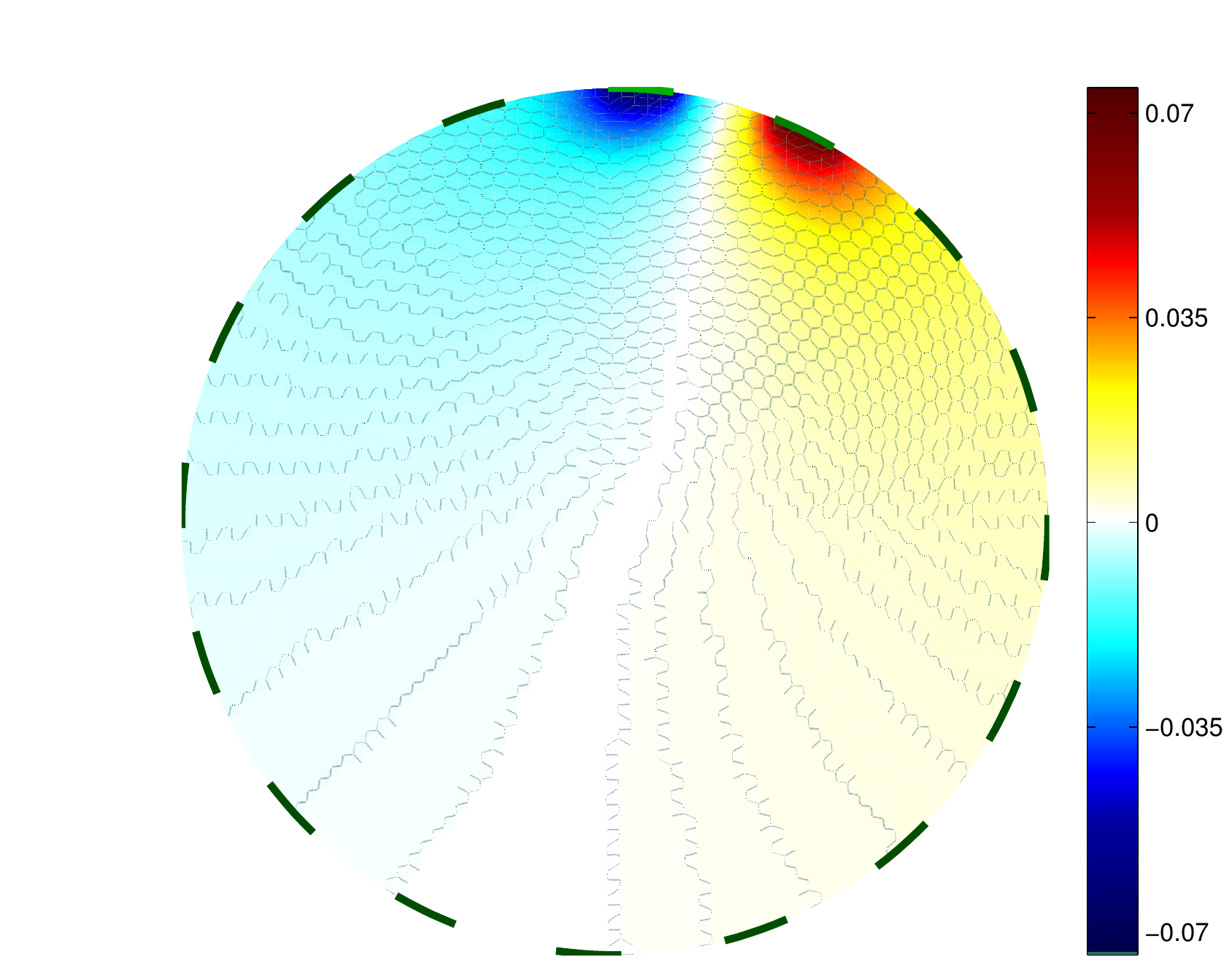}
\includegraphics[scale=0.27]{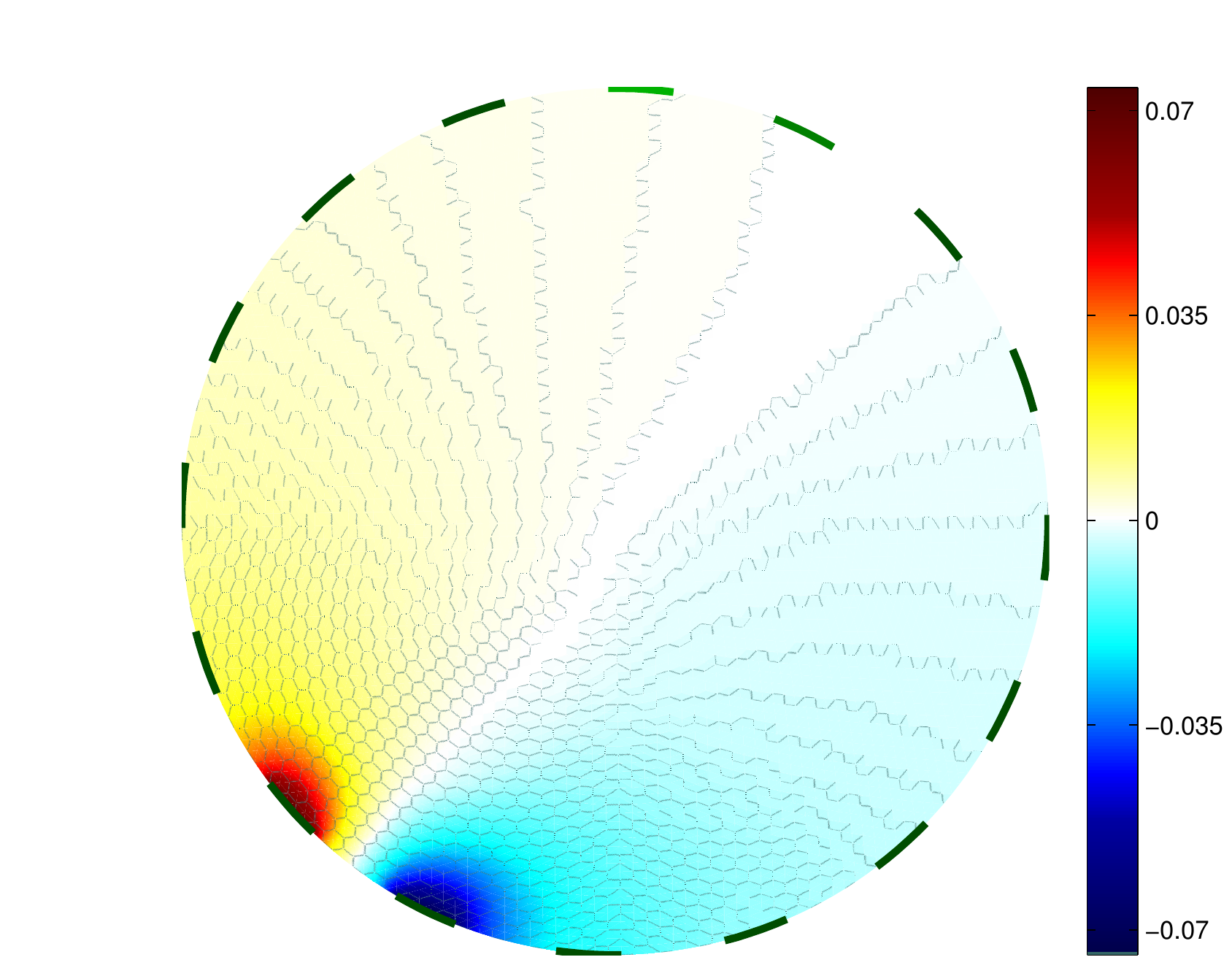}
 \caption{Top-left: true conductivity. Top-middle and top-right: two of the 15 voltage patterns generated with the true conductivity. } \label{Fig_L_EIT0}
\end{center}
\end{figure}

In this subsection we consider again the EIT problem but assume that we are now interested in recovering a conductivity with sharp discontinuities. Let us then consider the conductivity and the electrode configuration displayed in \Fref{Fig_L_EIT0} (left). We use this field as the true conductivity that we aim at estimating with the ensemble Kalman method by using the level-set parametrization of expression (\ref{103}). Synthetic data are generated (and avoiding inverse crimes) as described in subsection \ref{num_EIT}. Some true voltages are displayed in \Fref{Fig_L_EIT0} (middle and right panels).  Level-set approaches for shape identification in EIT problems have been studied, for example, in \cite{LS_EIT1,LS_EIT2}. However, these approaches are variational and use a level-set equation for the evolution of the shape. 

\begin{figure}[htbp]
\begin{center}
\includegraphics[scale=0.175]{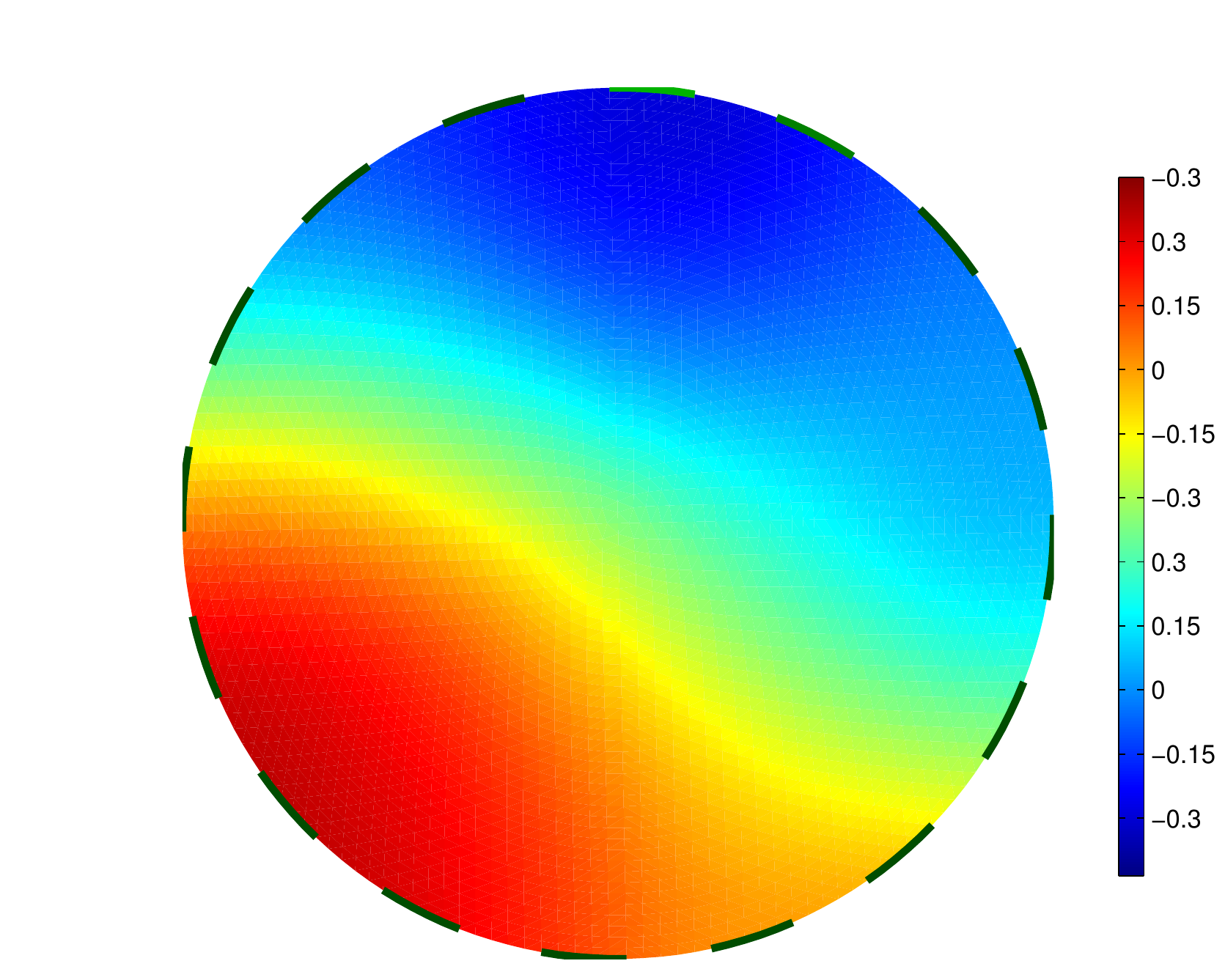}
\includegraphics[scale=0.175]{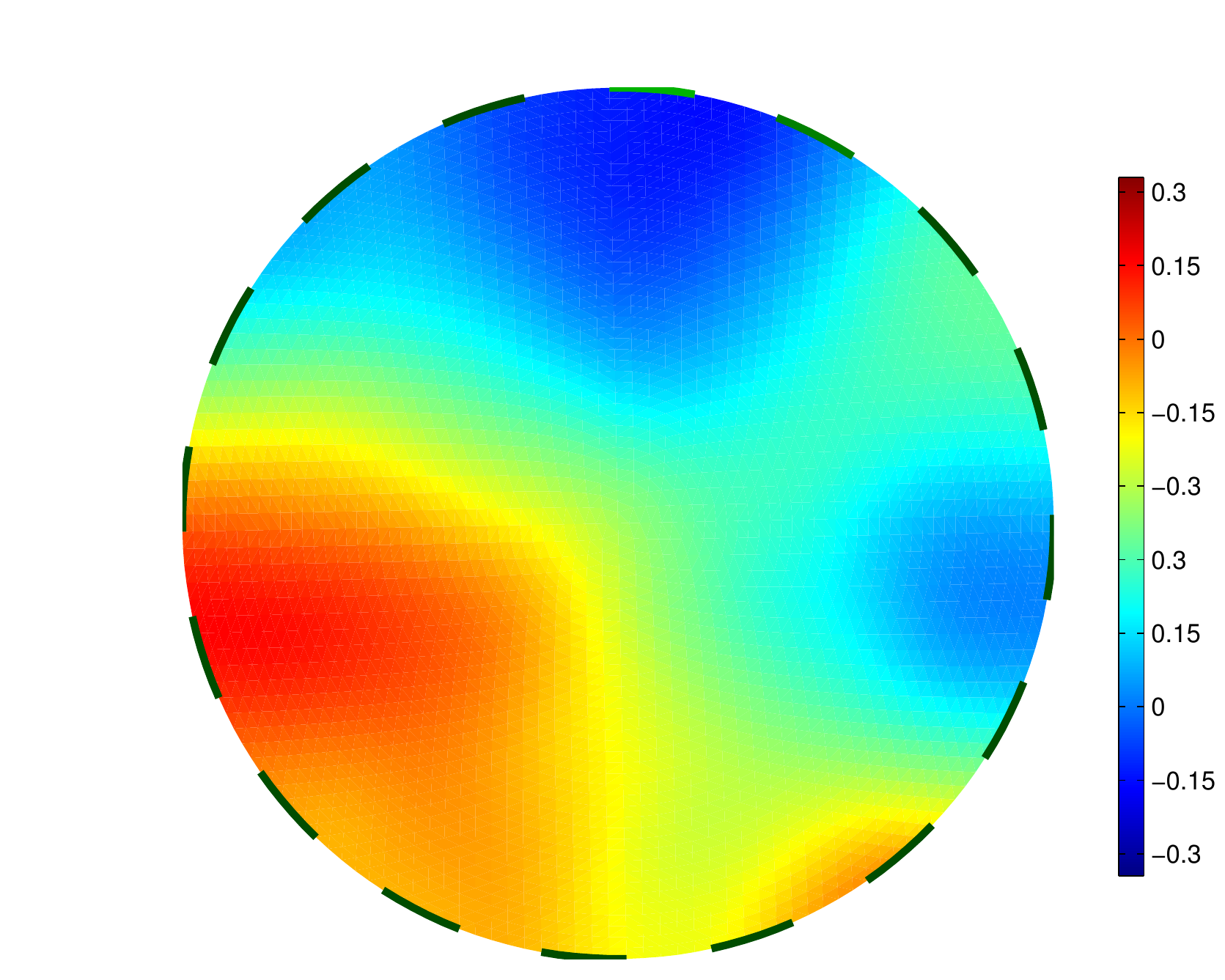}
\includegraphics[scale=0.175]{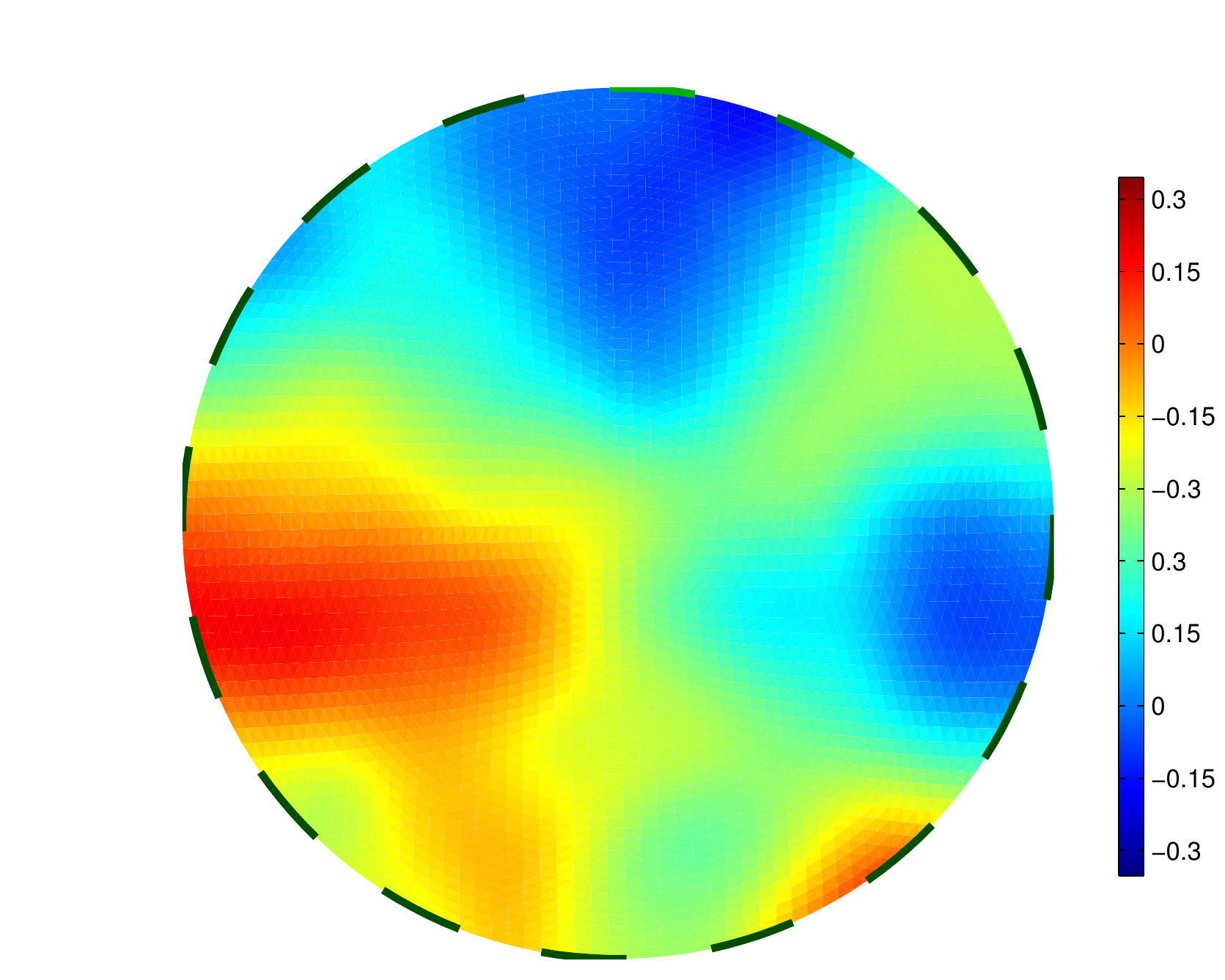}
\includegraphics[scale=0.175]{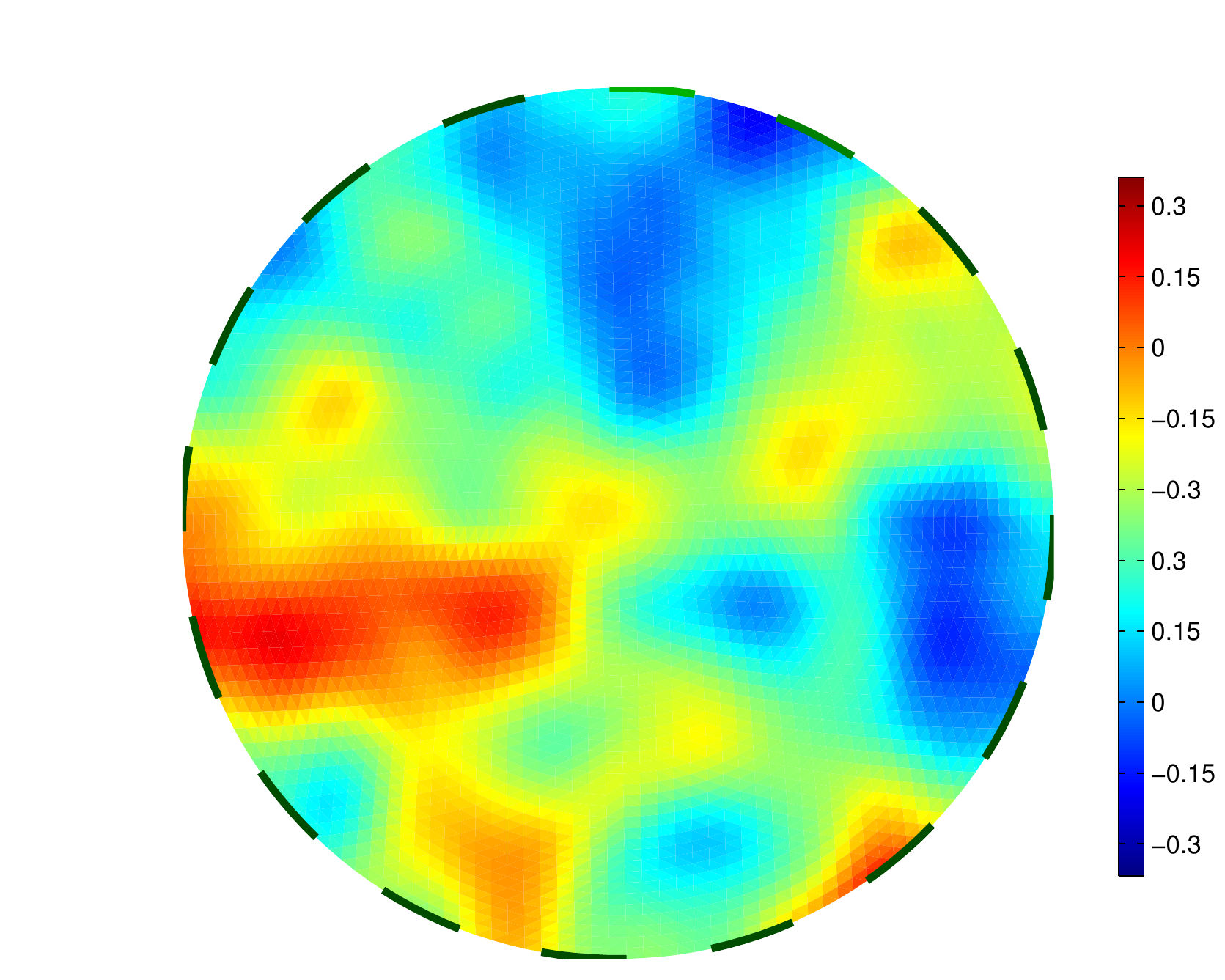}
\includegraphics[scale=0.175]{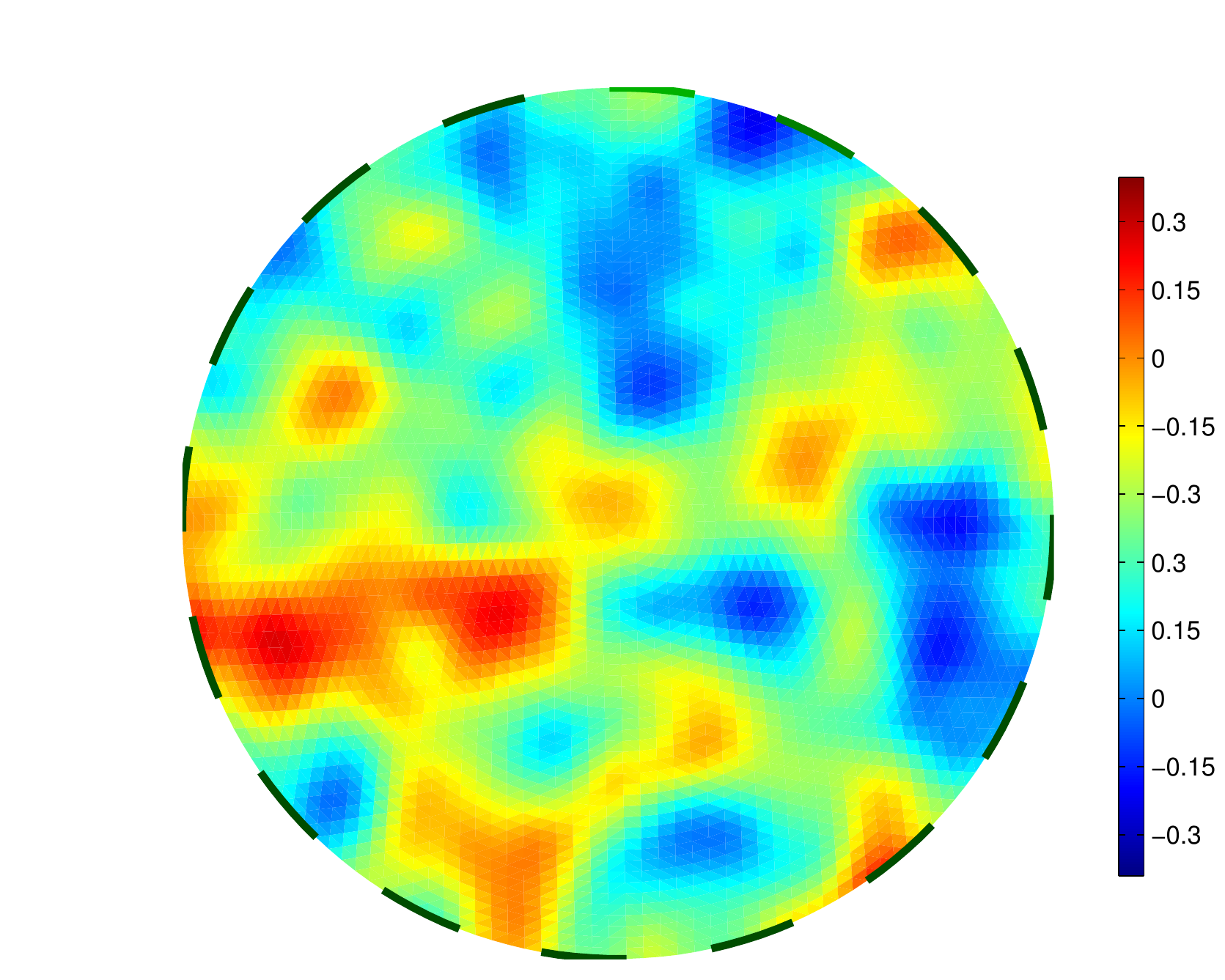}\\
\includegraphics[scale=0.18]{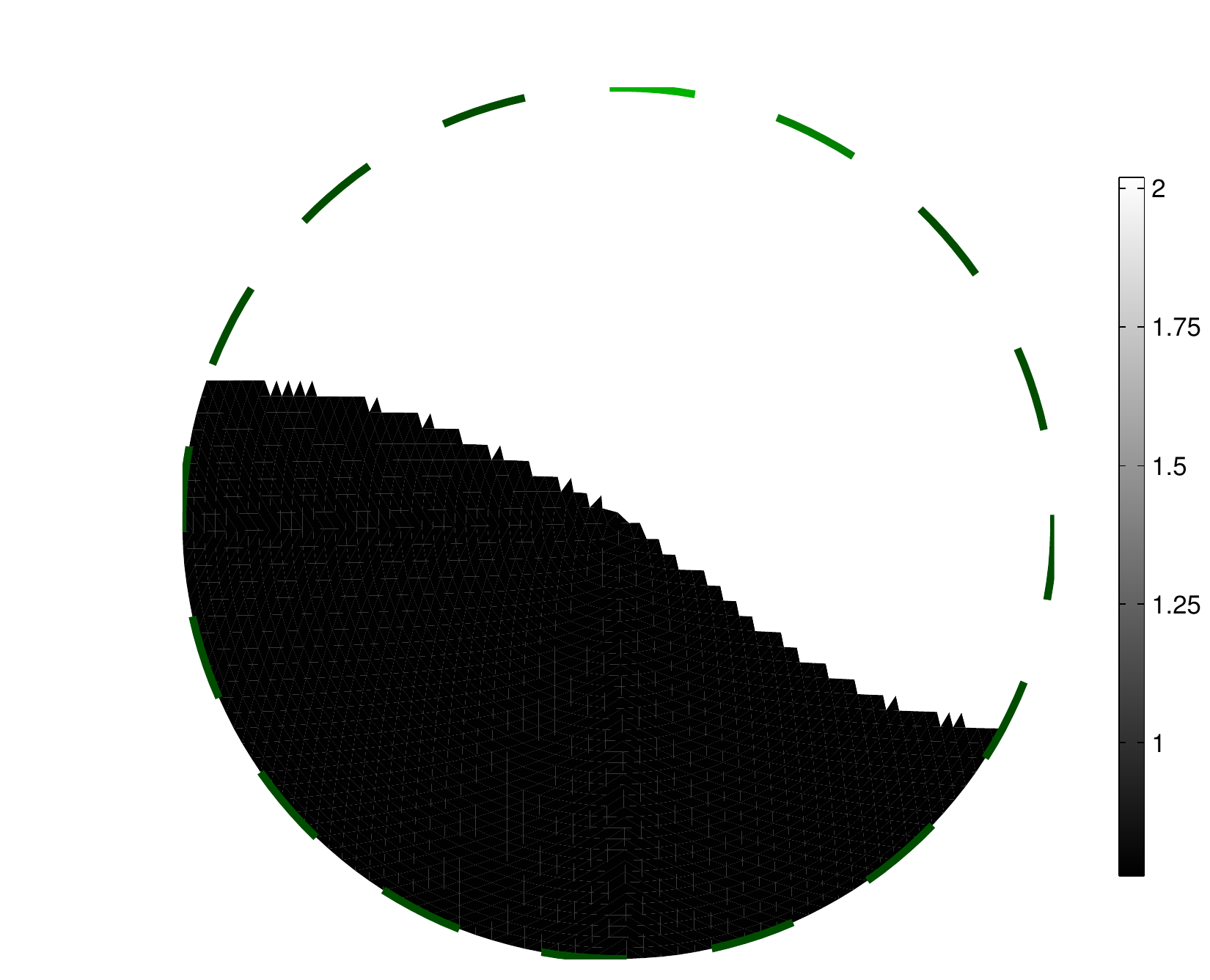}
\includegraphics[scale=0.18]{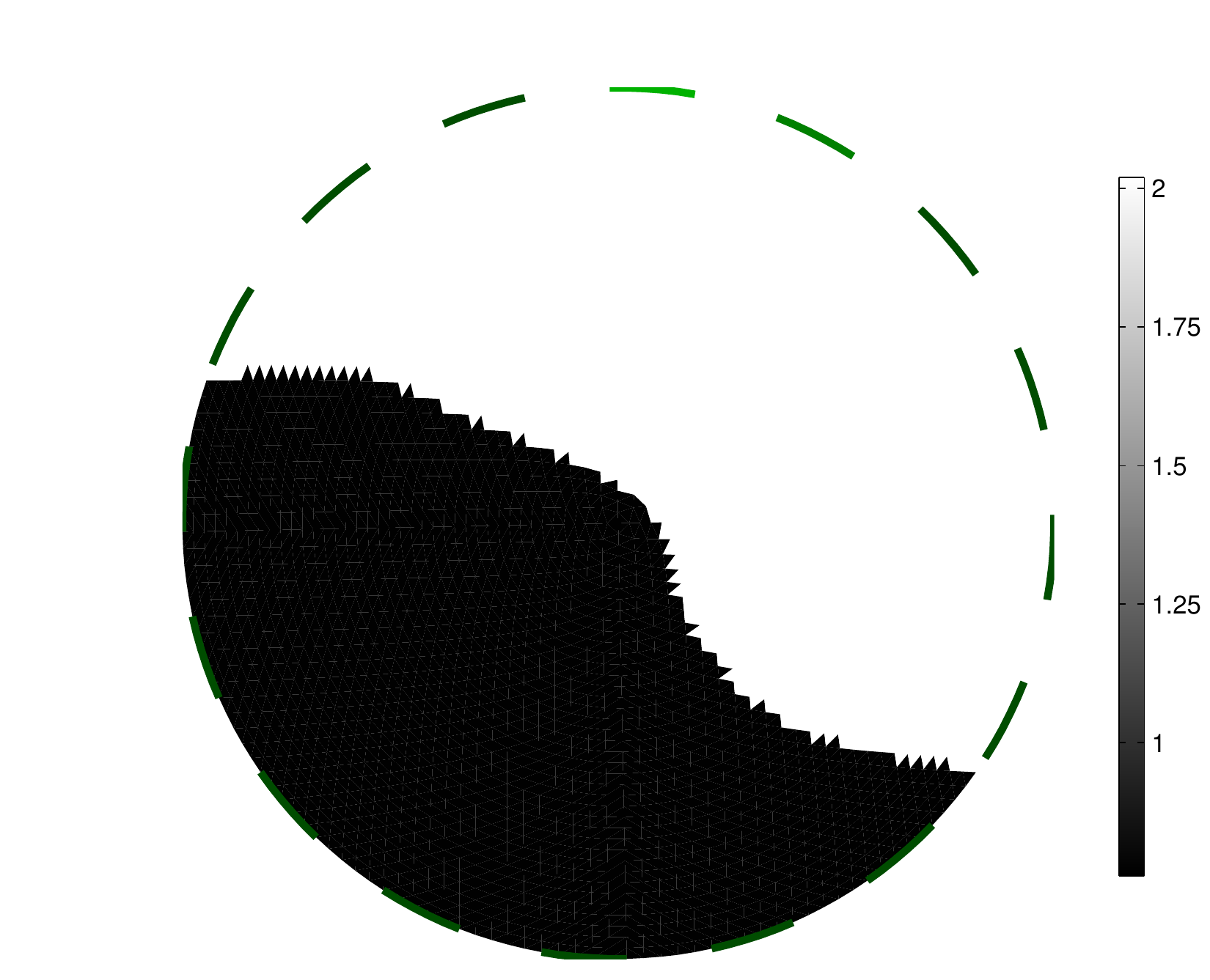}
\includegraphics[scale=0.18]{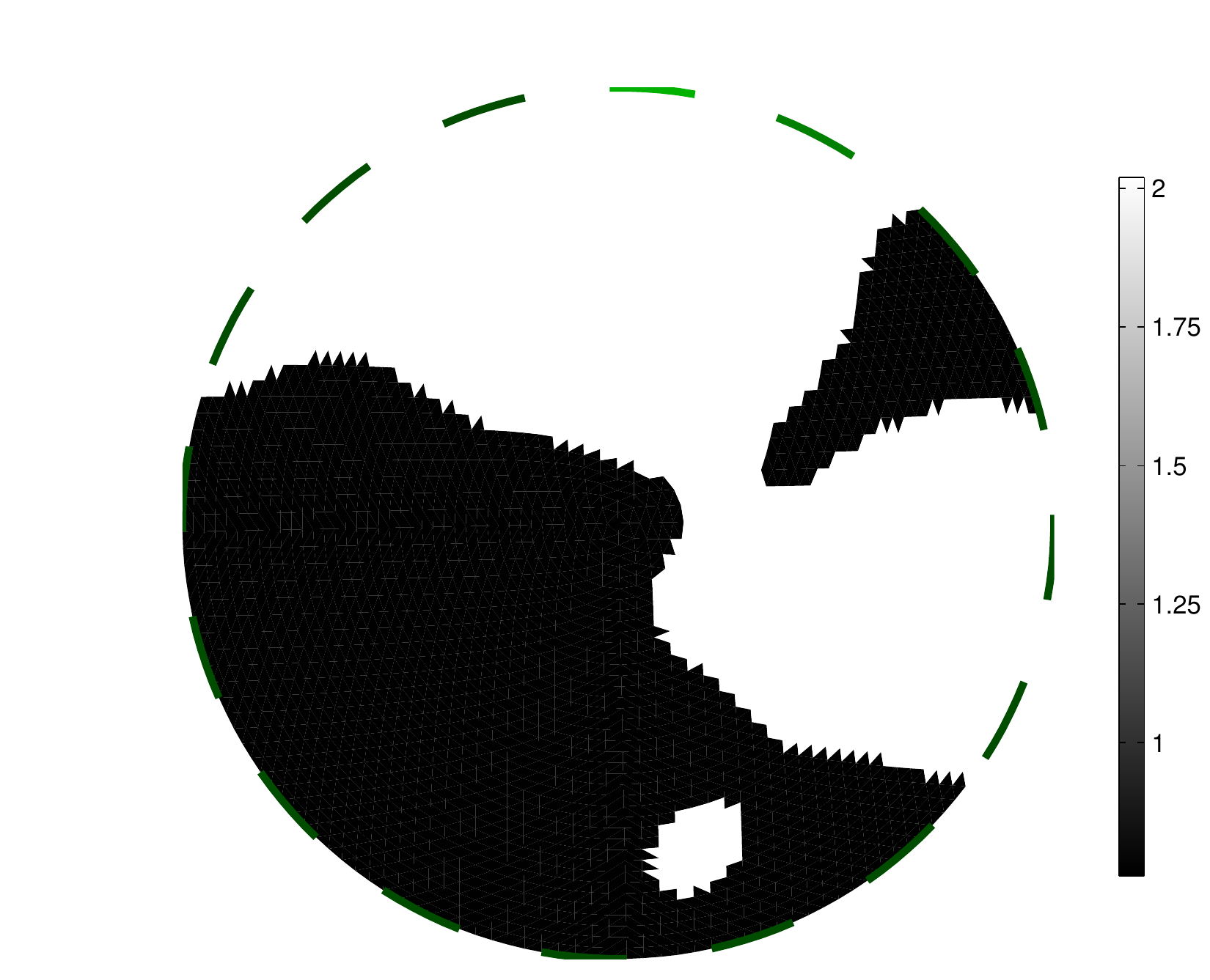}
\includegraphics[scale=0.18]{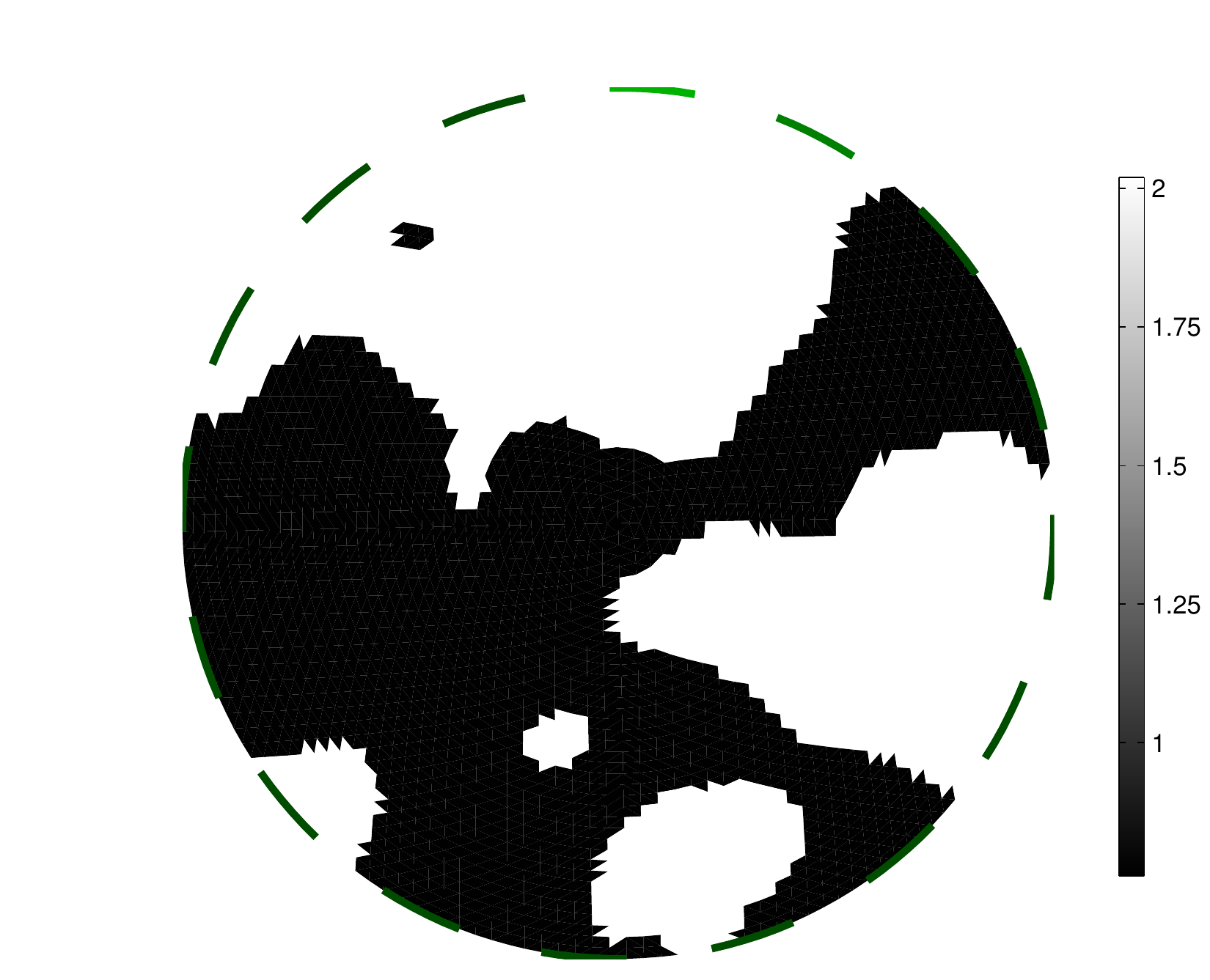}
\includegraphics[scale=0.18]{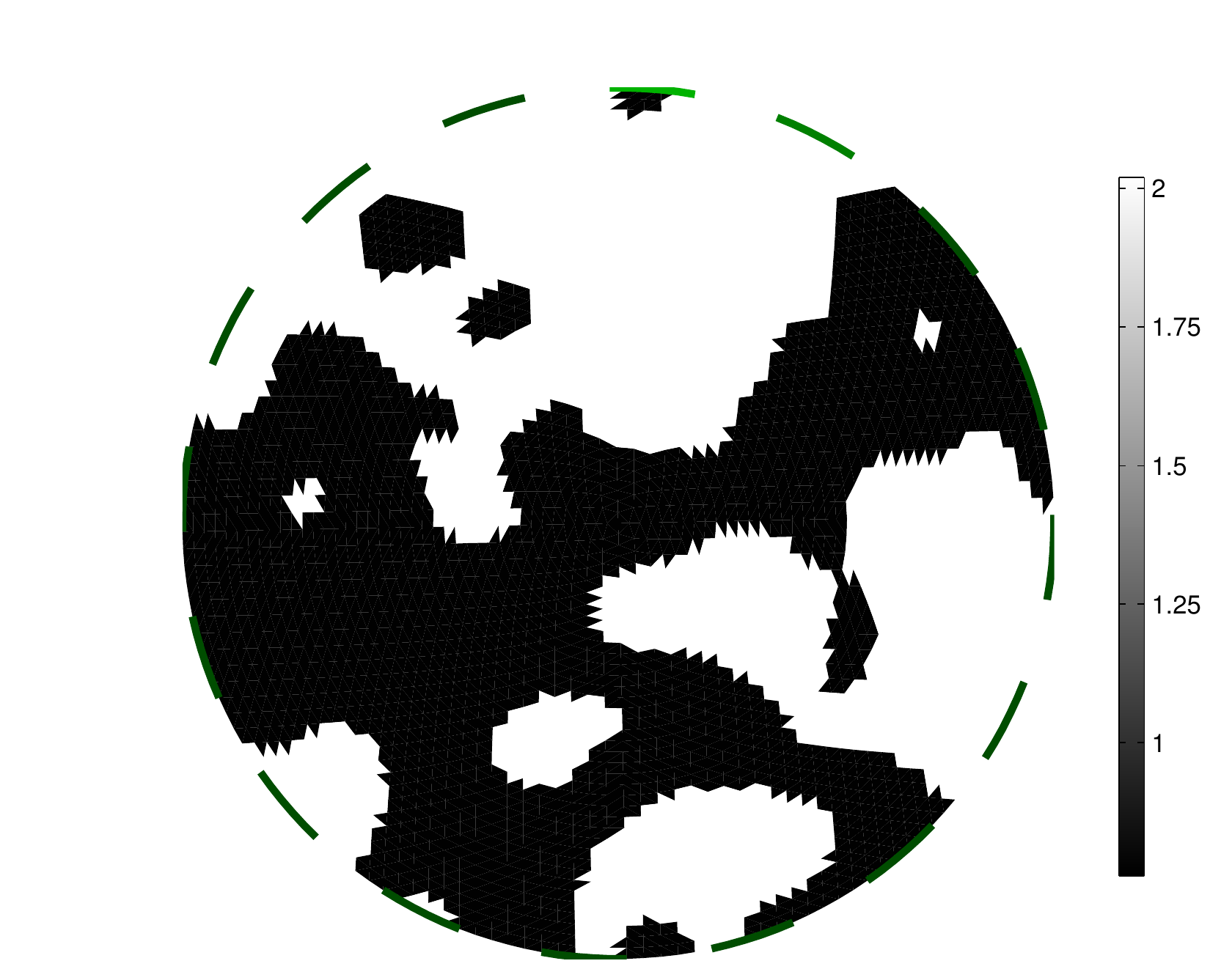}
 \caption{Top: Level-set functions sampled of $N(0,C)$ with values of $L$ in (\ref{eq:cova1}) (from left to right) $L=0.2, 0.1, 0.06,0.04, 0.03$ (these samples have the same realization of KL coefficients). Bottom: Conductivities computed from these sample level-set functions (by means of (\ref{103})).} \label{Fig_L_EIT1}
\end{center}
\end{figure}

\begin{figure}[htbp]
\begin{center}
\includegraphics[scale=0.175]{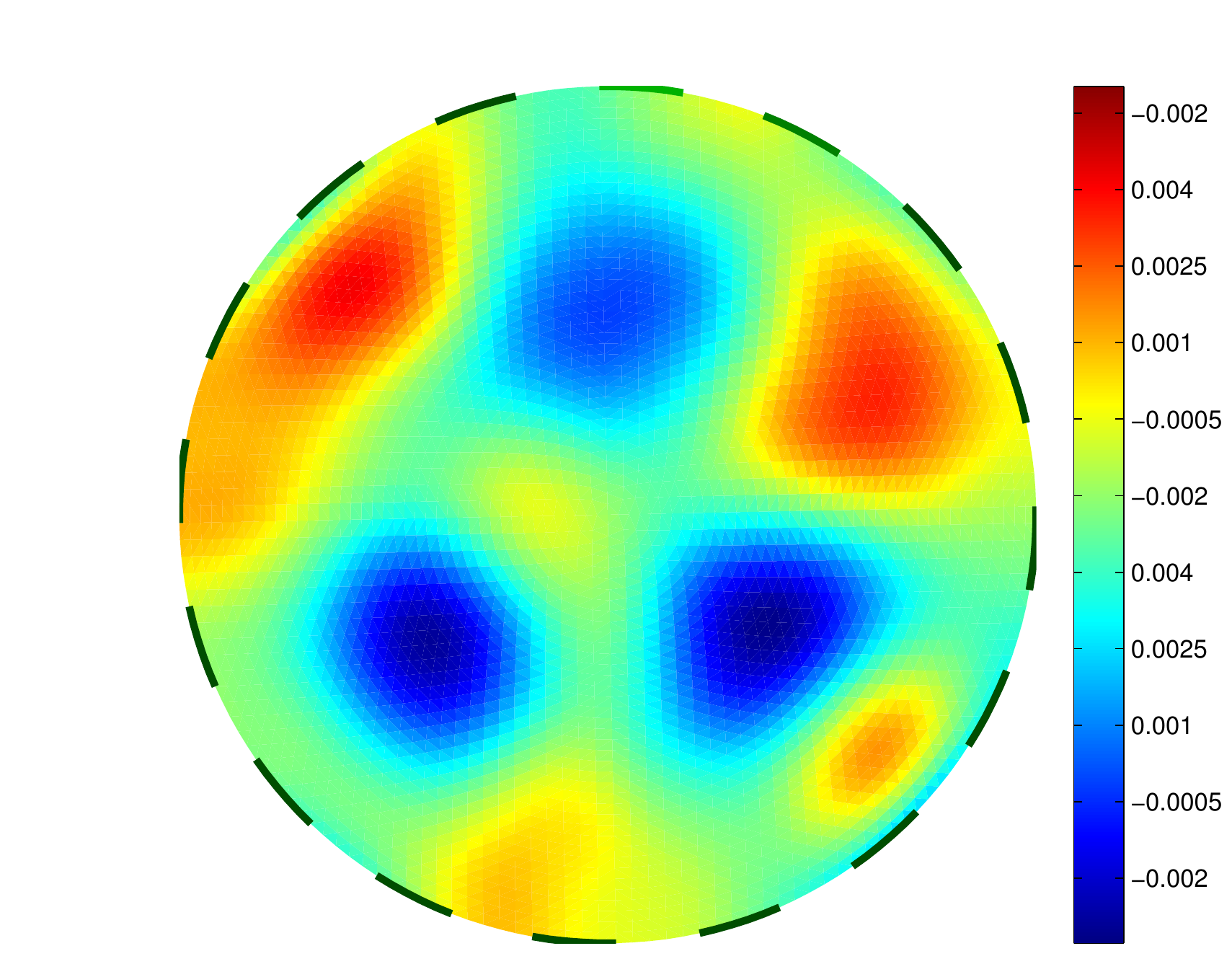}
\includegraphics[scale=0.175]{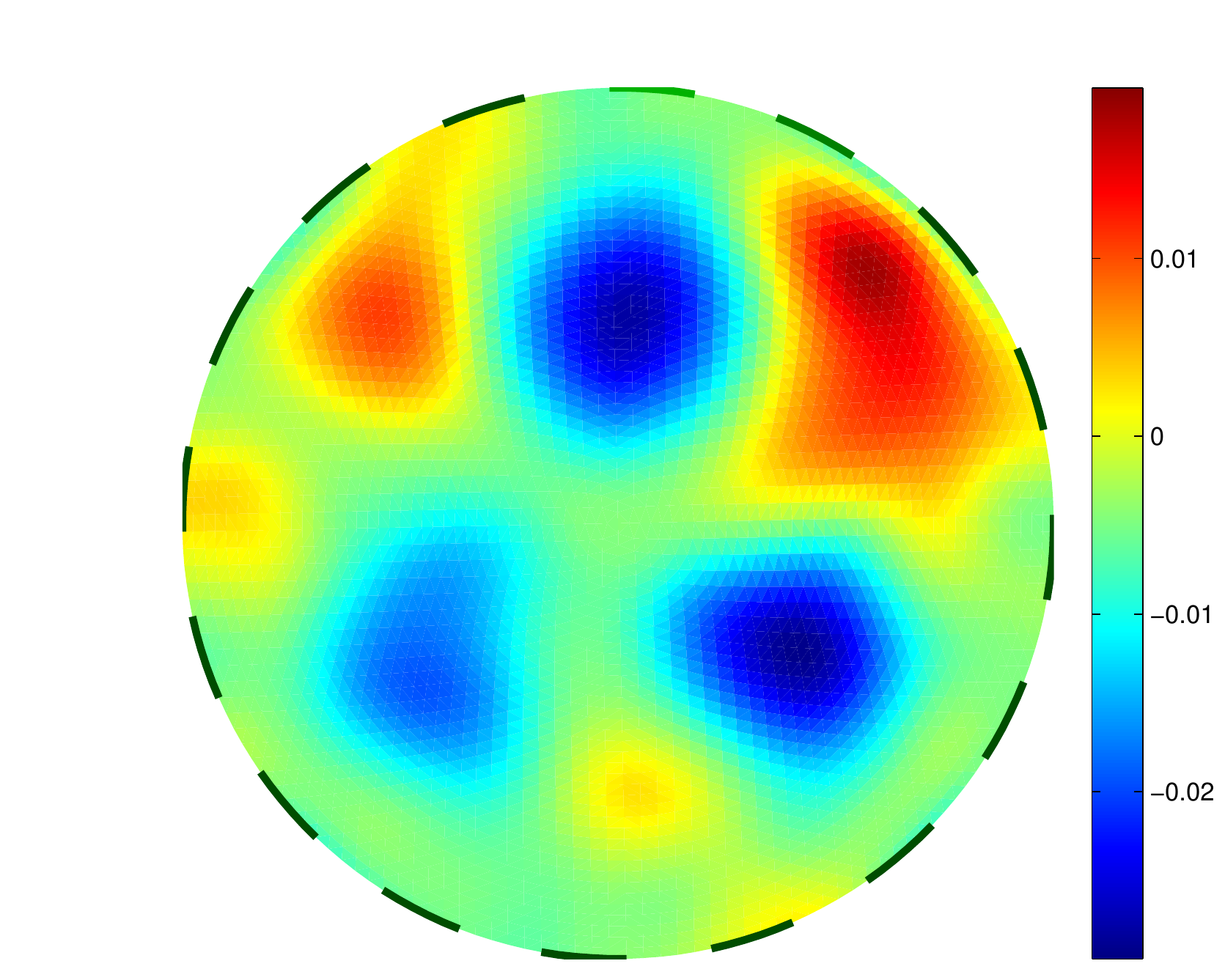}
\includegraphics[scale=0.175]{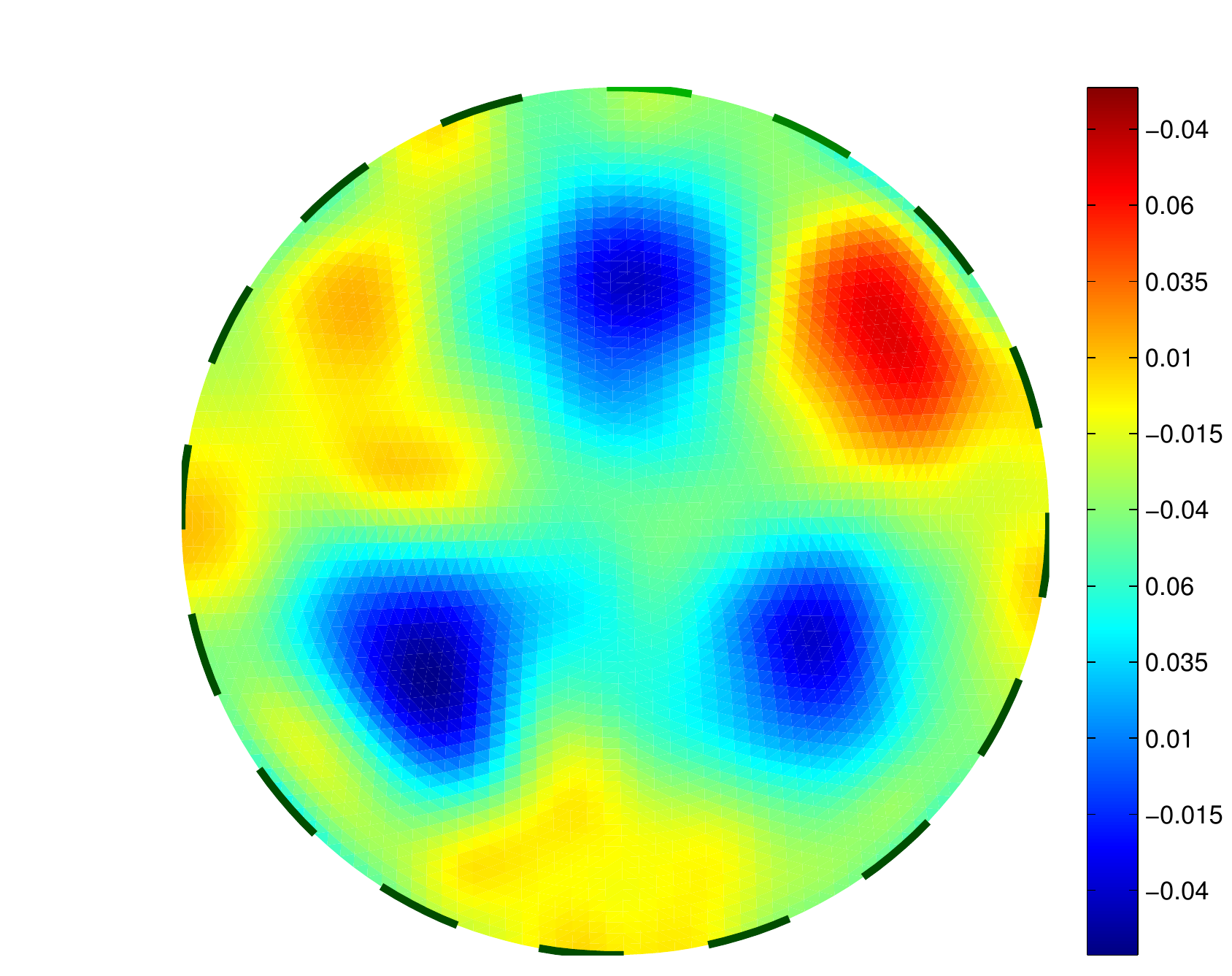}
\includegraphics[scale=0.175]{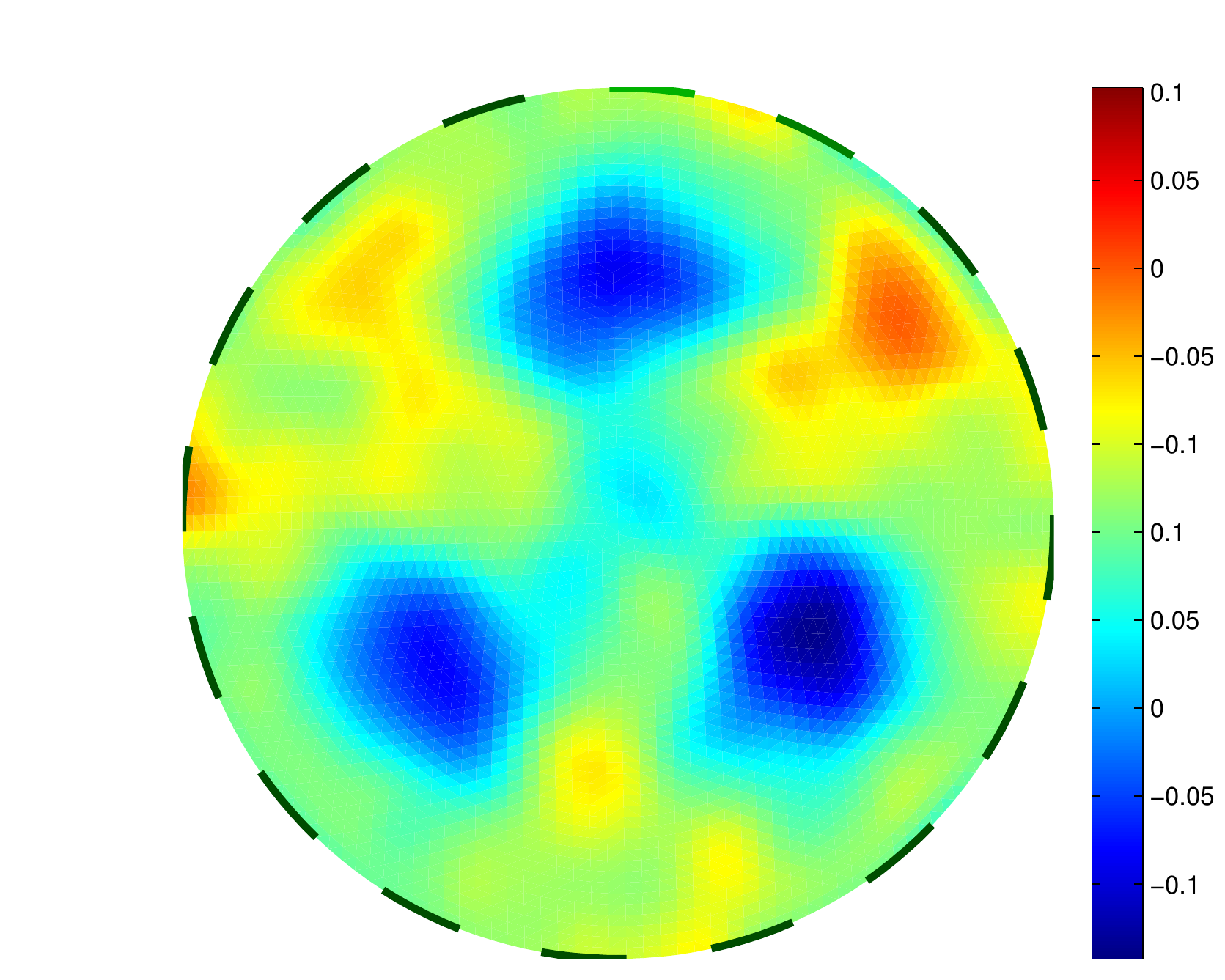}
\includegraphics[scale=0.175]{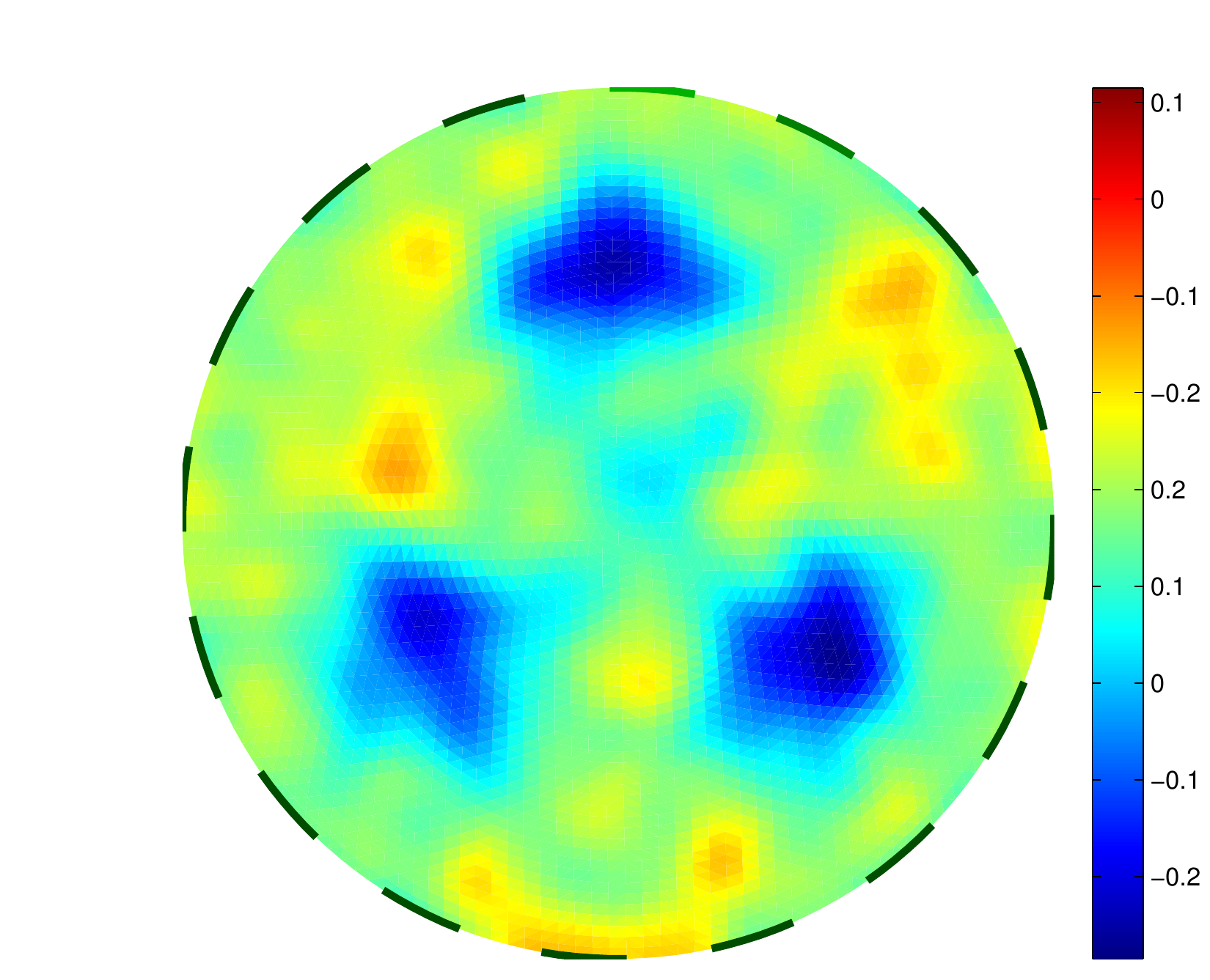}\\
\includegraphics[scale=0.18]{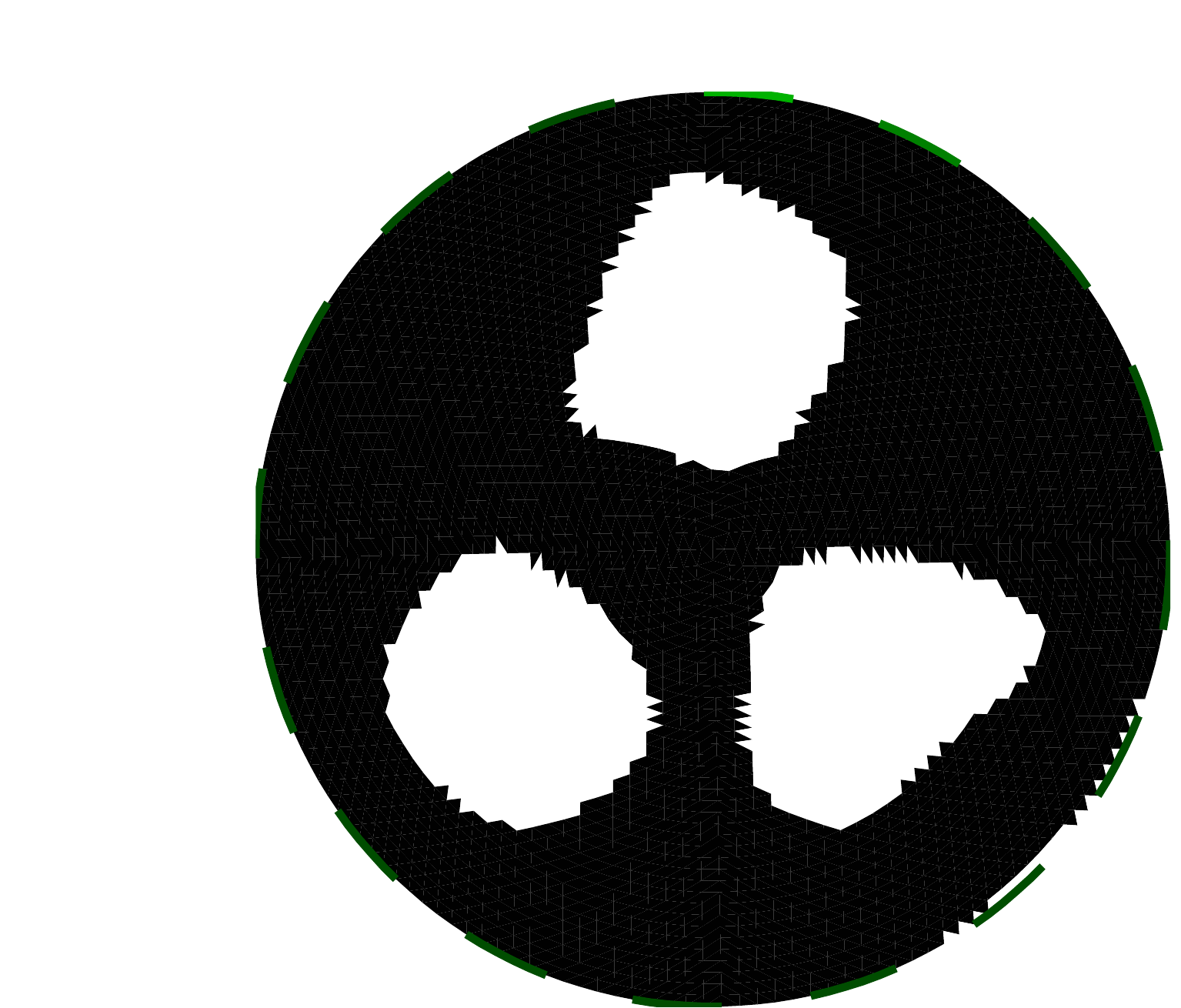}~
\includegraphics[scale=0.18]{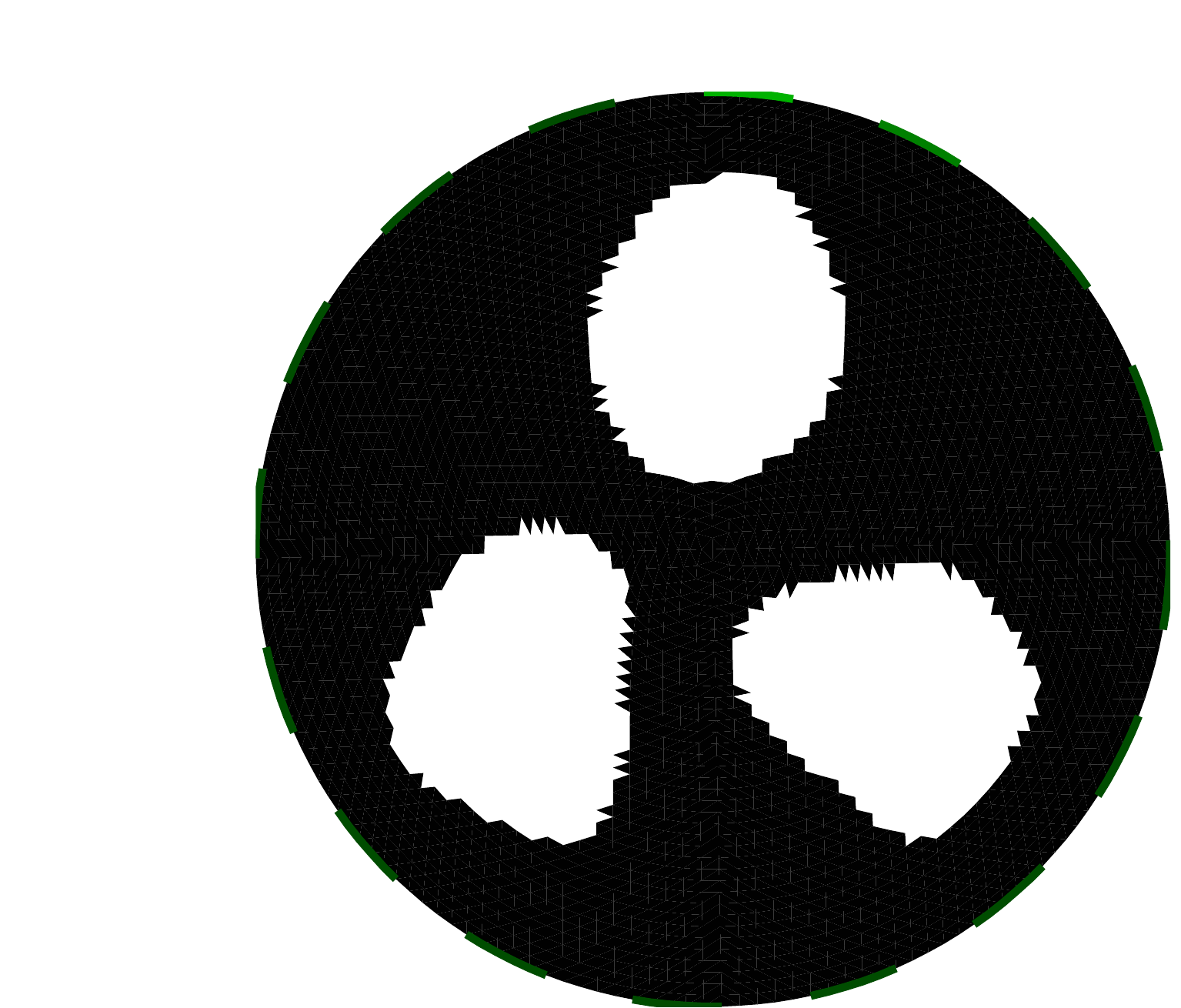}~
\includegraphics[scale=0.18]{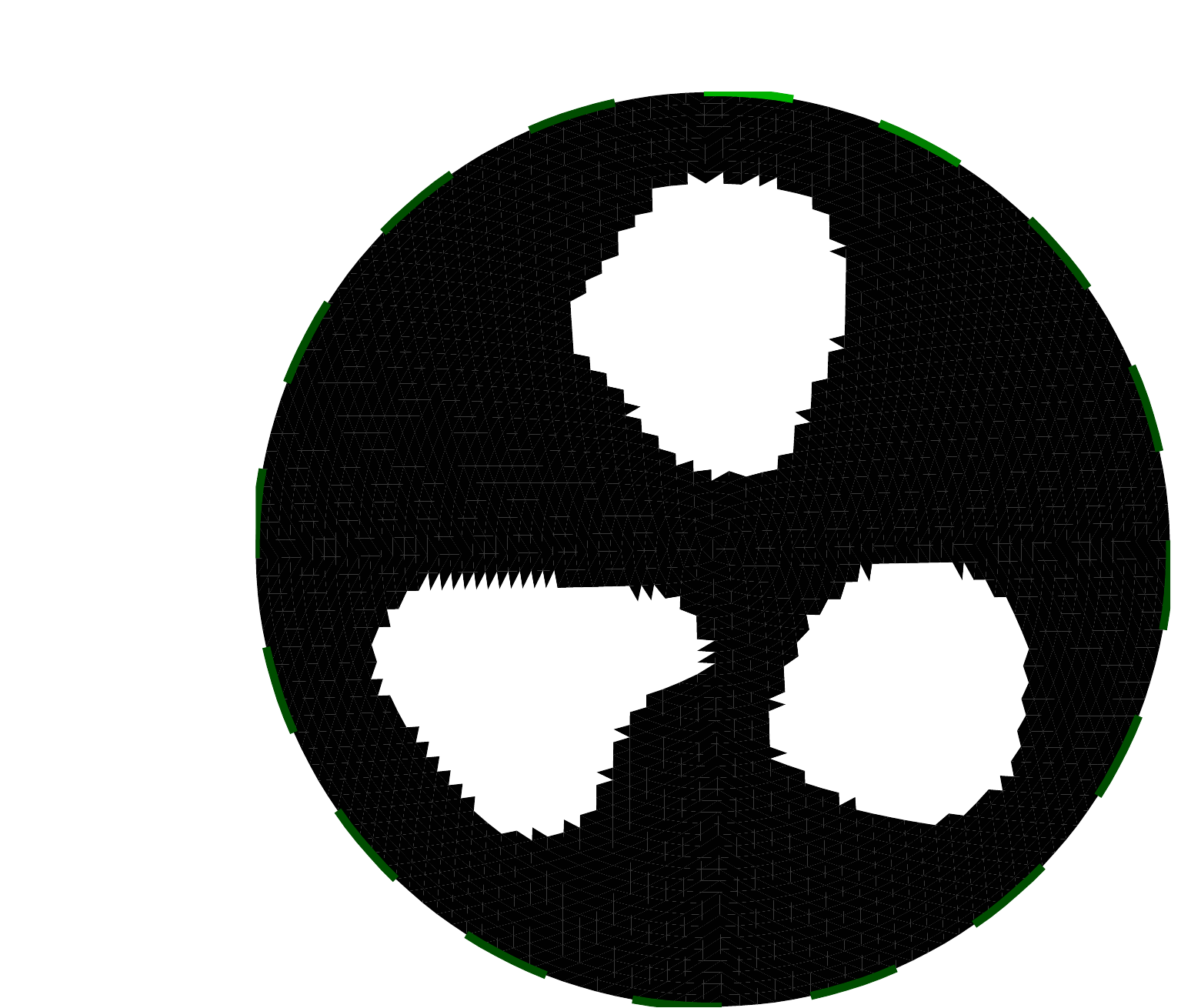}~
\includegraphics[scale=0.18]{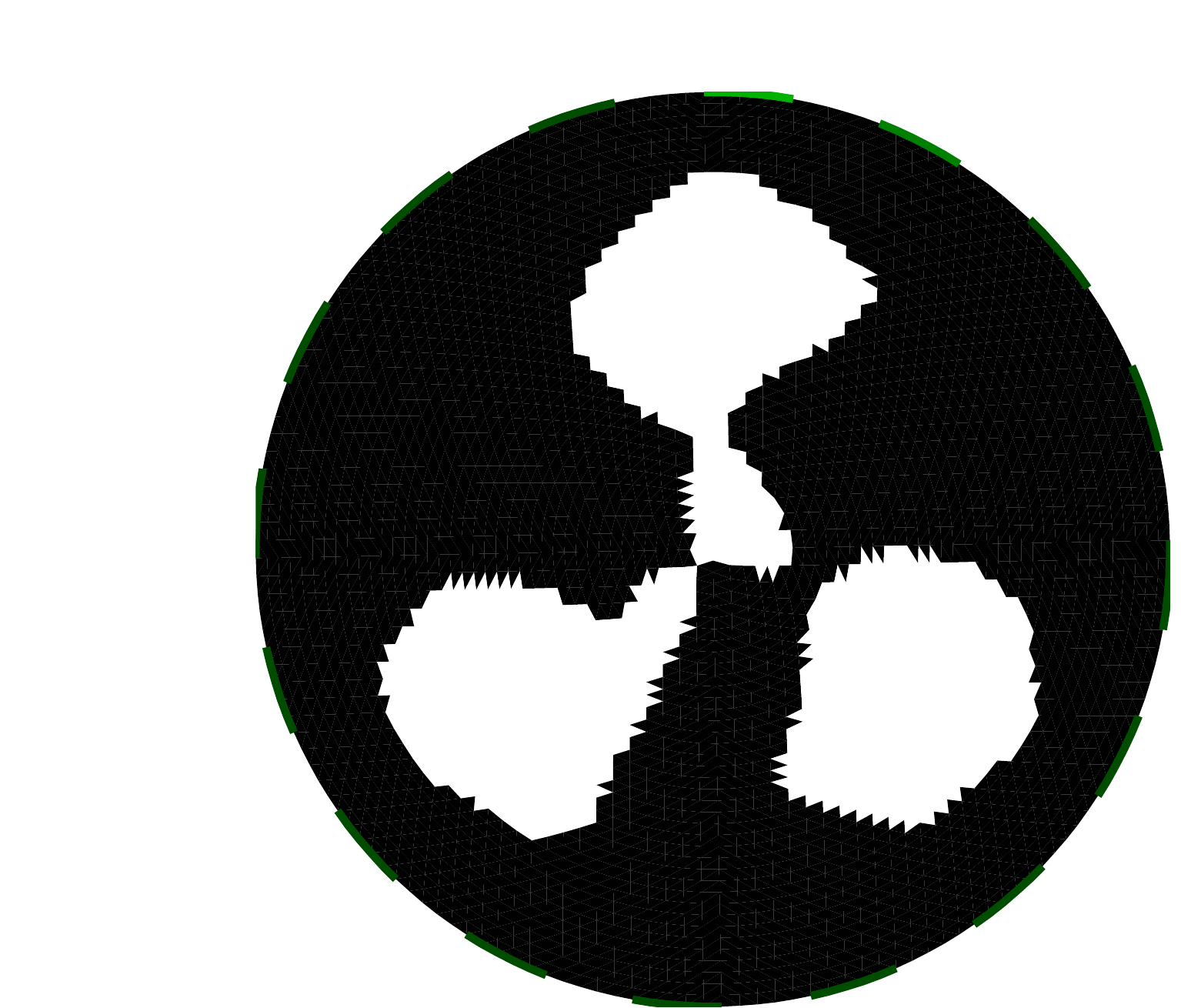}~
\includegraphics[scale=0.18]{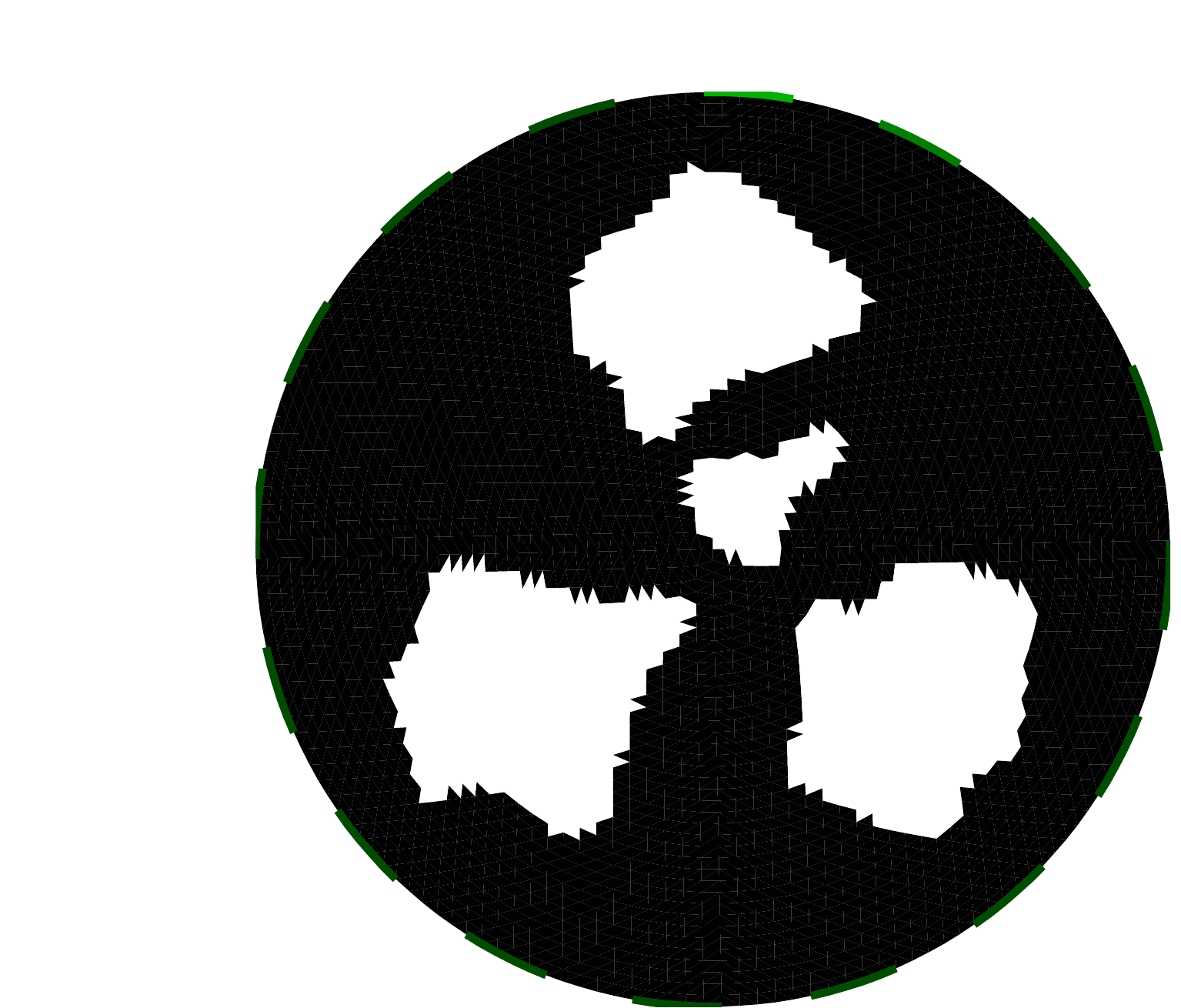}\\
\includegraphics[scale=0.18]{True_LS_EIT}
\caption{Top: Estimates of level-set obtained from 5 experiments with different initial ensembles with samples of $N(0,C)$ with values of $L$ in (\ref{eq:cova1}) (from left to right) $L=0.2, 0.1, 0.06,0.04, 0.03$.  Middle: Conductivities computed from these estimated level-set functions (by means of (\ref{103})). Bottom: True conductivity. } \label{Fig_L_EIT2}
\end{center}
\end{figure}

\begin{figure}[htbp]
\begin{center}
\includegraphics[scale=0.20]{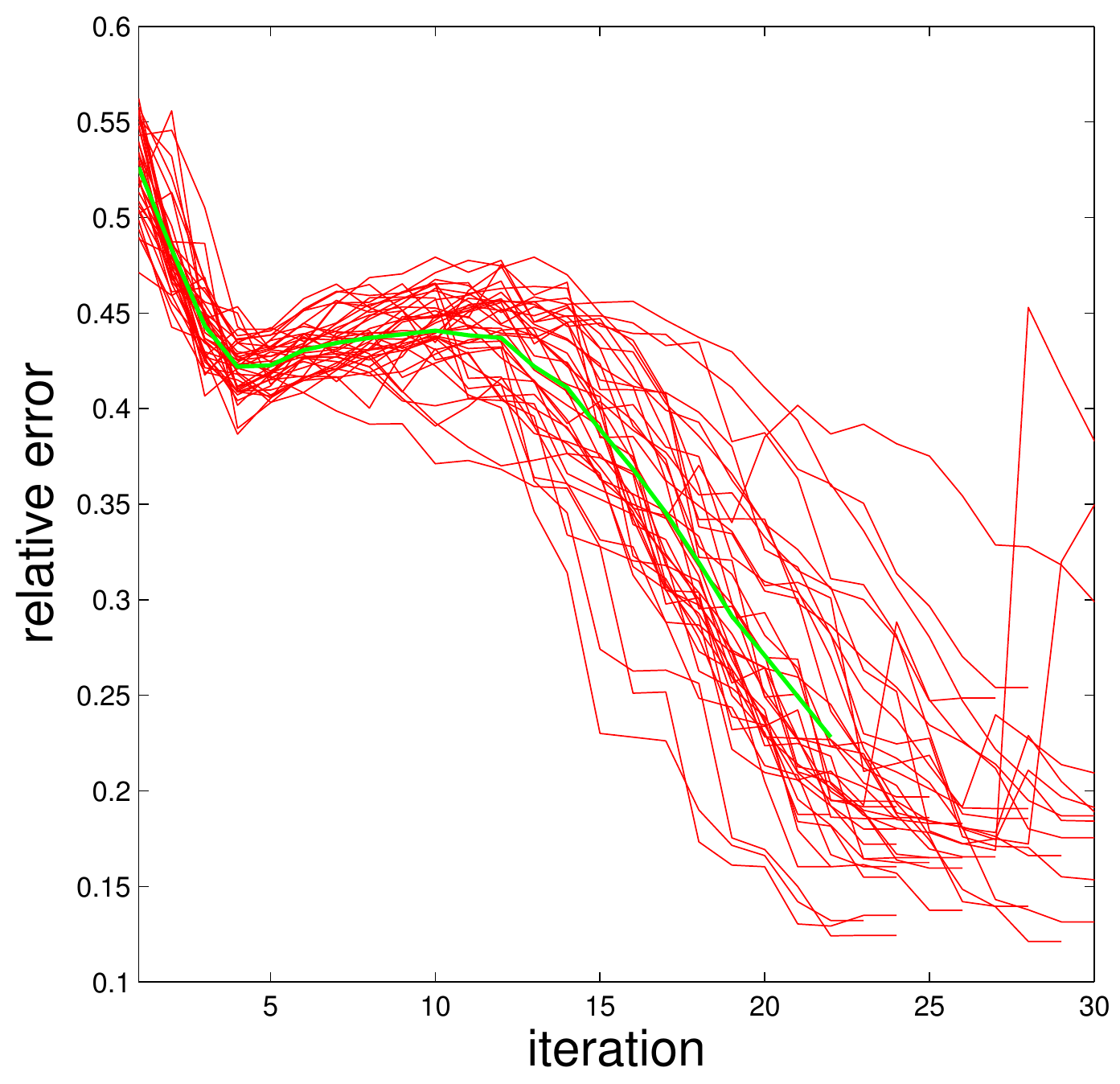}
\includegraphics[scale=0.20]{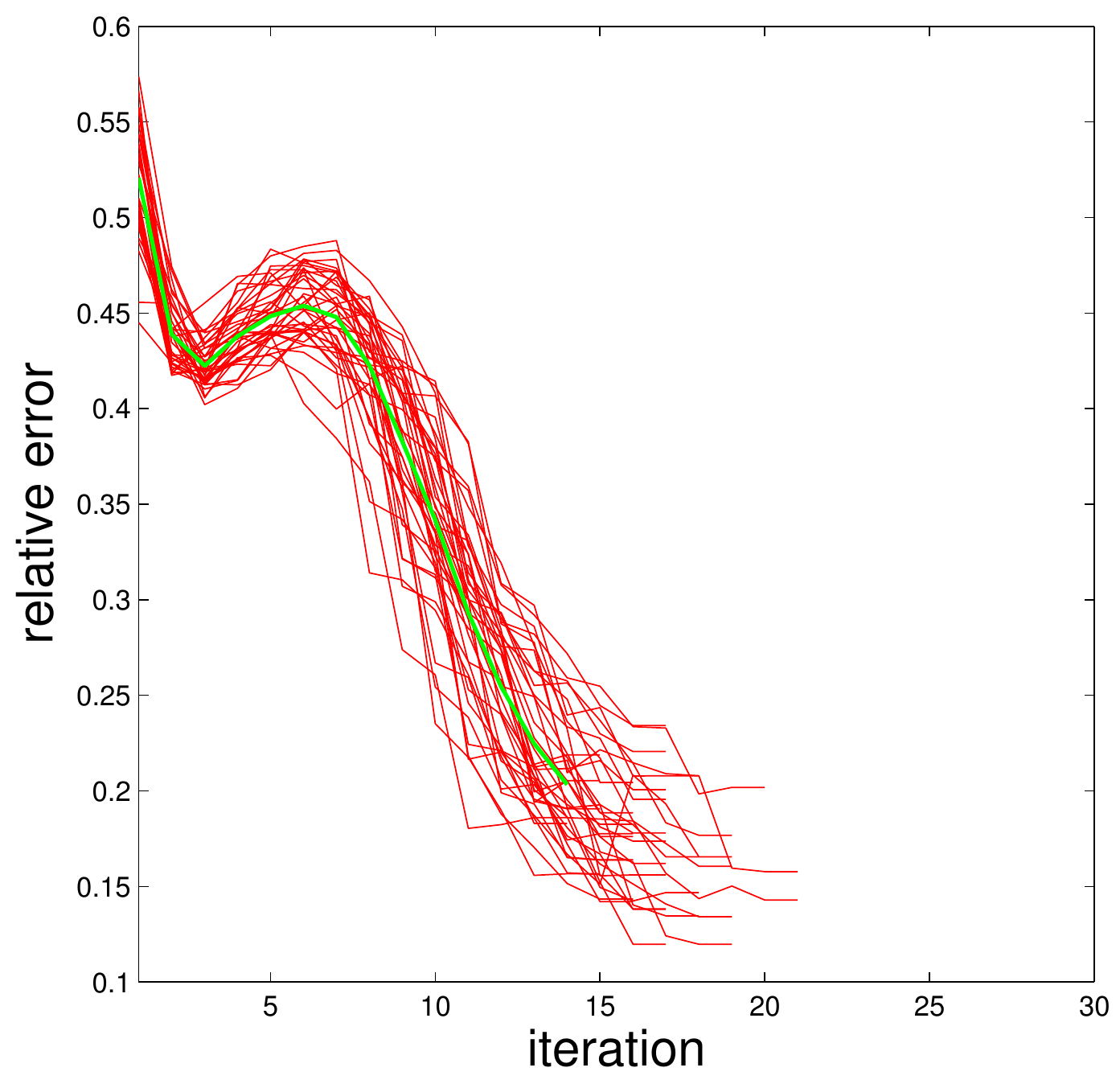}
\includegraphics[scale=0.20]{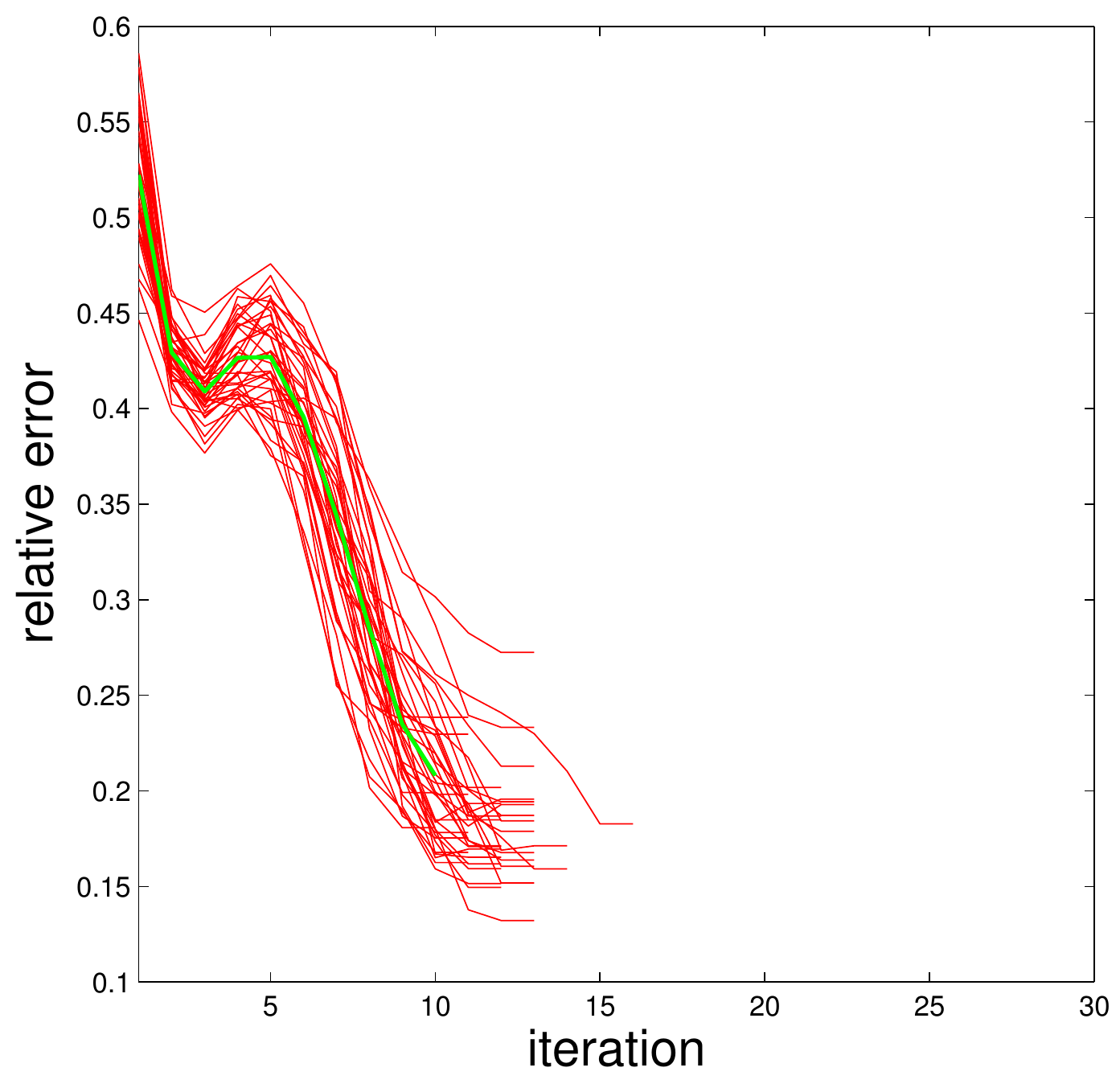}
\includegraphics[scale=0.20]{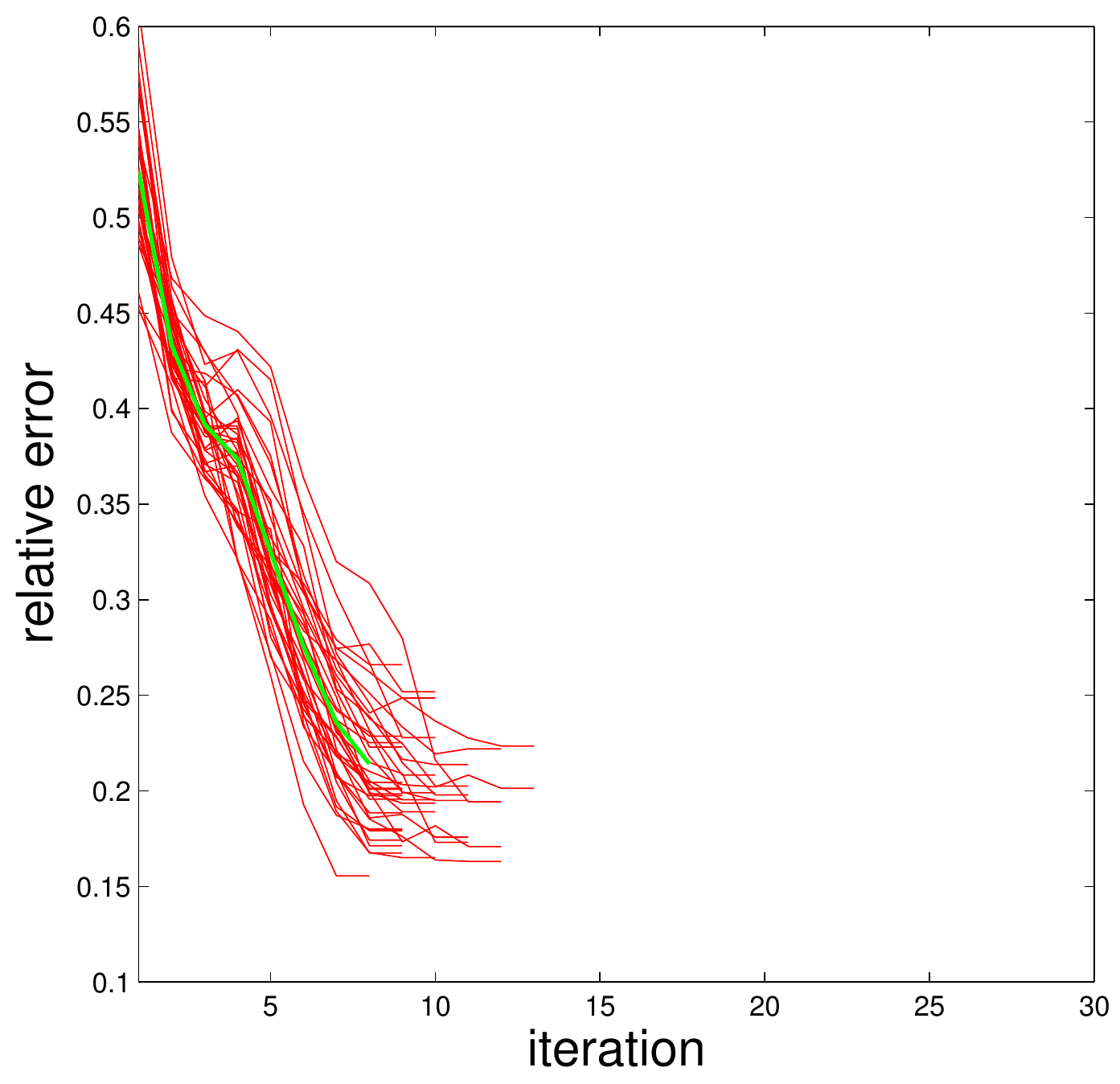}
\includegraphics[scale=0.20]{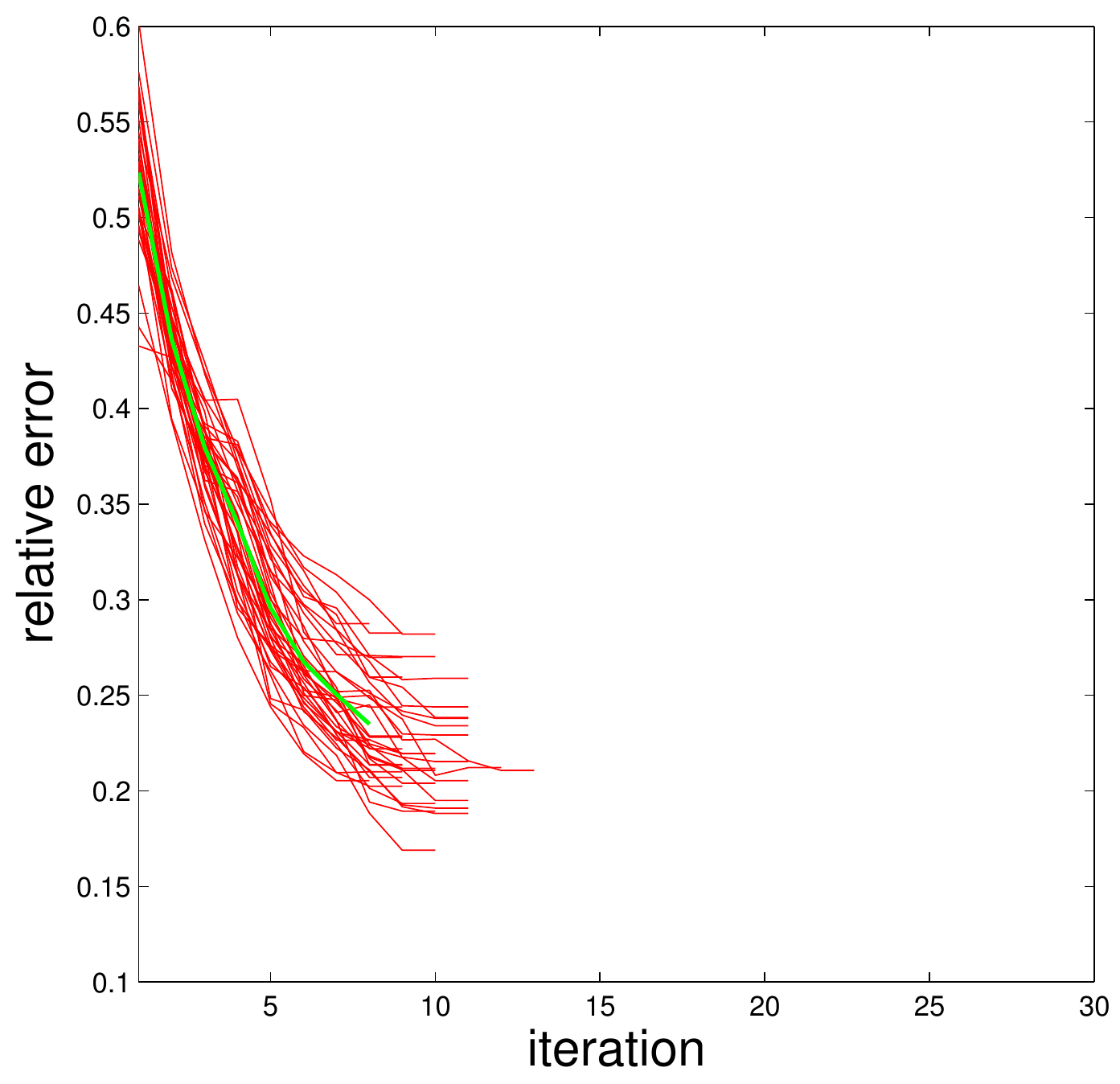}\\
\includegraphics[scale=0.20]{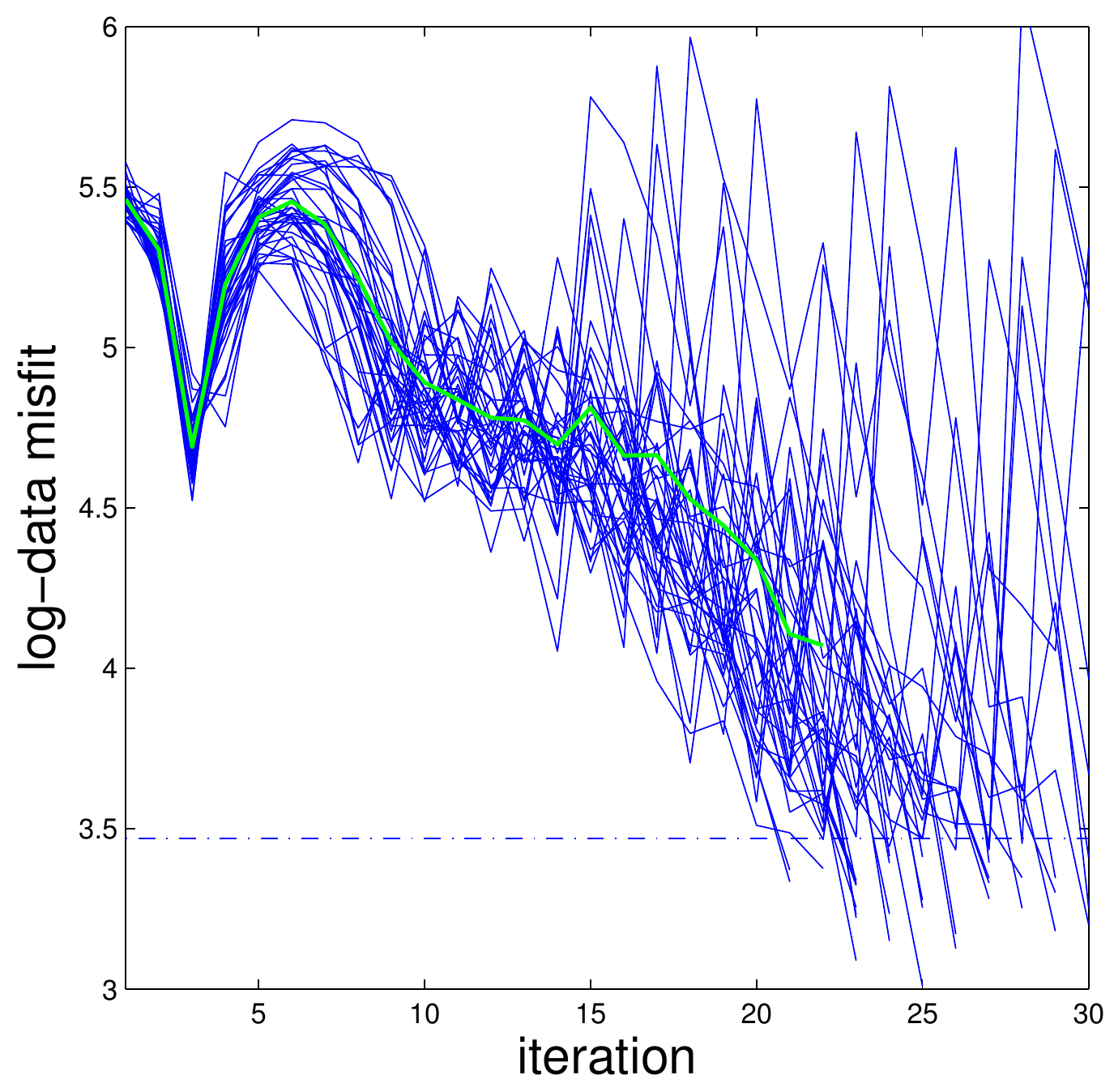}
\includegraphics[scale=0.20]{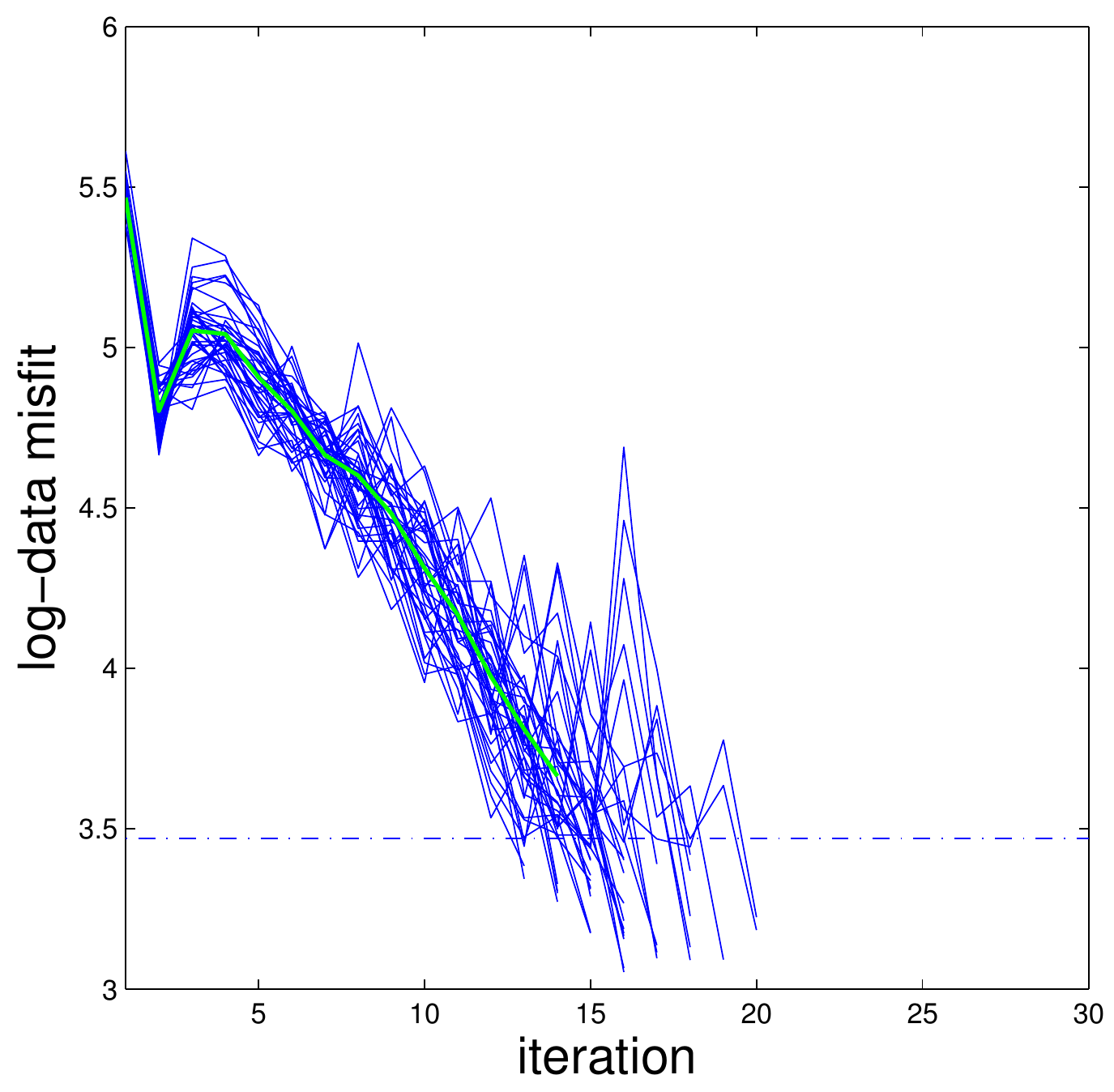}
\includegraphics[scale=0.20]{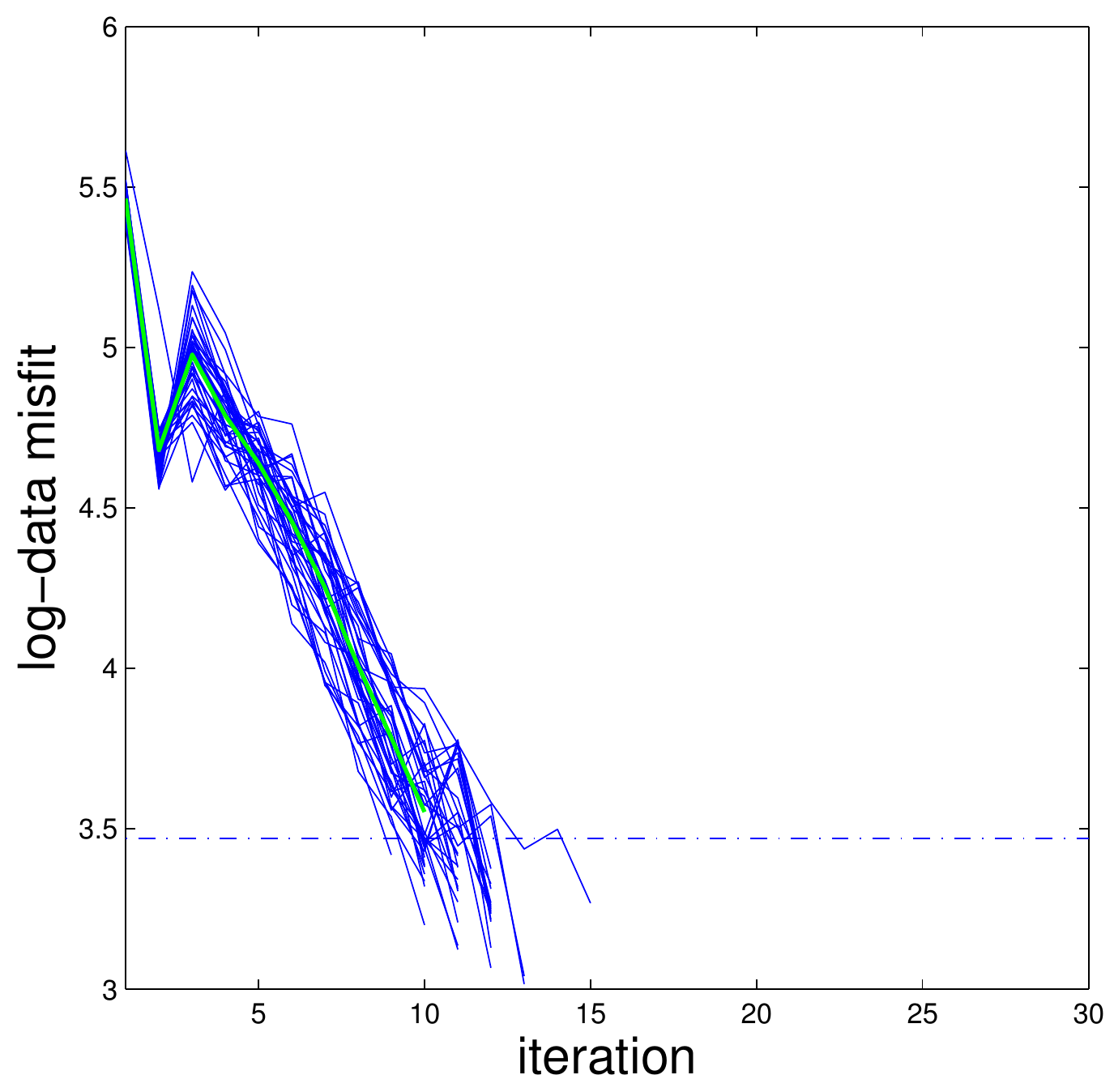}
\includegraphics[scale=0.20]{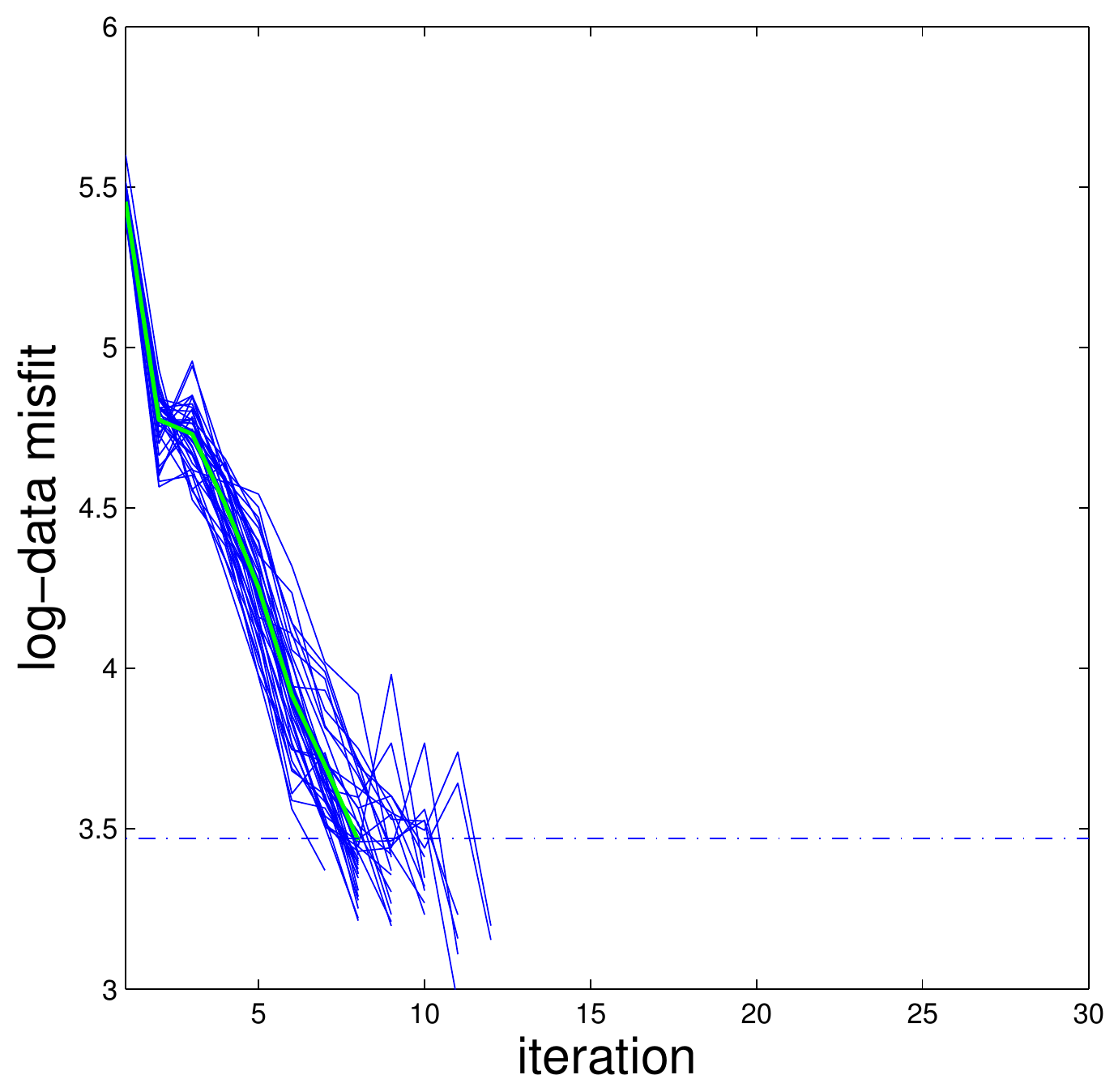}
\includegraphics[scale=0.20]{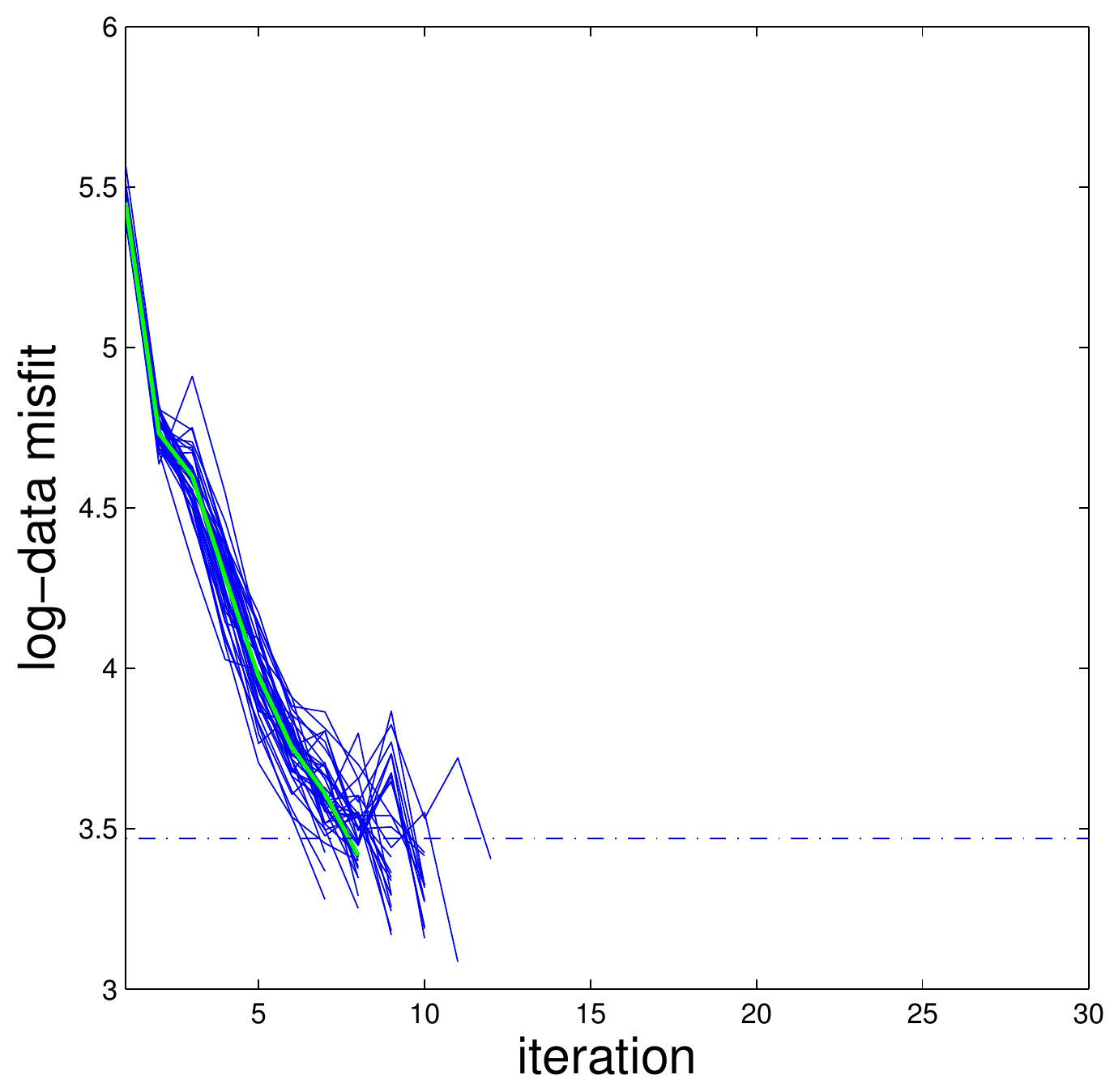}

 \caption{Error with respect to the truth (top) and log data misfit (bottom) from 40 experiments with different initial ensembles generated from samples of $N(0,C)$ with values of $L$ in (\ref{eq:cova1}) (from left to right) $L=0.2, 0.1, 0.06,0.04, 0.03$. }\label{Fig_L_EIT3}

\end{center}
\end{figure}

As in the preceding subsection, we postulate an artificial Gaussian measure $N(0,C)$ that we use to generate an initial ensemble. We consider again $C$ described by (\ref{eq:cova1}). However, note that for this subsection, such Gaussian is used to generate the ensemble of level-set function rather than conductivities. Nonetheless, as in the experiments of subsection \ref{num_EIT}, we expect the parameter $L$ to have an influence on the initial ensemble and this the estimate of the level-set and the associated conductivity obtained with (\ref{103}). In \Fref{Fig_L_EIT1} we present some samples of the level-set function (top) and the corresponding conductivity (bottom) obtained with different values of $L$ ($L=0.2, 0.1, 0.06,0.04, 0.03$) but with the same set of KL coefficients. We clearly observe that the correlation length of the level-set is reflected in the spatial correlation of the interface between the regions of different conductivities.  In \Fref{Fig_L_EIT2} we show estimates of level-set function (top) and conductivities (bottom) obtained, by means of Algorithm \ref{Al1}, with an initial ensemble of $N_{e}=200$ elements generated from the Gaussian distribution described above with the aforementioned values of $L$. We clearly observe that there is an optimal choice of $L$ which yields estimates that visually agree better with the truth displayed in \Fref{Fig_L_EIT2} (bottom). In this case, $0.06\leq L\leq 0.1$ provides (on average) the lowest error w.r.t.truth. In \Fref{Fig_L_EIT3} we show the relative error (top) and the log-data misfit (bottom) from different experiments corresponding to multiple initial ensembles generated with the Gaussian $N(0,C)$ for the aforementioned choices of $L$. For larger correlation lengths $L\ge 0.2$ the data misfit displays larger fluctuations that seem to arise from the appearance of high conductivity close to the electrodes. In addition, very small values of $L$ yields conductivities where the interface between different values has short correlation length. In the current framework, knowledge of the optimal correlation length can be used for the generation of initial ensemble. However, for more general/realistic cases where such parameter is unknown, the estimation of $L$ should be conducted within the method; this is beyond the scope of the present manuscript.

\section{Conclusions}\label{Conclusions}

Our numerical investigation indicates that for sufficiently large ensemble size, the proposed regularizing ensemble Kalman method inherits the regularizing properties of the LM scheme \cite{Hanke}. In other words, when $N_{e}$ is large enough, we observe that the proposed selection of the regularization parameter and early termination of the scheme prevent the lack of stability typical of unregularized schemes. Fortunately, the aforementioned size needed for stability is reasonable and often used in standard ensemble implementations for large-scale applications. For the Darcy problem, for example, we found that with $N_{e}=150$ stable and reasonably accurate estimates can be obtained, in average, within 12 iterations of the scheme when $\rho=0.7$ (and with $M=100$ measurements). The computational cost for computing these estimates is 1800 forward model simulations which is comparable with other gradient-based techniques where also dozens of iterations are need for the convergence, but at each iteration, the explicit computation of the Fr\'echet derivative (if possible at all) has often the cost of $M$ linearized forward model simulations (recall $M$ is the total number of measurements). From the results of the preceding section we can clearly appreciate that the computational cost and accuracy of the proposed method with large enough $N_{e}$ is very competitive with variational iterative regularization methods. However, as we have consistently reiterated throughout this manuscript, the main advantage of the proposed method is that no derivatives of the forward map are needed for the proposed ensemble regularizing scheme.

For the EIT problem, we found that a reasonable ensemble size $N_{e}=100$ is sufficient for the method to display stability in the case of 16 electrodes considered in the subsection \ref{EIT} and for wide class of tunable parameters. However, we observe that the selection of the initial ensemble is crucial for the stability and accuracy of the proposed scheme. In concrete, there is a range for the optimal selection of the spatial correlation length of the members of the initial ensemble. As we expected, good approximations of the truth are obtained if the samples that we use for the generation of the initial ensemble have a correlation length similar to the one of the true conductivity. While prior knowledge of the regularity and spatial features from the truth could be incorporated for the selection of the initial ensemble, we recognize that in more general cases such information may not be available. Further research should address the estimation of the parameters that control the regularity/spatial features of the initial ensemble.

The numerical experiments from the preceding section show that the proposed method is computationally flexible and suitable for a wide class of parameterizations of PDE parameters. We demonstrated that the regularizing ensemble Kalman method can be successfully used for the identification of shapes where the geometry is parameterized with a level-set function. The advantage of having a regularized ensemble method has been showcased with the examples of \Sref{Sec:Applications} where we observe that, with a reasonable ensemble size, the proposed method provides stable computations of the unknown regions of high conductivity. The regularization properties and our selection of the initial ensemble generated from a Gaussian distribution is key for developing a robust level-set scheme that avoids solving the level-set equation for the evolution of the shape; our scheme induces the motion of the shape in a controlled fashion.

Standard ensemble Kalman methods have been numerously applied for solving identification problems where the forward map arises from a large-scale forward model.  However, these standard implementations are often unstable in particular when small ensemble sizes (like the ones used here) are used. We have shown that the proposed regularizing ensemble method, when the ensemble size is large enough, not only addresses the ill-poseness proper of the inverse problem but also the one associated to the small size effect. There is a broad class of ensemble methods where the Fr\'echet derivative of the forward map is approximated with an ensemble \cite{EmeRey}. While most of these methods have been developed to address concrete applications in a statistical setting, they could be potentially applied for the solution of generic parameter identification problems if these methods are empowered with the regularization/stabilization needed to compute the solution of such inverse problems. This manuscript offers numerical evidence that importing ideas from iterative regularization methods can potentially use with ensemble methods to stabilize/regularize the computations of inverse estimates.\\

\noindent{\bf Acknowledgments} 
The author would like to thank Andrew Stuart for his helpful comments on the manuscript. The author also wishes to thank Matthew Dunlop for his help with using EIDORS.

\section{Appendix}
\subsection{Proof of Lemma \ref{lema:Rep1}}
Note that if $X$ is a finite dimensional, $DG(u_{n})$, $D^{\ast}G(u_{n})$ and $C$ in (\ref{eq:m2}) are matrices. Moreover, since in this case $C$ is invertible (recall $C$ is positive definite) simple matrix algebra can be applied to show the desired equivalence \cite{Tarantola}. When $X$ is an infinite-dimensional space the argument for matrices can no longer be applied. In this case, however, we can use a representer-based argument. First we note that, since $Y=\mathbb{R}^{M}$, $DG(u_{n}):X\to Y$ has the form $DG(u_{n})v=[DG_{1}(u_{n})v,\dots ,D_{M}G(u_{n})v]$ with $D_{j}G(u_{n}):X\to \mathbb{R}$ for all $j\in \{1,\dots,M\}$. Moreover, we note that 
\begin{eqnarray}\label{eq:ape1}
wDG_{m}(u_{n})v = \langle v,D_{m}G^{*}(u_{n})w\rangle_{X}
\end{eqnarray}
for all $w\in \mathbb{R}$ and $v\in X$. Therefore,  $$\langle  w,DG(u_{n})v \rangle_{\mathbb{R}^{M}} = \sum_{m=1}^{M}w_{m}D_{m}G(u_{n})v = \Big\langle v, \sum_{m=1}^{M}D_{m}G^{*}(u_{n})w_{m}\Big\rangle_{X}$$ and so $DG^{\ast}(u_{n})w\equiv \sum_{m=1}^{M}D_{m}G^{*}(u_{n})w_{m}$ for all $w\in \mathbb{R}^{M}$. Let us denote $\{e_{i}\}_{i=1}^{M}$ the canonical basis in $\mathbb{R}^{M}$. Note that
\begin{eqnarray}DG(u_{n})\,C\,DG^{*}(u_{n})e_{i}=DG(u_{n})\,C\,DG_{i}^{*}(u_{n})1\nonumber\\
=[D_{1}G(u_{n})\,C\,DG_{i}^{*}(u_{n})1,\dots ,D_{M}G(u_{n})\,C\,DG_{i}^{*}(u_{n})1] 
\end{eqnarray}
 and so
\begin{eqnarray}\label{eq:ape2}
e_{j}^{T}DG(u_{n})\,C\,DG^{*}(u_{n})e_{i}=D_{j}G(u_{n})\,C\,D_{i}G^{*}(u_{n})
\end{eqnarray}
Let $r_{m}\in X$ be the unique ``representer'' such that 
\begin{eqnarray}\label{rep}
D_{m}G(u_{n})v =\langle C^{-1/2}r _{m},C^{-1/2}v\rangle_{X}
\end{eqnarray}
 for all $v\in X$. On the other hand, from (\ref{eq:ape1}) we have $D_{m}G(u_{n})v = \langle v,DG_{m}^{*}(u_{n})1\rangle_{X}$ for all $m\in \{1,\dots,M\}$, then $C^{-1}r_{m}= DG_{m}^{*}(u_{n})1$ and so, formally, 
\begin{eqnarray}\label{eq:ape3}
r_{m}=CDG_{m}^{*}(u_{n})1
\end{eqnarray}
We may now consider the following representation in terms of representers
\begin{eqnarray}\label{eq:ape4}
v=\sum_{m=1}^{M}\beta_{m}r_{m} +b
\end{eqnarray}
where $b\perp r_{j}$ in $X$ for all $j\in \{1,\dots,M\}$ (i.e. $\langle C^{-1/2} r _{m},C^{-1/2}b\rangle_{X}$). From (\ref{eq:ape3})  and (\ref{eq:ape4}) we have
\begin{eqnarray}\label{eq:ape5}
\fl DG(u_{n})v=\sum_{m=1}^{M}\beta_{m}DG(u_{n})r_{m} +DG(u_{n})b=\sum_{j=1}^{M}\beta_{m}DG(u_{n}) CD_{m}^{*}G(u_{n})1 +DG(u_{n})b
\end{eqnarray}
Note that $D_{m}G(u_{n})b=\langle C^{-1/2}r_{m},C^{-1/2}b\rangle_{X} =0$. Therefore,  we arrive at
\begin{eqnarray}\label{eq:ape6}
DG(u_{n})v=\sum_{j=1}^{M}\beta_{m}DG(u_{n}) C D_{m}G^{*}(u_{n})1= DG(u_{n}) C DG^{*}(u_{n})\beta
\end{eqnarray}
where $\beta=(\beta_{1},\dots, \beta_{M})$. Similarly, from (\ref{eq:ape2})- (\ref{eq:ape4}) we obtain
\begin{eqnarray}\label{eq:ape7}
\fl \vert\vert C^{-1/2}v\vert\vert_{X}^2=\sum_{jm}\beta_{j}
\beta_{m}\langle C^{-1/2} r_{m},C^{-1/2} r_{j}\rangle_{X} + \langle C^{-1/2} b,C^{-1/2}b\rangle_{X}\nonumber\\
\fl =\sum_{jm}\beta_{j}
\beta_{m}D_{m}G(u_{n})r_{j}+\vert\vert C^{-1/2}b \vert\vert_{X}^2  =\beta^{T}DG(u_{n}) C DG^{*}(u_{n})\beta+\vert\vert C^{-1/2}b \vert\vert_{X}^2 
\end{eqnarray}
Combining (\ref{eq:ape6}) and (\ref{eq:ape7}) we write 
\begin{eqnarray*}
\fl J(v)\equiv \vert\vert \Gamma^{-1/2}(y^{\eta}-G(u_{n}))-DG(u_{n})v\vert\vert_{Y}^2+\alpha\vert\vert C^{-1/2} v\vert\vert_{X}^{2}
\end{eqnarray*}
as
\begin{eqnarray}\label{eq:ape7B}
\fl J_{LM}(\beta,b)\equiv \vert\vert \Gamma^{-\frac{1}{2}}(y-G(u_{n})-DG(u_{n}) C DG^{*}(u_{n})\beta)\vert\vert^2\nonumber\\ +\alpha \beta^{T}DG(u_{n}) C DG^{*}(u_{n})\beta+\alpha \vert\vert C^{-1/2}b \vert\vert_{X}^2 
\end{eqnarray}
From the linearly independece of $\{D_{m}G(u_{n})\}_{m=1}^{M}$, it is not difficult to see that the matrix $DG(u_{n}) C DG^{*}(u_{n})$ is invertible. Then, from standard arguments it follows that the unique minimizer of (\ref{eq:ape7B}) in $\mathbb{R}^{M}\times (\textrm{span}\{r_{m}\}_{m=1}^{M})^{\perp}$ is
\begin{eqnarray}\label{eq:ape8}
\beta= [DG(u_{n}) C DG^{*}(u_{n}) + \alpha \Gamma ]^{-1}( y-G(u_{n})),\qquad b=0
\end{eqnarray}
Expression  (\ref{eq:m6}) then follows from writing $u_{n+1}-u_{n}=v$ and using (\ref{eq:ape3})-(\ref{eq:ape4}) and (\ref{eq:ape8}). $\Box$

\section*{References}
\bibliography{Ensemble_bib}

\begin{thebibliography}{10}

\bibitem{EnKF_Review}
S.~I. Aanonsen, G.~Nævdal, D.~S. Oliver, A.~C. Reynolds, and B~Valles.
\newblock The ensemble kalman filter in reservoir engineering--a review.
\newblock {\em SPE Journal}, 14.

\bibitem{EIDORS}
Andy Adler and William R~B Lionheart.
\newblock Uses and abuses of eidors: an extensible software base for eit.
\newblock {\em Physiological Measurement}, 27(5):S25, 2006.

\bibitem{Bear}
J.~Bear.
\newblock {\em Dynamics of Fluids in Porous Media}.
\newblock Dover Pulications, New York, 1972.

\bibitem{EIT_revew}
Liliana Borcea.
\newblock Electrical impedance tomography.
\newblock {\em Inverse Problems}, 18(6):R99, 2002.

\bibitem{BurgerGB}
M.~Burger.
\newblock A framework for the construction of level set methods for shape
  optimization and reconstruction.
\newblock {\em Interfaces and Free Boundaries}, 5:301--329, 2002.

\bibitem{BurgerSurvey}
M.~Burger and S.~Osher.
\newblock A survey on level set methods for inverse problems and optimal
  design.
\newblock {\em Eur. J. Appl. Math.}, 16:263--301, 2005.

\bibitem{Carrera}
J.~Carrera and S.~P. Neuman.
\newblock Estimation of aquifer parameters under transient and steady state
  conditions: 3. application to synthetic and field data.
\newblock {\em Water Resources Research}, 22.

\bibitem{LS_EIT1}
Eric~T. Chung, Tony~F. Chan, and Xue-Cheng Tai.
\newblock Electrical impedance tomography using level set representation and
  total variational regularization.
\newblock {\em Journal of Computational Physics}, 205(1):357 -- 372, 2005.

\bibitem{geos}
C.~V. Deutsch.
\newblock {\em Geostatistical Reservoir Modeling}.
\newblock Oxford University Press, Oxford, 2002.

\bibitem{EmeRey}
AlexandreA. Emerick and AlbertC. Reynolds.
\newblock Investigation of the sampling performance of ensemble-based methods
  with a simple reservoir model.
\newblock {\em Computational Geosciences}, 17(2):325--350, 2013.

\bibitem{Ill-posed}
H~W Engl, K~Kunisch, and A~Neubauer.
\newblock Convergence rates for {T}ikhonov regularisation of non-linear
  ill-posed problems.
\newblock {\em Inverse Problems}, 5(4):523, 1989.

\bibitem{Engl}
H.W. Engl, M.~Hanke, and A.~Neubauer.
\newblock {\em Regularization of inverse problems}, volume 375.
\newblock Springer, 1996.

\bibitem{evensen2009data}
G.~Evensen.
\newblock {\em Data assimilation: the ensemble Kalman filter}.
\newblock Springer Verlag, 2009.

\bibitem{Hanke}
M.~Hanke.
\newblock A regularizing {L}evenberg-{M}arquardt scheme, with applications to
  inverse groundwater filtration problems.
\newblock {\em Inverse Problems}, 13:79--95, 1997.

\bibitem{Repre}
M.~A. Iglesias and C.~Dawson.
\newblock The representer method for state and parameter estimation in
  single-phase {D}arcy flow.
\newblock {\em Computer Methods in Applied Mechanics and Engineering},
  196(1):4577--4596, 2007.

\bibitem{EnKF_US}
Marco~A Iglesias, Kody J~H Law, and Andrew~M Stuart.
\newblock Ensemble kalman methods for inverse problems.
\newblock {\em Inverse Problems}, 29(4):045001, 2013.

\bibitem{Level_set_US}
Marco~A Iglesias, Yulong Lu, and Andrew Stuart.
\newblock A bayesian level-set method for geometric inverse problems.
\newblock {\em submitted}, 2015.

\bibitem{Iglesias3}
Marco~A Iglesias and Dennis McLaughlin.
\newblock Level-set techniques for facies identification in reservoir modeling.
\newblock {\em Inverse Problems}, 27(3):035008, 2011.

\bibitem{Yo}
MarcoA. Iglesias.
\newblock Iterative regularization for ensemble data assimilation in reservoir
  models.
\newblock {\em Computational Geosciences}, pages 1--36, 2014.

\bibitem{LM}
MarcoA. Iglesias and Clint Dawson.
\newblock The regularizing levenberg‚Äìmarquardt scheme for history
  matching of petroleum reservoirs.
\newblock {\em Computational Geosciences}, pages 1--21, 2013.

\bibitem{Evaluation}
MarcoA. Iglesias, KodyJ.H. Law, and AndrewM. Stuart.
\newblock Evaluation of gaussian approximations for data assimilation in
  reservoir models.
\newblock {\em Computational Geosciences}, pages 1--35, 2013.

\bibitem{Iterative}
B.~Katltenbacher, A.~Neubauer, and O.~Scherzer.
\newblock {\em {I}terative {R}egularization {M}ethods for {N}onlinear
  {I}ll-{P}osed {P}roblems}.
\newblock Radon Series on Computational and Applied Mathematics, de Gruyter,
  Berlin, 1st edition, 2008.

\bibitem{Kirsch}
Andreas Kirsch.
\newblock {\em An Introduction to the Mathematical Theory of Inverse Problems}.
\newblock Springer-Verlag New York, Inc., New York, NY, USA, 1996.

\bibitem{IterativeEnKF}
G.~Li and A.~Reynolds.
\newblock Iterative {E}nsemble {K}alman {F}ilters for {D}ata {A}ssimilation.
\newblock {\em in Proccedings of the SPE Annual Technical Conference and
  Exhibition, Anaheim, California, USA, 11-14 Noviembre}, {S}{P}{E} 109808,
  2007.

\bibitem{Besov}
M.Dashti, S.~Harris, and A.M.Stuart.
\newblock Besov priors for bayesian inverse problems.
\newblock {\em Inverse Problems and Imaging}, 6:183--200, 2012.

\bibitem{Osher}
S.~Osher and J.A. Sethian.
\newblock Fronts propagating with curvature dependent speed: Algorithms based
  on hamilton-jacobi formulations.
\newblock {\em Journal of Computational Physics}, 79(1):12--49, 1988.

\bibitem{LS_EIT2}
Peyman Rahmati and Andy Adler.
\newblock A level set based regularization framework for eit image
  reconstruction.
\newblock {\em Journal of Physics: Conference Series}, 434(1):012083, 2013.

\bibitem{WM}
Lassi Roininen, Janne M.~J. Huttunen, and Sari Lasanen.
\newblock Whittle-mate ́rn priors for bayesian statistical inversion with
  applications in electrical impedance tomography.
\newblock {\em Inverse Problems and Imaging}, 8(2), 2014.

\bibitem{Mixed}
T.F. Russell and M.F. Wheeler.
\newblock Finite element and finite difference methods for continuous flows in
  porous media.
\newblock {\em In: R.E. Ewing, Editor, Mathematics of Reservoir Simulation},
  SIAM,Philadelphia, PA.

\bibitem{Santosa}
F.~Santosa.
\newblock A level-set approach for inverse problems involving obstacles.
\newblock {\em ESAIM Contr�le Optim. Calc. Var}, 1:17--33, 1996.

\bibitem{cheney}
Erkki Somersalo, Margaret Cheney, and David Isaacson.
\newblock Existence and uniqueness for electrode models for electric current
  computed tomography.
\newblock {\em SIAM Journal on Applied Mathematics}, 52(4):pp. 1023--1040,
  1992.

\bibitem{Andrew}
A.M. Stuart.
\newblock Inverse problems: a {B}ayesian perspective.
\newblock In {\em Acta Numerica}, volume~19. 2010.

\bibitem{Tarantola}
A.~Tarantola.
\newblock {\em Inverse {P}roblems {T}heory: {M}ethods for {D}ata {F}itting and
  {M}odel {P}arameter {E}stimation}.
\newblock Elsevier,, New York, 1st edition, 1987.

\end{thebibliography}
\bibliographystyle{plain}

\end{document}